\documentclass{article}

\usepackage{amssymb,amsfonts,amsmath}
\usepackage{cite,enumerate,float}
\usepackage{color}

\usepackage{physics}

\usepackage{mathrsfs}

\usepackage{tikz}
\usepackage{tikz-3dplot}
\usepackage{xparse}
\usetikzlibrary{arrows}
\usetikzlibrary{arrows.meta}
\usetikzlibrary{positioning}
\usetikzlibrary{shapes,snakes}
\usetikzlibrary{fit}
\usetikzlibrary{decorations.pathmorphing,decorations.pathreplacing,decorations.markings, matrix,fit,backgrounds,positioning}
\usetikzlibrary{calc}
\usepackage{inconsolata}
\usepackage{ytableau}

\usepackage{caption}
\usepackage{subcaption}

\usepackage{xcolor}
\definecolor{burgundy}{rgb}{0.5, 0.0, 0.13}
\definecolor{olive}{rgb}{0.50, 0.50, 0.0}
\usepackage[linktocpage=true,colorlinks=true,linkcolor=burgundy,citecolor=black!20!blue,urlcolor=black!40!violet]{hyperref}

\usepackage[most]{tcolorbox}

\usepackage{pdflscape}
\usepackage{array, longtable}
\newcolumntype{C}{>{$}c<{$}}

\def\be{\begin{eqnarray}}
\def\ee{\end{eqnarray}}

\def\p{\partial}

\def\Tr{{\rm Tr}\,}

\definecolor{red}{rgb}{1,0,0}
\definecolor{orange}{rgb}{1,0.5,0}
\definecolor{violet}{rgb}{0.7,0,1}

\def\CA{{\cal A}}

\def\CI {{\cal I}}
\def\CJ {{\cal J}}

\def\CM {{\cal M}}
\def\CN {{\cal N}}

\def\CP {{\cal P}}
\def\CR {{\cal R}}

\def\CW {{\cal W}}

\def\CI {{{\cal I}}}

\def\CQ {{\cal Q}}

\def\IC{\mathbb{C}}

\def\ff{\mathfrak{f}}

\def\fg{\mathfrak{g}}

\def\fl{\mathfrak{l}}

\def\fs{\mathfrak{s}}

\def\fs{\mathfrak{s}}

\def\fA{\mathfrak{A}}

\def\fS{\mathfrak{S}}

\def\bPsi{{\boldsymbol{\Psi}}}

\def\lm{\limits}

\def\myY{\mathsf{Y}}

\DeclareSymbolFont{bbsymbol}{U}{bbold}{m}{n}
\DeclareMathSymbol{\bbzero}{\mathbin}{bbsymbol}{"30}
\DeclareMathSymbol{\bbone}{\mathbin}{bbsymbol}{"31}
\DeclareMathSymbol{\bbtwo}{\mathbin}{bbsymbol}{"32}
\DeclareMathSymbol{\bbthree}{\mathbin}{bbsymbol}{"33}
\DeclareMathSymbol{\bbfour}{\mathbin}{bbsymbol}{"34}
\DeclareMathSymbol{\bbfive}{\mathbin}{bbsymbol}{"35}
\DeclareMathSymbol{\bbsix}{\mathbin}{bbsymbol}{"36}
\DeclareMathSymbol{\bbseven}{\mathbin}{bbsymbol}{"37}
\DeclareMathSymbol{\bbeight}{\mathbin}{bbsymbol}{"38}
\DeclareMathSymbol{\bbnine}{\mathbin}{bbsymbol}{"39}

\newcommand{\ceil}[1]{\left\lceil #1 \right\rceil}

\newcommand\sqbox[1]{{
	\setbox0=\hbox{\mbox{$\Box$}}
	\setbox1=\hbox{\mbox{\raisebox{0.35ex}{\small #1}}}
	\mbox{\raisebox{-0.2ex}{\rlap{\hbox to \wd0{\hss{\box1}\hss}}\box0}}
	}}
	\newcommand\ssqbox[1]{{
	\setbox0=\hbox{\mbox{$\scriptstyle\Box$}}
	\setbox1=\hbox{\mbox{\raisebox{0.35ex}{\tiny #1}}}
	\mbox{\raisebox{-0.2ex}{\rlap{\hbox to \wd0{\hss{\box1}\hss}}\box0}}
	}}
	
	\def\myblue{white!40!blue}
	\def\mygray{white!40!gray}
	
	\def\myred{white!40!red}
	\def\mygreen{white!40!green}
	\def\mypurple{white!40!purple}

	\tcbset{boxrule=0.6pt,colback=white!97!blue}

	\NewDocumentCommand{\GTrow}{O{0} m}{%
		\begingroup
		\def\spacing{1.05}%
		\def\row{#2}%
		\pgfmathsetmacro{\N}{0}%
		\foreach \i in \row { \pgfmathparse{\N+1} \global\let\N=\pgfmathresult }%
		\foreach \val [count=\i from 1] in \row {
			\pgfmathsetmacro{\x}{(\i - (\N+1)/2)*\spacing}
			\node at (\x, #1) {\val};
		}
		\endgroup
	}

	\DeclareMathOperator{\sgn}{\mathop{sgn}}
	\DeclareMathOperator{\rk}{\mathop{rk}}
	\newcommand{\Eul}{\text{Eul}}
	\renewcommand{\d}{\partial}
	\renewcommand{\bar}{\overline}

	\usepackage{tabularray} 
	
	\hoffset= - 1.0in         

	\numberwithin{equation}{section}
\begin{document}

\hfill MIPT/TH-26/25

\hfill ITEP/TH-35/25


\vskip 1.5in
\begin{center}
	
	{\bf\Large Quiver Yangian algebras associated to Dynkin diagrams of A-type and their rectangular representations}
	
	\vskip 0.2in
	\renewcommand{\thefootnote}{\alph{footnote}}
	{
		A. Gavshin$^{1,2,3,}$\footnote[1]{e-mail: gavshin.an@phystech.edu}	
	}
	\vskip 0.2in
	\renewcommand{\thefootnote}{\roman{footnote}}
	{\small{
			\textit{$^1$MIPT, 141701, Dolgoprudny, Russia}
			\vskip 0 cm
			\textit{$^2$NRC ``Kurchatov Institute'', 123182, Moscow, Russia}
			%
			%
			\vskip 0 cm
			\textit{$^3$ITEP, Moscow, Russia}
	}}
\end{center}

\vskip 0.2in
\baselineskip 16pt

\begin{abstract}
	The connection between simple Lie algebras and their Yangian algebras has a long history. In this work, we construct finite-dimensional representations of Yangian algebras $\myY(\fs\fl_{n})$ using the quiver approach. Starting from quivers associated to Dynkin diagrams of type A, we construct a family of quiver Yangians. We show that the quiver description of these algebras enables an effective construction of representations with a single non-zero Dynkin label. For these representations, we provide an explicit construction using the equivariant integration over the corresponding quiver moduli spaces. The resulting states admit a crystal description and can be identified with the Gelfand-Tsetlin bases for $\fs\fl_{n}$ algebras. Finally, we show that the resulting Yangians possess notable algebraic properties, and the algebras are isomorphic to their alternative description known as the second Drinfeld realization.
\end{abstract}

\ytableausetup{boxsize = 0.5em}
\tableofcontents

\section{Introduction}\label{sec: Introduction}

Yangian algebras were originally introduced by Drinfeld \cite{Drinfeld:1985rx} for simple finite-dimensional Lie algebras $\fg$. The Yangians were defined as a canonical deformation of the universal enveloping algebra U$(\fg[z])$ for the corresponding current algebra $\fg[z]$. However, this definition turned out to be inconvenient for describing the representations of these algebras.

Later Drinfeld \cite{Drinfeld:1987sy} addressed these limitations and introduced an alternative definition of these algebras that resembled the Chevalley bases of the corresponding Lie algebras. 
The new construction, also known as the Drinfeld second realization, highlighted the highest-weight structure of the finite-dimensional irreducible representations and completely classified them in terms of Drinfeld polynomials \cite{Drinfeld:1987sy, chari1994guide}.

Another fruitful perspective on Yangian algebras has emerged from the context of superstring theory. In particular, the study of BPS states in systems of D-branes wrapping toric Calabi-Yau three-folds gives rise to algebraic structures known as BPS algebras \cite{Harvey:1995fq, Harvey:1996gc, Kontsevich:2008fj, Kontsevich:2010px, Galakhov:2018lta, rapcak2021branesquiversbpsalgebras}. These algebras can be classified by quivers that encode the effective gauge theories describing the low-energy dynamics of brane systems. In these settings, the quivers are constructed from toric diagrams via brane tiling \cite{Yamazaki:2008bt, Li:2020rij}. These Yangians are commonly referred to as quiver Yangian algebras in the literature.

One of the richest and most interesting directions of research is to study the representation theory of these algebras. Although the full theory is still under development, numerous important results have already been established. In the paper, we focus on the special class of Yangian representations --- called crystal representations --- whose states can be described by the statistical model of crystal melting \cite{Ooguri:2009ijd, Aganagic:2010qr, Yamazaki:2010fz}. The requirement that crystals transform into valid neighboring crystals under the action of the Yangian allows one to bootstrap the corresponding matrix elements of the generators in terms of meromorphic functions \cite{Li:2020rij}. An alternative and more computationally convenient approach is to use equivariant localization techniques \cite{Rapcak:2018nsl, Galakhov:2020vyb}, where the matrix elements are computed as integrals over quiver moduli spaces \cite{NakajimaALE}. Duistermaat-Heckman integration formulae \cite{Pestun:2016qko, guillemin2013supersymmetry} ensure that the calculations localize to contributions from the neighborhoods of fixed points on the quiver varieties.  

There has been a particular interest in affine quiver Yangians, especially $\myY(\widehat{\fg\fl}_{n})$ and its supersymmetric counterparts $\myY(\widehat{\fg\fl}_{n|m})$ that also fall into the latter category. These affine algebras have been studied extensively, have wide-ranging applications, and have rich algebraic structures. The simplest example of these algebras is the Yangian $\myY(\widehat{\fg\fl}_{1})$ \cite{Prochazka:2015deb, Tsymbaliuk_2017, maulik2018quantumgroupsquantumcohomology}, which admits the well-known Fock representation. This representation is a crystal representation whose states are parameterized by Young diagrams, and it plays a central role in the theory of quantum integrable systems such as Calogero-Moser-Sutherland systems \cite{calogero1971solution, MOSER1975197} and WLZZ models \cite{Wang:2022fxr, Mironov:2023pnd, Mironov:2023zwi}. The eigenfunctions of Calogero Hamiltonians are classical Schur/Jack polynomials \cite{Galakhov_2024, nakajima1996jackpolynomialshilbertschemes, van2012calogero}, which are tightly related to superintegrability \cite{Mishnyakov:2024cgl, Mironov_2022}. The study of similar families of polynomials is itself a broad and active research area \cite{Azheev:2025wti}. Notably, the Yangian structure allows natural generalizations of these families. For instance, the MacMahon representation of $\myY(\widehat{\fg\fl}_{1})$ yields 3-dimensional analogues of Young diagrams and Jack polynomials \cite{Morozov_2023}. The lift to the $\myY(\widehat{\fg\fl}_{n})$ algebras gives rise to colored Young diagrams and Uglov polynomials \cite{Galakhov:2024mbz, Mishnyakov:2024cgl}.

In the supersymmetric case, affine super Yangians $\myY(\widehat{\fg\fl}_{n|m})$ generalize these structures further. The most prominent example is the algebra $\myY(\widehat{\fg\fl}_{1|1})$ \cite{Galakhov:2023mak, Galakhov:2024foa, Galakhov_2024}, which admits semi-Fock crystal representation with states enumerated by super-Young diagrams. In a similar way, the corresponding family of orthogonal polynomials  emerges, known as super-Schur/Jack polynomials \cite{Galakhov:2023mak, Galakhov_2024mpsp}.

One can extend the ideas above and consider the quiver Yangians associated with arbitrary affine Dynkin diagrams \cite{Bao:2023ece, Li:2023zub}. While affine Yangians exhibit rich structure and wide applications, their analysis is often complicated by the presence of infinite-dimensional representations and intricate algebraic relations. In this paper, we focus instead on the Yangian algebras associated with finite Dynkin diagrams. These algebras admit finite-dimensional irreducible representations and possess simpler, more tractable structures, which makes them suitable for explicit construction.

More specifically, we explore the Yangian algebras $\myY(\fs\fl_{n})$ associated with Dynkin diagrams of A type. We use the quiver approach combined with equivariant integration techniques to explicitly describe their representations. The case of $\myY(\fs\fl_{2})$ was previously described using similar methods in \cite{Galakhov_2024}; see also \cite{Bykov:2019cst, Yang:2024ubh}. While the approach is robust and, in theory, can extend beyond Dynkin classification, it naturally describes only the highest-weight representations with a single non-zero highest weight. These representations are labeled by rectangular Young diagrams; therefore, we refer to them as \textit{rectangular}. Notably, the states of these representations still have crystal structure. Moreover, they allow natural parametrization in terms of the Gelfand-Tsetlin bases \cite{molev2002gelfandtsetlinbasesclassicallie, Gelfand:1950ihs}. This correspondence we gradually introduce through the text.

The paper is organized as follows. Section \ref{sec: QYA} provides an overview of the construction of quiver Yangian algebras. We start with general data and gradually focus on the Dynkin quivers of A type. In section \ref{sec: Y(sl(n))} we describe the algebras $\myY(\fs\fl_{n})$ and their representations in detail. We begin with a warm-up example, the Yangian $\myY(\fs\fl_{3})$, in section \ref{subsec: Y(sl(3))}, proceed to the representations of the $\myY(\fs\fl_{4})$ algebra in section \ref{subsec: Y(sl(4))}, and generalize the previous results to an arbitrary $\myY(\fs\fl_{n})$ algebra in section \ref{subsec: Y(sl(n)) details}. Finally, in section \ref{subsec: Comments} we show that the constructed Yangians are in fact isomorphic to Drinfeld Yangians \cite{Drinfeld:1987sy} and give our comments on the construction.


\section{Quiver Yangian algebras \texorpdfstring{$\mathsf{Y}(\CQ)$}{}}\label{sec: QYA}

\subsection{Quiver Data}\label{subsec: Quiver Data}

\begin{figure}[ht!]
	\centering
	\begin{tikzpicture}[scale=0.8, every path/.style={>={Stealth[scale=0.6]}}]
		\foreach \x in {1,..., 4}
		{
			\begin{scope}[rotate = \x * 90]
				\draw[postaction={decorate},decoration={markings, 
					mark= at position 0.6 with {\arrow{>}}}] (-1, 1) to (-1,-1);
				\draw[postaction={decorate},decoration={markings, 
					mark= at position 0.6 with {\arrow{>}}}] (-1, -1) to[out=105,in=-105] (-1,1);
				\draw[postaction={decorate},decoration={markings, 
					mark= at position 0.6 with {\arrow{>}}}] (0, 0) -- (-1, 1);
				\draw[postaction={decorate},decoration={markings, 
					mark= at position 0.6 with {\arrow{>}}}] (-1, 1) -- (-2, 0);
				\draw[postaction={decorate},decoration={markings, 
					mark= at position 0.6 with {\arrow{>}}}] (-1, 1) -- (0, 2);
				\draw[postaction={decorate},decoration={markings, 
					mark= at position 0.3 with {\arrowreversed{>}}, mark= at position 0.4 with {\arrowreversed{>}}, mark= at position 0.6 with {\arrow{>}}, mark= at position 0.7 with {\arrow{>}}}] (-2, 0) to[out=65, in=-155] (0, 2);
			\end{scope}
		}
		\draw[fill=\myblue] (-1, 1) circle (0.1) (1, 1) circle (0.1) (-1,  -1) circle (0.1) (1, -1) circle (0.1) (0, 2) circle (0.1) (2, 0) circle (0.1) (0, -2) circle (0.1) (-2, 0) circle (0.1);
		\draw[fill=burgundy] (-0.1,-0.1) -- (-0.1,0.1) -- (0.1,0.1) -- (0.1,-0.1) -- cycle;
	\end{tikzpicture}
	\caption{Example of a quiver}\label{fig:Quiv_exmpl}
\end{figure}
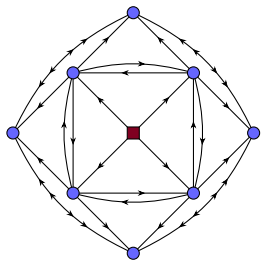

We start with a pair $(\CQ, \CW)$ of a quiver diagram $\CQ = (Q_{0}, Q_{1})$ where $Q_{0}$ is the set of vertices of the quiver $\CQ$, and $Q_{1}$ is the set of the arrows between these vertices. The function $\CW$ is a superpotential. We also introduce some useful notations:
\begin{itemize}
	\item $\{a\to b\}$ -- a set of arrows flowing from node $a$ to node $b$.
	\item $|a\to b|$ -- a number of arrows flowing from node $a$ to node $b$.
\end{itemize}
If we are to specify the head and the tail of an arrow $I$, we denote it as $I\colon a \to b$.
For a pair of nodes $a, b\in Q_{0}$ we define their chirality by:
\begin{equation}
	\chi_{ab} = |a \to b| - |b \to a|\,.
\end{equation}
The quiver is called \textit{non-chiral} when $\chi_{ab} = 0$ for any pair of nodes; otherwise the quiver is called \textit{chiral}. For some notes related to chiral quivers, see, for example, \cite{Galakhov:2021xum} or \cite{Galakhov:2021vbo}. From now on we consider only non-chiral quivers.

The quiver formalism has previously been applied to describe BPS algebras arising from D-brane systems on Calabi–Yau three-folds \cite{Li:2020rij, Galakhov:2020vyb, Galakhov_2024}. In this work, we employ essentially the same framework, adapted to the study of Dynkin quivers, following the general ideas and motivations presented, for example, in \cite{Li:2023zub, Bao:2023ece}.

We now discuss the quiver data in more detail.
\begin{enumerate}
	\item The quiver nodes can be divided into gauge nodes and framing nodes. We draw them using different shapes in the diagram: the gauge nodes are denoted as round nodes, whereas framing nodes are denoted as square nodes (see fig. \ref{fig:Quiv_exmpl}). With each node $a\in Q_{0}$ we associate a vector space $V_{a} = \mathbb{C}^{d_{a}}$ and a gauge or a flavor group $GL(d_{a}, \mathbb{C})$ with a dimension parameter $d_{a}\geqslant 0$. We also assume that there is a field $\Phi_{a} \in {\rm Hom}(V_{a}, V_{a})$ associated with the node.
	\item We work with the quivers with only one type of arrows. To each arrow $I\in \{a \to b\}$ we assign a chiral bifundamental field $q_{I}$ that is charged by $\bar{GL(d_{a})}\times GL(d_{b})$. If a node $a$ or $b$ is a framing node, $q_{I}$ is just an (anti-)fundamental field. We suppose it to be a linear map from $V_{a}$ to $V_{b}$, therefore, $q_{I:a\to b} \in {\rm Hom}(V_{a}, V_{b})$. We also assume that additional $U(1)$ flavor symmetries are associated with framing nodes and assign (equivariant) weight $h_{I} \in \mathbb{C}$ to each field. Moreover, we introduce R-charges $R_{I}$ for each arrow of the quiver. There are $|Q_{1}|$ weights, or R-charges, that are, in general, unconstrained.
	\item The superpotential $\CW$ is a holomorphic gauge invariant function of fields. We impose that the superpotential can be decomposed into a sum of monomials in $q_{I}$. Gauge invariance requires the monomials to form closed loops $\{L\}$ in the quiver $\CQ$. The superpotential is flavor invariant as well. This leads us to the loop constraints that, together with vertex constraints, define the relations between equivariant weights:
	\begin{equation}\label{Loop and Vertex Constraints}
		\begin{aligned}
			&\text{Loop constraint: } \quad\sum_{I\in L} h_{I} = 0\,, \quad \forall L\in \CW\\
			&\text{Vertex constraint: } \quad \sum_{I \in a}\sgn_{a}(I)\,h_{I} = 0\,, \quad \forall a \in Q_{0}\,,
		\end{aligned}
	\end{equation}
	where $\sgn_{a}(I)$ is equal to $+1$ if the arrow flows towards the vertex $a$, $-1$ if the arrow flows outwards from the vertex $a$, and $0$ otherwise. The reasons behind the vertex constraints we address later in the text.
	The R-charges are also constrained by the superpotential due to the fact that its R-charge is fixed:
	\begin{equation}\label{Loop Constraints R-charges}
		\sum_{I\in L} R_{I} = 2\,, \quad \forall L\in \CW\,.
	\end{equation}
\end{enumerate}

\subsection{Yangian Algebra}\label{subsec: Yangian Algebra}

Having defined the quiver data and superpotential, we now proceed to define the quadratic relations of the Yangian algebra following \cite{Li:2020rij, Li:2023zub}.

The generators of the algebra are divided into $|Q_{0}|$ triplets $(e^{(a)}_{n}, f^{(a)}_{n}, \psi^{(a)}_{n})$ where $n\in \mathbb{Z}$, $a\in Q_{0}$. They can be organized into generating functions for convenience\footnote{The expressions hold for so-called un-shifted Yangians. One could introduce the shift \cite{Galakhov:2021xum, Noshita:2021dgj, Kodera:2016faj} in the mode expression of $\psi^{(a)}(z)$. However, we do not cover this in our consideration.}:
\begin{equation}\label{General generating functions}
	e^{(a)}(z)=\sum\lm_{n=0}^{\infty}\frac{e_n^{(a)}}{z^{n+1}},\quad f^{(a)}(z)=\sum\lm_{n=0}^{\infty}\frac{f_n^{(a)}}{z^{n+1}},\quad \psi^{(a)}(z)=1 + \sum\lm_{n=0}^{\infty}\frac{\psi_n^{(a)}}{z^{n+1}}\,.
\end{equation}
Generators $e^{(a)}(z)$ and $f^{(a)}(z)$ have parity:
\begin{equation}
	|a|=(|a\to a|+1)\;{\rm mod}\;2\,,
\end{equation}
that counts the number of loops at vertex $a$. Generators $\psi^{(a)}(z)$ always have parity 0.

For each pair of vertices, we introduce the \textit{bonding factors} with the help of equivariant weights of arrows between these vertices:
\begin{equation}\label{bonding}
	\varphi_{a,b}(z)=\frac{\prod\lm_{I\in\{a\to b\}}(z+h_I)}{\prod\lm_{J\in\{b\to a\}}(z-h_J)}\,.
\end{equation}
The generating functions defined above help us to compactly write the quadratic relations of the algebra:
\begin{equation}\label{Yangian Relations}
	\begin{aligned}
		&e^{(a)}(z)e^{(b)}(w)\simeq (-1)^{|a||b|}\varphi_{a,b}(z-w)\,e^{(b)}(w)e^{(a)}(z)\,,\\
		&\psi^{(a)}(z)e^{(b)}(w)\simeq \varphi_{a,b}(z-w)e^{(b)}(w)\,\psi^{(a)}(z)\,,\\
		&f^{(a)}(z)f^{(b)}(w)\simeq (-1)^{|a||b|}\varphi_{a,b}(z-w)^{-1}f^{(b)}(w)f^{(a)}(z)\,,\\
		&\psi^{(a)}(z)f^{(b)}(w)\simeq \varphi_{a,b}(z-w)^{-1}f^{(b)}(w)\psi^{(a)}(z)\,,\\
		&\psi^{(a)}(z)\psi^{(b)}(w)=\psi^{(b)}(w)\psi^{(a)}(z)\,,\\
		&\left[e^{(a)}(z),f^{(b)}(w)\right\}\simeq-\delta_{ab}\frac{\psi^{(a)}(z)-\psi^{(b)}(w)}{z-w}\,,
	\end{aligned}
\end{equation}
where $[*,*\}$ denotes a super-commutator, and $\simeq$ denotes an equivalence of series up to $z^{n}w^{k\geq 0}$ and $z^{n\geq 0}w^{k}$.
One could unfold these relations in terms of modes. We refer the interested reader to a more detailed review \cite{Li:2020rij}.

In general, Yangian algebras have higher-order algebraic relations that are called Serre relations. However, cubic and higher-order relations require additional consideration and depend on the choice of the algebra. Some useful notes on Serre relations are given in \cite{Li:2020rij}. A conjecture in the case of a general quiver is given in \cite{Negut:2022pka}. The affine cases $\myY(\widehat{\fg\fl}_{m|n})$ are covered in \cite{Bezerra:2019kfl, Bezerra_2020}. The Serre relations for $\myY(\fs\fl_{n})$ algebras are mentioned later in the text.

\subsection{Basic Properties}\label{subsec: Basic Properties}

Let us now mention some basic properties of the algebra $\myY(\CQ)$ defined above. They are useful in understanding the structure of the relations \eqref{Yangian Relations} and address the motivation behind the vertex constraints \eqref{Loop and Vertex Constraints}. All the properties that we discuss are tightly related to automorphisms of the algebra and were analyzed, for example, in \cite{Li:2020rij, Prochazka:2015deb}.
\begin{itemize}
	\item The algebra admits $\mathbb{Z}_{2}$ symmetry:
	\begin{equation}\label{Z_2 symmetry}
		e^{(a)}(z) \leftrightarrow f^{(a)}(z), \qquad \psi^{(a)}(z) \to \psi^{(a)}(z)^{-1}\,,
	\end{equation}
	that introduces a $\mathbb{Z}_{2}$ grading of the generators.
	\item The relations also have rescaling symmetry\footnote{We subtly changed the form of this symmetry compared to \cite{Li:2020rij, Galakhov:2021xum} for our purposes. However, the concept behind it remains unchanged.}:
	\begin{equation}\label{Z symmetry}
		\begin{aligned}
			&h_{I} \to \sigma h_{I}\,, \qquad z \to \sigma z\,,\\
			&e^{(a)}(z) \to \sigma^{-1}e^{(a)}(z)\,, \quad f^{(a)}(z) \to \sigma^{-1}f^{(a)}(z)\,, \quad \psi^{(a)}(z) \to \sigma^{-1}\psi^{(a)}(z)\,,\\
			&e_{n}^{(a)} \to \sigma^{n}e_{n}^{(a)}\,, \quad f_{n}^{(a)} \to \sigma^{n}f_{n}^{(a)}\,, \quad \psi_{n}^{(a)} \to \sigma^{n}\psi_{n}^{(a)}\,.
		\end{aligned}
	\end{equation}
	The transformation rewritten in terms of modes shows that there is a natural $\mathbb{Z}$ grading. This grading is also known as level grading, or spin filtration \cite{Li:2020rij}.
	\item In fact, there is another scaling symmetry:
	\begin{equation}\label{normalization symmetry}
		\begin{aligned}
			&e^{(a)}(z) \to \alpha e^{(a)}(z)\,, \quad f^{(a)}(z) \to \alpha^{-1}f^{(a)}(z)\,, \quad \psi^{(a)}(z) \to \psi^{(a)}(z)\,,\\
			&e^{(a)}_{n} \to \alpha e_{n}^{(a)}\,, \quad f^{(a)}_{n} \to \alpha^{-1}f_{n}^{(a)}\,, \quad \psi_{n}^{(a)} \to \psi_{n}^{(a)}\,.
		\end{aligned}
	\end{equation}
	This symmetry has a simpler form and can be fixed by choosing a bilinear form with respect to which the generators are conjugate, in other words, a normalization.
	\item The last symmetry that we mention in this text has a more complicated structure. It changes the spectral parameter $z$ of the algebra and therefore is known as \textit{spectral shift}. The symmetry is parameterized by a complex number $s_{a} \in \mathbb{C}$ and takes the following form:
	\begin{equation}\label{spectral shift}
		\begin{aligned}
			&e^{(a)}(u) \to e^{(a)}(u - s_{a})\,, \qquad e^{(a)}_{j} \to \sum_{k = 0}^{j}\binom{j}{k}s_{a}^{j-k}e^{(a)}_{k}\,,\\
			&f^{(a)}(u) \to f^{(a)}(u - s_{a})\,, \qquad f^{(a)}_{j} \to \sum_{k = 0}^{j}\binom{j}{k}s_{a}^{j-k}f^{(a)}_{k}\,,\\
			&\psi^{(a)}(u) \to \psi^{(a)}(u - s_{a})\,, \qquad \psi^{(a)}_{j} \to \sum_{k = 0}^{j}\binom{j}{k}s_{a}^{j-k}\psi^{(a)}_{k}\,.
		\end{aligned}
	\end{equation}
	One could easily confirm that the relations \eqref{Yangian Relations} are invariant under this transformation in terms of generating functions. The binomial identities are useful to prove it in terms of modes \cite{Prochazka:2015deb}.
\end{itemize}

We related the equivariant parameters $h_{I}$ to the global symmetry of the theory. We, therefore, imposed the loop constraints \eqref{Loop and Vertex Constraints}. One could notice that these constraints are invariant under the following shift of equivariant parameters:
\begin{equation}\label{gauge shift}
	h_{I} \quad \to \quad h_{I} + \varepsilon_{a}\sgn_{a}(I) - \varepsilon_{b}\sgn_{b}(I)\,,
\end{equation}
where $\varepsilon_{a}$ are parameters of the transformation and $I\in \{a\to b\}$. This means, we have redundancy in the parametrization. From physical point of view, it represents the fact that some flavor symmetries are gauge \cite{Li:2020rij}. Indeed, if we mix the global symmetry with a gauge symmetry associated with a vertex $a$, we get \eqref{gauge shift}.

Now, let us modify the generators in each vertex $a\in Q_{0}$:
\begin{equation}
	\widetilde{e}^{(a)}(z) = e^{(a)}\bigl(z + \varepsilon_{a}\bigr)\,, \quad \widetilde{f}^{(a)}(z) = f^{(a)}\bigl(z + \varepsilon_{a}\bigr)\,, \quad \widetilde{\psi}^{(a)}(z) = \psi^{(a)}\bigl(z + \varepsilon_{a}\bigr)\,.
\end{equation}
One could ask if these new generators still form a quiver Yangian algebra. In order to satisfy the relations \eqref{Yangian Relations}, we should also modify the bond factors \eqref{bonding} in the following way:
\begin{equation}
	\varphi_{a,b}(u)=\frac{\prod\lm_{I\in\{a\to b\}}(u+h_I)}{\prod\lm_{J\in\{b\to a\}}(u-h_J)} \qquad \to \qquad \widetilde{\varphi}_{a, b}(u) = \frac{\prod\lm_{I\in\{a\to b\}}\bigl(u+h_I + \varepsilon_{a} - \varepsilon_{b}\bigr)}{\prod\lm_{J\in\{b\to a\}}\bigl(u-h_J + \varepsilon_{a} - \varepsilon_{b}\bigr)}\,.
\end{equation}
Using appropriate spectral shifts \eqref{spectral shift}, one could compensate \eqref{gauge shift}. This procedure only shuffles the generators by linear combinations. This means, one could regard the spectral shift \eqref{spectral shift} as gauge symmetry.

One could deal with this gauge freedom directly; however, the more common choice \cite{Li:2020rij, Galakhov:2021xum, Bao_2025} is to fix it by imposing the vertex constraint \eqref{Loop and Vertex Constraints}.

\subsection{States}\label{subsec: States}

Having defined the algebra relations, we now focus on its representations. The first step is to describe the states. Mathematically, the states are fixed points on the moduli space of quiver representations with respect to equivariant action originating from the toric geometry \cite{maulik2018quantumgroupsquantumcohomology, chuang2022hilbertschemesnonreduceddivisors, Galakhov:2020vyb}. First, we review the effective description of the fixed points that considers the states as sets of paths or posets\footnote{Partially ordered sets} on the quiver. This approach is covered, for example, in reviews \cite{Li:2023zub, Szendroi:2007nu, kirillov2016quiver}. Notably, for the toric Calabi-Yau threefolds we end up with 3-dimensional molten crystals \cite{Ooguri:2009ijd, Yamazaki:2010fz, Li:2020rij, Galakhov_2024}. This discussion applies (with new features) to toric CY fourfolds \cite{Galakhov_2024csd, franco20234dcrystalmeltingtoric, Bao_2024}. Even more general crystal models were constructed recently using the Jeffrey-Kirwan residue formulas \cite{jeffrey1994localizationnonabeliangroupactions} in \cite{Bao_2025}.

\subsubsection*{Sets of Paths}

To define the states directly in terms of quiver data, we briefly review some key notions of the path algebras of a quiver. The algebra $\mathbb{C}\CQ$ is generated by arrows in the quiver diagram, and the multiplication is defined by concatenation. A path $\CP$ is a sequence of $n$ arrows $I_{s}$:
\begin{equation}\label{path example}
	I_{1}\cdot I_{2}\cdot\ldots\cdot I_{n}\,,
\end{equation}
where each target of a previous arrow points at the source of the next arrow, which means that a path is connected.

In the presence of a superpotential, one should impose equivalence relations on the paths given by $F$-term relations:
\begin{equation}\label{general quiver F-terms}
	F_{I} = \d_{I}\CW = 0\,.
\end{equation}
For example, let us consider the superpotential:
\begin{equation}
	\CW = \Tr(I_{1}I_{3} - I_{1}I_{4}I_{2})\,.
\end{equation}
The relation $\d_{I_{1}}\CW = 0$ defines the path equivalence:
\begin{equation}
	I_{1} \simeq I_{4}I_{2}\,.
\end{equation}
This algebra is known as Jacobian algebra of a quiver:
\begin{equation}\label{Jacobian algebra}
	\CJ_{(Q, W)} = \mathbb{C}\CQ\,/\,F_{I}\,.
\end{equation}

We have already mentioned the framing nodes in section \ref{subsec: Quiver Data}. A choice of framing nodes and the arrows connecting them with the remaining quiver is called ``framing". The states and therefore representations of the quiver depend on the choice of framings. We assume that we have a single framing node in the quiver and its framing dimension is $d_{\ff} = 1$, and label the arrow connecting the framing node to a gauge node as $R$. This choice is known as a canonical framing. For more consistent discussion of various framings, see, for example, \cite{Li:2023zub, Galakhov:2022uyu, Galakhov:2021xum, Chen:2025xoe, rapcak2021branesquiversbpsalgebras}.

Framing projects the Jacobian algebra to a subalgebra that we denote by $\CJ^{\sharp}$. The paths grow from the framing field $R$:
\begin{equation}\label{Jacobian algebra framed}
	\CJ^{\sharp}_{(Q, W)} = \CJ_{(Q, W)}\cdot R\,.
\end{equation}
Analogously to the case of toric Calabi-Yau threefolds, we call each path an ``atom", and a state of the representation should be viewed as a set of these atoms, up to the homotopic equivalence defined by \eqref{general quiver F-terms}, which is called ``crystal"\footnote{We note that in our notation it is not a traditional molten crystal but rather its complement. One should keep in mind that the sets of paths in general can have complicated structure. We discuss later in the text why, in our case, we still consider the states as crystals.}. In these settings, one can think of a framing as a choice of a root atom.

There are two sets of parameters that are related to the arrows and therefore the paths on a quiver, namely, equivariant weights and R-charges. The constraints on equivariant weights can be resolved, and the parameters form a $k$-dimensional equivariant linear space:
\begin{equation}\label{general arrow weights}
	h_{I} = \sum_{j = 1}^{k}x^{j}_{I}\epsilon_{j}\,,
\end{equation}
where $x^{j}_{I}$ is the $j$-th coordinate of the $h_{I}$ in the equivariant space. For example, in the case of toric ${\bf CY}_{3}$, there are two independent equivariant parameters \cite{Li:2020rij}.

Naturally, the equivariant weights and R-changes can also be defined for the paths $\CP$ \eqref{path example}:
\begin{equation}\label{general path weights}
	\begin{aligned}
		&x_{\CP}^{j} = x^{j}_{I_{1}} + x^{j}_{I_{2}} + \ldots + x_{I_{n}}^{j}, \qquad R_{\CP} = R_{I_{1}} + R_{I_{2}} + \ldots + R_{I_{n}}\,,\\
		&h_{\CP} = h_{I_{1}} + h_{I_{2}} + \ldots + h_{I_{n}}\,.
	\end{aligned}
\end{equation}
The starting point of a path is fixed by $R$. Together with the homotopic equivalence, it means that the atoms are specified by their heads; therefore, they can be represented as points on a $k + 1$-dimensional space:
\begin{equation}
	\CP \quad \leadsto \quad (\vec{x}_{\CP}, R_{\CP})\,.
\end{equation}
We will refer to this space as the \textit{R-equivariant space} $\Xi$. For later convenience, we denote a path with the ending node $a$ as $\sqbox{$a$}$. The ending vertex is called a \textit{color} of an atom. When two atoms $\sqbox{$a$}$ and $\sqbox{$b$}$ are connected by an arrow of the quiver $I \in \{a \to b\}$, we connect these points in $\Xi$ by an arrow. These arrows are also known as the chemical bonds. Therefore the crystals are graphs in this space.

We also impose an extra condition on the crystals known as the \textit{no-overlap condition} \cite{Bao_2025}:
\begin{tcolorbox}
	Any crystal $\Lambda$ contains $\sum\limits_{a\in Q_{0}}d_{a}$ different paths for any dimension vector $\vec{d}$.
\end{tcolorbox}
The condition means that the atoms in a crystal do not overlap in the space $\Xi$. This holds for the known molten crystals for toric Calabi-Yau manifolds.

Earlier we introduced the vertex constraints \eqref{Loop and Vertex Constraints} on the weights, and in section \ref{subsec: Basic Properties} we stated that they fix gauge redundancies. However, one should be careful while imposing them or any other gauge-fixing relations. There are examples \cite{Bao_2025} where one cannot satisfy the no-overlap condition if they impose the vertex constraints. As we will see later in the paper, this holds for our cases as well. Since we can work with the gauge freedom using the spectral shift \eqref{spectral shift}, the correct procedure is to impose the vertex constraints only after the computations.

\subsubsection*{Fixed Points on Quiver Varieties}

Now, we proceed further and describe the states from a more mathematical point of view. We refer the reader to \cite{Galakhov:2020vyb, maulik2018quantumgroupsquantumcohomology, chuang2022hilbertschemesnonreduceddivisors, Rapcak:2018nsl, Rapcak:2020ueh} for more concise reviews. Instead, we limit ourselves with technical details of the construction.

With each arrow $I\in \{a \to b\}$, we have associated a matrix $q_{I}$ that acts from $V_{a}$ to $V_{b}$. The F-relations \eqref{general quiver F-terms} are translated into the matrix equations
\begin{equation}\label{general F-terms}
	F_{I} = \d_{q_{I}}\CW = 0\,,
\end{equation}
and cut off a complex algebraic variety. Factorizing it modulo gauge transformations, we end up with the moduli space of quiver representations \cite{NakajimaALE, king1994moduli, Nakajima98}:
\begin{equation}
	\mathscr{M}_{\vec{d}} = \left\{\begin{array}{c}
		\text{stable}\\
		F_{I} = 0
	\end{array}\right\}/\prod_{a \in Q_{0}}GL(d_{a}, \mathbb{C})\,.
\end{equation}
We do not address the definition of stability here. We emphasize, however, that the effective approach that was mentioned above allows us to describe only the so-called cyclic stability chamber of the moduli spaces\footnote{We refer the interested reader to \cite{Galakhov:2024foa} for some details of the construction of state models in different stability chambers.}.

Now, we narrow this space further using the equivariant group generated by the vector fields:
\begin{equation}\label{equiv vector fields}
	v_{I:a\to b}=\Tr(\Phi_b q_I-q_I\Phi_a-h_I q_I)\frac{\p}{\p q_I}\,.
\end{equation}
The states in this construction are associated with the fixed points of \eqref{equiv vector fields}. This procedure is analogous to assigning the equivariant weights to the paths. Notably, the crystals are represented as a set of matrices $q_{I}$ with a given dimension vector $\vec{d}$. The paths on the quiver \eqref{path example} are naturally identified with the fixed points as follows:
\begin{equation}
	\CP = I_{1}\cdot I_{2}\cdot\ldots\cdot I_{n} \quad \leadsto \quad q_{\CP} = q_{I_{1}}\cdot q_{I_{2}}\cdot \ldots \cdot q_{I_{n}}\,.
\end{equation}
For homotopic paths we have $q_{\CP} = q_{\CP'}$, which is ensured by the F-relations \eqref{general F-terms}.
The atoms become vectors in the spaces $V_{a}$, where $a\in Q_{0}$:
\begin{equation}\label{cyc}
	V_a={\rm Span}\left\{q_{I_n}\cdot\ldots\cdot q_{I_2}\cdot q_{I_1}\cdot R \right\}\,.
\end{equation}

One could directly solve the F-term relations, fix the gauge freedom, and find the corresponding fixed points for each possible choice of the dimension vector $\vec{d}$. However, in our case, we can construct the fixed points by using the information from the Jacobian algebra $\CJ^{\sharp}$ of the quiver. This approach is also described in \cite{Galakhov_2024}.

First, we fix the dimension vector and identify the crystal structure $\Lambda$ of the state in the $k + 1$-dimensional space. We also choose the numeration of the vectors in $V_{a}$:
\begin{equation}
	V_a=\bigoplus_{\ssqbox{$a$}\in \Lambda}\IC|\sqbox{$a$}\rangle, \quad a\in Q_0\,,
\end{equation}
where $\sqbox{$a$}$ specifies a \textit{color} of an atom.

Second, we construct the vacuum expectation values of the matrices $\Phi_{a}$:
\begin{equation}\label{gen Phi}
	\Phi_a={\rm diag}\left\{\phi_{\ssqbox{$a$}_1},\phi_{\ssqbox{$a$}_2},\ldots,\phi_{\ssqbox{$a$}_{d_a}}\right\},\quad \phi_{\ssqbox{$a$}_i}= \sum_{j = 1}^{k} x_{\ssqbox{$a$}_i}^{j} \epsilon_{j}\,.
\end{equation}
The eigenvalues of these matrices encode the weights \eqref{general path weights} of the corresponding atoms.

Finally, we construct explicitly the matrices of $q_I$ in fixed points as follows ($i=1,\ldots,d_b$, $j=1,\ldots,d_a$):
\begin{equation}\label{morphisms}
	\left(q_{I:a\to b}\right)_{ij}=\left\{\begin{array}{ll}
		1, & \mbox{ if } (\vec{x}_{\ssqbox{$a$}_j} +\vec{x}_I,R_{\ssqbox{$a$}_j}+R_I)=(\vec{x}_{\ssqbox{$b$}_i},R_{\ssqbox{$b$}_i})\,;\\
		0, &\mbox{ otherwise}\,.
	\end{array}\right.
\end{equation}

\subsection{Yangian Representations}\label{subsec: Yangian Representation}

Now we proceed to define the representation itself. Vectors in our modules are described in section \ref{subsec: States}. Given a state $\Lambda$, one should define how the generators \eqref{General generating functions} act on it. We use the model where a raising operator $e^{(a)}(z)$ adds a single atom of a color $a$ to a crystal, and a lowering operator $f^{(a)}(z)$ removes an atom from a crystal.

When we add/remove an atom, we can accidentally violate the molten crystal structure. It leads us to a natural constraint that a new set of atoms $\Lambda \pm \sqbox{$a$}$ is also a crystal. The sets $\text{Add}(\Lambda)$ and $\text{Rem}(\Lambda)$ consist of all possible atoms that can be added to (removed from) a given crystal $\Lambda$ correspondingly.

Finally, we end up with the ansatz for Yangian action on crystal states:
\begin{equation}\label{General Ansatz}
	\begin{aligned}
		&\psi^{(a)}(z)|\Lambda\rangle = {\bf \Psi}_{\Lambda}^{(a)}(z)|\Lambda\rangle\,, \\
		&e^{(a)}(z)|\Lambda\rangle = \sum_{\ssqbox{$a$}\in \text{Add}(\Lambda)}\dfrac{{\bf E}_{\Lambda, \Lambda + \ssqbox{$a$}}}{z - \phi_{\ssqbox{$a$}}}|\Lambda + \sqbox{$a$}\rangle\,,\\
		&f^{(a)}(z)|\Lambda\rangle = \sum_{\ssqbox{$a$}\in \text{Rem}(\Lambda)}\dfrac{{\bf F}_{\Lambda, \Lambda - \ssqbox{$a$}}}{z - \phi_{\ssqbox{$a$}}}|\Lambda - \sqbox{$a$}\rangle\,.
	\end{aligned}
\end{equation}
The eigenvalues of operators $\psi^{(a)}(z)$ are given by the formula:
\begin{equation}\label{general eigenvalues}
	\bPsi^{(a)}_{\Lambda}(z)=\prod\lm_{I\in\{a\to a\}}\left(-\frac{1}{h_I}\right)\times\frac{\prod\lm_{K\in\{a\to\ff\}}(-z-h_K)}{\prod\lm_{J\in\{\ff\to a\}}(z-h_I)}\times \prod\lm_{b\in Q_0}\prod\lm_{\ssqbox{$b$}\in \Lambda}\varphi_{a,b}(z-\phi_{\ssqbox{$b$}})\,.
\end{equation}
In fact, the poles of functions $\bPsi^{(a)}_{\Lambda}(z)$ coincide with the set $\text{Add}(\Lambda) \cup \text{Rem}(\Lambda)$; we refer the interested reader to \cite{Li:2023zub} for a detailed review. If a pole is equal to a weight of an atom in the crystal $\Lambda$, it belongs to $\text{Rem}(\Lambda)$ and to $\text{Add}(\Lambda)$ otherwise.

One could derive that for \eqref{General Ansatz} to be a representation of the Yangian \eqref{Yangian Relations}, matrix elements have to satisfy so-called \textit{hysteresis} relations \cite{Li:2020rij, Galakhov_2024}:
\begin{equation}\label{general hysteresis}
	\begin{aligned}
		&{\bf E}_{\Lambda+\ssqbox{$a$},\Lambda+\ssqbox{$a$}+\ssqbox{$b$}}{\bf F}_{\Lambda+\ssqbox{$a$}+\ssqbox{$b$},\Lambda+\ssqbox{$b$}}=(-1)^{|a||b|}{\bf F}_{\Lambda+\ssqbox{$a$},\Lambda}{\bf E}_{\Lambda,\Lambda+\ssqbox{$b$}}\,,\\
		&\frac{{\bf E}_{\Lambda,\Lambda+\ssqbox{$a$}}{\bf E}_{\Lambda+\ssqbox{$a$},\Lambda+\ssqbox{$a$}+\ssqbox{$b$}}}{{\bf E}_{\Lambda,\Lambda+\ssqbox{$b$}}{\bf E}_{\Lambda+\ssqbox{$b$},\Lambda+\ssqbox{$a$}+\ssqbox{$b$}}}\varphi_{a,b}(\phi_{\ssqbox{$a$}}-\phi_{\ssqbox{$b$}})=(-1)^{|a||b|}\,,\\
		&\frac{{\bf F}_{\Lambda+\ssqbox{$a$}+\ssqbox{$b$},\Lambda+\ssqbox{$a$}}{\bf F}_{\Lambda+\ssqbox{$a$},\Lambda}}{{\bf F}_{\Lambda+\ssqbox{$a$}+\ssqbox{$b$},\Lambda+\ssqbox{$b$}}{\bf F}_{\Lambda+\ssqbox{$b$},\Lambda}}\varphi_{a,b}(\phi_{\ssqbox{$a$}}-\phi_{\ssqbox{$b$}})=(-1)^{|a||b|}\,,\\
		&{\bf E}_{\Lambda,\Lambda+\ssqbox{$a$}}{\bf F}_{\Lambda+\ssqbox{$a$},\Lambda}=\mathop{\rm res}\lm_{z=\phi_{\ssqbox{$a$}}}\Psi^{(a)}_{\Lambda}(z)\,.
	\end{aligned}
\end{equation}

\subsection{Equivariant Matrix Coefficients}\label{subsec: Equivariant Matrix Coefficients}

For \eqref{General Ansatz} to be a representation, we need an explicit formula for the matrix elements ${\bf E}_{\Lambda, \Lambda + \ssqbox{$a$}}$ and ${\bf F}_{\Lambda + \ssqbox{$a$}, \Lambda}$. The ansatz allows a slight ambiguity in the definition of the coefficients that was briefly addressed in \cite{Galakhov_2024}. This is caused by the algebra symmetries that we listed in section \ref{subsec: Basic Properties}.

One of the options is a \textit{square-root} representation \cite{Galakhov_2024}. We refer the interested reader to \cite{Li:2020rij} for a more detailed review. In this case the coefficients are chosen as follows:
\begin{equation}\label{sq_rt}
	{\bf E}^{(root)}_{\Lambda,\Lambda+\ssqbox{$a$}}={\bf F}^{(root)}_{\Lambda+\ssqbox{$a$},\Lambda}\sim\sqrt{\mathop{\rm res}\lm_{z=\phi_{\ssqbox{$a$}}}\Psi^{(a)}_{\Lambda}(z)}\,.
\end{equation}
The states in this representation are normalized to unity.
However, due to a lack of a canonical way to choose the branch of the square root function, we do not use \eqref{sq_rt} to calculate the matrix elements.

Instead, we exploit the geometric approach that relies on an equivariant integration over quiver representation moduli spaces \cite{Nakajima_lect,Rapcak:2018nsl,Rapcak:2020ueh}. For more physical motivations, see, for example \cite{Galakhov:2020vyb, Galakhov_2024}. Here, we briefly review this construction.

One of the key ideas is that equivariant localization \cite{guillemin2013supersymmetry, Witten:1982im, Cordes:1994fc} allows us to extract all the information about the system properties in the vicinities of the fixed points. The equivariant fields \eqref{equiv vector fields} take diagonal form in a neighborhood of a fixed point and grade the corresponding tangent space $\CN$:
\begin{equation}
	v=\sum\lm_i w_i \, z_i\frac{\p}{\p z_i}, \qquad \CN = \bigoplus_{i}\mathbb{C}|w_{i}\rangle\,,
\end{equation}
where $z_{i}$ are the local coordinates on $\CN$. Note that some of the weights $w_{i}$ can be zero.

In these settings, natural characteristics of a fixed point are the corresponding Euler classes. One should be careful, however, in directly applying the algorithm. The varieties $\mathscr{M}$ that we work with can happen to be singular. It leads to jumps in the dimensions of the tangent spaces. For more details on the regularization in these cases, we refer to \cite{Galakhov:2020vyb, Galakhov_2024}. After that, we end up with the modified version of the Euler classes:
\begin{equation}\label{reg_Eul}
	{\rm Eul}\;\CN\,=\,(-1)^{\left\lfloor\frac{1}{2}\#\{i:\, w_i=0\}\right\rfloor}\prod\lm_{i:\,w_i\neq 0}w_i\,.
\end{equation}

For each fixed point we define the corresponding Euler class:
\begin{equation}
	{\rm Eul}_{\Lambda}={\rm Eul}\; \mathsf{T}_{\Lambda}\mathscr{M}\,.
\end{equation}

To proceed further and define the matrix elements ${\bf E}_{\Lambda,\Lambda'}$, ${\bf F}_{\Lambda',\Lambda}$ where $\Lambda'=\Lambda+\Box$ we consider two corresponding quiver representations $q_{I},~ q'_{I}$. The representations are homomorphic if there exists a set of maps $\tau_a$, $a\in Q_0$ making the following diagrams commutative:
\begin{equation}\label{incidence}
	\begin{array}{c}
		\begin{tikzpicture}
			\node(A) at (0,0) {$V_a$};
			\node(B) at (2,0) {$V_b$};
			\node(C) at (0,-1.2) {$V_a'$};
			\node(D) at (2,-1.2) {$V_b'$};
			\path (A) edge[->] node[above] {$\scriptstyle q_{I:a\to b}$} (B) (C) edge[->] node[above] {$\scriptstyle q_{I:a\to b}'$} (D) (A) edge[->] node[left] {$\scriptstyle \tau_a$} (C) (B) edge[->] node[right] {$\scriptstyle \tau_b$} (D);
		\end{tikzpicture}
	\end{array},\quad q_{I:a\to b}'\cdot \tau_a=\tau_b\cdot q_{I:a\to b},\quad \forall I\in Q_1\,.
\end{equation}

We call an \emph{incidence locus} $\CI$ a surface in the Cartesian product of two representations, $q_I$ and $q_I'$, where this homomorphism exists. The tangent space to the incidence locus $\mathsf{T}_{\Lambda,\Lambda'}\CI\subset \mathsf{T}_{\Lambda}\mathscr{M}\oplus \mathsf{T}_{\Lambda'}\mathscr{M}$ is naturally an equivariantly weighted space, which means that we are able to define the corresponding Euler class:
\begin{equation}
	{\rm Eul}_{\Lambda,\Lambda'}={\rm Eul}\;\mathsf{T}_{\Lambda,\Lambda'}\CI\,.
\end{equation}

The matrix coefficients are constructed as Fourier-Mukai transforms \cite{huybrechts2006fourier} from $\mathsf{T}_{\Lambda}\mathscr{M}$ to $\mathsf{T}_{\Lambda'}\mathscr{M}$ and inverse with a kernel given by the structure sheaf of $\CI$. The construction boils down to the canonical pullback-pushforward integration \cite{Rapcak:2020ueh, Rapcak:2018nsl, Nakajima_lect}. All the integrals are equivariant, and after applying the canonical Berline-Vergne-Atiyah-Bott localization formula \cite{zbMATH03814600, Atiyah:1984px, Alekseev:2000fe} we acquire the result as a ratio of Euler classes:
\begin{equation}\label{equiv_matrix_coeffs}
	\begin{aligned}
		&{\bf E}^{(equiv)}_{\Lambda, \Lambda + \ssqbox{$a$}} = \dfrac{\Eul_{\Lambda}}{\Eul_{\Lambda, \Lambda + \ssqbox{$a$}}}\,, \\
		&{\bf F}^{(equiv)}_{\Lambda + \ssqbox{$a$}, \Lambda} = \dfrac{\Eul_{\Lambda + \ssqbox{$a$}}}{\Eul_{\Lambda, \Lambda + \ssqbox{$a$}}}\,.
	\end{aligned}
\end{equation}
Equivariant integration gives us the normalization of the states \cite{Galakhov:2020vyb}:
\begin{equation}
	\langle \Lambda |\Lambda \rangle = {\rm Eul}_{\Lambda}\,.
\end{equation}
One could expect that changing the norm to 1 returns the root representation. Although the statement is not proven in the general case, there are a few examples \cite{Galakhov_2024} where the representations coincide. This paper provides even more examples to the point. The coefficients are related using the formula \eqref{normalization symmetry}:
\begin{equation}\label{Rescaling Equiv-to-Root}
	{\bf E}^{(root)}_{\Lambda, \Lambda + \Box} = {\bf E}^{(equiv)}_{\Lambda, \Lambda + \Box}\sqrt{\dfrac{{\rm Eul}_{\Lambda + \Box}}{{\rm Eul}_{\Lambda}}}, \quad {\bf F}^{(root)}_{\Lambda + \Box, \Lambda} = {\bf F}^{(equiv)}_{\Lambda + \Box, \Lambda}\sqrt{\dfrac{{\rm Eul}_{\Lambda}}{{\rm Eul}_{\Lambda + \Box}}}\,.
\end{equation}
For the rest of the paper, if not specified, the notation ${\bf E}_{\Lambda, \Lambda + \ssqbox{$a$}}$ stands for ${\bf E}^{(equiv)}_{\Lambda, \Lambda + \ssqbox{$a$}}$.

\subsection{Dynkin Diagrams}\label{subsec: Quiver Dynkin Diagrams}

In the preceding sections, we reviewed the general machinery for the construction of Yangian algebras and their representations using the quiver data. We now specialize to the class of quivers derived from Dynkin diagrams. This bridge between quivers and Dynkin diagrams is provided by the McKay correspondence \cite{McKay:1981aap, Bridgeland:2001xf, mclean2023mckay, cirafici2013curvecountinginstantonsmckay}. The original correspondence \cite{McKay:1981aap} connects the smooth resolutions of orbifold singularities $\mathbb{C}^{2}/\Gamma$, for finite groups $\Gamma \subset SL(2, \mathbb{C})$, and Dynkin graphs of ADE type\footnote{More generally, the McKay correspondence refers to the study of smooth resolutions of orbifold singularities $\mathbb{C}^{n}/\Gamma$, for finite groups $\Gamma \subset SL(n, \mathbb{C})$ \cite{mclean2023mckay}; $\mathbb{C}^{n}$ can be replaced by a complex $n$-dimensional variety \cite{Bridgeland:2001xf}.}. The generalization to the non-simply laced cases is given by the Slodowy correspondence \cite{Slodowy_1980, stekolshchik2005notescoxetertransformationsmckay}, see also \cite{Bershadsky_1996, Cecotti:2012gh} for the geometric construction.

Here we present an algorithm to acquire the quivers using a given Dynkin diagram. We refer to \cite{cirafici2013curvecountinginstantonsmckay, Li:2023zub, Bao:2023ece} for more details.

A McKay quiver $\CQ_{\Gamma}$, for a finite group $\Gamma$, is defined as follows. Irreducible representations $\rho_{a}$ of $\Gamma$ label the set of vertices $Q_{0}$ of the quiver\footnote{Including the trivial $\rho_{0}$ for the affine cases}. One can determine the number of arrows $n_{ab}$ between arbitrary nodes $a, b$ using the decomposition:
\begin{equation}\label{McKay arrows}
	\CR \otimes \rho_{a} = \bigoplus_{b \in Q_{0}}n_{ab}\rho_{b}\,,
\end{equation}
where $\CR\subset \mathbb{C}^{2}$ is the fundamental representation of $\Gamma$ induced by inclusion $\Gamma \subset SL(2, \mathbb{C})$. This results in ``doubling" of the original Dynkin graph. We also extend this quiver further by adding a self-loop $C_{a}$ to each vertex. These quivers are sometimes referred to as triple quivers \cite{Bao:2023ece, Li:2023zub, Ginzburg:2006fu}. We depict this construction, including non-simply laced cases\footnote{We avoid the superalgebras in this consideration. The question of constructing a quiver and assigning a superpotential for a given McKay quiver is more subtle in this case. We refer the interested reader to \cite{Bao:2023ece, Galakhov:2021vbo} for more detailed reviews.}, and assign the corresponding superpotential in \eqref{McKay quivers}:
\begin{equation}\label{McKay quivers}
	\begin{tblr}{c|c|c}
		\mbox{Dynkin graph} & \mbox{Quiver} & \mbox{Superpotential } \delta W_{ab} \\ \hline
		\begin{tikzpicture}[baseline={(0, -0.1)}]
			\draw[postaction=decorate, decoration={markings}, thick] (0,0) to (1.9,0);
			\draw[fill=\myblue] (0, 0) circle (0.1) (1.9, 0) circle (0.1);
			\node[below] at (0, -0.1) {$\scriptstyle a$};
			\node[below] at (1.9, -0.1) {$\scriptstyle b$};
		\end{tikzpicture} & \begin{tikzpicture}[baseline={(0,-0.1)}]
			\draw[postaction=decorate, decoration={markings, mark= at position 0.8 with {\arrow{stealth}}}] (0,0) to[out=60,in=0] (0,0.6) to[out=180,in=120] (0,0);
			\draw[postaction=decorate, decoration={markings, mark= at position 0.8 with {\arrow{stealth}}}] (2,0) to[out=60,in=0] (2,0.6) to[out=180,in=120] (2,0);
			\draw[postaction=decorate, decoration={markings, mark= at position 0.7 with {\arrow{stealth}}}] (0,0) to[out=10,in=170] node[pos=0.5,above] {$\scriptstyle X_{ab}$} (2,0);
			\draw[postaction=decorate, decoration={markings, mark= at position 0.7 with {\arrow{stealth}}}] (2,0) to[out=-170,in=-10] node[pos=0.5,below] {$\scriptstyle X_{ba}$} (0,0);
			\draw[fill=\myblue] (0, 0) circle (0.1) (2, 0) circle (0.1);
			\node[below] at (0, -0.1) {$\scriptstyle a$};
			\node[below] at (2, -0.1) {$\scriptstyle b$};
			\node[above] at (0,0.6) {$\scriptstyle C_a$};
			\node[above] at (2,0.6) {$\scriptstyle C_b$};
		\end{tikzpicture} & \Tr(X_{ab}X_{ba}C_{a} - X_{ba}X_{ab}C_{b})\\
		\begin{tikzpicture}[scale=1, baseline={(0, -0.1)}]
			\draw[thick] (0,0.08) -- (1.9,0.08);
			\draw[thick] (0,-0.08) -- (1.9,-0.08);
			\draw[thick] (0.9,0.2) to[out=-80, in=170] (1.15,0) to[out=-170, in=80] (0.9,-0.2);
			\draw[fill=\myblue] (0,0) circle (0.1);
			\draw[fill=\myblue] (1.9,0) circle (0.1);
			\node[below] at (0, -0.1) {$\scriptstyle a$};
			\node[below] at (1.9, -0.1) {$\scriptstyle b$};
		\end{tikzpicture} & \begin{tikzpicture}[baseline={(0,-0.1)}]
			\draw[postaction=decorate, decoration={markings, mark= at position 0.8 with {\arrow{stealth}}}] (0,0) to[out=60,in=0] (0,0.6) to[out=180,in=120] (0,0);
			\draw[postaction=decorate, decoration={markings, mark= at position 0.8 with {\arrow{stealth}}}] (2,0) to[out=60,in=0] (2,0.6) to[out=180,in=120] (2,0);
			\draw[postaction=decorate, decoration={markings, mark= at position 0.7 with {\arrow{stealth}}}] (0,0) to[out=10,in=170] node[pos=0.5,above] {$\scriptstyle X_{ab}$} (2,0);
			\draw[postaction=decorate, decoration={markings, mark= at position 0.7 with {\arrow{stealth}}}] (2,0) to[out=-170,in=-10] node[pos=0.5,below] {$\scriptstyle X_{ba}$} (0,0);
			\draw[fill=\myblue] (0, 0) circle (0.1) (2, 0) circle (0.1);
			\node[below] at (0, -0.1) {$\scriptstyle a$};
			\node[below] at (2, -0.1) {$\scriptstyle b$};
			\node[above] at (0,0.6) {$\scriptstyle C_a$};
			\node[above] at (2,0.6) {$\scriptstyle C_b$};
		\end{tikzpicture} & \Tr\Bigl(X_{ab}X_{ba}C_{a}^{|\CA_{ab}|} - X_{ba}X_{ab}C_{b}^{|\CA_{ba}|}\Bigr)
	\end{tblr}
\end{equation}
In \eqref{McKay quivers} $\CA_{ab}$ denotes the corresponding Cartan matrix. The full unframed superpotential is recovered by summing up all the possible pairs in the quiver: $\CW = \sum\limits_{(a, b)}\delta W_{ab}$. By construction:
\begin{equation}
	|Q_{0}| = \rk(\fg)\,,
\end{equation}
where $\fg$ is the corresponding Lie algebra. Also, we would like to emphasize that although the quiver diagram in \eqref{McKay quivers} seems to be the same for non-simply laced cases, the resulting theory is different. Indeed, the superpotentials have different forms, which change the weight assignment to the arrows of the quiver.

Although the procedure can be applied in theory for an arbitrary Dynkin diagram, we restrict ourselves to the A-type diagrams. When supplemented with suitable framings and superpotentials, the construction defines the Yangian algebras $\myY(\fs\fl_{n})$ and their representations. An arbitrary Dynkin diagram A$_{n - 1}$ is given by the picture:
\begin{equation}
	\begin{tikzpicture}
		\draw[postaction=decorate, decoration={markings}, thick] (0,0) to (1,0);
		\draw[postaction=decorate, decoration={markings}, thick] (1,0) to (2,0);
		\draw[postaction=decorate, decoration={markings}, thick] (4,0) to (5,0);
		\draw[postaction=decorate, decoration={markings}, thick] (2,0) to (2.5,0);
		\draw[postaction=decorate, decoration={markings}, thick] (3.5,0) to (4,0);
		\draw[fill=black] (3, 0) circle (0.02) (2.8, 0) circle (0.02) (3.2, 0) circle (0.02);
		\draw[fill=\myblue] (0, 0) circle (0.1) (1, 0) circle (0.1) (2, 0) circle (0.1) (4, 0) circle (0.1) (5,0) circle (0.1);
		\node[below] at (0, -0.1) {$\scriptstyle 1$};
		\node[below] at (1, -0.1) {$\scriptstyle 2$};
		\node[below] at (5, -0.1) {$\scriptstyle n - 1$};
	\end{tikzpicture}
\end{equation}

Applying the algorithm, we end up with the following family of quivers to describe rectangular representations of $\myY(\fs\fl_{n})$ algebras:
\begin{equation}\label{family of quivers and superpotentials for sl(n)}
	\CQ_{n,p,\lambda}=\left\{\begin{array}{c}
		\begin{tikzpicture}
			\draw[postaction=decorate, decoration={markings, mark= at position 0.7 with {\arrow{stealth}}}] (0,0) to[out=20,in=160] node[pos=0.5,above] {$\scriptstyle A_1$} (1.5,0);
			\draw[postaction=decorate, decoration={markings, mark= at position 0.7 with {\arrow{stealth}}}] (1.5,0) to[out=200,in=340] node[pos=0.5,below] {$\scriptstyle B_1$} (0,0);
			\draw[postaction=decorate, decoration={markings, mark= at position 0.8 with {\arrow{stealth}}}] (0,0) to[out=60,in=0] (0,0.6) to[out=180,in=120] (0,0);
			\node[above] at (0,0.6) {$\scriptstyle C_1$};
			\begin{scope}[shift={(1.5,0)}]
				\draw[postaction=decorate, decoration={markings, mark= at position 0.8 with {\arrow{stealth}}}] (0,0) to[out=60,in=0] (0,0.6) to[out=180,in=120] (0,0);
				\node[above] at (0,0.6) {$\scriptstyle C_2$};
			\end{scope}
			\begin{scope}
				\draw[postaction=decorate, decoration={markings, mark= at position 0.7 with {\arrow{stealth}}}, dashed] (0,-1.2) to[out=100,in=260] (0,0);
				\draw[postaction=decorate, decoration={markings, mark= at position 0.7 with {\arrow{stealth}}}, dashed] (0,0) to[out=280,in=80] (0,-1.2);
				\begin{scope}[shift={(0,-1.2)}, dashed]
					\draw[fill=gray!50!white] (-0.08,-0.08) -- (-0.08,0.08) -- (0.08,0.08) -- (0.08,-0.08) -- cycle;
				\end{scope}
			\end{scope}
			\begin{scope}[shift={(1.5, 0)}]
				\draw[postaction=decorate, decoration={markings, mark= at position 0.7 with {\arrow{stealth}}}, dashed] (0,-1.2) to[out=100,in=260] (0,0);
				\draw[postaction=decorate, decoration={markings, mark= at position 0.7 with {\arrow{stealth}}}, dashed] (0,0) to[out=280,in=80] (0,-1.2);
				\begin{scope}[shift={(0,-1.2)}, dashed]
					\draw[fill=gray!50!white] (-0.08,-0.08) -- (-0.08,0.08) -- (0.08,0.08) -- (0.08,-0.08) -- cycle;
				\end{scope}
			\end{scope}
			\draw[fill=\myblue] (0,0) circle (0.08) (1.5,0) circle (0.08);
			\begin{scope}[shift = {(4,0)}]
				\begin{scope}[shift={(-1.5, 0)}]
					\draw[postaction=decorate, decoration={markings, mark= at position 0.7 with {\arrow{stealth}}}, dashed] (0,-1.2) to[out=100,in=260] (0,0);
					\draw[postaction=decorate, decoration={markings, mark= at position 0.7 with {\arrow{stealth}}}, dashed] (0,0) to[out=280,in=80] (0,-1.2);
					\begin{scope}[shift={(0,-1.2)}, dashed]
						\draw[fill=gray!50!white] (-0.08,-0.08) -- (-0.08,0.08) -- (0.08,0.08) -- (0.08,-0.08) -- cycle;
					\end{scope}
				\end{scope}
				\begin{scope}[shift={(1.5, 0)}]
					\draw[postaction=decorate, decoration={markings, mark= at position 0.7 with {\arrow{stealth}}}, dashed] (0,-1.2) to[out=100,in=260] (0,0);
					\draw[postaction=decorate, decoration={markings, mark= at position 0.7 with {\arrow{stealth}}}, dashed] (0,0) to[out=280,in=80] (0,-1.2);
					\begin{scope}[shift={(0,-1.2)}, dashed]
						\draw[fill=gray!50!white] (-0.08,-0.08) -- (-0.08,0.08) -- (0.08,0.08) -- (0.08,-0.08) -- cycle;
					\end{scope}
				\end{scope}
				\draw[postaction=decorate, decoration={markings, mark= at position 0.7 with {\arrow{stealth}}}] (0,-1.2) to[out=100,in=260] node[pos=0.3, left] {$\scriptstyle R_{p}$} (0,0);
				\draw[postaction=decorate, decoration={markings, mark= at position 0.7 with {\arrow{stealth}}}] (0,0) to[out=280,in=80] node[pos=0.7, right] {$\scriptstyle S_{p}$} (0,-1.2);
				\draw[postaction=decorate, decoration={markings, mark= at position 0.7 with {\arrow{stealth}}}] (0,0) to[out=20,in=160] node[pos=0.5,above] {$\scriptstyle A_p$} (1.5,0);
				\draw[postaction=decorate, decoration={markings, mark= at position 0.7 with {\arrow{stealth}}}] (1.5,0) to[out=200,in=340] node[pos=0.5,below] {$\scriptstyle B_p$} (0,0);
				\begin{scope}[shift={(-1.5,0)}]
					\draw[postaction=decorate, decoration={markings, mark= at position 0.7 with {\arrow{stealth}}}] (0,0) to[out=20,in=160] node[pos=0.5,above] {$\scriptstyle A_{p-1}$} (1.5,0);
					\draw[postaction=decorate, decoration={markings, mark= at position 0.7 with {\arrow{stealth}}}] (1.5,0) to[out=200,in=340] node[pos=0.5,below] {$\scriptstyle B_{p-1}$} (0,0);
				\end{scope}
				\draw[postaction=decorate, decoration={markings, mark= at position 0.8 with {\arrow{stealth}}}] (0,0) to[out=60,in=0] (0,0.6) to[out=180,in=120] (0,0);
				\node[above] at (0,0.6) {$\scriptstyle C_p$};
				\begin{scope}[shift={(1.5,0)}]
					\draw[postaction=decorate, decoration={markings, mark= at position 0.8 with {\arrow{stealth}}}] (0,0) to[out=60,in=0] (0,0.6) to[out=180,in=120] (0,0);
					\node[above] at (0,0.6) {$\scriptstyle C_{p+1}$};
				\end{scope}
				\begin{scope}[shift={(-1.5,0)}]
					\draw[postaction=decorate, decoration={markings, mark= at position 0.8 with {\arrow{stealth}}}] (0,0) to[out=60,in=0] (0,0.6) to[out=180,in=120] (0,0);
					\node[above] at (0,0.6) {$\scriptstyle C_{p-1}$};
				\end{scope}
				\draw[fill=\myblue] (-1.5,0) circle (0.08) (0,0) circle (0.08) (1.5,0) circle (0.08);
				\begin{scope}[shift={(0,-1.2)}]
					\draw[fill=burgundy] (-0.08,-0.08) -- (-0.08,0.08) -- (0.08,0.08) -- (0.08,-0.08) -- cycle;
				\end{scope}
			\end{scope}
			\begin{scope}[shift = {(8,0)}]
				\begin{scope}[shift={(-1.5, 0)}]
					\draw[postaction=decorate, decoration={markings, mark= at position 0.7 with {\arrow{stealth}}}, dashed] (0,-1.2) to[out=100,in=260] (0,0);
					\draw[postaction=decorate, decoration={markings, mark= at position 0.7 with {\arrow{stealth}}}, dashed] (0,0) to[out=280,in=80] (0,-1.2);
					\begin{scope}[shift={(0,-1.2)}, dashed]
						\draw[fill=gray!50!white] (-0.08,-0.08) -- (-0.08,0.08) -- (0.08,0.08) -- (0.08,-0.08) -- cycle;
					\end{scope}
				\end{scope}
				\begin{scope}
					\draw[postaction=decorate, decoration={markings, mark= at position 0.7 with {\arrow{stealth}}}, dashed] (0,-1.2) to[out=100,in=260] (0,0);
					\draw[postaction=decorate, decoration={markings, mark= at position 0.7 with {\arrow{stealth}}}, dashed] (0,0) to[out=280,in=80] (0,-1.2);
					\begin{scope}[shift={(0,-1.2)}, dashed]
						\draw[fill=gray!50!white] (-0.08,-0.08) -- (-0.08,0.08) -- (0.08,0.08) -- (0.08,-0.08) -- cycle;
					\end{scope}
				\end{scope}
				\begin{scope}[shift={(-1.5,0)}]
					\draw[postaction=decorate, decoration={markings, mark= at position 0.7 with {\arrow{stealth}}}] (0,0) to[out=20,in=160] node[pos=0.5,above] {$\scriptstyle A_{n-2}$} (1.5,0);
					\draw[postaction=decorate, decoration={markings, mark= at position 0.7 with {\arrow{stealth}}}] (1.5,0) to[out=200,in=340] node[pos=0.5,below] {$\scriptstyle B_{n-2}$} (0,0);
				\end{scope}
				\draw[postaction=decorate, decoration={markings, mark= at position 0.8 with {\arrow{stealth}}}] (0,0) to[out=60,in=0] (0,0.6) to[out=180,in=120] (0,0);
				\node[above] at (0,0.6) {$\scriptstyle C_{n-1}$};
				\begin{scope}[shift={(-1.5,0)}]
					\draw[postaction=decorate, decoration={markings, mark= at position 0.8 with {\arrow{stealth}}}] (0,0) to[out=60,in=0] (0,0.6) to[out=180,in=120] (0,0);
					\node[above] at (0,0.6) {$\scriptstyle C_{n-2}$};
				\end{scope}
				\draw[fill=\myblue] (-1.5,0) circle (0.08) (0,0) circle (0.08);
			\end{scope}
			\draw[fill=black] (1.75,0) circle (0.03) (2,0) circle (0.03) (2.25,0) circle (0.03) (5.75,0) circle (0.03) (6,0) circle (0.03) (6.25,0) circle (0.03);
		\end{tikzpicture}\\
		W=\Tr\left[A_1C_1B_1+\sum\lm_{a=2}^{n-2}\left(A_aC_aB_a-B_{a-1}C_aA_{a-1}\right)-B_{n-2}C_{n-1}A_{n-2}+{\color{burgundy}C_p^{\lambda}R_{p}S_{p}}\right]
	\end{array}\right\}\,.
\end{equation}

Note that for all $a \in \CQ_{n,p,\lambda}$, we have $|a| = 0$.
The fields $R_{a}, S_{a}$ are indexed depending on what node they are attached to. The unframed quiver is independent of fields $R_{a}$, $S_{a}$ and of parameters $p$, $\lambda$. We emphasize that technically, we should consider the perturbations of all the framing nodes in equivariant calculations:
\begin{equation}
	R_{a} = \bar{R}_{a} + \delta R_{a}\,, \quad S_{a} = \bar{S}_{a} + \delta S_{a}\,,
\end{equation}
including the fields where $a \ne p$. This means that we added more than 1 framing with the dimensions $\ff_{a} = 1$. This could affect the construction of the states described in section \ref{subsec: States}. However, we will see later that $\bar{R}_{a} = 0,~ \bar{S}_{a} = 0$ for $a \ne p$. The latter ensures that there are no additional branches starting at $R_{a},~ a\ne p$.

The $F$-terms cut off toric varieties in this case. In fact, this is no longer true for any other Dynkin diagram \cite{Bao:2023ece, Bao_2025, Chen:2025xoe} and complicates the analysis in the general case. 

Our framing choice, depending explicitly on $p$ and $\lambda$, cuts out a specific representation of $\fs\fl_n$, classified by the following Young diagram:
\begin{equation}\label{gen rectangular reps}
	\Upsilon_{p,\lambda}=\begin{array}{c}
		\begin{tikzpicture}[scale=0.4]
			\foreach \i in {0,1,2,3}
			{
				\draw[thick] (0,\i) -- (4,\i);
			}
			\foreach \i in {0,1,2,3,4}
			{
				\draw[thick] (\i,0) -- (\i,3);
			}
			\draw (0,0) to[out=270, in=180] (0.25,-0.25) -- (1.75,-0.25) to[out=0,in=90] (2,-0.5) to[out=90,in=180] (2.25,-0.25) -- (3.75,-0.25) to[out=0,in=270] (4,0);
			\draw (4,0) to[out=0,in=270] (4.25,0.25) -- (4.25,1.25) to[out=90,in=180] (4.5,1.5) to[out=180,in=270] (4.25,1.75) -- (4.25,2.75) to[out=90,in=0] (4,3);
			\node[below] at (2,-0.5) {$\scriptstyle \lambda$ \scriptsize columns};
			\node[right] at (4.5,1.5) {$\scriptstyle p$ \scriptsize rows};
		\end{tikzpicture}
	\end{array}\,,
\end{equation}
The quivers have an additional $\mathbb{Z}_{2}$-symmetry:
\begin{equation}\label{sl(n) quiver symmetry Z_2}
	\mathcal{Q}_{n, p, \lambda} \to \mathcal{Q}_{n, n - p, \lambda}\,,
\end{equation}
that corresponds to complex conjugation of the representations $\Upsilon_{p,\lambda} \to \Upsilon_{n - p,\lambda}$. Due to this fact, we can omit the description of the representations with $p > \ceil{\frac{n-1}{2}}$. That leaves us with the representations:
\begin{equation}\label{gen rectangular reps after cutoff}
	\Upsilon_{p,\lambda}\,, \quad \text{where } \quad p \leqslant \ceil{\frac{n-1}{2}}\,.
\end{equation}

\subsection*{Counting the number of equivariant parameters}

Since we have specified  the quivers, we need to check the independence of their equivariant parameters. The dependencies are linear and defined by the system of equations \eqref{Loop and Vertex Constraints}. Let us denote the matrix corresponding to this system as $\CM$. Basic linear algebra tells us that the number of independent variables can be calculated as follows:
\begin{equation}\label{Number of equivariant parameters}
	\#(\epsilon) = |Q_{1}| - \rk(\CM)\,.
\end{equation}
For the given quivers we have:
\begin{equation}
	|Q_{1}| = (n - 1) + 2(n-2) = 3n - 5\,.
\end{equation}
The number of all constraints is:
\begin{equation}
	\#(\text{Constraints}) = \underset{\text{Vertex}}{\underbrace{n - 1}} + \underset{\text{Loop}}{\underbrace{2(n - 2)}} = 3n - 5\,.
\end{equation}
The next step is to determine the rank. One can notice that the loop constraints and the vertex constraints are separable in our case \eqref{family of quivers and superpotentials for sl(n)}. This is because the loop constraints involve the weights $h_{C_{a}}$ whereas the vertex constraints do not.
We can check that the vertex constraints are not independent:
\begin{equation}
	\sum_{a = 1}^{n - 1} \text{Vertex constraint}_{a} = 0\,.
\end{equation}
For example, in the case of $n = 4$:
\begin{equation}
	\begin{aligned}
		&\text{Vertex constraint}_{1} = h_{B_{1}} - h_{A_{1}}\,,\\
		&\text{Vertex constraint}_{2} = h_{A_{1}} - h_{B_{1}} + h_{B_{2}} - h_{A_{2}}\,,\\
		&\text{Vertex constraint}_{3} = h_{A_{2}} - h_{B_{2}}\,,\\
		&\sum_{a = 1}^{3} \text{Vertex constraint}_{a} = h_{B_{1}} - h_{A_{1}} + h_{A_{1}} - h_{B_{1}} + h_{B_{2}} - h_{A_{2}} + h_{A_{2}} - h_{B_{2}} = 0\,.
	\end{aligned}
\end{equation}
The matrix of the loop constraints has a full rank. It means that the rank of the whole matrix $\CM$ reads:
\begin{equation}
	\rk(\CM) = n-2 + 2(n - 2) = 3n - 6\,.
\end{equation}
Finally, we end with the result:
\begin{equation}
	\fbox{$\#(\epsilon) = 3n - 5 - (3n -6) = 1$}\,.
\end{equation}

\subsection*{Analysis of the result}

Let us examine the weight assignment in the presence of the constraints \eqref{Loop and Vertex Constraints} more closely. First, one can solve the loop constraints as follows:
\begin{equation}
	\begin{aligned}
		&h_{C_{1}} = \dots = h_{C_{a}} = \dots = h_{C_{n-1}} = \epsilon\,,\\
		&h_{A_{a}} = h_{a} - \dfrac{\epsilon}{2}\,, \quad
		h_{B_{a}} = - h_{a} - \frac{\epsilon}{2}\,,
	\end{aligned}
\end{equation}
where $\epsilon$ and $h_{a}$ are $n - 1$ parameters. For example, for the fields $C_{a}$ and $C_{a + 1}$ we have:
\begin{equation}
	\begin{aligned}
		&h_{A_{a}} + h_{C_{a + 1}} + h_{B_{a}} = 0,\,\\
		&h_{A_{a}} + h_{C_{a}} + h_{B_{a}} = 0\,.
	\end{aligned}
\end{equation}
Subtracting the one from the other we get $h_{C_{a}} = h_{C_{a + 1}}$.

We now examine the vertex constraints. There are two types of these constraints. The constraints that correspond to the vertices from $2$ to $n - 2$, and the other constraints that correspond to the vertices labeled by $1$ and $n - 1$. There are $n - 3$ constraints of the first type, but all of them have the same structure. For a vertex $a$, the constraint takes the following form:
\begin{equation}
	\begin{aligned}
		&h_{A_{a - 1}} - h_{B_{a - 1}} + h_{B_{a}} - h_{A_{a}} = 0\,\\
		&h_{a - 1} - \frac{\epsilon}{2} + h_{a - 1} + \frac{\epsilon}{2} - h_{a} - \frac{\epsilon}{2} - h_{a} + \frac{\epsilon}{2} = 0\,\\
		&h_{a - 1} = h_{a}\,.
	\end{aligned}
\end{equation}
This means that all the parameters $h_{a}$ are equal to each other. We denote them as $h$.

The ``edge" conditions resolve as:
\begin{equation}\label{h = 0}
	h = -h \qquad \rightarrow \qquad h = 0\,,
\end{equation}
which, as expected, mirrors the result of our calculations above that we have a single independent parameter $\epsilon$.

In other words, we demonstrated that our R-equivariant space $\Xi$ happens to be two-dimensional. However, as we will see, for example, in section \ref{subsec: Y(sl(4))} atoms in the crystals overlap if we impose the vertex constraints. In section \ref{subsec: States} we mentioned that, in this case, one should impose the loop constraints, construct a representation, and only after that impose the vertex constraints, ensuring that the theory is still gauge invariant.

We will use the suggested procedure, although with a slight modification. In our work we do not need to keep track of all the $n - 1$ parameters $(h_{a}, \epsilon)$ that are left after the loop constraints. Therefore, we just loose the vertex constraints and impose only the constraints of the first type. That leaves us with only two equivariant parameters $h, \epsilon$. Together with R-charge, the space $\Xi$ becomes three-dimensional. In the final formulas we substitute $h = 0$.

\section{Yangian Algebras $\myY(\fs\fl_{n})$}\label{sec: Y(sl(n))}

In this section, we consider a few examples and applications of the construction introduced in section \ref{sec: QYA} in detail.

\subsection{$\myY(\fs\fl_{3})$ example}\label{subsec: Y(sl(3))}

\subsubsection{The Algebra}

In this section we discuss one of the simplest examples, the $\mathsf{Y}(\fs\fl_{3})$ algebra, and its rectangular representations that are related to $(\lambda, 0)$ or $(0, \lambda)$ representations\footnote{We encode the representations by the Dynkin labels of their highest weights for simplicity. For example, the notation $(\lambda, 0)$ is equivalent to $\Upsilon_{1, \lambda}$ for the algebra $\myY(\fs\fl_{3})$. However, when we consider the representations of the algebra $\myY(\fs\fl_{4})$, $\Upsilon_{1, \lambda}$ corresponds to $(\lambda, 0, 0)$. Therefore, the notation $\Upsilon_{p, \lambda}$ depends on the algebra indirectly.} of $\fs\fl_{3}$ algebra. These representations are complex conjugates of each other; see \eqref{gen rectangular reps after cutoff}. Therefore, it is sufficient to focus on $(\lambda, 0)$ representations. The quiver depicted in fig. \ref{sl(3) quiver} can be used to describe these representations.

\begin{figure}[h]
	\begin{center}
		\begin{tikzpicture}
			\draw[postaction=decorate, decoration={markings, mark= at position 0.7 with {\arrow{stealth}}}] (0,0) to[out=20,in=160] node[pos=0.5,above] {$\scriptstyle A$} (1.5,0);
			\draw[postaction=decorate, decoration={markings, mark= at position 0.7 with {\arrow{stealth}}}] (1.5,0) to[out=200,in=340] node[pos=0.5,below] {$\scriptstyle B$} (0,0);
			\draw[postaction=decorate, decoration={markings, mark= at position 0.8 with {\arrow{stealth}}}] (0,0) to[out=60,in=0] (0,0.6) to[out=180,in=120] (0,0);
			\node[above] at (0,0.6) {$\scriptstyle C_1$};
			\begin{scope}[shift={(1.5,0)}]
				\draw[postaction=decorate, decoration={markings, mark= at position 0.8 with {\arrow{stealth}}}] (0,0) to[out=60,in=0] (0,0.6) to[out=180,in=120] (0,0);
				\node[above] at (0,0.6) {$\scriptstyle C_2$};
			\end{scope}
			\begin{scope}[shift = {(0,0)}]
				\draw[postaction=decorate, decoration={markings, mark= at position 0.7 with {\arrow{stealth}}}] (0,-1.2) to[out=100,in=260] node[pos=0.3, left] {$\scriptstyle R_{1}$} (0,0);
				\draw[postaction=decorate, decoration={markings, mark= at position 0.7 with {\arrow{stealth}}}] (0,0) to[out=280,in=80] node[pos=0.7, right] {$\scriptstyle S_{1}$} (0,-1.2);
				\begin{scope}[shift={(0,-1.2)}]
					\draw[fill=burgundy] (-0.08,-0.08) -- (-0.08,0.08) -- (0.08,0.08) -- (0.08,-0.08) -- cycle;
				\end{scope}
			\end{scope}
			\begin{scope}[shift={(1.5, 0)}]
				\draw[postaction=decorate, decoration={markings, mark= at position 0.7 with {\arrow{stealth}}}, dashed] (0,-1.2) to[out=100,in=260] (0,0);
				\draw[postaction=decorate, decoration={markings, mark= at position 0.7 with {\arrow{stealth}}}, dashed] (0,0) to[out=280,in=80] (0,-1.2);
				\begin{scope}[shift={(0,-1.2)}, dashed]
					\draw[fill=gray!30!white] (-0.08,-0.08) -- (-0.08,0.08) -- (0.08,0.08) -- (0.08,-0.08) -- cycle;
				\end{scope}
			\end{scope}
			\draw[fill=\myblue] (0,0) circle (0.08);
			\draw[fill=\mygray] (1.5,0) circle (0.08);
			\node[left] at (0, 0) {$\scriptstyle n_{1}$};
			\node[right] at (1.5, 0) {$\scriptstyle n_{2}$};
		\end{tikzpicture}
	\end{center}
	\caption{The quiver that describes $(\lambda, 0)$ reps of $\mathsf{Y}(\fs\fl_{3})$}\label{sl(3) quiver}
\end{figure}
The dimensions of the vector spaces that correspond to the nodes are $n_{1}$ and $n_{2}$.
We also add the following superpotential:
\begin{equation}\label{sl(3) superpotential}
	\mathcal{W} = \Tr(C_{1}BA - C_{2}AB + C_{1}^{\lambda}R_{1}S_{1} + R_{2}S_{2})\,.
\end{equation}
The weights and R-charges of fields read:
\begin{equation}\label{sl(3) weights}
	\begin{tblr}{c|c|c|c|c|c|c}
		\mbox{Fields} & C_{1} & C_{2} & A & B & R_{1} & S_{1}\\
		\hline
		\mbox{Weights} & \epsilon & \epsilon & -\frac{\epsilon}{2} + h & -\frac{\epsilon}{2} - h & 0 & -\lambda \epsilon \\
		\mbox{R-charges} & 0 & 0 & 1 & 1 & 0 & 2
	\end{tblr}\,.
\end{equation}

Yangian algebra $\mathsf{Y}(\fs\fl_{3})$ contains two families of generators $(e^{(a)}_{k}, f^{(a)}_{k}, \psi^{(a)}_{k})$, $a\in \{1, 2\}$, which can be assembled in generating functions:
\begin{equation}\label{sl(3) generating functions}
	e^{(a)}(z) = \sum_{n = 0}^{\infty}\dfrac{e^{(a)}_{n}}{z^{n + 1}}, \quad f^{(a)}(z) = \sum_{n = 0}^{\infty}\dfrac{f^{(a)}_{n}}{z^{n + 1}}, \quad \psi^{(a)}(z) = 1 + \sum_{n = 0}^{\infty}\dfrac{\psi^{(a)}_{n}}{z^{n + 1}}\,.
\end{equation}
Using quiver \ref{sl(3) quiver} we construct the bonding factors using the formula \eqref{bonding}:
\begin{equation}\label{sl(3) bond factors}
	\begin{aligned}
		&\varphi_{1, 1}(z) = \dfrac{z + \epsilon}{z - \epsilon}\,, \quad
		\varphi_{1, 2}(z) = \dfrac{z - \frac{\epsilon}{2} + h}{z + \frac{\epsilon}{2} + h}\,,\\
		&\varphi_{2, 1}(z) = \dfrac{z - \frac{\epsilon}{2} - h}{z + \frac{\epsilon}{2} - h}\,, \quad 
		\varphi_{2, 2}(z) = \dfrac{z + \epsilon}{z - \epsilon}\,,
	\end{aligned}
\end{equation}
The generating functions satisfy the equations \eqref{Yangian Relations}:
\begin{equation}\label{sl(3) relations}
	\begin{aligned}
		&e^{(a)}(z)e^{(b)}(w) \simeq \varphi_{a, b}(z- w)e^{(b)}(w)e^{(a)}(z)\,,\\
		&\psi^{(a)}(z)e^{(b)}(w) \simeq \varphi_{a, b}(z - w)e^{(b)}(w)\psi^{(a)}(z)\,,\\
		&f^{(a)}(z)f^{(b)}(w) \simeq \varphi_{a, b}(z - w)^{-1} f^{(b)}(w)f^{(a)}(z)\,,\\
		&\psi^{(a)}(z)f^{(b)}(w) \simeq \varphi_{a, b}(z - w)^{-1} f^{(b)}(w)\psi^{(a)}(z)\,,\\
		&\psi^{(a)}(z)\psi^{(b)}(w) = \psi^{(b)}(w)\psi^{(a)}(z)\,,\\
		&\Bigl[e^{(a)}(z), f^{(b)}(w)\Bigr] \simeq -\delta_{ab}\dfrac{\psi^{(a)}(z) - \psi^{(b)}(w)}{z - w}\,,
	\end{aligned}
\end{equation}
together with the Serre relations \cite{Guay:2018wel}:
\begin{equation}\label{sl(3) Serre relations}
	\begin{aligned}
		&\sum_{\sigma \in \fS_{2}}\Bigl[e^{(a)}\bigl(u_{\sigma(1)}\bigr)\,, \Bigl[e^{(a)}\bigl(u_{\sigma(2)}\bigr), e^{(b)}(v)\Bigr]\Bigr] = 0\,,\\
		&\sum_{\sigma \in \fS_{2}}\Bigl[f^{(a)}\bigl(u_{\sigma(1)}\bigr)\,, \Bigl[f^{(a)}\bigl(u_{\sigma(2)}\bigr), f^{(b)}(v)\Bigr]\Bigr] = 0\,,
	\end{aligned}
\end{equation}
where $a \ne b$, and $\fS_{2}$ is the symmetric group on two letters. These relations can be translated into the mode relations. From the relations with $\varphi_{1, 1}$ we get:
\begin{equation}\label{sl(3) phi_11 relations}
	\begin{aligned}
		&[e^{(1)}_{n + 1}, e^{(1)}_{m}] - [e^{(1)}_{n}, e^{(1)}_{m + 1}] = \epsilon \{e^{(1)}_{n}, e^{(1)}_{m}\}\,,\\
		&[\psi^{(1)}_{n + 1}, e^{(1)}_{m}] - [\psi^{(1)}_{n}, e^{(1)}_{m + 1}] = \epsilon \{\psi^{(1)}_{n}, e^{(1)}_{m}\}\,,\\
		&[f^{(1)}_{n + 1}, f^{(1)}_{m}] - [f^{(1)}_{n}, f^{(1)}_{m + 1}] = -\epsilon \{f^{(1)}_{n}, f^{(1)}_{m}\}\,,\\
		&[\psi^{(1)}_{n + 1}, f^{(1)}_{m}] - [\psi^{(1)}_{n}, f^{(1)}_{m + 1}] = -\epsilon \{\psi^{(1)}_{n}, f^{(1)}_{m}\}\,;
	\end{aligned}
\end{equation}
From $\varphi_{1, 2}$ we get:
\begin{equation}\label{sl(3) phi_12 relations}
	\begin{aligned}
		&[e^{(1)}_{n + 1}, e^{(2)}_{m}] - [e^{(1)}_{n}, e^{(2)}_{m + 1}] = -\frac{\epsilon}{2} \{e^{(1)}_{n}, e^{(2)}_{m}\} - h[e^{(1)}_{n}, e^{(2)}_{m}]\,,\\
		&[\psi^{(1)}_{n + 1}, e^{(2)}_{m}] - [\psi^{(1)}_{n}, e^{(2)}_{m + 1}] = -\frac{\epsilon}{2} \{\psi^{(1)}_{n}, e^{(2)}_{m}\} - h[\psi^{(1)}_{n}, e^{(2)}_{m}]\,,\\
		&[f^{(1)}_{n + 1}, f^{(2)}_{m}] - [f^{(1)}_{n}, f^{(2)}_{m + 1}] = \frac{\epsilon}{2} \{f^{(1)}_{n}, f^{(2)}_{m}\} + h[f^{(1)}_{n}, f^{(2)}_{m}]\,,\\
		&[\psi^{(1)}_{n + 1}, f^{(2)}_{m}] - [\psi^{(1)}_{n}, f^{(2)}_{m + 1}] = \frac{\epsilon}{2} \{\psi^{(1)}_{n}, f^{(2)}_{m}\} + h[\psi^{(1)}_{n}, f^{(2)}_{m}]\,;
	\end{aligned}
\end{equation}
And from the relations with $\varphi_{2, 1},~ \varphi_{2, 2}$ we get the mode relations similar to \eqref{sl(3) phi_12 relations} and \eqref{sl(3) phi_11 relations} correspondingly. This fact is even more transparent when we set $h = 0$, because after that $\varphi_{1, 1}(z) = \varphi_{2, 2}(z)$ and $\varphi_{1, 2}(z) = \varphi_{2, 1}(z)$.

The latter two lines of \eqref{sl(3) relations} are translated into:
\begin{equation}
	\begin{aligned}
		&[\psi_{n}^{(a)}, \psi_{m}^{(b)}] = 0\,,\\
		&[e^{(a)}_{n}, f^{(b)}_{m}] = \delta_{ab}\psi_{n + m}\,,
	\end{aligned}
\end{equation}
In the discussion above we skipped so-called ``boundary conditions" that are also derived from \eqref{sl(3) relations}:
\begin{equation}\label{sl(3) boundary conditions}
	\begin{aligned}
		&[\psi_{0}^{(1)}, e_{m}^{(1)}] = 2 e_{m}^{(1)},\quad &[\psi_{0}^{(1)}, f_{m}^{(1)}] = -2 f_{m}^{(1)}\,, \\
		&[\psi^{(2)}_{0}, e_{m}^{(1)}] = -e_{m}^{(1)},\quad &[\psi_{0}^{(2)}, f_{m}^{(1)}] = f_{m}^{(1)}\,, \\
		&[\psi_{0}^{(1)}, e^{(2)}_{m}] = -e_{m}^{(2)},\quad &[\psi_{0}^{(1)}, f^{(2)}_{m}] = f_{m}^{(2)}\,,\\
		&[\psi_{0}^{(2)}, e_{m}^{(2)}] = 2 e_{m}^{(2)},\quad &[\psi_{0}^{(2)}, f_{m}^{(2)}] = -2 f_{m}^{(2)}\,, \\
	\end{aligned}
\end{equation}
Finally, for $a \ne b$ the Serre relations \eqref{sl(3) Serre relations} take the form \cite{Drinfeld:1987sy}:
\begin{equation}\label{sl(3) Serre relations in modes}
	\begin{aligned}
		&\sum_{\sigma \in \fS_{2}}\big[e_{n_{\sigma(1)}}^{(a)}\big[e_{n_{\sigma(2)}}^{(a)}, e^{(b)}_{m}\big]\big] = 0\,,\\
		&\sum_{\sigma \in \fS_{2}}\big[f_{n_{\sigma(1)}}^{(a)}\big[f_{n_{\sigma(2)}}^{(a)}, f^{(b)}_{m}\big]\big] = 0\,.
	\end{aligned}
\end{equation}
In this section we provided the calculations with the parameter $h$. However, we will not focus on it later in the text and will automatically assume $h = 0$ in any calculations except the crystal structures.

\subsubsection{The States}

From the superpotential \eqref{sl(3) superpotential}, we can derive the corresponding F-term equations:
\begin{equation}\label{sl(3) F-term relations}
	\begin{aligned}
		&\d_{A}\CW = C_{1}B - BC_{2} = 0\,,\quad 
		&\d_{B}\CW = AC_{1} - C_{2}A = 0\,, \\
		&\d_{C_{1}}\mathcal{W} = BA = 0\,, 
		&\d_{C_{2}}\mathcal{W} = AB = 0\,, \\
		&\d_{S_{1}}\mathcal{W} = C_{1}^{\lambda}R_{1} = 0\,,
		&\d_{S_{2}}\CW = R_{2} = 0\,.
	\end{aligned}
\end{equation}
Fixed points on this variety are defined by the following equations:
\begin{equation}\label{sl(3) fixed points}
	\begin{aligned}
		&[\Phi_{1}, C_{1}] = \epsilon~C_{1}\,, &[\Phi_{2}, C_{2}] = \epsilon~C_{2}\,,\\
		&\Phi_{2}A - A\Phi_{1} = \Bigl(-\frac{\epsilon}{2} + h\Bigr)~A \,, &\Phi_{1}B - B\Phi_{2} = \Bigl(-\frac{\epsilon}{2} - h\Bigr)~B \,, \\
		&\Phi_{1} R_{1} = 0\,, & S_{1}\Phi_{1} = -\lambda\epsilon~S_{1}\,.
	\end{aligned}
\end{equation}
It is helpful to depict fixed points in terms of sets of paths on the quiver \ref{sl(3) quiver}. Using the equivalence relations given by the F-terms, one can prove that:
\begin{equation}
	C_{1}^{a}B = BC_{2}^{a}\,, \quad AC_{1}^{b} = C_{2}^{b}A\,.
\end{equation}
It can be shown by exploiting the first two relations from \eqref{sl(3) F-term relations}. For example:
\begin{equation}
	C_{1}^{2}B = C_{1}(C_{1}B) \overset{F}{=} (C_{1}B)C_{2} \overset{F}{=} BC_{2}C_{2} = BC_{2}^{2}\,.
\end{equation}
Suppose that the field $B$ appears in a path. Then:
\begin{equation}
	BC_{2}^{a}AC_{1}^{b}R_{1} = BC_{2}^{a}C_{2}^{b}AR_{1} = BC_{2}^{a + b}AR_{1} = C_{1}^{a + b}BAR_{1} = C_{1}^{a + b}(BA)R_{1} = 0\,.
\end{equation}
Therefore all the paths have the following form:
\begin{equation}
	C_{1}^{a}R_{1}\,, \qquad AC_{1}^{a}R_{1} = C_{2}^{a}AR_{1}\,.
\end{equation}
Note that the equivariant weights and the R-charges of the paths are:
\begin{equation}
	\begin{aligned}
		&h(C_{1}^{a}R_{1}) = a\epsilon\,, \qquad h(AC_{1}^{a}R_{1}) = -\frac{\epsilon}{2} + h +  a\epsilon\,,\\
		&R(C_{1}^{a}R_{1}) = 0\,, \qquad R(AC_{1}^{a}R_{1}) = 1\,.
	\end{aligned}
\end{equation}
Having these two building blocks and denoting them as points, we depict the random fixed point of the quiver in the space of weights and R-charges $\Xi$:
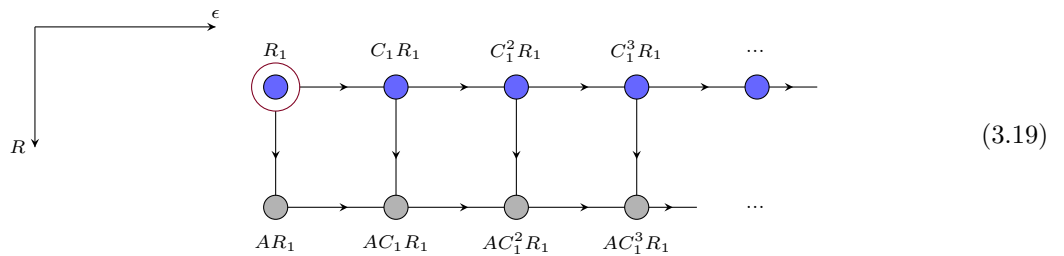
\begin{equation}\label{sl(3) random fixed point}
	\begin{array}{c}
		\begin{tikzpicture}[rotate = 0, scale=0.8]
			\draw[postaction={decorate},
			decoration={markings, mark= at position 0.6 with {\arrow{stealth}}}] (0, 0) -- (2,0);
			\draw[postaction={decorate},
			decoration={markings, mark= at position 0.6 with {\arrow{stealth}}}] (2, 0) -- (4,0);
			\draw[postaction={decorate},
			decoration={markings, mark= at position 0.6 with {\arrow{stealth}}}] (4, 0) -- (6,0);
			\draw[postaction={decorate},
			decoration={markings, mark= at position 0.6 with {\arrow{stealth}}}] (6, 0) -- (8,0);
			\draw[postaction={decorate},
			decoration={markings, mark= at position 0.6 with {\arrow{stealth}}}] (8, 0) -- (9,0);
			\draw[postaction={decorate},
			decoration={markings, mark= at position 0.6 with {\arrow{stealth}}}] (0, -2) -- (2,-2);
			\draw[postaction={decorate},
			decoration={markings, mark= at position 0.6 with {\arrow{stealth}}}] (2, -2) -- (4,-2);
			\draw[postaction={decorate},
			decoration={markings, mark= at position 0.6 with {\arrow{stealth}}}] (4, -2) -- (6,-2);
			\draw[postaction={decorate},
			decoration={markings, mark= at position 0.6 with {\arrow{stealth}}}] (6, -2) -- (7,-2);
			\draw[postaction={decorate},
			decoration={markings, mark= at position 0.6 with {\arrow{stealth}}}] (0, 0) -- (0,-2);
			\draw[postaction={decorate},
			decoration={markings, mark= at position 0.6 with {\arrow{stealth}}}] (2, 0) -- (2,-2);
			\draw[postaction={decorate},
			decoration={markings, mark= at position 0.6 with {\arrow{stealth}}}] (4, 0) -- (4,-2);
			\draw[postaction={decorate},
			decoration={markings, mark= at position 0.6 with {\arrow{stealth}}}] (6, 0) -- (6,-2);
			\draw[fill=white, draw=burgundy] circle (0.4);
			\draw[fill=\myblue] (0, 0) circle (0.2) (2, 0) circle (0.2) (4, 0) circle (0.2) (6, 0) circle (0.2) (8, 0) circle (0.2);
			\draw (0, 0.6) node {$\scriptstyle R_{1}$};
			\draw (2, 0.6) node {$\scriptstyle C_{1}R_{1}$};
			\draw (4, 0.6) node {$\scriptstyle C_{1}^{2}R_{1}$};
			\draw (6, 0.6) node {$\scriptstyle C_{1}^{3}R_{1}$};
			\draw (8, 0.6) node {$\scriptstyle \dots$};
			\begin{scope}[shift = {(0, 0)}]
				\draw[fill=\mygray] (0, -2) circle (0.2) (2, -2) circle (0.2) (4, -2) circle (0.2) (6, -2) circle (0.2);
				\draw (0, -2.6) node {$\scriptstyle AR_{1}$};
				\draw (2, -2.6) node {$\scriptstyle AC_{1}R_{1}$};
				\draw (4, -2.6) node {$\scriptstyle AC_{1}^{2}R_{1}$};
				\draw (6, -2.6) node {$\scriptstyle AC_{1}^{3}R_{1}$};
				\draw (8, -2.0) node {$\scriptstyle \dots$};
			\end{scope}
			\draw[postaction={decorate},
			decoration={markings, mark= at position 1 with {\arrow{stealth}}}] (-4, 1) to (-1, 1);
			\draw[postaction={decorate},
			decoration={markings, mark= at position 1 with {\arrow{stealth}}}] (-4, 1) to (-4, -1);
			\node[above] at (-1, 1) {$\scriptstyle \epsilon$};
			\node[left] at (-4, -1) {$\scriptstyle R$};
		\end{tikzpicture}
	\end{array}
\end{equation}
The figure requires the fulfillment of the F-terms \eqref{sl(3) F-term relations}. Otherwise the point $AC_{1}R_{1}$ would be split into two different points: $AC_{1}R_{1}$ and $C_{2}AR_{1}$, which would mean that we could not form a square. It also applies to any rectangle that we can see on the figure.
The F-relation $C_{1}^{\lambda}R_{1} = 0$ serves as a cut-off.

Note also that the picture \eqref{sl(3) random fixed point} shows that it is necessary that $n_{1} \geqslant n_{2}$. A state grows from the field $R_{1}$. We now assume that there is a path where:
\begin{equation}
	C_{2}^{m}AR_{1} \ne 0, \quad \text{and } \quad C_{1}^{m}R_{1} = 0\,.
\end{equation}
Then we can apply the equivalence relation above and derive:
\begin{equation}
	C_{2}^{m}AR_{1} = AC_{1}^{m}R_{1} = 0\,,
\end{equation}
which was not true due to our suggestion.

\subsubsection*{State Counting}

The next step is to check the dimension of the representation $(\lambda, 0)$. The dimension of the finite-dimensional representations of Yangian algebras $\fs\fl_{n}$ coincide with the dimensions of highest-weight representations of $\fs\fl_{n}$, as was shown in \cite{Drinfeld:1987sy}. Therefore, the dimension of our representations on the fixed points should match with:
\begin{equation}
	\dim_{\fs\fl_{3}}|(\lambda, 0)| = \frac{1}{2}(\lambda + 1)(\lambda + 2)\,.
\end{equation}

To count the number of fixed points in $(\lambda, 0)$ representations, we reimagine the states. First, we introduce an empty box (a planar rectangular graph in this case):
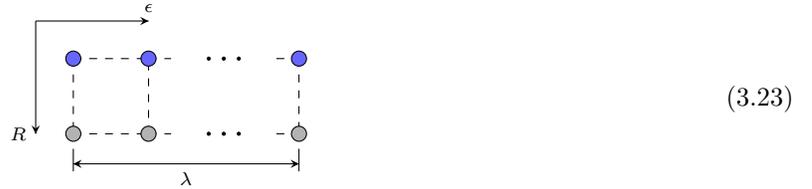
\begin{equation}\label{sl(3) empty box}
	\begin{array}{c}
		\begin{tikzpicture}[rotate = 0, scale=1]
			\draw[dashed] (0, 0) -- (1, 0) (0, 0) -- (0, -1) (0, -1) -- (1, -1) (1, -1) -- (1, 0);
			\draw[dashed] (1, 0) -- (1.3, 0) (1, -1) -- (1.3, -1) (2.7, 0) -- (3, 0) (2.7, -1) -- (3, -1) (3, 0) -- (3, -1);
			\draw[fill=\myblue] (0, 0) circle (0.1) (1, 0) circle (0.1) (3, 0) circle (0.1);
			\draw[fill=\mygray] (0, -1) circle (0.1) (1, -1) circle (0.1) (3, -1) circle (0.1);
			\draw[fill=black] (2, 0) circle (0.02) (1.8, 0) circle (0.02) (2.2, 0) circle (0.02) (2, -1) circle (0.02) (1.8, -1) circle (0.02) (2.2, -1) circle (0.02);
			\draw (0, -1.2) -- (0, -1.5) (3, -1.2) -- (3, -1.5);
			\draw[postaction={decorate},
			decoration={markings, mark= at position 0 with {\arrowreversed{stealth}}, mark= at position 1 with {\arrow{stealth}}}] (0, -1.4) -- (3, -1.4);
			\node at (1.5, -1.6) {$\scriptstyle \lambda$};
			\draw[postaction={decorate},
			decoration={markings, mark= at position 1 with {\arrow{stealth}}}] (-0.5, 0.5) to (1, 0.5);
			\draw[postaction={decorate},
			decoration={markings, mark= at position 1 with {\arrow{stealth}}}] (-0.5, 0.5) to (-0.5, -1);
			\node[above] at (1, 0.5) {$\scriptstyle \epsilon$};
			\node[left] at (-0.5, -1) {$\scriptstyle R$};
		\end{tikzpicture}
	\end{array}
\end{equation}
The graph is embedded in the R-equivariant space $\Xi$ and represents the \textit{allowed positions} of atoms in the crystals. This is also our vacuum state where $n_{1} = 0,~n_{2} = 0$. When we increase the dimensions $n_{1}$ and $n_{2}$, ``put atoms into the box", we obtain the remaining states. The atoms, the filled spaces, we will denote as bigger nodes in contrast to the empty nodes.

The states can be uniquely parameterized by $n_{1}$ and $n_{2}$ since in each node there are only vectors of one type:
\begin{equation}
	\begin{aligned}
		&\text{First node: } \quad &C_{1}^{a}R_{1}\,,\\
		&\text{Second node: } \quad &AC_{1}^{a}R_{1}\,.
	\end{aligned}
\end{equation}

It leads us to the picture:
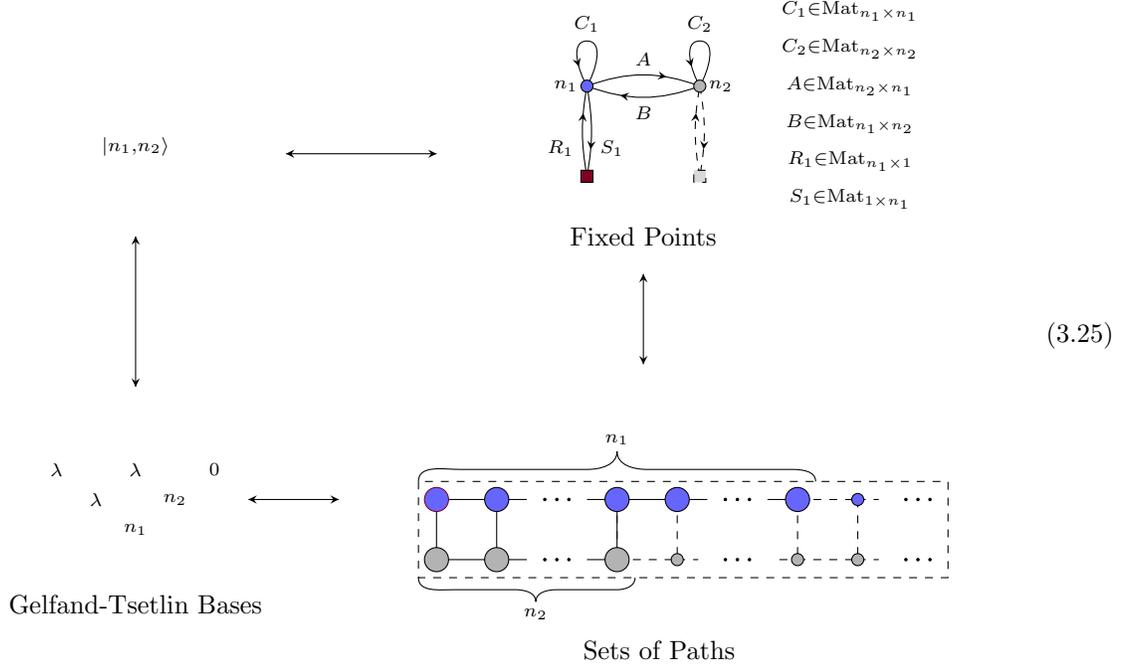
\begin{equation}\label{sl(3) states picture}
	\begin{array}{c}
		\begin{tikzpicture}
			\begin{scope}[shift={(2, 0.5)}]
				\draw[postaction=decorate, decoration={markings, mark= at position 0.7 with {\arrow{stealth}}}] (0,0) to[out=20,in=160] node[pos=0.5,above] {$\scriptstyle A$} (1.5,0);
				\draw[postaction=decorate, decoration={markings, mark= at position 0.7 with {\arrow{stealth}}}] (1.5,0) to[out=200,in=340] node[pos=0.5,below] {$\scriptstyle B$} (0,0);
				\draw[postaction=decorate, decoration={markings, mark= at position 0.8 with {\arrow{stealth}}}] (0,0) to[out=60,in=0] (0,0.6) to[out=180,in=120] (0,0);
				\node[above] at (0,0.6) {$\scriptstyle C_1$};
				\begin{scope}[shift={(1.5,0)}]
					\draw[postaction=decorate, decoration={markings, mark= at position 0.8 with {\arrow{stealth}}}] (0,0) to[out=60,in=0] (0,0.6) to[out=180,in=120] (0,0);
					\node[above] at (0,0.6) {$\scriptstyle C_2$};
				\end{scope}
				\begin{scope}[shift = {(0,0)}]
					\draw[postaction=decorate, decoration={markings, mark= at position 0.7 with {\arrow{stealth}}}] (0,-1.2) to[out=100,in=260] node[pos=0.3, left] {$\scriptstyle R_{1}$} (0,0);
					\draw[postaction=decorate, decoration={markings, mark= at position 0.7 with {\arrow{stealth}}}] (0,0) to[out=280,in=80] node[pos=0.7, right] {$\scriptstyle S_{1}$} (0,-1.2);
					\begin{scope}[shift={(0,-1.2)}]
						\draw[fill=burgundy] (-0.08,-0.08) -- (-0.08,0.08) -- (0.08,0.08) -- (0.08,-0.08) -- cycle;
					\end{scope}
				\end{scope}
				\begin{scope}[shift={(1.5, 0)}]
					\draw[postaction=decorate, decoration={markings, mark= at position 0.7 with {\arrow{stealth}}}, dashed] (0,-1.2) to[out=100,in=260] (0,0);
					\draw[postaction=decorate, decoration={markings, mark= at position 0.7 with {\arrow{stealth}}}, dashed] (0,0) to[out=280,in=80] (0,-1.2);
					\begin{scope}[shift={(0,-1.2)}, dashed]
						\draw[fill=gray!30!white] (-0.08,-0.08) -- (-0.08,0.08) -- (0.08,0.08) -- (0.08,-0.08) -- cycle;
					\end{scope}
				\end{scope}
				\draw[fill=\myblue] (0,0) circle (0.08);
				\draw[fill=\mygray] (1.5,0) circle (0.08);
				\node[left] at (0, 0) {$\scriptstyle n_{1}$};
				\node[right] at (1.5, 0) {$\scriptstyle n_{2}$};
				\begin{scope}[shift={(0, 1)}]
					\node at (3.5, 0) {$\scriptstyle C_{1} \in \text{Mat}_{n_{1}\times n_{1}}$};
					\node at (3.5, -0.5) {$\scriptstyle C_{2} \in \text{Mat}_{n_2\times n_{2}}$};
					\node at (3.5, -1) {$\scriptstyle A \in \text{Mat}_{n_{2}\times n_{1}}$};
					\node at (3.5, -1.5) {$\scriptstyle B \in \text{Mat}_{n_{1}\times n_{2}}$};
					\node at (3.5, -2) {$\scriptstyle R_{1} \in \text{Mat}_{n_{1}\times 1}$};
					\node at (3.5, -2.5) {$\scriptstyle S_{1} \in \text{Mat}_{1 \times n_{1}}$};
				\end{scope}
				\node at (0.75, -2) {Fixed Points};
			\end{scope}
			\begin{scope}[shift={(-4, -0.3)}]
				\node at (0, 0) {$\scriptstyle |n_{1}, n_{2}\rangle$};
			\end{scope}
			\begin{scope}[shift={(0, -5)}, scale=0.8]
				\draw (0, 0) -- (1, 0) (1, 0) -- (1.5, 0) (2.5, 0) -- (3, 0) (3, 0) -- (4, 0) (4, 0) -- (4.5, 0) (5.5, 0) -- (6, 0);
				\draw (0, -1) -- (1, -1) (1, -1) -- (1.5, -1) (2.5, -1) -- (3, -1);
				\draw (0, 0) -- (0, -1) (1, 0) -- (1, -1) (3, 0) -- (3, -1);
				\draw[dashed] (3, -1) -- (4, -1) (4, -1) -- (4.35, -1) (5.7, -1) -- (6, -1) (6, -1) -- (7, -1) (7, -1) -- (7.35, -1) (6, 0) -- (7, 0) (7, 0) -- (7.35, 0);
				\draw[dashed] (3, 0) -- (3, -1) (4, 0) -- (4, -1) (6, 0) -- (6, -1) (7, 0) -- (7, -1);
				\draw[fill = \myblue, draw=burgundy] (0, 0) circle (0.2);
				\draw[fill = \myblue] (1, 0) circle (0.2) (3, 0) circle (0.2) (4, 0) circle (0.2) (6, 0) circle (0.2);
				\draw[fill = \mygray] (0, -1) circle (0.2) (1, -1) circle (0.2) (3, -1) circle (0.2);
				\draw[fill=\myblue] (7, 0) circle (0.1);
				\draw[fill=\mygray] (4, -1) circle (0.1) (6, -1) circle (0.1) (7, -1) circle (0.1);
				\draw[fill=black] (2, 0) circle (0.02) (1.8, 0) circle (0.02) (2.2, 0) circle (0.02) (2, -1) circle (0.02) (1.8, -1) circle (0.02) (2.2, -1) circle (0.02) (5, 0) circle (0.02) (4.8, 0) circle (0.02) (5.2, 0) circle (0.02) (5, -1) circle (0.02) (4.8, -1) circle (0.02) (5.2, -1) circle (0.02) (8, 0) circle (0.02) (7.8, 0) circle (0.02) (8.2, 0) circle (0.02) (8, -1) circle (0.02) (7.8, -1) circle (0.02) (8.2, -1) circle (0.02);
				\draw[dashed] (-0.3, 0.3) to (-0.3, -1.3) -- (-0.3, -1.3) to (8.5, -1.3) -- (8.5, -1.3) to (8.5, 0.3) -- (8.5, 0.3) to (-0.3, 0.3);
				\draw (-0.3,-1.3) to[out=270, in=180] (0.3,-1.5) -- (1.3,-1.5) to[out=0,in=90] (1.65,-1.7) to[out=90,in=180] (2,-1.5) -- (3,-1.5) to[out=0,in=270] (3.3,-1.3);
				\draw (-0.3, 0.3) to[out=90, in=180] (0.4, 0.5) -- (2.5, 0.5) to[out=0,in=270] (3,0.8) to[out=270,in=180] (3.5,0.5) -- (5.9, 0.5) to[out=0,in=90] (6.3, 0.3);
				\node at (1.65, -1.9) {$\scriptstyle n_{2}$};
				\node at (3, 1) {$\scriptstyle n_{1}$};
				\node at (3.7, -2.5) {Sets of Paths};
			\end{scope}
			\begin{scope}[shift={(-4, -5)}]
				\GTrow[0.4]{$\scriptstyle \lambda$, $\scriptstyle \lambda$, $\scriptstyle 0$}
				\GTrow[0]{$\scriptstyle \lambda$, $\scriptstyle n_{2}$}
				\GTrow[-0.4]{$\scriptstyle n_{1}$}
				\node at (0, -1.4) {Gelfand-Tsetlin Bases};
			\end{scope}
			\begin{scope}[shift={(-2, -0.4)}]
				\draw[postaction=decorate, decoration={markings, mark= at position 0 with {\arrowreversed{stealth}}, mark= at position 1 with {\arrow{stealth}}}] (0, 0) to (2, 0);
			\end{scope}
			\begin{scope}[shift={(-4, -1.5)}]
				\draw[postaction=decorate, decoration={markings, mark= at position 0 with {\arrowreversed{stealth}}, mark= at position 1 with {\arrow{stealth}}}] (0, 0) to (0, -2);
			\end{scope}
			\begin{scope}[shift = {(-2.5, -5)}]
				\draw[postaction=decorate, decoration={markings, mark= at position 0 with {\arrowreversed{stealth}}, mark= at position 1 with {\arrow{stealth}}}] (0, 0) to (1.2, 0);
			\end{scope}
			\begin{scope}[shift={(2.75, -2)}]
				\draw[postaction=decorate, decoration={markings, mark= at position 0 with {\arrowreversed{stealth}}, mark= at position 1 with {\arrow{stealth}}}] (0, 0) to (0, -1.2);
			\end{scope}
		\end{tikzpicture}
	\end{array}
\end{equation}
Starting from the vacuum \eqref{sl(3) empty box}, we can add the paths into the first row without restrictions. It means $n_{1} \in [0, \lambda]$. Then, when we are to add the paths into the second line, we must keep in mind that $n_{2} \leqslant n_{1}$. Therefore $n_{2} \in [0, n_{1}]$, and the number of states reads:
\begin{equation}
	\dim_{f.p.}|(\lambda, 0)| = \sum_{n_{2} = 0}^{\lambda}\sum_{n_{1} = n_{2}}^{\lambda}1 = \sum_{n_{1} = 0}^{\lambda}\sum_{n_{2} = 0}^{n_{1}}1 = \sum_{n_{1} = 0}^{\lambda}(n_{1} + 1) = \sum_{k = 1}^{\lambda + 1}k = \dfrac{1}{2}(\lambda + 1)(\lambda + 2) = \dim_{\fs\fl_{3}}|(\lambda, 0)|\,.
\end{equation}

\subsubsection*{Gelfand-Tsetlin Bases}
  
In the \eqref{sl(3) states picture} we mentioned the Gelfand-Tsetlin bases \cite{Gelfand:1950ihs, molev2002gelfandtsetlinbasesclassicallie}. In this case they parametrize the states equivalently as $(n_{1}, n_{2})$. However, this description will be useful within the general approach. We address this correspondence later in sections \ref{subsec: Y(sl(4))} and \ref{subsec: Y(sl(n)) details}.

\subsubsection{The Representations $\Upsilon_{1,\lambda}$}

Having defined the states, we can proceed and construct the Yangian representations. First, we adapt the ansatz \eqref{General Ansatz} to our case. It takes the following form:
\begin{equation}\label{Y(sl(3)) anzatz}
	\begin{aligned}
		&e^{(1)}_{[\lambda, 0]}(z)|n_{1}, n_{2}\rangle = \dfrac{{\bf E}^{[\lambda, 0]}_{(n_{1}, n_{2})\to (n_{1} + 1, n_{2})}}{z - n_{1}\epsilon}|n_{1} + 1, n_{2}\rangle\,, \\
		&e^{(2)}_{[\lambda, 0]}(z)|n_{1}, n_{2}\rangle = \dfrac{{\bf E}^{[\lambda, 0]}_{(n_{1}, n_{2})\to (n_{1}, n_{2} + 1)}}{z + \frac{\epsilon}{2} - n_{2}\epsilon}|n_{1}, n_{2} + 1\rangle\,,\\
		&f^{(1)}_{[\lambda, 0]}(z)|n_{1}, n_{2}\rangle = \dfrac{{\bf F}^{[\lambda, 0]}_{(n_{1}, n_{2})\to (n_{1} - 1, n_{2})}}{z - (n_{1} - 1)\epsilon}|n_{1} - 1, n_{2}\rangle\,, \\
		&f^{(2)}_{[\lambda, 0]}(z)|n_{1}, n_{2}\rangle = \dfrac{{\bf F}^{[\lambda, 0]}_{(n_{1}, n_{2})\to (n_{1}, n_{2} - 1)}}{z + \frac{\epsilon}{2} - (n_{2} - 1)\epsilon}|n_{1}, n_{2} - 1\rangle\,,\\
		&\psi^{(s)[\lambda, 0]}(z)|n_{1}, n_{2}\rangle = \Psi^{(s)[\lambda, 0]}_{(n_{1}, n_{2})}|n_{1}, n_{2}\rangle\,.
	\end{aligned}
\end{equation}
The eigenvalues of $\psi^{(s)[\lambda, 0]}(z)$ read as follows:
\begin{equation}\label{sl(3) eigenfunctions}
	\begin{aligned}
		\Psi^{(1)[\lambda, 0]}_{(n_{1}, n_{2})}(z) &= -\frac{1}{\epsilon}\dfrac{(z - \lambda\epsilon)}{z}\prod_{k = 1}^{n_{1}}\varphi_{1, 1}\bigl(z - (k-1)\epsilon\bigr)\prod_{k = 1}^{n_{2}}\varphi_{1, 2}\Bigl(z + \frac{\epsilon}{2} - (k-1)\epsilon\Bigr) =\\ 
		&= -\frac{1}{\epsilon}\dfrac{(z - \lambda\epsilon)\bigl(z - (n_{2} - 1)\epsilon\bigr)}{(z - n_{1}\epsilon)\bigl(z - (n_{1} - 1)\epsilon\bigr)}\,, \\
		\Psi^{(2)[\lambda, 0]}_{(n_{1}, n_{2})}(z) &= -\frac{1}{\epsilon}\prod_{k = 1}^{n_{1}}\varphi_{2, 1}\bigl(z - (k - 1)\epsilon\bigr)\prod_{k = 1}^{n_{2}}\varphi_{2, 2}\Bigl(z + \frac{\epsilon}{2} - (k - 1)\epsilon\Bigr) = \\
		&= -\frac{1}{\epsilon}\frac{\Bigl(z + \frac{3\epsilon}{2}\Bigr)\bigl(z + \frac{\epsilon}{2} - n_{1}\epsilon\bigr)}{\Bigl(z + \frac{\epsilon}{2} - n_{2}\epsilon\Bigr)\Bigl(z + \frac{\epsilon}{2} - (n_{2} - 1)\epsilon\Bigr)}\,.
	\end{aligned}
\end{equation}
Again, we have set $h = 0$ in these formulas.

\subsubsection*{Amplitudes}

For later convenience we introduce useful ``ICO" matrices:
\begin{equation}\label{ICO base matrices}
	\begin{aligned}
		\text{If } n \leqslant m: \quad &I(n, m) := \begin{tikzpicture}[baseline]
			\matrix[matrix of math nodes,
			left delimiter={(},
			right delimiter={)},
			nodes in empty cells,
			row sep=3pt, column sep=6pt,
			ampersand replacement=\&] (A) {
				1 \& 0 \& \cdots \& 0 \& 0 \\
				0 \& 1 \& \cdots \& 0 \& 0 \\
				\vdots \& \vdots \& \ddots \& \vdots \& \vdots \\
				0 \& 0 \& \cdots \& 1 \& 0 \\
			};
			\node[draw=black, dashed, thick, fit=(A-1-1)(A-4-4), inner sep=2pt] {};
			\draw[decorate,decoration={brace,mirror,amplitude=5pt}, thick]
			(A-4-5.south west) -- (A-4-5.south east)
			node[midway, below=6pt] {\scriptsize $m - n$};
			\node at ([yshift=-2em]A-4-3.south) {\scriptsize $n \times m$};
		\end{tikzpicture}\,,\quad
		n > m\,: \quad 
		&I(n, m) := \begin{tikzpicture}[baseline]
			\matrix[matrix of math nodes,
			left delimiter={(},
			right delimiter={)},
			nodes in empty cells,
			row sep=3pt, column sep=6pt,
			ampersand replacement=\&] (A) {
				1 \& 0 \& \cdots \& 0 \\
				0 \& 1 \& \cdots \& 0 \\
				\vdots \& \vdots \& \ddots \& \vdots \\
				0 \& 0 \& \cdots \& 1 \\
				0 \& 0 \& \cdots \& 0 \\
			};
			\node[draw=black, dashed, thick, fit=(A-1-1)(A-4-4), inner sep=2pt] {};
			\draw[decorate,decoration={brace,mirror,amplitude=5pt}, thick]
			(A-5-4.south east) -- (A-5-4.north east)
			node[midway, right=8pt] {\scriptsize $m - n$};
			\node at ([yshift=-1.5em]A-5-3.south) {\scriptsize $n \times m$};
		\end{tikzpicture}\,, \\
		&C(n) = \begin{tikzpicture}[baseline]
			\matrix[matrix of math nodes,
			left delimiter={(},
			right delimiter={)},
			nodes in empty cells,
			row sep=3pt, column sep=6pt,
			ampersand replacement=\&] (A) {
				0 \& 0 \& 0 \& \dots \& 0\\
				1 \& 0 \& 0 \& \dots \& 0\\
				0 \& 1 \& 0 \& \dots \& 0\\
				\vdots \& \vdots \& \ddots \& \ddots \& \vdots\\
				0 \& 0 \& \dots \& 1 \& 0 \\
			};
			\node at ([yshift=-1.5em]A-5-3.south) {\scriptsize $n \times n$};
		\end{tikzpicture}\,,
		&O(n, m) = \begin{tikzpicture}[baseline]
			\matrix[matrix of math nodes,
			left delimiter={(},
			right delimiter={)},
			nodes in empty cells,
			row sep=3pt, column sep=6pt,
			ampersand replacement=\&] (A) {
				0 \& 0 \& \cdots \& 0 \\
				0 \& 0 \& \cdots \& 0 \\
				\vdots \& \vdots \& \ddots \& \vdots \\
				0 \& 0 \& \cdots \& 0 \\
			};
			\node at ([yshift=-1.5em]A-4-3.south) {\scriptsize $n \times m$};
		\end{tikzpicture}\,,
	\end{aligned}
\end{equation}
where $I(n, m)$ is a generalization of an identity matrix and $O(n, m)$ is just a null matrix. These matrices will be building blocks of the matrices that correspond to the fixed points.

The vacuum expectation values of the fields take the form:
\begin{equation}\label{sl(3) vacuum matrices}
	\begin{aligned}
		&\bar{C}_{1} = C(n_{1})\,, \quad \bar{C}_{2} = C(n_{2})\,, \quad \bar{A} = I(n_{2}, n_{1})\,, \quad \bar{B} = O(n_{1}, n_{2})\,,\\ 
		&\bar{R}_{1} = I(n_{1}, 1)\,, \quad \bar{S}_{1} = O(1, n_{1})\,, \quad \bar{R}_{2} = O(n_{2}, 1)\,, \quad \bar{S}_{2} = O(1, n_{2})\,.
	\end{aligned}
\end{equation}
In the next examples we do not highlight the fields that have zero expectation values.

Implementing the algorithm described in section \ref{subsec: Equivariant Matrix Coefficients}, we can evaluate the corresponding Euler classes:
\begin{equation}\label{sl(3) Euler classes}
	\begin{aligned}
		&\Eul^{[\lambda, 0]}_{(n_{1}, n_{2})} = (\epsilon)^{2n_{1}}\frac{\lambda!}{(\lambda - n_{1})!}(n_{1} - n_{2})!n_{2}!\prod_{k = 1}^{n_{2}}\Bigl(\frac{1}{2}(2k-3)\epsilon\Bigr)^{2}, \\
		&\Eul^{[\lambda, 0]}_{(n_{1}, n_{2})\to (n_{1}+1, n_{2})} = (-\epsilon)(\epsilon)^{2n_{1}}\frac{\lambda!}{(\lambda - n_{1})!}(n_{1} - n_{2})!n_{2}!\prod_{k = 1}^{n_{2}}\Bigl(\frac{1}{2}(2k-3)\epsilon\Bigr)^{2}\,, \\
		&\Eul^{[\lambda, 0]}_{(n_{1}, n_{2})\to (n_{1}, n_{2} + 1)} = \frac{(2n_{2} - 1)}{2}\epsilon~(\epsilon)^{2n_{1}}\dfrac{\lambda!}{(\lambda - n_{1})!}(n_{1} - n_{2} - 1)!n_{2}!\prod_{k = 1}^{n_{2}}\Bigl(\frac{1}{2}(2k-3)\epsilon\Bigr)^{2}\,.
	\end{aligned}
\end{equation}
Next we derive the corresponding matrix coefficients using the formulas \eqref{equiv_matrix_coeffs}:
\begin{equation}\label{sl(3) yangian coefficients}
	\begin{aligned}
		&{\bf E}^{[\lambda, 0]}_{(n_{1}, n_{2})\to (n_{1} + 1, n_{2})} = -\dfrac{1}{\epsilon}, \\
		&{\bf E}^{[\lambda, 0]}_{(n_{1}, n_{2}) \to (n_{1}, n_{2} + 1)} = \dfrac{(n_{1} - n_{2})}{n_{2}\epsilon - \frac{\epsilon}{2}}\,, \\
		&{\bf F}^{[\lambda, 0]}_{(n_{1}, n_{2})\to (n_{1} - 1, n_{2})} = -\epsilon(n_{1} - n_{2})(\lambda - n_{1} + 1)\,, \\
		&{\bf F}^{[\lambda, 0]}_{(n_{1}, n_{2}) \to (n_{1}, n_{2} - 1)} = n_{2}\Bigl(n_{2} - \frac{3}{2}\Bigr)\epsilon\,.
	\end{aligned}
\end{equation}

As in the case of $\myY(\fs\fl_{2})$ \cite{Galakhov_2024}, we expect that the lowest operators of the algebra form the representation of the $\fs\fl_{3}$ algebra itself. In order to see that, one should rescale the states to have unit norm. The corresponding rescaling of the coefficients is given in the formula \eqref{Rescaling Equiv-to-Root}.

The coefficients then take the following form:
\begin{equation}\label{sl(3) normalized coefficients}
	\begin{aligned}
		&e^{(1)}_{0}|n_{1}, n_{2}\rangle = \sqrt{(\lambda - n_{1})(n_{1} - n_{2} + 1)}|n_{1} + 1, n_{2}\rangle\,,\\
		&e^{(2)}_{0}|n_{1}, n_{2}\rangle = \sqrt{(n_{1} - n_{2})(n_{2} + 1)}|n_{1}, n_{2} + 1\rangle\,,\\
		&f^{(1)}_{0}|n_{1}, n_{2}\rangle = \sqrt{(n_{1} - n_{2})(\lambda - n_{1} + 1)}|n_{1} - 1, n_{2}\rangle\,,\\
		&f^{(2)}_{0}|n_{1}, n_{2}\rangle = \sqrt{n_{2}(n_{1} - n_{2} + 1)}|n_{1}, n_{2} - 1\rangle\,.
	\end{aligned}
\end{equation}
These expressions take the same form as the coefficients of the representations $(\lambda, 0)$ of $\fs\fl_{3}$ Lie algebra that was given in \cite{Gelfand:1950ihs}.

Finally, we check the hysteresis relations \eqref{general hysteresis} to ensure that the algebra is self-consistent. For our Yangian $\myY(\fs\fl_{3})$, these relations are satisfied and take the form:
\begin{equation}\label{sl(3) hysteresis relations}
	\begin{aligned}
		&{\bf E}^{[\lambda, 0]}_{(n_{1}, n_{2}) \to (n_{1} + 1, n_{2})}{\bf F}^{[\lambda, 0]}_{(n_{1} +1, n_{2})\to (n_{1}, n_{2})} = (n_{1} - n_{2} + 1)(\lambda - n_{1}) = \mathop{\rm res}\lm_{z= n_{1}\epsilon}\Psi^{(1)}_{\Lambda}(z)\,,\\
		&{\bf E}^{[\lambda, 0]}_{(n_{1}, n_{2}) \to (n_{1}, n_{2} + 1)}{\bf F}^{[\lambda, 0]}_{(n_{1}, n_{2} + 1)\to (n_{1}, n_{2})} = (n_{1} - n_{2})(n_{2} + 1) = \mathop{\rm res}\lm_{z= -\frac{\epsilon}{2} + n_{2}\epsilon}\Psi^{(2)}_{\Lambda}(z)\,,\\
		&\frac{{\bf E}^{[\lambda, 0]}_{(n_{1}, n_{2}) \to (n_{1} + 1, n_{2})}{\bf E}^{[\lambda, 0]}_{(n_{1} +1, n_{2})\to (n_{1} + 1, n_{2} + 1)}}{{\bf E}^{[\lambda, 0]}_{(n_{1}, n_{2}) \to (n_{1}, n_{2} + 1)}{\bf E}^{[\lambda, 0]}_{(n_{1}, n_{2} + 1)\to (n_{1}+1, n_{2}+1)}} = \dfrac{n_{1} + 1 - n_{2}}{n_{1} - n_{2}} = \varphi_{2, 1}\Bigl(n_{2}\epsilon - \frac{\epsilon}{2} - n_{1}\epsilon\Bigr)\,.
	\end{aligned}
\end{equation}

\subsection{$\myY(\fs\fl_{4})$ example}\label{subsec: Y(sl(4))}

Our next example is the algebra $\myY(\fs\fl_{4})$. In this case we have two different types of finite-dimensional representations. The first are symmetric representations $(\lambda, 0, 0)$; their description is similar to the previous case. The second are the $(0, \lambda, 0)$ representations. Their description requires a more involved approach. These representations serve as a crucial stepping point in the later generalization of our results to an arbitrary $\myY(\fs\fl_{n})$ algebra.

\subsubsection{The Algebra}

The quivers that we use to describe the $\mathsf{Y}(\fs\fl_{4})$ algebra and its representations take the forms presented in fig. \ref{sl(4) quivers}.
\begin{figure}[h]
	\centering
	\begin{subfigure}{0.45\textwidth}
		\begin{center}
			\begin{tikzpicture}
				\draw[postaction=decorate, decoration={markings, mark= at position 0.7 with {\arrow{stealth}}}] (0,0) to[out=20,in=160] node[pos=0.5,above] {$\scriptstyle A_{1}$} (1.5,0);
				\draw[postaction=decorate, decoration={markings, mark= at position 0.7 with {\arrow{stealth}}}] (1.5,0) to[out=200,in=340] node[pos=0.5,below] {$\scriptstyle B_{1}$} (0,0);
				\draw[postaction=decorate, decoration={markings, mark= at position 0.7 with {\arrow{stealth}}}] (1.5,0) to[out=20,in=160] node[pos=0.5,above] {$\scriptstyle A_{2}$} (3,0);
				\draw[postaction=decorate, decoration={markings, mark= at position 0.7 with {\arrow{stealth}}}] (3,0) to[out=200,in=340] node[pos=0.5,below] {$\scriptstyle B_{2}$} (1.5,0);
				\draw[postaction=decorate, decoration={markings, mark= at position 0.8 with {\arrow{stealth}}}] (0,0) to[out=60,in=0] (0,0.6) to[out=180,in=120] (0,0);
				\node[above] at (0,0.6) {$\scriptstyle C_1$};
				\begin{scope}[shift={(1.5,0)}]
					\draw[postaction=decorate, decoration={markings, mark= at position 0.8 with {\arrow{stealth}}}] (0,0) to[out=60,in=0] (0,0.6) to[out=180,in=120] (0,0);
					\node[above] at (0,0.6) {$\scriptstyle C_2$};
				\end{scope}
				\begin{scope}[shift={(3,0)}]
					\draw[postaction=decorate, decoration={markings, mark= at position 0.8 with {\arrow{stealth}}}] (0,0) to[out=60,in=0] (0,0.6) to[out=180,in=120] (0,0);
					\node[above] at (0,0.6) {$\scriptstyle C_3$};
				\end{scope}
				\begin{scope}[shift = {(0,0)}]
					\draw[postaction=decorate, decoration={markings, mark= at position 0.7 with {\arrow{stealth}}}] (0,-1.2) to[out=100,in=260] node[pos=0.3, left] {$\scriptstyle R_{1}$} (0,0);
					\draw[postaction=decorate, decoration={markings, mark= at position 0.7 with {\arrow{stealth}}}] (0,0) to[out=280,in=80] node[pos=0.7, right] {$\scriptstyle S_{1}$} (0,-1.2);
					\begin{scope}[shift={(0,-1.2)}]
						\draw[fill=burgundy] (-0.08,-0.08) -- (-0.08,0.08) -- (0.08,0.08) -- (0.08,-0.08) -- cycle;
					\end{scope}
				\end{scope}
				\begin{scope}[shift={(1.5, 0)}]
					\draw[postaction=decorate, decoration={markings, mark= at position 0.7 with {\arrow{stealth}}}, dashed] (0,-1.2) to[out=100,in=260] (0,0);
					\draw[postaction=decorate, decoration={markings, mark= at position 0.7 with {\arrow{stealth}}}, dashed] (0,0) to[out=280,in=80] (0,-1.2);
					\begin{scope}[shift={(0,-1.2)}, dashed]
						\draw[fill=gray!50!white] (-0.08,-0.08) -- (-0.08,0.08) -- (0.08,0.08) -- (0.08,-0.08) -- cycle;
					\end{scope}
				\end{scope}
				\begin{scope}[shift={(3, 0)}]
					\draw[postaction=decorate, decoration={markings, mark= at position 0.7 with {\arrow{stealth}}}, dashed] (0,-1.2) to[out=100,in=260] (0,0);
					\draw[postaction=decorate, decoration={markings, mark= at position 0.7 with {\arrow{stealth}}}, dashed] (0,0) to[out=280,in=80] (0,-1.2);	\def\myblue{white!40!blue}
					\def\mygray{white!40!gray}
					\begin{scope}[shift={(0,-1.2)}, dashed]
						\draw[fill=gray!50!white] (-0.08,-0.08) -- (-0.08,0.08) -- (0.08,0.08) -- (0.08,-0.08) -- cycle;
					\end{scope}
				\end{scope}
				\draw[fill=\myblue] (0,0) circle (0.08);
				\draw[fill=\mygray] (1.5,0) circle (0.08) (3, 0) circle (0.08);
			\end{tikzpicture}
			\subcaption{Describes $\Upsilon_{1, \lambda}$ reps}\label{quiver sl(4) (l, 0, 0)}
		\end{center}
	\end{subfigure}
	\hfill
	\begin{subfigure}{0.45\textwidth}
		\begin{center}
			\begin{tikzpicture}
				\draw[postaction=decorate, decoration={markings, mark= at position 0.7 with {\arrow{stealth}}}] (0,0) to[out=20,in=160] node[pos=0.5,above] {$\scriptstyle A_{1}$} (1.5,0);
				\draw[postaction=decorate, decoration={markings, mark= at position 0.7 with {\arrow{stealth}}}] (1.5,0) to[out=200,in=340] node[pos=0.5,below] {$\scriptstyle B_{1}$} (0,0);
				\draw[postaction=decorate, decoration={markings, mark= at position 0.7 with {\arrow{stealth}}}] (1.5,0) to[out=20,in=160] node[pos=0.5,above] {$\scriptstyle A_{2}$} (3,0);
				\draw[postaction=decorate, decoration={markings, mark= at position 0.7 with {\arrow{stealth}}}] (3,0) to[out=200,in=340] node[pos=0.5,below] {$\scriptstyle B_{2}$} (1.5,0);
				\draw[postaction=decorate, decoration={markings, mark= at position 0.8 with {\arrow{stealth}}}] (0,0) to[out=60,in=0] (0,0.6) to[out=180,in=120] (0,0);
				\node[above] at (0,0.6) {$\scriptstyle C_1$};
				\begin{scope}[shift={(1.5,0)}]
					\draw[postaction=decorate, decoration={markings, mark= at position 0.8 with {\arrow{stealth}}}] (0,0) to[out=60,in=0] (0,0.6) to[out=180,in=120] (0,0);
					\node[above] at (0,0.6) {$\scriptstyle C_2$};
				\end{scope}
				\begin{scope}[shift={(3,0)}]
					\draw[postaction=decorate, decoration={markings, mark= at position 0.8 with {\arrow{stealth}}}] (0,0) to[out=60,in=0] (0,0.6) to[out=180,in=120] (0,0);
					\node[above] at (0,0.6) {$\scriptstyle C_3$};
				\end{scope}
				\begin{scope}[shift = {(1.5,0)}]
					\draw[postaction=decorate, decoration={markings, mark= at position 0.7 with {\arrow{stealth}}}] (0,-1.2) to[out=100,in=260] node[pos=0.3, left] {$\scriptstyle R_{2}$} (0,0);
					\draw[postaction=decorate, decoration={markings, mark= at position 0.7 with {\arrow{stealth}}}] (0,0) to[out=280,in=80] node[pos=0.7, right] {$\scriptstyle S_{2}$} (0,-1.2);
					\begin{scope}[shift={(0,-1.2)}]
						\draw[fill=burgundy] (-0.08,-0.08) -- (-0.08,0.08) -- (0.08,0.08) -- (0.08,-0.08) -- cycle;
					\end{scope}
				\end{scope}
				\begin{scope}
					\draw[postaction=decorate, decoration={markings, mark= at position 0.7 with {\arrow{stealth}}}, dashed] (0,-1.2) to[out=100,in=260] (0,0);
					\draw[postaction=decorate, decoration={markings, mark= at position 0.7 with {\arrow{stealth}}}, dashed] (0,0) to[out=280,in=80] (0,-1.2);
					\begin{scope}[shift={(0,-1.2)}, dashed]
						\draw[fill=gray!50!white] (-0.08,-0.08) -- (-0.08,0.08) -- (0.08,0.08) -- (0.08,-0.08) -- cycle;
					\end{scope}
				\end{scope}
				\begin{scope}[shift={(3, 0)}]
					\draw[postaction=decorate, decoration={markings, mark= at position 0.7 with {\arrow{stealth}}}, dashed] (0,-1.2) to[out=100,in=260] (0,0);
					\draw[postaction=decorate, decoration={markings, mark= at position 0.7 with {\arrow{stealth}}}, dashed] (0,0) to[out=280,in=80] (0,-1.2);
					\begin{scope}[shift={(0,-1.2)}, dashed]
						\draw[fill=gray!50!white] (-0.08,-0.08) -- (-0.08,0.08) -- (0.08,0.08) -- (0.08,-0.08) -- cycle;
					\end{scope}
				\end{scope}
				\draw[fill=\myblue] (1.5,0) circle (0.08);
				\draw[fill=\mygray] (0,0) circle (0.08) (3, 0) circle (0.08);
			\end{tikzpicture}
			\subcaption{Describes $\Upsilon_{2,\lambda}$ reps}\label{quiver sl(4) (0, l, 0)}
		\end{center}
	\end{subfigure}
	\caption{$\myY(\fs\fl_{4})$ quivers}\label{sl(4) quivers}
\end{figure}

The dimensions of the vector spaces associated with the nodes are $n_{1},~n_{2},$ and $n_{3}$ correspondingly, and the superpotential takes the form:
\begin{equation}\label{sl(4) superpotential}
	\CW = \Tr\bigl(A_{1}C_{1}B_{1} + A_{2}C_{2}B_{2} - B_{1}C_{2}A_{1} - B_{2}C_{3}A_{2} + C_{1}^{\lambda_{1}}R_{1}S_{1} + C_{2}^{\lambda_{2}}R_{2}S_{2} + C_{3}^{\lambda_{3}}R_{3}S_{3}\bigr)\,,
\end{equation}
where $(\lambda_{1}, \lambda_{2}, \lambda_{3})$ is equal to $(\lambda, 0, 0)$ or $(0, \lambda, 0)$ correspondingly.

We can determine the weights and the R-charges of the fields:
\begin{equation}\label{sl(4) weights}
	\begin{tblr}{c|c|c|c|c|c|c|c|c|c|c|c|c|c}
		\mbox{Fields} & C_{1} & C_{2} & C_{3} & A_{1} & B_{1} & A_{2} & B_{2} &  R_{1} & S_{1} & R_{2} & S_{2} & R_{3} & S_{3} \\
		\hline
		\mbox{Weights} & \epsilon & \epsilon & \epsilon & -\frac{\epsilon}{2} + h & -\frac{\epsilon}{2} - h & -\frac{\epsilon}{2} + h & -\frac{\epsilon}{2} - h & 0 & -\lambda_{1} \epsilon & 0 & -\lambda_{2}\epsilon & 0 & -\lambda_{3}\epsilon\\
		\mbox{R-charges} & 0 & 0 & 0 & 1 & 1 & 1 & 1 & 0 & 2 & 0 & 2 & 0 & 2
	\end{tblr}\,.
\end{equation}

The Yangian algebra $\myY(\fs\fl_{4})$ involves three families of generators $(e^{(a)}_{k}, f^{(a)}_{k}, \psi^{(a)}_{k})$, $a\in \{1, 2, 3\}$, which we assemble in generating functions:
\begin{equation}\label{sl(4) generating functions}
	e^{(a)}(z) = \sum_{n = 0}^{\infty}\dfrac{e^{(a)}_{n}}{z^{n + 1}}\,, \quad f^{(a)}(z) = \sum_{n = 0}^{\infty}\dfrac{f^{(a)}_{n}}{z^{n + 1}}\,, \quad \psi^{(a)}(z) = 1 + \sum_{n = 0}^{\infty}\dfrac{\psi^{(a)}_{n}}{z^{n + 1}}\,.
\end{equation}
Using the unframed quiver, we construct the bonding factors:
\begin{equation}\label{sl(4) bonding factors}
	\begin{aligned}
		&\varphi_{1, 1}(z) = \varphi_{2, 2}(z) = \varphi_{3, 3}(z) = \frac{z + \epsilon}{z - \epsilon}\,,\\
		&\varphi_{1, 2}(z) = \varphi_{2, 1}(z) = \varphi_{3, 2}(z) = \varphi_{2, 3}(z) = \dfrac{z - \frac{\epsilon}{2}}{z + \frac{\epsilon}{2}}\,.
	\end{aligned}
\end{equation} 
The generating functions satisfy the equation identical to \eqref{sl(3) relations}. We also add the Serre relations:
\begin{equation}\label{sl(4) Serre relations}
	\begin{aligned}
		&\text{If } |a - b| = 1\colon \quad \sum_{\sigma \in \fS_{2}}\Big[e^{(a)}\bigl(u_{\sigma(1)}\bigr),\Bigl[e^{(a)}\bigl(u_{\sigma(2)}\bigr), e^{(b)}(v)\Bigr]\Big] = 0\,,\\
		&\text{If } |a - b| = 2\colon \quad [e^{(a)}(u), e^{(b)}(v)] = 0\,;\\
		&\text{If } |a - b| = 1\colon \quad \sum_{\sigma \in \fS_{2}}\Big[f^{(a)}\bigl(u_{\sigma(1)}\bigr),\Bigl[f^{(a)}\bigl(u_{\sigma(2)}\bigr), f^{(b)}(v)\Bigr]\Big] = 0\,,\\
		&\text{If } |a - b| = 2\colon \quad [f^{(a)}(u), f^{(b)}(v)] = 0\,.
	\end{aligned}
\end{equation}
We do not unfold the relations in terms of modes here, leaving it to the general case \eqref{sl(n) Yangian relations in modes}, \eqref{sl(n) boundary relations}, \eqref{sl(n) serre relations in modes}.

Having the superpotential \eqref{sl(4) superpotential}, we find the corresponding F-term relations:
\begin{equation}\label{sl(4) F-terms}
	\begin{aligned}
		&C_{1}B_{1} - B_{1}C_{2} = 0\,,
		&A_{1}C_{1} - C_{2}A_{1} = 0\,,\\
		&C_{2}B_{2} - B_{2}C_{3} = 0\,,
		&A_{2}C_{2} - C_{3}A_{2} = 0\,, \\
		&B_{1}A_{1} = A_{2}B_{2} = 0\,, 
		&B_{2}A_{2} - A_{1}B_{1} = 0\,,\\
		&C_{1}^{\lambda_{1}}R_{1} = 0\,, \qquad C_{2}^{\lambda_{2}}R_{2} = 0\,, \qquad &C_{3}^{\lambda_{3}}R_{3} = 0\,.
	\end{aligned}
\end{equation}

The fixed points on the variety above are defined by:
\begin{equation}\label{sl(4) fixed points}
	\begin{aligned}
		&[\Phi_{1}, C_{1}] = \epsilon C_{1}\,, \quad [\Phi_{2}, C_{2}] = \epsilon C_{2}\,, \quad [\Phi_{3}, C_{3}] = \epsilon C_{3}\,,\\
		&\Phi_{2}A_{1} - A_{1}\Phi_{1} = \Bigl(-\frac{\epsilon}{2} + h\Bigr)A_{1}\,, \quad \Phi_{1}B_{1} - B_{1}\Phi_{2} = \Bigl(-\frac{\epsilon}{2} - h\Bigr)B_{1}\,,\\
		&\Phi_{3}A_{2} - A_{2}\Phi_{2} = \Bigl(-\frac{\epsilon}{2} + h\Bigr)A_{2}\,, \quad \Phi_{2}B_{2} - B_{2}\Phi_{3} = \Bigl(-\frac{\epsilon}{2} - h\Bigr)B_{2}\,,\\
		&\Phi_{1}R_{1} = 0\,, \quad \Phi_{2}R_{2} = 0\,, \quad \Phi_{3}R_{3} = 0\,, \\
		&S_{1}\Phi_{1} = -\lambda_{1}\epsilon S_{1}\,, \quad \Phi_{2}S_{2} = -\lambda_{2}\epsilon S_{2}\,, \quad S_{3}\Phi_{3} = -\lambda_{3}\epsilon S_{3}\,.
	\end{aligned}
\end{equation}

\subsubsection{The Representations $\Upsilon_{1,\lambda}$}

In this section we use the quiver \ref{quiver sl(4) (l, 0, 0)} to describe the representations $\Upsilon_{1, \lambda} = (\lambda, 0, 0)$ of the Yangian $\myY(\fs\fl_{4})$. The analysis of these representations is quite similar to the previous example described in section \ref{subsec: Y(sl(3))}, so we skip some details of the construction.

\subsubsection*{The States}
 
The paths correspond to the words of the following form $\mathcal{P}_{1} = m(C_{1}, C_{2}, C_{3}, A_{1}, A_{2}, B_{1}, B_{2})\cdot R_{1}$, where $m$ is a monomial.
As before, using F-terms, we can convert the C-fields $C_{3} \to C_{2}, C_{2} \to C_{1}$:
\begin{equation}
	\begin{aligned}
		&B_{2}C_{3}^{a} = C_{2}^{a}B_{2}\,,
		&B_{1}C_{2}^{a} = C_{1}^{a}B_{1}\,, \\
		&C_{3}^{a}A_{2} = A_{2}C_{2}^{a}\,, 
		&C_{2}^{a}A_{1} = A_{1}C_{1}^{a}\,,
	\end{aligned}
\end{equation}
so we can consider only the paths that contain the field $C_{1}$.

\begin{itemize}
	\item If the path $\mathcal{P}_{1}$ doesn't contain the field $A_{1}$, it takes the only possible form:
	\begin{equation}
		C_{1}^{a}R_{1}\,.
	\end{equation}
	\item If $A_{1} \in \mathcal{P}_{1}$, but $A_{2}, B_{1}, B_{2} \notin \mathcal{P}_{1}$, the path takes form:
	\begin{equation}
		A_{1}C_{1}^{a}R_{1}\,.
	\end{equation}
	\item If $A_{1}, A_{2} \in \mathcal{P}_{1}$, the path takes form:
	\begin{equation}
		A_{2}A_{1}C_{1}^{a}R_{1}\,.
	\end{equation}
	\item However, if $B_{1}\in \mathcal{P}_{1}$:
	\begin{equation}
		B_{1}A_{1}C_{1}^{a}R_{1} = 0\,,
	\end{equation}
	so these monomials do not acquire expectation values in the vacuum, and thus $B_{1} \notin \mathcal{P}_{1}$.
	\item The last case when $B_{2} \in \mathcal{P}_{1}$:
	\begin{equation}
		B_{2}A_{2}A_{1}C_{1}^{a}R_{1} = A_{1}B_{1}A_{1}C_{1}^{a}R_{1} = A_{1}(B_{1}A_{1})C_{1}^{a}R_{1} = 0\,,
	\end{equation}
	so $B_{2} \notin \mathcal{P}_{1}$.
\end{itemize}

Summarizing, we observe that a set of possible paths contains elements of three types:
\begin{equation}\label{sl(4) possible paths (l, 0, 0)}
	C_{1}^{a}R_{1}\,, \quad A_{1}C_{1}^{a}R_{1}\,, \quad A_{2}A_{1}C_{1}^{a}R_{1}\,,
\end{equation}
that corresponds to three vertices of the quiver.

The weights and the R-charges of these paths read:
\begin{equation}
	\begin{aligned}
		&h(C_{1}^{a}R_{1}) = a\epsilon\,, \quad h(A_{1}C_{1}^{a}R_{1}) = -\frac{\epsilon}{2} + h + a\epsilon\,, \quad h(A_{2}A_{1}C_{1}^{a}R_{1}) = -\epsilon + 2h + a\epsilon\,, \\
		&R(C_{1}^{a}R_{1}) = 0\,, \quad R(A_{1}C_{1}^{a}R_{1}) = 1\,, \quad R(A_{2}A_{1}C_{1}^{a}R_{1}) = 2\,.
	\end{aligned}
\end{equation}
The general state of the representation can be depicted in the space of equivariant weights and R-charges similar to \eqref{sl(3) random fixed point}. The relation $C_{1}^{\lambda}R_{1} = 0$ serves as a cut-off.
As before, we note that $n_{1} \geqslant n_{2} \geqslant n_{3}$; this is a direct result of the imposed $F$-terms.

\subsubsection*{State counting}

Here we want to count the fixed points of the representation. Due to the arguments similar to that given in the previous example, the dimension should equal the expression from the algebra $\fs\fl_{4}$:
\begin{equation}
	\dim_{\fs\fl_{4}}|(\lambda, 0, 0)| = \dfrac{(\lambda + 1)(\lambda + 2)(\lambda + 3)}{3!}\,.
\end{equation}

The empty box is still a planar rectangle and takes the form:
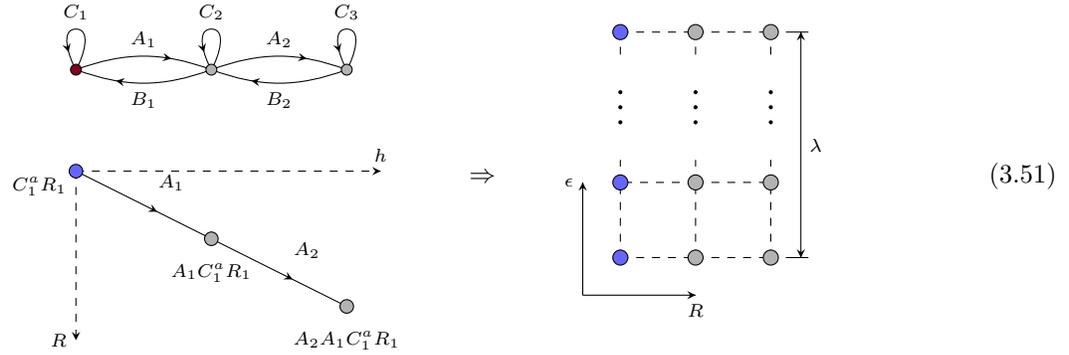
\begin{equation}\label{sl(4) (l, 0, 0) empty box}
	\begin{tikzpicture}[baseline = {(0, -1.5)}, scale = 0.9]
		\begin{scope}[shift = {(0,0)}]
			\draw[postaction=decorate, decoration={markings, mark= at position 0.7 with {\arrow{stealth}}}] (0,0) to[out=20,in=160] node[pos=0.5,above] {$\scriptstyle A_2$} (2,0);
			\draw[postaction=decorate, decoration={markings, mark= at position 0.7 with {\arrow{stealth}}}] (2,0) to[out=200,in=340] node[pos=0.5,below] {$\scriptstyle B_2$} (0,0);
			\begin{scope}[shift={(-2,0)}]
				\draw[postaction=decorate, decoration={markings, mark= at position 0.7 with {\arrow{stealth}}}] (0,0) to[out=20,in=160] node[pos=0.5,above] {$\scriptstyle A_{1}$} (2,0);
				\draw[postaction=decorate, decoration={markings, mark= at position 0.7 with {\arrow{stealth}}}] (2,0) to[out=200,in=340] node[pos=0.5,below] {$\scriptstyle B_{1}$} (0,0);
			\end{scope}
			\draw[postaction=decorate, decoration={markings, mark= at position 0.8 with {\arrow{stealth}}}] (0,0) to[out=60,in=0] (0,0.6) to[out=180,in=120] (0,0);
			\node[above] at (0,0.6) {$\scriptstyle C_2$};
			\begin{scope}[shift={(2,0)}]
				\draw[postaction=decorate, decoration={markings, mark= at position 0.8 with {\arrow{stealth}}}] (0,0) to[out=60,in=0] (0,0.6) to[out=180,in=120] (0,0);
				\node[above] at (0,0.6) {$\scriptstyle C_{3}$};
			\end{scope}
			\begin{scope}[shift={(-2,0)}]
				\draw[postaction=decorate, decoration={markings, mark= at position 0.8 with {\arrow{stealth}}}] (0,0) to[out=60,in=0] (0,0.6) to[out=180,in=120] (0,0);
				\node[above] at (0,0.6) {$\scriptstyle C_{1}$};
			\end{scope}
			\draw[fill=\mygray] (0,0) circle (0.08) (2,0) circle (0.08);
			\draw[fill=burgundy] (-2, 0) circle (0.08);
		\end{scope}
		\draw[dashed, postaction={decorate},
		decoration={markings, mark= at position 1 with {\arrow{stealth}}}] (-2, -1.5) -- (2.5, -1.5);
		\draw[dashed, postaction={decorate},
		decoration={markings, mark= at position 1 with {\arrow{stealth}}}] (-2, -1.5) -- (-2, -4);
		\node[above] at (2.5, -1.5) {$\scriptstyle h$};
		\node[left] at (-2, -4) {$\scriptstyle R$};
		\draw[postaction={decorate},
		decoration={markings, mark= at position 0.6 with {\arrow{stealth}}}] (0, -2.5) to node[above right=0.1] {$\scriptstyle A_{2}$} (2, -3.5);
		\draw[postaction={decorate},
		decoration={markings, mark= at position 0.6 with {\arrow{stealth}}}] (-2, -1.5) to node[above right=0.1] {$\scriptstyle A_{1}$} (0, -2.5);
		\draw[fill=\myblue] (-2, -1.5) circle (0.1);
		\draw[fill=\mygray] (0, -2.5) circle (0.1) (2, -3.5) circle (0.1);
		\node[above left] at (-2, -2) {$\scriptstyle C_{1}^{a}R_{1}$};
		\node at (0, -3) {$\scriptstyle A_{1}C_{1}^{a}R_{1}$};
		\node at (2, -4) {$\scriptstyle A_{2}A_{1}C_{1}^{a}R_{1}$};
	\end{tikzpicture}
	\qquad \Rightarrow \qquad
	\begin{tikzpicture}[baseline={(0, 1)}]
		\begin{scope}[rotate = 90]
			\draw[dashed] (0, 0) -- (1, 0) (0, 0) -- (0, -1) (0, -1) -- (1, -1) (1, -1) -- (1, 0) (0, -1) -- (0, -2) (1, -1) -- (1, -2) (3, -1) -- (3, -2) (0, -2) -- (1, -2);
			\draw[dashed] (1, 0) -- (1.3, 0) (1, -1) -- (1.3, -1) (2.7, 0) -- (3, 0) (2.7, -1) -- (3, -1) (3, 0) -- (3, -1) (1, -2) -- (1.3, -2) (2.7, -2) -- (3, -2);
			\draw[fill=\myblue] (0, 0) circle (0.1) (1, 0) circle (0.1) (3, 0) circle (0.1);
			\draw[fill=\mygray] (0, -1) circle (0.1) (1, -1) circle (0.1) (3, -1) circle (0.1);
			\draw[fill=\mygray] (0, -2) circle (0.1) (1, -2) circle (0.1) (3, -2) circle (0.1);
			\draw[fill=black] (2, 0) circle (0.02) (1.8, 0) circle (0.02) (2.2, 0) circle (0.02) (2, -1) circle (0.02) (1.8, -1) circle (0.02) (2.2, -1) circle (0.02) (2, -2) circle (0.02) (1.8, -2) circle (0.02) (2.2, -2) circle (0.02);
			\draw (0, -2.2) -- (0, -2.5) (3, -2.2) -- (3, -2.5);
			\draw[postaction={decorate},
			decoration={markings, mark= at position 0 with {\arrowreversed{stealth}}, mark= at position 1 with {\arrow{stealth}}}] (0, -2.4) -- (3, -2.4);
			\node at (1.5, -2.6) {$\scriptstyle \lambda$};
			\draw[postaction={decorate},
			decoration={markings, mark= at position 1 with {\arrow{stealth}}}] (-0.5, 0.5) to (1, 0.5);
			\draw[postaction={decorate},
			decoration={markings, mark= at position 1 with {\arrow{stealth}}}] (-0.5, 0.5) to (-0.5, -1);
			\node[left] at (1, 0.5) {$\scriptstyle \epsilon$};
			\node[below] at (-0.5, -1) {$\scriptstyle R$};
		\end{scope}
	\end{tikzpicture}
\end{equation}
Next we obtain the other states by increasing the dimensions $n_{1}, n_{2}$, and $n_{3}$. We can add the nodes into the first row before $n_{1} = \lambda$, so $n_{1} \in [0, \lambda]$. Similarly, the dimension $n_{2}$ is restricted by $n_{1}$ as follows: $n_{2} \in [0, n_{1}]$. Finally, the last dimension $n_{3}\in [0, n_{2}]$ and the number of the fixed points reads:
\begin{equation}
	\dim_{f.p.}|(\lambda, 0, 0)| = \sum_{n_{1} = n_{2}}^{\lambda}\sum_{n_{2} = n_{3}}^{\lambda}\sum_{n_{3} = 0}^{\lambda}1 = \sum_{n_{1} = 0}^{\lambda}\sum_{n_{2} = 0}^{n_{1}}\sum_{n_{3} = 0}^{n_{2}}1 = \sum_{n_{1} = 0}^{\lambda}\sum_{n_{2} = 0}^{n_{1}}(n_{2} + 1) = \dfrac{(\lambda + 1)(\lambda + 2)(\lambda + 3)}{3!} = \dim_{\fs\fl_{4}}|(\lambda, 0, 0)|\,.
\end{equation}

\subsubsection*{Gelfand-Tsetlin Bases}

In this representation each node of the quiver has the paths of only one type. Therefore the states can be uniquely parameterized by the dimensions $n_{1}, n_{2}, n_{3}$. The Gelfand-Tsetlin bases take the following form:
\begin{equation}\label{sl(4) (l, 0, 0) Gelfand-Tsetlin}
	|n_{1}, n_{2}, n_{3}\rangle = \quad
	\begin{tikzpicture}[baseline={(0,-0.0)}, every node/.style={inner sep=1pt}]
		\begin{scope}[shift={(0,0)}]
			\GTrow[1.2]{$\lambda$, $\lambda$, $\lambda$, $0$}
			\GTrow[0.6]{$\lambda$, $\lambda$, $n_{3}$}
			\GTrow[0]{$\lambda$, $n_{2}$}
			\GTrow[-0.6]{$n_{1}$}
		\end{scope}
	\end{tikzpicture}\,.
\end{equation}
We use the notation $|n_{1}, n_{2}, n_{3}\rangle$ for simplicity, as the descriptions are equivalent in this case.

\subsubsection*{The Representation}

The next step is to adopt the ansatz \eqref{General Ansatz}. It reads as follows:
\begin{equation}\label{sl(4) (l, 0, 0) ansatz}
	\begin{aligned}
		&e^{(1)}_{[\lambda, 0, 0]}(z)|n_{1}, n_{2}, n_{3}\rangle = \dfrac{{\bf E}^{[\lambda, 0, 0]}_{(n_{1}, n_{2}, n_{3})\to (n_{1} + 1, n_{2}, n_{3})}}{z - n_{1}\epsilon}|n_{1} + 1, n_{2}, n_{3}\rangle\,, \\
		&e^{(2)}_{[\lambda, 0, 0]}(z)|n_{1}, n_{2}, n_{3}\rangle = \dfrac{{\bf E}^{[\lambda, 0, 0]}_{(n_{1}, n_{2}, n_{3})\to (n_{1}, n_{2} + 1, n_{3})}}{z + \frac{\epsilon}{2} - n_{2}\epsilon}|n_{1}, n_{2} + 1, n_{3}\rangle\,,\\
		&e^{(3)}_{[\lambda, 0, 0]}(z)|n_{1}, n_{2}, n_{3}\rangle = \dfrac{{\bf E}^{[\lambda, 0, 0]}_{(n_{1}, n_{2}, n_{3})\to (n_{1}, n_{2}, n_{3} + 1)}}{z + \epsilon - n_{3}\epsilon}|n_{1}, n_{2}, n_{3} + 1\rangle\,,\\
		&f^{(1)}_{[\lambda, 0, 0]}(z)|n_{1}, n_{2}, n_{3}\rangle = \dfrac{{\bf F}^{[\lambda, 0, 0]}_{(n_{1}, n_{2}, n_{3})\to (n_{1} - 1, n_{2}, n_{3})}}{z - (n_{1} - 1)\epsilon}|n_{1} - 1, n_{2}, n_{3}\rangle\,, \\
		&f^{(2)}_{[\lambda, 0, 0]}(z)|n_{1}, n_{2}, n_{3}\rangle = \dfrac{{\bf F}^{[\lambda, 0, 0]}_{(n_{1}, n_{2}, n_{3})\to (n_{1}, n_{2} - 1, n_{3})}}{z + \frac{\epsilon}{2} - (n_{2} - 1)\epsilon}|n_{1}, n_{2} - 1, n_{3}\rangle\,,\\
		&f^{(3)}_{[\lambda, 0, 0]}(z)|n_{1}, n_{2}, n_{3}\rangle = \dfrac{{\bf F}^{[\lambda, 0, 0]}_{(n_{1}, n_{2}, n_{3})\to (n_{1}, n_{2}, n_{3} - 1)}}{z + \epsilon - (n_{3} - 1)\epsilon}|n_{1}, n_{2}, n_{3} - 1\rangle\,,\\
		&\psi^{(s)[\lambda, 0, 0]}(z)|n_{1}, n_{2}, n_{3}\rangle = \Psi^{(s)[\lambda, 0, 0]}_{(n_{1}, n_{2}, n_{3})}(z)|n_{1}, n_{2}, n_{3}\rangle\,.
	\end{aligned}
\end{equation}
The eigenfunctions take the form:
\begin{equation}\label{sl(4) (l, 0, 0) eigenfunctions}
	\begin{aligned}
		\Psi^{(1)[\lambda, 0, 0]}_{(n_{1}, n_{2}, n_{3})}(z) &= -\frac{1}{\epsilon}\dfrac{(z - \lambda\epsilon)}{z}\prod_{k = 1}^{n_{1}}\varphi_{1, 1}\bigl(z - (k-1)\epsilon\bigr)\prod_{k = 1}^{n_{2}}\varphi_{1, 2}\Bigl(z + \frac{\epsilon}{2} - (k-1)\epsilon\Bigr) =\\ 
		&= -\frac{1}{\epsilon}\dfrac{(z - \lambda\epsilon)\bigl(z - (n_{2} - 1)\epsilon\bigr)}{(z - n_{1}\epsilon)\bigl(z - (n_{1} - 1)\epsilon\bigr)}\,, \\
		\Psi^{(2)[\lambda, 0, 0]}_{(n_{1}, n_{2}, n_{3})}(z) &= -\frac{1}{\epsilon}\prod_{k = 1}^{n_{1}}\varphi_{2, 1}\bigl(z - (k - 1)\epsilon\bigr)\prod_{k = 1}^{n_{2}}\varphi_{2, 2}\Bigl(z + \frac{\epsilon}{2} - (k - 1)\epsilon\Bigr)\prod_{k = 1}^{n_{3}}\varphi_{2, 3}\Bigl(z + \epsilon - (k-1)\epsilon\Bigr) = \\
		&= -\frac{1}{\epsilon}\frac{\Bigl(z + \frac{\epsilon}{2} - n_{1}\epsilon\Bigr)\Bigl(z + \frac{3\epsilon}{2} - n_{3}\epsilon\Bigr)}{\Bigl(z + \frac{\epsilon}{2} - n_{2}\epsilon\Bigr)\Bigl(z + \frac{\epsilon}{2} - (n_{2} - 1)\epsilon\Bigr)}\,, \\
		\Psi^{(3)[\lambda, 0, 0]}_{(n_{1}, n_{2}, n_{3})}(z) &= \frac{1}{\epsilon}\prod_{k = 1}^{n_{2}}\varphi_{3, 2}\Bigl(z + \frac{\epsilon}{2} - (k-1)\epsilon\Bigr)\prod_{k = 1}^{n_{3}}\varphi_{3, 3}\Bigl(z + \epsilon - (k - 1)\epsilon\Bigr) = \\
		&= -\frac{1}{\epsilon}\frac{(z + 2\epsilon)(z + \epsilon - n_{2}\epsilon)}{\bigl(z + \epsilon - n_{3}\epsilon\bigr)\bigl(z + \epsilon - (n_{3} - 1)\epsilon\bigr)}\,.
	\end{aligned}
\end{equation}

\subsubsection*{Amplitudes}

The vacuum expectation values of the fields take the form:
\begin{equation}\label{sl(4) (l, 0, 0) vacuum fixed points}
	\begin{aligned}
		&\bar{C}_{1} = C(n_{1}), \quad \bar{C}_{2} = C(n_{2})\,, \quad \bar{C}_{3} = C(n_{3})\,,\\
		&\bar{A}_{1} = I(n_{2}, n_{1})\,, \quad \bar{A}_{2} = I(n_{3}, n_{2})\,, \quad \bar{R}_{1} = I(n_{1}, 1)\,,
	\end{aligned}
\end{equation}
whereas the remaining fields $B_{1}, B_{2}, S_{1}, R_{2}, S_{2}, R_{3},$ and $S_{3}$ acquire zero expectation values.

Applying the equivariant methods, we get the matrix elements of the representation:
\begin{equation}\label{sl(4) (l, 0, 0) equiv representation}
	\begin{aligned}
		&{\bf E}^{[\lambda, 0, 0]}_{(n_{1}, n_{2}, n_{3}) \to (n_{1}+1, n_{2}, n_{3})} = -\dfrac{1}{\epsilon}\,,\\
		&{\bf E}^{[\lambda, 0, 0]}_{(n_{1}, n_{2}, n_{3}) \to (n_{1}, n_{2}+1, n_{3})} = \dfrac{n_{1} - n_{2}}{n_{2}\epsilon - \frac{\epsilon}{2}}\,,\\
		&{\bf E}^{[\lambda, 0, 0]}_{(n_{1}, n_{2}, n_{3}) \to (n_{1}, n_{2}, n_{3} + 1)} = \dfrac{n_{2}-n_{3}}{n_{3}\epsilon - \epsilon}\,,\\
		&{\bf F}^{[\lambda, 0, 0]}_{(n_{1}, n_{2}, n_{3}) \to (n_{1}-1, n_{2}, n_{3})} = -\epsilon(n_{1}-n_{2})(\lambda-n_{1}+1)\,,\\
		&{\bf F}^{[\lambda, 0, 0]}_{(n_{1}, n_{2}, n_{3}) \to (n_{1}, n_{2}-1, n_{3})} = (n_{2}-n_{3})\Bigl(n_{2} - \frac{3}{2}\Bigr)\epsilon\,,\\
		&{\bf F}^{[\lambda, 0, 0]}_{(n_{1}, n_{2}, n_{3}) \to (n_{1}, n_{2}, n_{3}-1)} = n_{3}(n_{3}-2)\epsilon\,.
	\end{aligned}
\end{equation}
The formulas for the Euler classes are quite bulky, so we do not write them explicitly here.

Next, we normalize the coefficients using \eqref{Rescaling Equiv-to-Root}. The result reads as follows:
\begin{equation}\label{sl(4) Lie algebra coefficients}
	\begin{aligned}
		&e^{(1)}_{0}|n_{1}, n_{2}, n_{3}\rangle = \sqrt{(\lambda - 	n_{1})(n_{1} - n_{2} + 1)}|n_{1} + 1, n_{2}, n_{3}\rangle\,,\\
		&e^{(2)}_{0}|n_{1}, n_{2}, n_{3}\rangle = \sqrt{(n_{1} - n_{2})(n_{2} - n_{3} + 1)}|n_{1}, n_{2} + 1, n_{3}\rangle\,,\\
		&e^{(3)}_{0}|n_{1}, n_{2}, n_{3}\rangle = \sqrt{(n_{2} - n_{3})(n_{3} + 1)}|n_{1}, n_{2}, n_{3} + 1\rangle\,;\\
		&f^{(1)}_{0}|n_{1}, n_{2}, n_{3}\rangle = \sqrt{(n_{1} - 	n_{2})(\lambda - n_{1} + 1)}|n_{1} - 1, n_{2}, n_{3}\rangle\,,\\
		&f^{(2)}_{0}|n_{1}, n_{2}, n_{3}\rangle = \sqrt{(n_{2} - n_{3})(n_{1} - 	n_{2} + 1)}|n_{1}, n_{2} - 1, n_{3}\rangle\,,\\
		&f^{(3)}_{0}|n_{1}, n_{2}, n_{3}\rangle = \sqrt{n_{3}(n_{2} - n_{3} + 1)}|n_{1}, n_{2}, n_{3}-1\rangle\,.
	\end{aligned}
\end{equation}
As expected, these expressions coincide with the formulas introduced by Gelfand in \cite{Gelfand:1950ihs}.

Finally, we check the hysteresis relations \eqref{general hysteresis}. For the representation $\Upsilon_{1, \lambda}$ of the algebra $\myY(\fs\fl_{4})$, the relations are satisfied and take the form:
\begin{equation}\label{sl(4) (l, 0, 0) hysteresis}
	\begin{aligned}
		&{\bf E}^{[\lambda, 0, 0]}_{(n_{1}, n_{2}, n_{3}) \to (n_{1} + 1, n_{2}, n_{3})}{\bf F}^{[\lambda, 0, 0]}_{(n_{1} +1, n_{2}, n_{3})\to (n_{1}, n_{2}, n_{3})} = (n_{1} - n_{2} + 1)(\lambda - n_{1}) = \mathop{\rm res}\lm_{z= n_{1}\epsilon}\Psi^{(1)}_{\Lambda}(z)\,,\\
		&{\bf E}^{[\lambda, 0, 0]}_{(n_{1}, n_{2}, n_{3}) \to (n_{1}, n_{2} + 1, n_{3})}{\bf F}^{[\lambda, 0, 0]}_{(n_{1}, n_{2} + 1, n_{3})\to (n_{1}, n_{2}, n_{3})} = (n_{1} - n_{2})(n_{2} - n_{3} + 1) = \mathop{\rm res}\lm_{z= -\frac{\epsilon}{2} + n_{2}\epsilon}\Psi^{(2)}_{\Lambda}(z)\,,\\
		&{\bf E}^{[\lambda, 0, 0]}_{(n_{1}, n_{2}, n_{3}) \to (n_{1}, n_{2}, n_{3} + 1)}{\bf F}^{[\lambda, 0, 0]}_{(n_{1}, n_{2}, n_{3} - 1)\to (n_{1}, n_{2}, n_{3})} = (n_{2} - n_{3})(n_{3} + 1) = \mathop{\rm res}\lm_{z= -\epsilon + n_{3}\epsilon}\Psi^{(2)}_{\Lambda}(z)\,;\\
		&\frac{{\bf E}^{[\lambda, 0, 0]}_{(n_{1}, n_{2}, n_{3}) \to (n_{1} + 1, n_{2}, n_{3})}{\bf E}^{[\lambda, 0, 0]}_{(n_{1} +1, n_{2}, n_{3})\to (n_{1} + 1, n_{2} + 1, n_{3})}}{{\bf E}^{[\lambda, 0, 0]}_{(n_{1}, n_{2}, n_{3}) \to (n_{1}, n_{2} + 1, n_{3})}{\bf E}^{[\lambda, 0, 0]}_{(n_{1}, n_{2} + 1, n_{3})\to (n_{1}+1, n_{2}+1, n_{3})}} = \dfrac{n_{1} + 1 - n_{2}}{n_{1} - n_{2}} = \varphi_{2, 1}\Bigl(n_{2}\epsilon - \frac{\epsilon}{2} - n_{1}\epsilon\Bigr)\,,\\
		&\frac{{\bf E}^{[\lambda, 0, 0]}_{(n_{1}, n_{2}, n_{3}) \to (n_{1}, n_{2} + 1, n_{3})}{\bf E}^{[\lambda, 0, 0]}_{(n_{1}, n_{2} + 1, n_{3})\to (n_{1}, n_{2} + 1, n_{3} + 1)}}{{\bf E}^{[\lambda, 0, 0]}_{(n_{1}, n_{2}, n_{3}) \to (n_{1}, n_{2}, n_{3} + 1)}{\bf E}^{[\lambda, 0, 0]}_{(n_{1}, n_{2}, n_{3} + 1)\to (n_{1}, n_{2}+1, n_{3} + 1)}} = \dfrac{n_{2} + 1 - n_{3}}{n_{2} - n_{3}} = \varphi_{3, 2}\Bigl(n_{3} - \epsilon - n_{2}\epsilon + \frac{\epsilon}{2}\Bigr)\,.
	\end{aligned}
\end{equation}

\subsubsection{The Representations $\Upsilon_{2,\lambda}$}

The $F$-terms that serve as a cut-off for crystal growth in this case read:
\begin{equation}
	R_{1} = 0\,, \quad C_{2}^{\lambda}R_{2} = 0\,, \quad R_{3} = 0\,.
\end{equation}

The paths correspond to the words of the form $\mathcal{P}_{2} = m(C_{1}, C_{2}, C_{3}, A_{1}, A_{2}, B_{1}, B_{2})\cdot R_{2}$. Here we use the field $C_{2}$ instead of $C_{1}, C_{3}$ because this is more convenient.

\begin{itemize}
	\item We have the paths
	\begin{equation}
		C_{2}^{a}R_{2}, \quad A_{2}C_{2}^{a}R_{2}\,,
	\end{equation}
	similarly to the case above.
	\item If $B_{1} \in \mathcal{P}_{2}$ and $A_{1}, A_{2}, B_{2} \notin \mathcal{P}_{2}$ we have the paths:
	\begin{equation}
		B_{1}C_{2}^{a}R_{2}\,.
	\end{equation}
	This was not true for the case described above.
	\item We can proceed and get new paths using $A_{1}$ or $B_{2}$:
	\begin{equation}
		A_{1}\underline{B_{1}C_{2}^{a}R_{2}} = B_{2}\underline{A_{2}C_{2}^{a}R_{2}}\,,
	\end{equation}
	where we have applied the relation $A_{1}B_{1} = B_{2}A_{2}$. We emphasize that the vacuum expectation values of these monomials do not vanish. This is the main difference from the case described above.
	\item However, when we have more than two fields $A_{1}, A_{2}, B_{1}, B_{2}$ in a monomial $m$, the paths equal zero:
	\begin{equation}
		\begin{aligned}
			&A_{2}A_{1}B_{1}C_{2}^{a}R_{2} = (A_{2}B_{2})A_{2}C_{2}^{a}R_{2} = 0\,, \\
			&B_{1}B_{2}A_{2}C_{2}^{a}R_{2} = (B_{1}A_{1})B_{1}C_{2}^{a}R_{2} = 0\,.
		\end{aligned}
	\end{equation}
\end{itemize}

Therefore, the set of the paths $\mathcal{P}_{2}$ can contain only the elements of four types:
\begin{equation}\label{sl(4) possible paths (0, l, 0)}
	C_{2}^{a}R_{2}, \quad A_{2}C_{2}^{a}R_{2}, \quad B_{1}C_{2}^{a}R_{2}, \quad A_{1}B_{1}C_{2}^{a}R_{2} = B_{2}A_{2}C_{2}^{a}R_{2}\,.
\end{equation}

The equivariant weights and the R-charges of these atoms take the form:
\begin{equation}\label{sl(4) (0, l, 0) paths weights and r-charges}
	\begin{aligned}
		&h(C_{2}^{a}R_{2}) = a\epsilon\,, \quad h(A_{2}C_{2}^{a}R_{2}) = -\frac{\epsilon}{2} + h + a\epsilon\,, \quad h(B_{1}C_{2}^{a}R_{2}) = -\frac{\epsilon}{2} - h + a\epsilon\,, \quad h(A_{1}B_{1}C_{2}^{a}R_{2}) = -\epsilon + a\epsilon\,,\\
		&R(C_{2}^{a}R_{2}) = 0\,, \quad R(A_{2}C_{2}^{a}R_{2}) = 1\,, \quad R(B_{1}C_{2}^{a}R_{2}) = 1\,, \quad R(A_{1}B_{1}C_{2}^{a}R_{2}) = 2\,.
	\end{aligned}
\end{equation}
We emphasize that the additional weight $h$ helps us to distinguish the paths $A_{2}C_{2}^{a}R_{2}$ and $B_{1}C_{2}^{a}R_{2}$, which would overlap if we assumed $h = 0$ from the beginning. One could also note that the paths $C_{2}^{a}R_{2}$ and $A_{1}B_{1}C_{2}^{a}R_{2}$ live in the same vector space $V_{2}$ yet have different R-charges, therefore, they also do not overlap in the space $\Xi$.

We can highlight these facts on the R-equivariant space and draw the corresponding empty box \eqref{sl(4) (0, l, 0) paths on quiver}:
\tdplotsetmaincoords{70}{130}
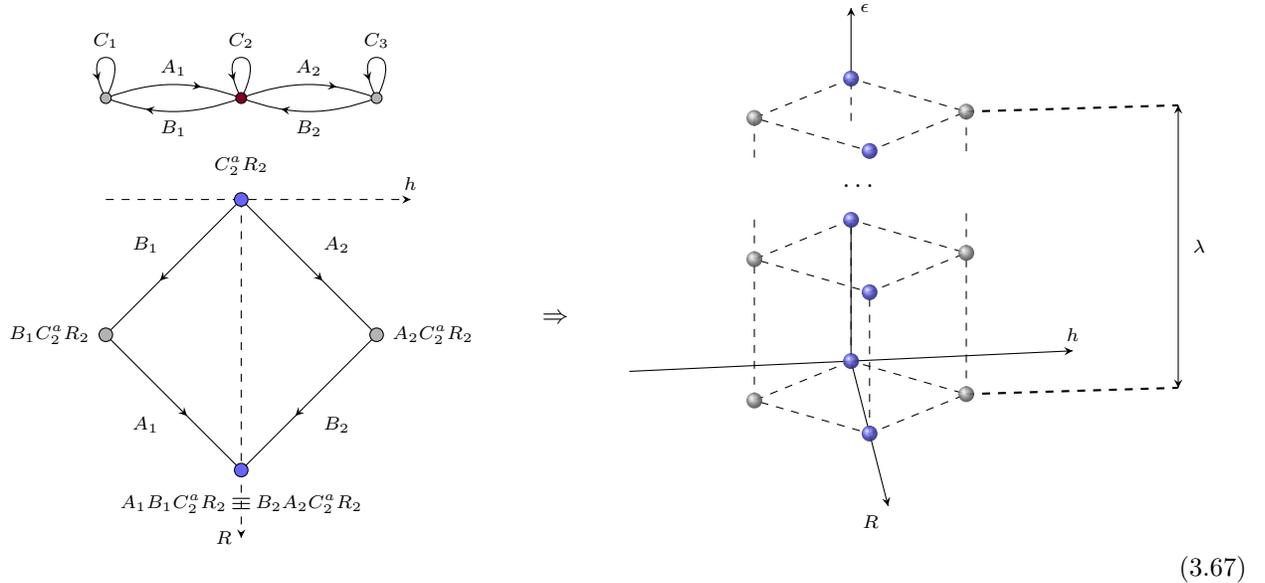
\begin{equation}\label{sl(4) (0, l, 0) paths on quiver}
	\begin{tikzpicture}[scale = 0.9, baseline={(0, -3)}]
		\begin{scope}[shift = {(0,0)}]
			\draw[postaction=decorate, decoration={markings, mark= at position 0.7 with {\arrow{stealth}}}] (0,0) to[out=20,in=160] node[pos=0.5,above] {$\scriptstyle A_2$} (2,0);
			\draw[postaction=decorate, decoration={markings, mark= at position 0.7 with {\arrow{stealth}}}] (2,0) to[out=200,in=340] node[pos=0.5,below] {$\scriptstyle B_2$} (0,0);
			\begin{scope}[shift={(-2,0)}]
				\draw[postaction=decorate, decoration={markings, mark= at position 0.7 with {\arrow{stealth}}}] (0,0) to[out=20,in=160] node[pos=0.5,above] {$\scriptstyle A_{1}$} (2,0);
				\draw[postaction=decorate, decoration={markings, mark= at position 0.7 with {\arrow{stealth}}}] (2,0) to[out=200,in=340] node[pos=0.5,below] {$\scriptstyle B_{1}$} (0,0);
			\end{scope}
			\draw[postaction=decorate, decoration={markings, mark= at position 0.8 with {\arrow{stealth}}}] (0,0) to[out=60,in=0] (0,0.6) to[out=180,in=120] (0,0);
			\node[above] at (0,0.6) {$\scriptstyle C_2$};
			\begin{scope}[shift={(2,0)}]
				\draw[postaction=decorate, decoration={markings, mark= at position 0.8 with {\arrow{stealth}}}] (0,0) to[out=60,in=0] (0,0.6) to[out=180,in=120] (0,0);
				\node[above] at (0,0.6) {$\scriptstyle C_{3}$};
			\end{scope}
			\begin{scope}[shift={(-2,0)}]
				\draw[postaction=decorate, decoration={markings, mark= at position 0.8 with {\arrow{stealth}}}] (0,0) to[out=60,in=0] (0,0.6) to[out=180,in=120] (0,0);
				\node[above] at (0,0.6) {$\scriptstyle C_{1}$};
			\end{scope}
			\draw[fill=\mygray] (-2,0) circle (0.08) (2,0) circle (0.08);
			\draw[fill=burgundy] (0, 0) circle (0.08);
		\end{scope}
		\draw[dashed, postaction={decorate},
		decoration={markings, mark= at position 1 with {\arrow{stealth}}}] (-2, -1.5) -- (2.5, -1.5);
		\draw[dashed, postaction={decorate},
		decoration={markings, mark= at position 1 with {\arrow{stealth}}}] (0, -1.5) -- (0, -6.5);
		\node[above] at (2.5, -1.5) {$\scriptstyle h$};
		\node[left] at (0, -6.5) {$\scriptstyle R$};
		\draw[postaction={decorate},
		decoration={markings, mark= at position 0.6 with {\arrow{stealth}}}] (0, -1.5) to node[above left=0.1] {$\scriptstyle B_{1}$} (-2, -3.5);
		\draw[postaction={decorate},
		decoration={markings, mark= at position 0.6 with {\arrow{stealth}}}] (0, -1.5) to node[above right=0.1] {$\scriptstyle A_{2}$} (2, -3.5);
		\draw[postaction={decorate},
		decoration={markings, mark= at position 0.6 with {\arrow{stealth}}}] (-2, -3.5) to node[below left=0.1] {$\scriptstyle A_{1}$} (0, -5.5);
		\draw[postaction={decorate},
		decoration={markings, mark= at position 0.6 with {\arrow{stealth}}}] (2, -3.5) to node[below right=0.1] {$\scriptstyle B_{2}$} (0, -5.5);
		\draw[fill=\myblue] (0, -1.5) circle (0.1) (0, -5.5) circle (0.1);
		\draw[fill=\mygray] (-2, -3.5) circle (0.1) (2, -3.5) circle (0.1);
		\node at (0, -1) {$\scriptstyle C_{2}^{a}R_{2}$};
		\node[right] at (2.1, -3.5) {$\scriptstyle A_{2}C_{2}^{a}R_{2}$};
		\node[left] at (-2.1, -3.5) {$\scriptstyle B_{1}C_{2}^{a}R_{2}$};
		\node at (0, -6) {$\equiv$};
		\node at (-1, -6) {$\scriptstyle A_{1}B_{1}C_{2}^{a}R_{2}$};
		\node at (1, -6) {$\scriptstyle B_{2}A_{2}C_{2}^{a}R_{2}$};
	\end{tikzpicture}
	\qquad \Rightarrow \qquad
	\begin{tikzpicture}[tdplot_main_coords, scale=2, baseline={(0, 0.5)}]
		\coordinate (A) at (0,0,0);
		\coordinate (B) at (1,0,0);
		\coordinate (C) at (1,1,0);
		\coordinate (D) at (0,1,0);
		\coordinate (A') at (0,0,1);
		\coordinate (B') at (1,0,1);
		\coordinate (C') at (1,1,1);
		\coordinate (D') at (0,1,1);
		\coordinate (E) at (0,0,2);
		\coordinate (F) at (1,0,2);
		\coordinate (G) at (1,1,2);
		\coordinate (H) at (0,1,2);
		\draw[thin] (A) -- (A'); 
		\draw[postaction={decorate},
		decoration={markings, mark= at position 1 with {\arrow{stealth}}}] (E) -- (0, 0, 2.5);
		\draw[dashed] ($(A')+(0, 0,0.7)$) -- (E);
		\draw[postaction={decorate},
		decoration={markings, mark= at position 1 with {\arrow{stealth}}}] (A) -- ($(C)+(1, 1, 0)$);
		\draw[postaction={decorate},
		decoration={markings, mark= at position 1 with {\arrow{stealth}}}] (1.1, -1, 0) -- (-1.1, 1, 0);
		\node[right] at (0, 0, 2.5) {$\scriptstyle \epsilon$};
		\node[below left] at (2, 2, 0) {$\scriptstyle R$};
		\node[above] at (-1.1, 1, 0) {$\scriptstyle h$};
		\draw[dashed] (A) -- (B) -- (C) -- (D) -- cycle;
		\draw[dashed] (A') -- (B') -- (C') -- (D') -- cycle;
		\draw[dashed] (E) -- (F) -- (G) -- (H) -- cycle;
		\draw[dashed] (A) -- (A');
		\draw[dashed] (B) -- (B');
		\draw[dashed] (C) -- (C');
		\draw[dashed] (D) -- (D');
		\draw[dashed] (B') -- (1,0,1.3);
		\draw[dashed] (D') -- (0,1,1.3);
		\draw[dashed] (F) -- (1,0,1.7);
		\draw[dashed] (H) -- (0,1,1.7);
		\node at (0.5,0.5, 1.5) {$\dots$};
		\draw[thick, dashed] ($(D)+(0,0,0)$) -- ($(D)+(0,0,0)+(-1,1,0)$);
		\draw[thick, dashed] ($(H)+(0,0,0)$) -- ($(H)+(0,0,0)+(-1,1,0)$);
		\draw[postaction={decorate},
		decoration={markings, mark= at position 0 with {\arrowreversed{stealth}}, mark= at position 1 with {\arrow{stealth}}}] ($(D)+(-1,1,0)+(0,0,0)$) -- 
		($(H)+(0,0,0)+(-1,1,0)$)
		node[midway, right=2pt] {$\scriptstyle \lambda$};
		\foreach \point in {B,D,B',D',F,H}
		{
			\shade[ball color=\mygray, opacity=0.9] (\point) circle (1.5pt);
		}
		\foreach \point in {A,C,A',C',E,G}
		{
			\shade[ball color=\myblue, opacity=0.9] (\point) circle (1.5pt);
		}
	\end{tikzpicture}
\end{equation}
The picture explicitly shows that the crystals of the representation $\Upsilon_{2,\lambda}$ naturally form a 3-dimensional graph. It justifies our assumption about the existence of an additional parameter $h$. 

We want to emphasize some key differences between this case and the previous cases that corresponded to symmetric representations. Although the previous crystals also lived in the 3-dimensional space, they were effectively 2-dimensional. All the nodes consisted of the atoms of only one type. Here, we have a visible 3-dimensional structure. The vector space $V_{2}$ include the atoms of two types: 
\begin{equation}\label{sl(4) (0, l, 0) structure of V_2}
	V_{2} = {\rm Span}\left\{C_{2}^{a}\cdot R_{2}\,,~~ A_{1} B_{1} C_{2}^{a}\cdot R_{2} \right\} = \nu_{1, 2}\oplus\nu_{2, 2}\,.
\end{equation}
That means that we do not have a unique way to add an atom to the second node. In order to solve this problem, we utilize the Gelfand-Tsetlin bases, which, as we demonstrate, naturally arise in our description. 

\subsubsection*{Example: Representation $(0,1,0)$}

We now demonstrate the statement above in the simplest example, the representation $(0, 1, 0)$ where $\lambda = 1$. We need to distinguish the vectors in the space $V_{2}$. The number $m_{1}$ will count the atoms of the type $C_{2}^{a}R_{2}$, whereas the number $m_{2}$ will count the others:
\begin{equation}\label{sl(4) (0, l, 0) the numeration}
	\begin{aligned}
		\dim\nu_{1, 2} = \#(C_{2}^{a}R_{2}) = m_{1}\,,\\
		\dim\nu_{2, 2} = \#(A_{1}B_{1}C_{2}^{a}R_{2}) = m_{2}\,.
	\end{aligned}
\end{equation}
The formula \eqref{sl(4) (0, l, 0) structure of V_2} gives the connection between the dimension $n_{2}$ of $V_{2}$ and the parameters $m_{1},~m_{2}$:
\begin{equation}
	n_{2} = m_{1} + m_{2}\,.
\end{equation}
The cut-off takes the form $C_{2}R_{2} = 0$ in this case; therefore, our empty box is of the form:
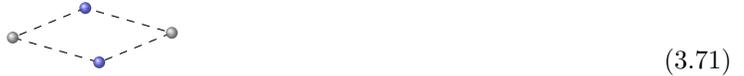
\begin{equation}
	\centering
	\begin{tikzpicture}[tdplot_main_coords, scale=1.5]
		\begin{scope}
			\coordinate (A) at (0,0,0);
			\coordinate (B) at (1,0,0);
			\coordinate (C) at (1,1,0);
			\coordinate (D) at (0,1,0);
			\draw[dashed] (A) -- (B) -- (C) -- (D) -- cycle;
			\foreach \point in {B,D}
			{
				\shade[ball color=\mygray, opacity=0.9] (\point) circle (1.5pt);
			}
			\foreach \point in {A,C}
			{
				\shade[ball color=\myblue, opacity=0.9] (\point) circle (1.5pt);
			}
		\end{scope}
	\end{tikzpicture}
\end{equation}
Next, we construct the crystals in this representation by gradually adding atoms into the box.

We start with the atom ${R_{2}}$. It corresponds to the action of the generator $e^{(2)}(z)$ of the Yangian:
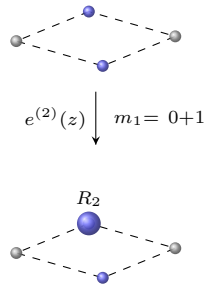
\begin{center}
	\begin{tikzpicture}[tdplot_main_coords, scale=1.5]
		\begin{scope}
			\coordinate (A) at (0,0,0);
			\coordinate (B) at (1,0,0);
			\coordinate (C) at (1,1,0);
			\coordinate (D) at (0,1,0);
			\draw[dashed] (A) -- (B) -- (C) -- (D) -- cycle;
			\foreach \point in {B,D}
			{
				\shade[ball color=\mygray, opacity=0.9] (\point) circle (1.5pt);
			}
			\foreach \point in {A,C}
			{
				\shade[ball color=\myblue, opacity=0.9] (\point) circle (1.5pt);
			}
		\end{scope}
		\begin{scope}[shift={(0,0,-2)}]
			\coordinate (A) at (0,0,0);
			\coordinate (B) at (1,0,0);
			\coordinate (C) at (1,1,0);
			\coordinate (D) at (0,1,0);
			\draw[dashed] (A) -- (B) -- (C) -- (D) -- cycle;
			\foreach \point in {B,D}
			{
				\shade[ball color=\mygray, opacity=0.9] (\point) circle (1.5pt);
			}
			\foreach \point in {A,C}
			{
				\shade[ball color=\myblue, opacity=0.9] (\point) circle (1.5pt);
			}
			\shade[ball color=\myblue, opacity=0.9] (A) circle (3.0pt);
			\node[above=3pt] at (0, 0, 0) {$\scriptstyle R_{2}$};
		\end{scope}
		\draw[postaction={decorate},
		decoration={markings, mark= at position 1 with {\arrow{stealth}}}] (0.5,0.5,-0.5) to node[left] {$\scriptstyle e^{(2)}(z)$} node[right] {$\scriptstyle ~m_{1} =~0 + 1$} (0.5,0.5,-1);
	\end{tikzpicture}
\end{center}
The path $R_{2}$ belongs to the first type of vectors in \eqref{sl(4) (0, l, 0) structure of V_2}; therefore, the parameter $m_{1}$ equals $1$ now.

Then we add a second atom. We can only add an atom to the first vector space $V_{1}$ or to the third one $V_{3}$ using the operators $e^{(1)}(z)$ or $e^{(3)}(z)$ correspondingly. This gives us two new states. These crystals take the form:
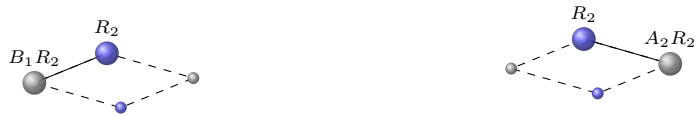
\begin{center}
	\begin{tikzpicture}[tdplot_main_coords, scale=1.5]
		\begin{scope}[shift={(1.5,-1.5,0)}]
			\coordinate (A) at (0,0,0);
			\coordinate (B) at (1,0,0);
			\coordinate (C) at (1,1,0);
			\coordinate (D) at (0,1,0);
			\draw[dashed] (A) -- (B) -- (C) -- (D) -- cycle;
			\draw (A) -- (B);
			\foreach \point in {B,D}
			{
				\shade[ball color=\mygray, opacity=0.9] (\point) circle (1.5pt);
			}
			\foreach \point in {A,C}
			{
				\shade[ball color=\myblue, opacity=0.9] (\point) circle (1.5pt);
			}
			\shade[ball color=\myblue, opacity=0.9] (A) circle (3.0pt);
			\shade[ball color=\mygray, opacity=0.9] (B) circle (3.0pt);
			\node[above=3pt] at (0, 0, 0) {$\scriptstyle R_{2}$};
			\node[above=3pt] at (1, 0, 0) {$\scriptstyle B_{1}R_{2}$};
		\end{scope}
		\begin{scope}[shift={(-1.5,1.5,0)}]
			\coordinate (A) at (0,0,0);
			\coordinate (B) at (1,0,0);
			\coordinate (C) at (1,1,0);
			\coordinate (D) at (0,1,0);
			\draw[dashed] (A) -- (B) -- (C) -- (D) -- cycle;
			\draw (A) -- (D);
			\foreach \point in {B,D}
			{
				\shade[ball color=\mygray, opacity=0.9] (\point) circle (1.5pt);
			}
			\foreach \point in {A,C}
			{
				\shade[ball color=\myblue, opacity=0.9] (\point) circle (1.5pt);
			}
			\shade[ball color=\myblue, opacity=0.9] (A) circle (3.0pt);
			\shade[ball color=\mygray, opacity=0.9] (D) circle (3.0pt);
			\node[above=3pt] at (0, 0, 0) {$\scriptstyle R_{2}$};
			\node[above=3pt] at (0, 1, 0) {$\scriptstyle A_{2}R_{2}$};
		\end{scope}
	\end{tikzpicture}
\end{center}

We can proceed further and get the remaining states. The whole structure of the states of this representation is depicted in fig. \ref{sl(4) (0, 1, 0) representation}. We emphasize that the dimension of the vector space $V_{2}$ \eqref{sl(4) (0, l, 0) structure of V_2} is shifted twice. The first shift corresponds to the parameter $m_{1}$, whereas the second shift corresponds to the parameter $m_{2}$. Nonetheless, both shifts are performed by the generator $e^{(2)}(z)$.
\begin{figure}[h]
	\centering
	\begin{tikzpicture}[tdplot_main_coords, scale=1.5]
		\begin{scope}[shift={(0,0,2)}]
			\coordinate (A) at (0,0,0);
			\coordinate (B) at (1,0,0);
			\coordinate (C) at (1,1,0);
			\coordinate (D) at (0,1,0);
			\draw[dashed] (A) -- (B) -- (C) -- (D) -- cycle;
			\foreach \point in {B,D}
			{
				\shade[ball color=\mygray, opacity=0.9] (\point) circle (1.5pt);
			}
			\foreach \point in {A,C}
			{
				\shade[ball color=\myblue, opacity=0.9] (\point) circle (1.5pt);
			}
		\end{scope}
		\begin{scope}[shift={(0,0,0)}]
			\coordinate (A) at (0,0,0);
			\coordinate (B) at (1,0,0);
			\coordinate (C) at (1,1,0);
			\coordinate (D) at (0,1,0);
			\draw[dashed] (A) -- (B) -- (C) -- (D) -- cycle;
			\foreach \point in {B,D}
			{
				\shade[ball color=\mygray, opacity=0.9] (\point) circle (1.5pt);
			}
			\foreach \point in {A,C}
			{
				\shade[ball color=\myblue, opacity=0.9] (\point) circle (1.5pt);
			}
			\shade[ball color=\myblue, opacity=0.9] (A) circle (3.0pt);
			\node[above=3pt] at (0, 0, 0) {$\scriptstyle R_{2}$};
		\end{scope}
		\begin{scope}[shift={(1.5,-1.5,-1.5)}]
			\coordinate (A) at (0,0,0);
			\coordinate (B) at (1,0,0);
			\coordinate (C) at (1,1,0);
			\coordinate (D) at (0,1,0);
			\draw[dashed] (A) -- (B) -- (C) -- (D) -- cycle;
			\draw (A) -- (B);
			\foreach \point in {B,D}
			{
				\shade[ball color=\mygray, opacity=0.9] (\point) circle (1.5pt);
			}
			\foreach \point in {A,C}
			{
				\shade[ball color=\myblue, opacity=0.9] (\point) circle (1.5pt);
			}
			\shade[ball color=\myblue, opacity=0.9] (A) circle (3.0pt);
			\shade[ball color=\mygray, opacity=0.9] (B) circle (3.0pt);
			\node[above=3pt] at (0, 0, 0) {$\scriptstyle R_{2}$};
			\node[above=3pt] at (1, 0, 0) {$\scriptstyle B_{1}R_{2}$};
		\end{scope}
		\begin{scope}[shift={(-1.5,1.5,-1.5)}]
			\coordinate (A) at (0,0,0);
			\coordinate (B) at (1,0,0);
			\coordinate (C) at (1,1,0);
			\coordinate (D) at (0,1,0);
			\draw[dashed] (A) -- (B) -- (C) -- (D) -- cycle;
			\draw (A) -- (D);
			\foreach \point in {B,D}
			{
				\shade[ball color=\mygray, opacity=0.9] (\point) circle (1.5pt);
			}
			\foreach \point in {A,C}
			{
				\shade[ball color=\myblue, opacity=0.9] (\point) circle (1.5pt);
			}
			\shade[ball color=\myblue, opacity=0.9] (A) circle (3.0pt);
			\shade[ball color=\mygray, opacity=0.9] (D) circle (3.0pt);
			\node[above=3pt] at (0, 0, 0) {$\scriptstyle R_{2}$};
			\node[above=3pt] at (0, 1, 0) {$\scriptstyle A_{2}R_{2}$};
		\end{scope}
		\begin{scope}[shift={(0,0,-3)}]
			\coordinate (A) at (0,0,0);
			\coordinate (B) at (1,0,0);
			\coordinate (C) at (1,1,0);
			\coordinate (D) at (0,1,0);
			\draw[dashed] (A) -- (B) -- (C) -- (D) -- cycle;
			\draw (A) -- (B) (A) -- (D);
			\foreach \point in {B,D}
			{
				\shade[ball color=\mygray, opacity=0.9] (\point) circle (1.5pt);
			}
			\foreach \point in {A,C}
			{
				\shade[ball color=\myblue, opacity=0.9] (\point) circle (1.5pt);
			}
			\shade[ball color=\myblue, opacity=0.9] (A) circle (3.0pt);
			\shade[ball color=\mygray, opacity=0.9] (B) circle (3.0pt);
			\shade[ball color=\mygray, opacity=0.9] (D) circle (3.0pt);
			\node[above=3pt] at (0, 0, 0) {$\scriptstyle R_{2}$};
			\node[above=3pt] at (0, 1, 0) {$\scriptstyle A_{2}R_{2}$};
			\node[above=3pt] at (1, 0, 0) {$\scriptstyle B_{1}R_{2}$};
		\end{scope}
		\begin{scope}[shift={(0,0,-5)}]
			\coordinate (A) at (0,0,0);
			\coordinate (B) at (1,0,0);
			\coordinate (C) at (1,1,0);
			\coordinate (D) at (0,1,0);
			\draw  (A) -- (B) -- (C) -- (D) -- cycle;
			\foreach \point in {B,D}
			{
				\shade[ball color=\mygray, opacity=0.9] (\point) circle (3.0pt);
			}
			\foreach \point in {A,C}
			{
				\shade[ball color=\myblue, opacity=0.9] (\point) circle (3.0pt);
			}
			\node[above=3pt] at (0, 0, 0) {$\scriptstyle R_{2}$};
			\node[above=3pt] at (0, 1, 0) {$\scriptstyle A_{2}R_{2}$};
			\node[above=3pt] at (1, 0, 0) {$\scriptstyle B_{1}R_{2}$};
			\node[below=3pt] at (1, 1, 0) {$\scriptstyle A_{1}B_{1}R_{2}$};
		\end{scope}
		\draw[postaction={decorate},
		decoration={markings, mark= at position 1 with {\arrow{stealth}}}] (0.5,0.5,1.5) to node[left] {$\scriptstyle e^{(2)}(z)$} node[right] {$\scriptstyle ~m_{1} =~0 + 1$} (0.5,0.5,1);
		\draw[postaction={decorate},
		decoration={markings, mark= at position 1 with {\arrow{stealth}}}] (0.5,-0.5,-1) to node[above left] {$\scriptstyle e^{(1)}(z)$} node[below right] {$\scriptstyle ~n_{1} =~0 + 1$} (1.5,-0.5,-1);
		\draw[postaction={decorate},
		decoration={markings, mark= at position 1 with {\arrow{stealth}}}] (-0.5,0.5,-1) to node[above right] {$\scriptstyle e^{(3)}(z)$} node[below left] {$\scriptstyle ~n_{3} =~0 + 1$} (-0.5,1.5,-1);
		\draw[postaction={decorate},
		decoration={markings, mark= at position 1 with {\arrow{stealth}}}] (1.5,-0.5,-2) to node[below left] {$\scriptstyle e^{(3)}(z)$} node[above right] {$\scriptstyle ~n_{3} =~0 + 1$} (1.5,0.5,-2);
		\draw[postaction={decorate},
		decoration={markings, mark= at position 1 with {\arrow{stealth}}}] (-0.5,1.5,-2) to node[below right] {$\scriptstyle e^{(1)}(z)$} node[above left] {$\scriptstyle ~n_{1} =~0 + 1$} (0.5,1.5,-2);
		\draw[postaction={decorate},
		decoration={markings, mark= at position 1 with {\arrow{stealth}}}] (0.5,0.5,-3.5) to node[left] {$\scriptstyle e^{(2)}(z)$} node[right] {$\scriptstyle ~m_{2} =~0 + 1$} (0.5,0.5,-4);
	\end{tikzpicture}
	\caption{The structure of the representation $\Upsilon_{2, 1}$ of the algebra $\myY(\fs\fl_{4})$}\label{sl(4) (0, 1, 0) representation}
\end{figure}

\subsubsection*{Multiplicities and Gelfand-Tsetlin bases}

In this section we parametrize the states of the representations $\Upsilon_{2,\lambda}$ and emphasize the importance of Gelfand-Tsetlin bases in our description. In all previous examples\footnote{Including also the representation $(0, 1, 0)$.}, we were able to assign the dimension vector $\vec{n}$ to a particular state of the representation. However, this is no longer true for the representations $\Upsilon_{2,\lambda}$, which means that a crystal cannot be parameterized just by three numbers $n_{1}, n_{2}, n_{3}$.

This problem is not new in the theory of quiver Yangian algebras. For example, the quiver for $\myY(\widehat{\fg\fl}_{1})$ has a single node with the dimension parameter $n$. It is a well-known fact \cite{rapcak2021branesquiversbpsalgebras, Galakhov:2020vyb, Galakhov:2022uyu} that there are $p(n)$ corresponding states in the Fock representation, where $p(n)$ is the number of Young diagrams of size $n$.

Notably, this is also tightly related to a classical problem in the representation theory of Lie algebras. Having a highest weight state $|\vec{\lambda}\rangle$, all the states can be acquired by the action of the raising operators of an algebra:
\begin{equation}
	\text{States: }\quad |\vec{\mu}\rangle = \prod_{i}e_{i}|\vec{\lambda}\rangle \quad \leftrightarrow \quad \vec{\mu} = \vec{\lambda} - \sum_{i}n_{i}\vec{\alpha}_{i}\,,
\end{equation}
where $\vec{\alpha}_{i}$ are simple roots. The numbers $n_{i}$ resemble the dimension vector $\vec{n}$ of the quiver by construction. The weights $\vec{\mu}$ form so-called weight system. Although the approach allows to systematically construct all the weights in the representation, it does not track multiplicities. This means that we can have different states $|\vec{\mu}_{1}\rangle \ne |\vec{\mu}_{2}\rangle$ with $\vec{\mu}_{1} = \vec{\mu}_{2}$. Depending on the needs of a researcher, there are various ways to overcome this problem. The one we are interested in is introducing the Gelfand-Tsetlin bases \cite{Gelfand:1950ihs} as the parametrization of the states.

In order to connect this to the crystal states, we start with the simplest example, where we encounter non-trivial multiplicities in our construction. Namely, we consider the particular states in the representation $(0, 2, 0)$.

The empty box is fixed by the relation $C_{2}^{2}R_{2} = 0$. We do not unwrap the whole set of states in the representation due to its high dimension. Instead, let us focus on the state $e^{(3)}(z_{3})e^{(1)}(z_{1})e^{(2)}(z_{2})|\varnothing\rangle$. Its dimension vector is $\vec{n} = (1, 1, 1)$, and the crystal diagram takes the following form:
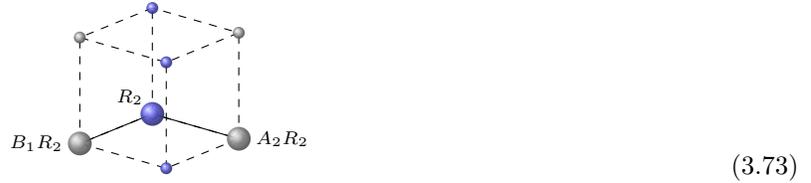
\begin{equation}\label{sl(4) (0, 2, 0) a state}
	\begin{tikzpicture}[tdplot_main_coords, scale=1.5]
		\begin{scope}[shift={(0,0,0)}]
			\coordinate (A) at (0,0,0);
			\coordinate (B) at (1,0,0);
			\coordinate (C) at (1,1,0);
			\coordinate (D) at (0,1,0);
			\coordinate (A') at (0,0,1);
			\coordinate (B') at (1,0,1);
			\coordinate (C') at (1,1,1);
			\coordinate (D') at (0,1,1);
			\draw[dashed] (A) -- (B) -- (C) -- (D) -- cycle;
			\draw[dashed] (A') -- (B') -- (C') -- (D') -- cycle;
			\draw[dashed] (A) -- (A') (B) -- (B') (D) -- (D') (C) -- (C');
			\draw (A) -- (B) (A) -- (D);
			\foreach \point in {B',D'}
			{
				\shade[ball color=\mygray, opacity=0.9] (\point) circle (1.5pt);
			}
			\foreach \point in {A',C'}
			{
				\shade[ball color=\myblue, opacity=0.9] (\point) circle (1.5pt);
			}
			\foreach \point in {B,D}
			{
				\shade[ball color=\mygray, opacity=0.9] (\point) circle (1.5pt);
			}
			\foreach \point in {A,C}
			{
				\shade[ball color=\myblue, opacity=0.9] (\point) circle (1.5pt);
			}
			\shade[ball color=\myblue, opacity=0.9] (A) circle (3.0pt);
			\shade[ball color=\mygray, opacity=0.9] (B) circle (3.0pt);
			\shade[ball color=\mygray, opacity=0.9] (D) circle (3.0pt);
			\node[above left] at (0, 0, 0) {$\scriptstyle R_{2}$};
			\node[right=3pt] at (0, 1, 0) {$\scriptstyle A_{2}R_{2}$};
			\node[left=3pt] at (1, 0, 0) {$\scriptstyle B_{1}R_{2}$};
		\end{scope}
	\end{tikzpicture}
\end{equation}
We are interested in expanding this crystal to get new states. In fact, there are two possibilities that we highlight in red for a moment:
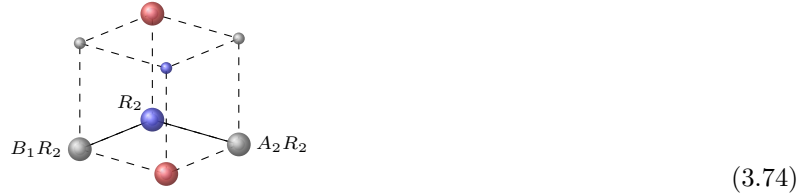
\begin{equation}\label{sl(4) (0, 2, 0) a state posiibilities}
	\begin{tikzpicture}[tdplot_main_coords, scale=1.5]
		\begin{scope}[shift={(0,0,0)}]
			\coordinate (A) at (0,0,0);
			\coordinate (B) at (1,0,0);
			\coordinate (C) at (1,1,0);
			\coordinate (D) at (0,1,0);
			\coordinate (A') at (0,0,1);
			\coordinate (B') at (1,0,1);
			\coordinate (C') at (1,1,1);
			\coordinate (D') at (0,1,1);
			\draw[dashed] (A) -- (B) -- (C) -- (D) -- cycle;
			\draw[dashed] (A') -- (B') -- (C') -- (D') -- cycle;
			\draw[dashed] (A) -- (A') (B) -- (B') (D) -- (D') (C) -- (C');
			\draw (A) -- (B) (A) -- (D);
			\foreach \point in {B',D'}
			{
				\shade[ball color=\mygray, opacity=0.9] (\point) circle (1.5pt);
			}
			\foreach \point in {A',C'}
			{
				\shade[ball color=\myblue, opacity=0.9] (\point) circle (1.5pt);
			}
			\foreach \point in {B,D}
			{
				\shade[ball color=\mygray, opacity=0.9] (\point) circle (1.5pt);
			}
			\foreach \point in {A,C}
			{
				\shade[ball color=\myblue, opacity=0.9] (\point) circle (1.5pt);
			}
			\shade[ball color=\myblue, opacity=0.9] (A) circle (3.0pt);
			\shade[ball color=\mygray, opacity=0.9] (B) circle (3.0pt);
			\shade[ball color=\mygray, opacity=0.9] (D) circle (3.0pt);
			\shade[ball color=\myred, opacity=0.9] (C) circle (3.0pt);
			\shade[ball color=\myred, opacity=0.9] (A') circle (3.0pt);
			\node[above left] at (0, 0, 0) {$\scriptstyle R_{2}$};
			\node[right=3pt] at (0, 1, 0) {$\scriptstyle A_{2}R_{2}$};
			\node[left=3pt] at (1, 0, 0) {$\scriptstyle B_{1}R_{2}$};
		\end{scope}
	\end{tikzpicture}
\end{equation}
Adding the atoms of the colors of the first and the third nodes is prohibited due to the F-terms. The direct calculation shows that $C_{1}B_{1}R_{2}$ indeed acquires zero expectation:
\begin{equation}
	\begin{aligned}
		&C_{1}B_{1}R_{2} = B_{1}C_{2}R_{2} = B_{1}(C_{2}R_{2}) = 0\,,
	\end{aligned}
\end{equation}
where we used the fact that $C_{2}R_{2} = 0$ because we do not have this atom in the crystal \eqref{sl(4) (0, 2, 0) a state}.

Now, we focus on the two states that we get from \eqref{sl(4) (0, 2, 0) a state}. Both crystals are obtained by the action of the generator $e^{(2)}(z)$. Moreover, the dimension vector of both states equals $\vec{n} = (1, 2, 1)$.
\begin{enumerate}
	\item We get the first one by shifting the parameter $m_{1}$ from $1$ to $2$. The crystal diagram takes the form:
	
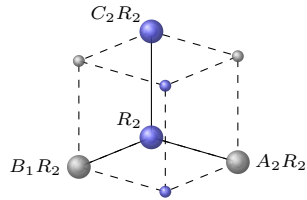
\begin{center}
		\begin{tikzpicture}[tdplot_main_coords, scale=1.5]
			\begin{scope}
				\coordinate (A) at (0,0,0);
				\coordinate (B) at (1,0,0);
				\coordinate (C) at (1,1,0);
				\coordinate (D) at (0,1,0);
				\coordinate (A') at (0,0,1);
				\coordinate (B') at (1,0,1);
				\coordinate (C') at (1,1,1);
				\coordinate (D') at (0,1,1);
				\draw[dashed]  (A) -- (B) -- (C) -- (D) -- cycle;
				\draw[dashed] (A') -- (B') -- (C') -- (D') -- cycle;
				\draw (A) -- (B) (A) -- (D) (A) -- (A');
				\draw[dashed] (B) -- (B') (D) -- (D') (C) -- (C');
				\foreach \point in {B,D,D',B'}
				{
					\shade[ball color=\mygray, opacity=0.9] (\point) circle (1.5pt);
				}
				\foreach \point in {A,C,A',C'}
				{
					\shade[ball color=\myblue, opacity=0.9] (\point) circle (1.5pt);
				}
				\shade[ball color=\myblue, opacity=0.9] (A) circle (3.0pt);
				\shade[ball color=\mygray, opacity=0.9] (B) circle (3.0pt);
				\shade[ball color=\mygray, opacity=0.9] (D) circle (3.0pt);
				\shade[ball color=\myblue, opacity=0.9] (A') circle (3.0pt);
			\end{scope}
			\node[above left] at (0, 0, 0) {$\scriptstyle R_{2}$};
			\node[above left] at (0, 0, 1) {$\scriptstyle C_{2}R_{2}$};
			\node[right=3pt] at (0, 1, 0) {$\scriptstyle A_{2}R_{2}$};
			\node[left=3pt] at (1, 0, 0) {$\scriptstyle B_{1}R_{2}$};
		\end{tikzpicture}
	\end{center}
	The corresponding vacuum expectation values of the fields read as follows:
	\begin{equation}
		\begin{aligned}
			&C_{1} = C_{3} = (0)\,, \\
			&C_{2} = \begin{pmatrix}
				0 & 0 \\ 1 & 0
			\end{pmatrix}\,,\\
			&A_{1} = (0, 0)^{T}, ~ A_{2} = (1, 0)\,,\\
			&B_{1} = (1, 0)\,, ~ B_{2} = (0, 0)^{T}\,,\\
			&R_{2} = (1, 0)^{T}\,, ~ S_{2} = (0, 0)\,.
		\end{aligned}
	\end{equation}
	\item By shifting the parameter $m_{2}$ from $0$ to $1$ we acquire the second state. Its diagram takes the form:
	
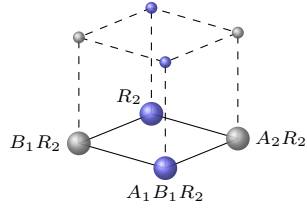
\begin{center}
		\begin{tikzpicture}[tdplot_main_coords, scale=1.5]
			\begin{scope}
				\coordinate (A) at (0,0,0);
				\coordinate (B) at (1,0,0);
				\coordinate (C) at (1,1,0);
				\coordinate (D) at (0,1,0);
				\coordinate (A') at (0,0,1);
				\coordinate (B') at (1,0,1);
				\coordinate (C') at (1,1,1);
				\coordinate (D') at (0,1,1);
				\draw  (A) -- (B) -- (C) -- (D) -- cycle;
				\draw[dashed] (A') -- (B') -- (C') -- (D') -- cycle;
				\draw[dashed] (A) -- (A') (B) -- (B') (D) -- (D') (C) -- (C');
				\foreach \point in {B',D'}
				{
					\shade[ball color=\mygray, opacity=0.9] (\point) circle (1.5pt);
				}
				\foreach \point in {A',C'}
				{
					\shade[ball color=\myblue, opacity=0.9] (\point) circle (1.5pt);
				}
				\foreach \point in {B,D}
				{
					\shade[ball color=\mygray, opacity=0.9] (\point) circle (3pt);
				}
				\foreach \point in {A,C}
				{
					\shade[ball color=\myblue, opacity=0.9] (\point) circle (3pt);
				}
				\node[above left] at (0, 0, 0) {$\scriptstyle R_{2}$};
				\node[right=3pt] at (0, 1, 0) {$\scriptstyle A_{2}R_{2}$};
				\node[left=3pt] at (1, 0, 0) {$\scriptstyle B_{1}R_{2}$};
				\node[below=3pt] at (1, 1, 0) {$\scriptstyle A_{1}B_{1}R_{2}$}; 
			\end{scope}
		\end{tikzpicture}
	\end{center}
	The fields acquire the following vacuum expectation values:
	\begin{equation}
		\begin{aligned}
			&C_{1} = C_{3} = (0)\,, \\
			&C_{2} = \begin{pmatrix}
				0 & 0 \\ 0 & 0
			\end{pmatrix}\,,\\
			&A_{1} = (0, 0)^{T}\,, ~ A_{2} = (0, 0)\,,\\
			&B_{1} = (1, 0)\,, ~ B_{2} = (0, 0)^{T}\,,\\
			&R_{2} = (1, 0)^{T}\,, ~ S_{2} = (0)\,.
		\end{aligned}
	\end{equation}
\end{enumerate}

We just provided the example of the representation, where one cannot parametrize the states by the dimension vector $\vec{n}$. We now lift the construction to an arbitrary $\Upsilon_{2,\lambda}$ representation. Effectively, we have four vector spaces:
\begin{equation}
	\{V_{1},~V_{2},~V_{3}\} = \{V_{1}, \nu_{1, 2}, \nu_{2, 2}, V_{3}\}\,.
\end{equation}
We suppose that their corresponding dimensions $n_{1},~m_{1},~m_{2},~n_{3}$ parametrize the crystal in the representation. We encode this in the following form:
\begin{equation}\label{sl(4) (0, l, 0) Gelfand-Tsetlin}
	\mu =\quad
	\begin{tikzpicture}[baseline={(0,-0.0)}, every node/.style={inner sep=1pt}]
		\GTrow[1.2]{$\lambda$, $\lambda$, $0$, $0$}
		\GTrow[0.6]{$\lambda$, $n_{3}$, $0$}
		\GTrow[0]{$m_{1}$, $m_{2}$}
		\GTrow[-0.6]{$n_{1}$}
	\end{tikzpicture}\,,
\end{equation}
where $n_{2} = m_{1} + m_{2}$. This resembles the Gelfand-Tsetlin bases for the representation $(0, \lambda, 0)$.
The generators $e^{(1)}(z)$ and $e^{(3)}(z)$, or their lowering analogues, shift the dimensions $n_{1},~n_{3}$. We denote the  corresponding patterns as follows:
\begin{equation}\label{sl(4) (0, l, 0) the patterns n1, n3}
	\mu_{3, 3 \pm 1}^{2} =\quad
	\begin{tikzpicture}[baseline={(0,-0.0)}, every node/.style={inner sep=1pt}]
		\GTrow[1.2]{$\lambda$, $\lambda$, $0$, $0$}
		\GTrow[0.6]{$\lambda$, $n_{3}\pm 1$, $0$}
		\GTrow[0]{$m_{1}$, $m_{2}$}
		\GTrow[-0.6]{$n_{1}$}
	\end{tikzpicture}\,, \qquad
	\mu_{1, 1 \pm 1}^{1} \quad = \begin{tikzpicture}[baseline={(0,-0.0)}, every node/.style={inner sep=1pt}]
		\GTrow[1.2]{$\lambda$, $\lambda$, $0$, $0$}
		\GTrow[0.6]{$\lambda$, $n_{3}$, $0$}
		\GTrow[0]{$m_{1}$, $m_{2}$}
		\GTrow[-0.6]{$n_{1} \pm 1$}
	\end{tikzpicture}\,.
\end{equation}
As for the second node, the operator $e^{(2)}(z)$, or $f^{(2)}(z)$, can shift the parameters $m_{1}$ and $m_{2}$. These options we encode in the following form:
\begin{equation}\label{sl(4) (0, l, 0) the patterns n2}
	\mu_{2, 2 \pm 1}^{1} =\quad
	\begin{tikzpicture}[baseline={(0,-0.0)}, every node/.style={inner sep=1pt}]
		\GTrow[1.2]{$\lambda$, $\lambda$, $0$, $0$}
		\GTrow[0.6]{$\lambda$, $n_{3}$, $0$}
		\GTrow[0]{$m_{1} \pm 1$, $m_{2}$}
		\GTrow[-0.6]{$n_{1}$}
	\end{tikzpicture}\,, \qquad
	\mu_{2, 2 \pm 1}^{2} \quad= \begin{tikzpicture}[baseline={(0,-0.0)}, every node/.style={inner sep=1pt}]
		\GTrow[1.2]{$\lambda$, $\lambda$, $0$, $0$}
		\GTrow[0.6]{$\lambda$, $n_{3}$, $0$}
		\GTrow[0]{$m_{1}$, $m_{2} \pm 1$}
		\GTrow[-0.6]{$n_{1}$}
	\end{tikzpicture}\,.
\end{equation}
One might find the notation of the patterns \eqref{sl(4) (0, l, 0) the patterns n1, n3} and \eqref{sl(4) (0, l, 0) the patterns n2} confusing at this moment. The symbol $\mu_{k, k \pm 1}^{j}$ encodes the fact that we shift the dimension of the $j$-th element of the $k$-th row. For our needs we sometimes stack this notation to encode the fact that we shift two dimensions. For example:
\begin{equation}
	\bigl(\mu_{2, 2 + 1}^{1}\bigr)_{1, 1 + 1}^{1} = \bigl(\mu_{1, 1 + 1}^{1}\bigr)_{2, 2 + 1}^{1} = \quad
	\begin{tikzpicture}[baseline={(0,-0.0)}, every node/.style={inner sep=1pt}]
		\GTrow[1.2]{$\lambda$, $\lambda$, $0$, $0$}
		\GTrow[0.6]{$\lambda$, $n_{3}$, $0$}
		\GTrow[0]{$m_{1} + 1$, $m_{2}$}
		\GTrow[-0.6]{$n_{1} + 1$}
	\end{tikzpicture}\,.
\end{equation}

\subsubsection*{State Counting}

Strictly speaking, the introduced patterns \eqref{sl(4) (0, l, 0) Gelfand-Tsetlin} have not proven to be the Gelfand-Tsetlin bases yet. In order to associate them with each other, one should verify the triangular inequalities:
\begin{equation}\label{sl(4) (0, l, 0) triangular inequalities}
	\begin{aligned}
		&0 \leqslant m_{2} \leqslant n_{3}\,, \quad n_{3} \leqslant m_{1} \leqslant \lambda\,,
		&m_{2} \leqslant n_{1} \leqslant m_{1}\,.
	\end{aligned}
\end{equation}
Note that by construction of the crystal box, all the dimensions are limited by $\lambda$:
\begin{equation}
	\begin{aligned}
		&0 \leqslant n_{1} \leqslant \lambda\,,\quad 0 \leqslant n_{3} \leqslant \lambda\,,
		&0 \leqslant m_{1} \leqslant \lambda\,, \quad 0 \leqslant m_{2} \leqslant \lambda\,,
	\end{aligned}
\end{equation}
which means that there are only four inequalities in \eqref{sl(4) (0, l, 0) triangular inequalities} that we need to prove. They can be organized as follows:
\begin{equation}
	m_{2} \leqslant n_{1} \leqslant m_{2}\,, \quad m_{2} \leqslant n_{3} \leqslant m_{2}\,.
\end{equation}
In fact, when we prove the first sequence, we automatically prove the second due to the quiver symmetry.

We start with the inequality $n_{1} \leqslant m_{1}$. Let us assume the contrary. It means there is an atom for which:
\begin{equation}
	C_{1}^{m}B_{1}R_{2} \ne 0\,, \quad C_{2}^{m}R_{2} = 0\,.
\end{equation}
However, we have equivalence relations \eqref{sl(4) F-terms}. Using them, we get:
\begin{equation}
	C_{1}^{m}B_{1}R_{2} = B_{1}C_{2}^{m}R_{2} = B_{1}(C_{2}^{m}R_{2}) = 0\,,
\end{equation}
which was not true in our suggestion. Therefore, $n_{1} \leqslant m_{1}$.

The second inequality can be proven in a similar way. The key calculation reads as follows:
\begin{equation}
	C_{2}^{m}A_{1}B_{1}R_{2} = A_{1}C_{1}^{m}B_{1}R_{2} = A_{1}(C_{1}^{m}B_{1}R_{2}) = 0\,.
\end{equation}
and shows that $m_{2} \leqslant n_{1}$. It concludes the proof of both sequences of inequalities and verifies that the patterns \eqref{sl(4) (0, l, 0) Gelfand-Tsetlin} are indeed the Gelfand-Tsetlin bases.

Next, we count the dimension of our crystal representation. The inequalities \eqref{sl(4) (0, l, 0) triangular inequalities} limit our dimensions. Without loss of generality, we assume also that $n_{1} \geqslant n_{3}$ since the quiver \ref{quiver sl(4) (0, l, 0)} has $\mathbb{Z}_{2}$ symmetry \eqref{sl(n) quiver symmetry Z_2}. This leads us to the sum:
\begin{equation}
	\dim_{f.p.}|(0, \lambda, 0)| = \sum_{m_{2} = 0}^{\lambda}\sum_{n_{3} = m_{2}}^{\lambda}\sum_{n_{1} = m_{2}}^{\lambda}\sum_{m_{1} = n_{1}}^{\lambda} 1 = \sum_{m_{2} = 0}^{\lambda}\sum_{n_{3} = m_{2}}^{\lambda}\sum_{n_{1} = m_{2}}^{\lambda}(\lambda - n_{1} + 1) = \dots\,.
\end{equation}
After the calculation, we end up with the result:
\begin{equation}\label{sl(4) (0, l, 0) dimension}
	 \dim_{f.p.}|(0, \lambda, 0)| = \dfrac{(\lambda + 1)(\lambda + 2)(\lambda + 2)(\lambda + 3)}{3!2!} =\dim_{\fs\fl_{4}}|(0, \lambda, 0)|\,,
\end{equation}
which is exactly the dimension of the $(0, \lambda, 0)$ representation of the algebra $\fs\fl_{4}$.

\subsubsection*{The representation}

We need to adjust the ansatz \eqref{General Ansatz} to the representations $\Upsilon_{2,\lambda}$. Our states are now denoted as Gelfand-Tsetlin bases \eqref{sl(4) (0, l, 0) Gelfand-Tsetlin}. The operators corresponding to the first node and the third node transform them into the bases \eqref{sl(4) (0, l, 0) the patterns n1, n3}. As for the operators $e^{(2)}(z)$ and $f^{(2)}(z)$, their action is now split in two and gives the patterns \eqref{sl(4) (0, l, 0) the patterns n2}. Therefore, we end up with the ansatz of the following form:
\begin{equation}\label{sl(4) (0, l, 0) ansatz}
	\begin{aligned}
		&e^{(1)}_{\Upsilon_{2,\lambda}}(z)|\mu\rangle = \dfrac{{\bf E}_{\Upsilon_{2,\lambda}}[\mu \to \mu_{1, 1 + 1}^{1}]}{z + \frac{\epsilon}{2} - n_{1}\epsilon}|\mu_{1, 1 + 1}^{1}\rangle\,, \\
		&e^{(2)}_{\Upsilon_{2,\lambda}}(z)|\mu\rangle = \dfrac{{\bf E}_{\Upsilon_{2,\lambda}}[\mu\to \mu_{2, 2 + 1}^{1}]}{z - m_{1}\epsilon}|\mu_{2, 2 + 1}^{1}\rangle + \dfrac{{\bf E}_{\Upsilon_{2,\lambda}}[\mu\to \mu_{2, 2 + 1}^{2}]}{z + \epsilon - m_{2}\epsilon}|\mu_{2, 2 + 1}^{2}\rangle\,,\\
		&e^{(3)}_{\Upsilon_{2,\lambda}}(z)|\mu\rangle = \dfrac{{\bf E}_{\Upsilon_{2,\lambda}}[\mu \to \mu_{3, 3 + 1}^{2}]}{z + \frac{\epsilon}{2} - n_{3}\epsilon}|\mu_{3, 3 + 1}^{2}\rangle\,,\\
		&f^{(1)}_{\Upsilon_{2,\lambda}}(z)|\mu\rangle = \dfrac{{\bf F}_{\Upsilon_{2,\lambda}}[\mu \to \mu_{1, 1 - 1}^{1}]}{z + \frac{\epsilon}{2} - (n_{1} - 1)\epsilon}|\mu_{1, 1 - 1}^{1}\rangle\,, \\
		&f^{(2)}_{\Upsilon_{2,\lambda}}(z)|\mu\rangle = \dfrac{{\bf F}_{\Upsilon_{2,\lambda}}[\mu \to \mu_{2, 2, -1}^{1}]}{z - (m_{1} - 1)\epsilon}|\mu_{2, 2-1}^{1}\rangle + \dfrac{{\bf F}_{\Upsilon_{2,\lambda}}[\mu \to \mu_{2, 2 -1}^{2}]}{z + \epsilon - (m_{2} - 1)\epsilon}|\mu_{2, 2-1}^{2}\rangle\,,\\
		&f^{(3)}_{\Upsilon_{2,\lambda}}(z)|\mu\rangle = \dfrac{{\bf F}_{\Upsilon_{2,\lambda}}[\mu \to \mu_{3, 3, -1}^{2}]}{z + \frac{\epsilon}{2} - (n_{3} - 1)\epsilon}|\mu_{3, 3 - 1}^{2}\rangle,\\
		&\psi^{(s)}_{\Upsilon_{2,\lambda}}(z)|\mu\rangle = \Psi^{(s)}_{\mu, \Upsilon_{2,\lambda}}(z)|\mu\rangle\,.
	\end{aligned}
\end{equation}
The eigenfunctions of the operators $\psi^{(s)}_{\Upsilon_{2,\lambda}}(z)$ read:
\begin{equation}\label{sl(4) (0, l, 0) eigenfunctions}
	\begin{aligned}
		\Psi^{(1)}_{\mu, \Upsilon_{2,\lambda}}(z) &= -\frac{1}{\epsilon}\prod_{k = 1}^{n_{1}}\varphi_{1, 1}\Bigl(z + \frac{\epsilon}{2} - (k-1)\epsilon\Bigr)\prod_{k = 1}^{m_{1}}\varphi_{1, 2}\Bigl(z - (k-1)\epsilon\Bigr)\prod_{k = 1}^{m_{2}}\Bigl(z + \epsilon - (k-1)\epsilon\Bigr) =\\ 
		&= -\dfrac{1}{\epsilon}\dfrac{\Bigl(z + \dfrac{\epsilon}{2} - m_{1}\epsilon\Bigr)\Bigl(z + \frac{\epsilon}{2} - (m_{2} - 1)\epsilon\Bigr)}{\Bigl(z + \dfrac{\epsilon}{2}- n_{1}\epsilon\Bigr)\Bigl(z + \frac{
			\epsilon}{2}- (n_{1} - 1)\epsilon\Bigr)}\,, \\
		\Psi^{(2)}_{\mu, \Upsilon_{2,\lambda}}(z) &= -\frac{1}{\epsilon}\dfrac{(z - \lambda\epsilon)}{z}\prod_{k = 1}^{n_{1}}\varphi_{2, 1}\Bigl(z + \frac{\epsilon}{2} - (k - 1)\epsilon\Bigr)\prod_{k = 1}^{m_{1}}\varphi_{2, 2}\Bigl(z - (k - 1)\epsilon\Bigr)\prod_{k = 1}^{m_{2}}\varphi_{2, 2}\Bigl(z + \epsilon - (k-1)\epsilon\Bigr)\cdot \\
		&\cdot\prod_{k = 1}^{n_{3}}\varphi_{2, 3}\Bigl(z + \dfrac{\epsilon}{2} - (k-1)\epsilon\Bigr) = -\frac{1}{\epsilon}\frac{(z - \lambda \epsilon)(z + 2\epsilon)(z + \epsilon - n_{1}\epsilon)(z + \epsilon - n_{3}\epsilon)}{(z - m_{1}\epsilon)(z - (m_{1} - 1)\epsilon)(z + \epsilon - m_{2}\epsilon)(z + \epsilon - (m_{2} - 1)\epsilon)}\,, \\
		\Psi^{(3)}_{\mu, \Upsilon_{2,\lambda}}(z) &= \frac{1}{\epsilon}\prod_{k = 1}^{n_{2}}\varphi_{3, 2}\Bigl(z + \frac{\epsilon}{2} - (k-1)\epsilon\Bigr)\prod_{k = 1}^{n_{3}}\varphi_{3, 3}\Bigl(z + \epsilon - (k - 1)\epsilon\Bigr) = \\
		&= -\frac{1}{\epsilon}\frac{\Bigl(z + \dfrac{\epsilon}{2} - m_{1}\epsilon\Bigr)\Bigl(z + \dfrac{\epsilon}{2} - (m_{2} - 1)\epsilon\Bigr)}{\Bigl(z + \dfrac{\epsilon}{2} - n_{3}\epsilon\Bigr)\Bigl(z + \dfrac{\epsilon}{2} - (n_{3} - 1)\epsilon\Bigr)}\,.
	\end{aligned}
\end{equation}
The function $\Psi_{\mu, \Upsilon_{2,\lambda}}^{(2)}(z)$ now has four poles, which reflects the fact that there are more possibilities to add or remove an atom from the crystal $|\mu\rangle$. 

\subsubsection*{Amplitudes}

In order to determine the coefficients in the ansatz \eqref{sl(4) (0, l, 0) ansatz}, we need to construct the vacuum expectation values of the fields that correspond to a fixed point in \eqref{sl(4) F-terms}. As we have mentioned, the space $V_{2}$ consists of two parts, $\nu_{1, 2}$ and $\nu_{2, 2}$. We remind the reader that the matrices $A_{1}, B_{1}, A_{2}, B_{2}, C_{2}$ are linear maps of the form:
\begin{equation}\label{sl(4) (0, l, 0) linear maps involving V2}
	\begin{aligned}
		&A_{1} \colon V_{1} \to V_{2}\,, \qquad B_{1} \colon V_{2} \to V_{1}\,,\\
		&A_{2} \colon V_{2} \to V_{3}\,, \qquad B_{2} \colon V_{3} \to V_{2}\,,\\
		&C_{2} \colon V_{2} \to V_{2}\,,
	\end{aligned}
\end{equation}
and all of them involve the vector space $V_{2}$. The structure of this space \eqref{sl(4) (0, l, 0) structure of V_2} is a direct sum of two subspaces. This results in the fact that the fields act on these subspaces independently, which means the matrices have a block structure.

To separate the two subspaces $\nu_{1, 2},~\nu_{2, 2}$, we fix the ordering of the vectors in $V_{2}$. First, we enumerate all the vectors of the form $C_{2}^{a}R_{2}$, then the remaining ones. As an example, let us consider another state of the representation $(0, 2, 0)$, where we have four vectors in $V_{2}$:
\begin{equation}\label{sl(4) (0, l, 0) example of enumeration}
	\begin{minipage}{0.35\textwidth}
		\begin{equation*}
			\begin{aligned}
				&V_{1} = {\rm Span}\{B_{1}R_{2}\,,~B_{1}C_{2}R_{2}\}\\
				&\nu_{1, 2} = {\rm Span}\{R_{2}\,,~C_{2}R_{2}\}\\
				&\nu_{2, 2} = {\rm Span}\{A_{1}B_{1}R_{2}\,,~A_{1}B_{1}C_{2}R_{2}\}\\
				&V_{3} = {\rm Span}\{A_{2}R_{2}\,,~A_{2}C_{2}R_{2}\}
			\end{aligned}
		\end{equation*}
	\end{minipage}
	\qquad
	\begin{minipage}{0.3\textwidth}
		\begin{tikzpicture}[tdplot_main_coords, scale=1.5]
			\begin{scope}
				\coordinate (A) at (0,0,0);
				\coordinate (B) at (1,0,0);
				\coordinate (C) at (1,1,0);
				\coordinate (D) at (0,1,0);
				\coordinate (A') at (0,0,1);
				\coordinate (B') at (1,0,1);
				\coordinate (C') at (1,1,1);
				\coordinate (D') at (0,1,1);
				\draw  (A) -- (B) -- (C) -- (D) -- cycle;
				\draw (A') -- (B') -- (C') -- (D') -- cycle;
				\draw (A) -- (A') (B) -- (B') (D) -- (D') (C) -- (C');
				\foreach \point in {B',D'}
				{
					\shade[ball color=\mygray, opacity=0.9] (\point) circle (3pt);
				}
				\foreach \point in {A',C'}
				{
					\shade[ball color=\myblue, opacity=0.9] (\point) circle (3pt);
				}
				\foreach \point in {B,D}
				{
					\shade[ball color=\mygray, opacity=0.9] (\point) circle (3pt);
				}
				\foreach \point in {A,C}
				{
					\shade[ball color=\myblue, opacity=0.9] (\point) circle (3pt);
				}
				\node[above left] at (0, 0, 0) {$\scriptstyle 1$};
				\node[right=3pt] at (0, 1, 0) {$\scriptstyle 1$};
				\node[left=3pt] at (1, 0, 0) {$\scriptstyle 1$};
				\node[below right] at (1, 1, 0) {$\scriptstyle 3$};
				\node[above left] at (0, 0, 1) {$\scriptstyle 2$};
				\node[right=3pt] at (0, 1, 1) {$\scriptstyle 2$};
				\node[left=3pt] at (1, 0, 1) {$\scriptstyle 2$};
				\node[below right] at (1, 1, 1) {$\scriptstyle 4$};
			\end{scope}
		\end{tikzpicture}
	\end{minipage}
	\quad
	\begin{minipage}{0.3\textwidth}
		\begin{tikzpicture}[tdplot_main_coords, scale=1.5]
			\begin{scope}
				\coordinate (A) at (0,0,0);
				\coordinate (B) at (1,0,0);
				\coordinate (C) at (1,1,0);
				\coordinate (D) at (0,1,0);
				\coordinate (A') at (0,0,1);
				\coordinate (B') at (1,0,1);
				\coordinate (C') at (1,1,1);
				\coordinate (D') at (0,1,1);
				\draw  (A) -- (B) -- (C) -- (D) -- cycle;
				\draw (A') -- (B') -- (C') -- (D') -- cycle;
				\draw (A) -- (A') (B) -- (B') (D) -- (D') (C) -- (C');
				\foreach \point in {B',D'}
				{
					\shade[ball color=\mygray, opacity=0.9] (\point) circle (3pt);
				}
				\foreach \point in {A',C'}
				{
					\shade[ball color=\myblue, opacity=0.9] (\point) circle (3pt);
				}
				\foreach \point in {B,D}
				{
					\shade[ball color=\mygray, opacity=0.9] (\point) circle (3pt);
				}
				\foreach \point in {A,C}
				{
					\shade[ball color=\myblue, opacity=0.9] (\point) circle (3pt);
				}
				\node[above left] at (0, 0, 0) {$\scriptstyle R_{2}$};
				\node[right=3pt] at (0, 1, 0) {$\scriptstyle A_{2}R_{2}$};
				\node[left=3pt] at (1, 0, 0) {$\scriptstyle B_{1}R_{2}$};
				\node[below=3pt] at (1, 1, 0) {$\scriptstyle A_{1}B_{1}R_{2}$};
				\node[above left] at (0, 0, 1) {$\scriptstyle C_{2}R_{2}$};
				\node[right=3pt] at (0, 1, 1) {$\scriptstyle A_{2}C_{2}R_{2}$};
				\node[left=3pt] at (1, 0, 1) {$\scriptstyle B_{1}C_{2}R_{2}$};
			\end{scope}
		\end{tikzpicture}
	\end{minipage}
\end{equation}
In \eqref{sl(4) (0, l, 0) example of enumeration} we present the structure of the vector spaces, the enumeration of the atoms of the state, and the crystal itself.

We can now define the vacuum expectation values of the fields for an arbitrary GT base \eqref{sl(4) (0, l, 0) Gelfand-Tsetlin} using the ICO matrices introduced earlier \eqref{ICO base matrices}:
\begin{equation}\label{sl(4) (0, l, 0) vacuum fixed points}
	\begin{aligned}
		&\bar{C}_{1} = C(n_{1})\,, \quad 
		\bar{C}_{2} = \left(
		\begin{array}{c|c}
			C(m_{1}) & O(m_{1}, m_{2}) \\
			\hline
			O(m_{2}, m_{1}) & C(m_{2})
		\end{array}
		\right)\,, \quad
		\bar{C}_{3} = C(n_{3})\,, \\
		&\bar{A}_{1} = \left(\begin{array}{c}
			O(m_{1}, n_{1}) \\ \hline I(m_{2}, n_{1})
		\end{array}\right)\,, \quad
		\bar{A}_{2} = \left(\begin{array}{c|c}
			I(n_{3}, m_{1}) & O(n_{3}, m_{2})
		\end{array}\right)\,, \\
		&\bar{B}_{1} = \left(\begin{array}{c|c}
			I(n_{1}, m_{1}) & O(n_{1}, m_{2})
		\end{array}\right)\,, \quad
		\bar{B}_{2} = \left(\begin{array}{c}
			O(m_{1}, n_{3}) \\ \hline I(m_{2}, n_{3})
		\end{array}\right)\,,\\
		&\bar{R}_{1} = O(n_{1}, 1)\,, \quad
		\bar{R}_{2} = I(m_{1} + m_{2}, 1)\,, \quad
		\bar{R}_{3} = O(n_{3}, 1)\,.
	\end{aligned}
\end{equation}

Substituting the fixed point into the algorithm discussed in section \ref{subsec: Equivariant Matrix Coefficients}, we end up with the matrix elements of the following form:
\begin{equation}\label{sl(4) (0, l, 0) equivariant coefficients}
	\begin{aligned}
		&{\bf E}_{\Upsilon_{2,\lambda}}[\mu \to \mu_{1, 1 + 1}^{1}] = \frac{(m_{1} - n_{1})}{n_{1}\epsilon - \frac{\epsilon}{2}}\,,\\
		&{\bf E}_{\Upsilon_{2,\lambda}}[\mu \to \mu_{2, 2 + 1}^{1}] = -\frac{1}{\epsilon}\,,\\
		&{\bf E}_{\Upsilon_{2,\lambda}}[\mu \to \mu_{2, 2 + 1}^{2}] = -\frac{(n_{1} - m_{2})(n_{3} - m_{2})}{(m_{1} - m_{2})(m_{1} - m_{2} + 1)}\frac{1}{\epsilon}\,,\\
		&{\bf E}_{\Upsilon_{2,\lambda}}[\mu \to \mu_{3, 3 + 1}^{2}] = \frac{m_{1} - n_{3}}{n_{3}\epsilon - \frac{\epsilon}{2}}\,,\\
		&{\bf F}_{\Upsilon_{2,\lambda}}[\mu \to \mu_{1, 1 - 1}^{1}] = (n_{1} - m_{2})\Bigl(n_{1} - \frac{3}{2}\Bigr)\epsilon\,,\\
		&{\bf F}_{\Upsilon_{2,\lambda}}[\mu \to \mu_{2, 2 - 1}^{1}] = -\frac{(m_{1} + 1)(\lambda - m_{1} + 1)(m_{1} - n_{1})(m_{1} - n_{3})}{(m_{1} - m_{2} + 1)(m_{1} - m_{2})}\epsilon\,,\\
		&{\bf F}_{\Upsilon_{2,\lambda}}[\mu \to \mu_{2, 2 - 1}^{2}] = - m_{2}(\lambda - m_{2} + 2)\epsilon\,,\\
		&{\bf F}_{\Upsilon_{2,\lambda}}[\mu \to \mu_{3, 3 - 1}^{2}] = (n_{3} - m_{2})\Bigl(n_{3} - \frac{3}{2}\Bigr)\epsilon\,.
	\end{aligned}
\end{equation}
Using the normalization \eqref{Rescaling Equiv-to-Root}, we verify that at the zero level the construction gives us representations $(0, \lambda, 0)$ of the algebra $\fs\fl_{4}$:
\begin{equation}\label{sl(4) (0, l, 0) Lie algebra coefficients}
	\begin{aligned}
		e^{(1)}_{0}|\mu\rangle &= \sqrt{(m_{1} - n_{1})(n_{1} - m_{2} + 1)}|\mu_{1, 1 + 1}^{1}\rangle\\
		e_{0}^{(2)}|\mu\rangle &= 
		\sqrt{\frac{(m_{1} + 2)(\lambda - m_{1})(m_{1} - n_{1} + 1)(m_{1} - n_{3} + 1)}{(m_{1} - m_{2} + 2)(m_{1} - m_{2} + 1)}}|\mu_{2, 2 + 1}^{1}\rangle +\\ 
		&+ \sqrt{\frac{(m_{2} + 1)(\lambda - m_{2} + 1)(n_{1} - m_{2})(n_{3} - m_{2})}{(m_{1} - m_{2})(m_{1} - m_{2} + 1)}}|\mu_{2, 2 + 1}^{2}\rangle\\
		e^{(3)}_{0}|\mu\rangle &= \sqrt{(m_{1} - n_{3})(n_{3} - m_{2} + 1)}|\mu_{3, 3 + 1}^{2}\rangle,\\
		f^{(1)}_{0}|\mu\rangle &= \sqrt{(m_{1} - n_{1} + 1)(n_{1} - m_{2})}|\mu_{1, 1 - 1}^{1}\rangle\\
		f_{0}^{(2)}|\mu\rangle &= 
		\sqrt{\frac{(m_{1} + 1)(\lambda - m_{1} + 1)(m_{1} - n_{1})(m_{1} - n_{3})}{(m_{1} - m_{2} + 1)(m_{1} - m_{2})}}|\mu_{2, 2 - 1}^{1}\rangle +\\ 
		&+ \sqrt{\frac{m_{2}(\lambda - m_{2} + 2)(n_{1} - m_{2} + 1)(n_{3} - m_{2} + 1)}{(m_{1} - m_{2} + 1)(m_{1} - m_{2} + 2)}}|\mu_{2, 2 - 1}^{2}\rangle\\
		f^{(3)}_{0}|\mu\rangle &= \sqrt{(m_{1} - n_{3} + 1)(n_{3} - m_{2})}|\mu_{3, 3 - 1}^{2}\rangle,
	\end{aligned}
\end{equation}

Lastly, we need to check the hysteresis relations of the Yangian \eqref{general hysteresis}. They are also satisfied in this case and take the following form:
\begin{equation}\label{sl(4) (0, l, 0) hysteresis}
	\begin{aligned}
		&{\bf E}_{\Upsilon_{2,\lambda}}[\mu \to \mu_{1, 1 + 1}^{1}]{\bf F}_{\Upsilon_{2, \lambda}}[\mu_{1, 1 + 1}^{1} \to \mu] = (m_{1} - n_{1})(n_{1} - m_{2} + 1) = \mathop{\rm res}_{z = n_{1}\epsilon - \frac{\epsilon}{2}}\Psi_{\mu, \Upsilon_{2,\lambda}}^{(1)}(z)\,,\\
		&{\bf E}_{\Upsilon_{2,\lambda}}[\mu \to \mu_{2, 2 + 1}^{1}]{\bf F}_{\Upsilon_{2, \lambda}}[\mu_{2, 2 + 1}^{1} \to \mu] = \frac{(m_{1} + 2)(\lambda - m_{1})(m_{1} - n_{1} + 1)(m_{1} - n_{3} + 1)}{(m_{1} - m_{2} + 2)(m_{1} - m_{2} + 1)} = \mathop{\rm res}_{z = m_{1}\epsilon}\Psi_{\mu, \Upsilon_{2,\lambda}}^{(1)}(z)\,,\\
		&{\bf E}_{\Upsilon_{2,\lambda}}[\mu \to \mu_{2, 2 + 1}^{2}]{\bf F}_{\Upsilon_{2, \lambda}}[\mu_{2, 2 + 1}^{2} \to \mu] = \frac{(m_{2} + 1)(\lambda - m_{2} + 1)(n_{1} - m_{2})(n_{3} - m_{2})}{(m_{1} - m_{2})(m_{1} - m_{2} + 1)} = \mathop{\rm res}_{z = m_{2}\epsilon - \epsilon}\Psi_{\mu, \Upsilon_{2,\lambda}}^{(1)}(z)\,,\\
		&{\bf E}_{\Upsilon_{2,\lambda}}[\mu \to \mu_{3, 3 + 1}^{2}]{\bf F}_{\Upsilon_{2, \lambda}}[\mu_{3, 3 + 1}^{2} \to \mu] = (m_{1} - n_{3})(n_{3} - m_{2} + 1) = \mathop{\rm res}_{z = n_{3}\epsilon - \frac{\epsilon}{2}}\Psi_{\mu, \Upsilon_{2,\lambda}}^{(3)}(z)\,,\\
		&\frac{{\bf E}_{\Upsilon_{2, \lambda}}[\mu \to \mu_{2, 2 + 1}^{1}]{\bf E}_{\Upsilon_{2, \lambda}}\Bigl[\mu_{2, 2 + 1}^{1} \to \bigl(\mu_{2, 2 + 1}^{1}\bigr)_{1, 1 + 1}^{1}\Bigr]}{{\bf E}_{\Upsilon_{2, \lambda}}[\mu \to \mu_{1, 1 + 1}^{1}]{\bf E}_{\Upsilon_{2, \lambda}}\Bigl[\mu_{1, 1 + 1}^{1} \to \bigl(\mu_{1, 1 + 1}^{1}\bigr)_{2, 2 + 1}^{1}\Bigr]} = \frac{m_{1} + 1 - n_{1}}{m_{1} - n_{1}} = \varphi_{1, 2}\Bigl(n_{1}\epsilon - \frac{\epsilon}{2} - m_{1}\epsilon\Bigr)\,,\\
		&\frac{{\bf E}_{\Upsilon_{2, \lambda}}[\mu \to \mu_{1, 1 + 1}^{1}]{\bf E}_{\Upsilon_{2, \lambda}}\Bigl[\mu_{1, 1 + 1}^{1} \to \bigl(\mu_{1, 1 + 1}^{1}\bigr)_{2, 2 + 1}^{2}\Bigr]}{{\bf E}_{\Upsilon_{2, \lambda}}[\mu \to \mu_{2, 2 + 1}^{2}]{\bf E}_{\Upsilon_{2, \lambda}}\Bigl[\mu_{2, 2 + 1}^{2} \to \bigl(\mu_{2, 2 + 1}^{2}\bigr)_{1, 1 + 1}^{1}\Bigr]} = \frac{n_{1} + 1 - m_{2}}{n_{1} - m_{2}} = \varphi_{2, 1}\Bigl(m_{2}\epsilon - \epsilon - n_{1}\epsilon + \frac{\epsilon}{2}\Bigr)\,,\\
		&\frac{{\bf E}_{\Upsilon_{2, \lambda}}[\mu \to \mu_{2, 2 + 1}^{1}]{\bf E}_{\Upsilon_{2, \lambda}}\Bigl[\mu_{2, 2 + 1}^{1} \to \bigl(\mu_{2, 2 + 1}^{1}\bigr)_{3, 3 + 1}^{2}\Bigr]}{{\bf E}_{\Upsilon_{2, \lambda}}[\mu \to \mu_{3, 3 + 1}^{2}]{\bf E}_{\Upsilon_{2, \lambda}}\Bigl[\mu_{3, 3 + 1}^{2} \to \bigl(\mu_{3, 3 + 1}^{2}\bigr)_{2, 2 + 1}^{1}\Bigr]} = \frac{m_{1} + 1 - n_{3}}{m_{1} - n_{3}} = \varphi_{3, 2}\Bigl(n_{3}\epsilon - \frac{\epsilon}{2} - m_{1}\epsilon\Bigr)\,,\\
		&\frac{{\bf E}_{\Upsilon_{2, \lambda}}[\mu \to \mu_{3, 3 + 1}^{2}]{\bf E}_{\Upsilon_{2, \lambda}}\Bigl[\mu_{3, 3 + 1}^{2} \to \bigl(\mu_{3, 3 + 1}^{2}\bigr)_{2, 2 + 1}^{2}\Bigr]}{{\bf E}_{\Upsilon_{2, \lambda}}[\mu \to \mu_{2, 2 + 1}^{2}]{\bf E}_{\Upsilon_{2, \lambda}}\Bigl[\mu_{2, 2 + 1}^{2} \to \bigl(\mu_{2, 2 + 1}^{2}\bigr)_{3, 3 + 1}^{2}\Bigr]} = \frac{n_{3} + 1 - m_{2}}{n_{3} - m_{2}} = \varphi_{2, 3}\Bigl(m_{2}\epsilon - \epsilon - n_{3}\epsilon + \frac{\epsilon}{2}\Bigr)\,,\\
		&\frac{{\bf E}_{\Upsilon_{2, \lambda}}[\mu \to \mu_{2, 2 + 1}^{1}]{\bf E}_{\Upsilon_{2, \lambda}}\Bigl[\mu_{2, 2 + 1}^{1} \to \bigl(\mu_{2, 2 + 1}^{1}\bigr)_{2, 2 + 1}^{2}\Bigr]}{{\bf E}_{\Upsilon_{2, \lambda}}[\mu \to \mu_{2, 2 + 1}^{2}]{\bf E}_{\Upsilon_{2, \lambda}}\Bigl[\mu_{2, 2 + 1}^{2} \to \bigl(\mu_{2, 2 + 1}^{2}\bigr)_{2, 2 + 1}^{1}\Bigr]} = \frac{m_{1} - m_{2}}{m_{1} - m_{2} + 2} = \varphi_{2, 2}\bigl(m_{2}\epsilon - \epsilon - m_{1}\epsilon\bigr)\,.
	\end{aligned}
\end{equation}

\subsection{$\myY(\fs\fl_{n})$ details}\label{subsec: Y(sl(n)) details}

We are now ready to lift the construction and describe the key elements of the $\myY(\fs\fl_{n})$ algebras representations for an arbitrary $n$.

\subsubsection{The Algebra}

We work with the quiver $\CQ_{n,p,\lambda}$ \eqref{family of quivers and superpotentials for sl(n)}. The dimensions of the nodes are denoted as $n_{1}, \dots, n_{r}$ where $r = n - 1$. These parameters, as usual, can be assembled into the dimension vector $\vec{n}$. We double the superpotential here for convenience:
\begin{equation}\label{sl(n) superpotential}
	\CW=\Tr\left[A_1C_1B_1+\sum\lm_{a=2}^{n-2}\left(A_aC_aB_a-B_{a-1}C_aA_{a-1}\right)-B_{n-2}C_{n-1}A_{n-2}+{\color{burgundy}C_p^{\lambda}R_{p}S_{p}}\right]\,.
\end{equation}

The equivariant weights and R-charges that are assigned to the fields take the following form:
\begin{equation}\label{sl(n) weights}
	\begin{tblr}{c|c|c|c|c|c}
		\mbox{Fields} & C_{i} & A_{i} & B_{i} & R_{p} & S_{p}\\
		\hline
		\mbox{Weights} & \epsilon & -\frac{\epsilon}{2} + h & -\frac{\epsilon}{2} - h & 0 & -\lambda \epsilon \\ 
		\mbox{R-charges} & 0 & 1 & 1 & 0 & 2
	\end{tblr}\,.
\end{equation}

The generating functions of the Yangian read:
\begin{equation}\label{sl(n) generating functions}
	e^{(a)}(z) = \sum_{n = 0}^{\infty}\dfrac{e^{(a)}_{n}}{z^{n + 1}}\,, \quad f^{(a)}(z) = \sum_{n = 0}^{\infty}\dfrac{f^{(a)}_{n}}{z^{n + 1}}\,, \quad
	\psi^{(a)}(z) = 1 + \sum_{n = 0}^{\infty}\dfrac{\psi_{n}^{(a)}}{z^{n + 1}}\,.
\end{equation}
These functions satisfy the relations \eqref{Yangian Relations}, where the bonding factors for $|a - b |\leqslant 1$ are given by:
\begin{equation}\label{sl(n) bonding factors}
	\varphi_{a, b}(z) = \delta_{b, a + 1}\dfrac{2z - \epsilon}{2z + \epsilon} + \delta_{b, a}\dfrac{z + \epsilon}{z - \epsilon} + \delta_{b, a - 1}\dfrac{2z - \epsilon}{2z + \epsilon}\,.
\end{equation}
When $|a - b| > 1$, we set $\varphi_{a, b}(z) \equiv 1$. Again, we substitute $h = 0$ in the calculations related to the algebra.

We also impose the Serre relations \cite{Drinfeld:1987sy, Guay:2018wel}:
\begin{equation}\label{sl(n) serre relations}
	\begin{aligned}
		&\sum_{\sigma \in \fS_{m}}\Bigl[e^{(a)}\bigl(u_{\sigma(1)}\bigr), \Bigl[e^{(a)}\bigl(u_{\sigma(2)}\bigr), \dots, \Bigl[e^{(a)}\bigl(u_{\sigma(m)}\bigr), e^{(b)}(v)\Bigl]\dots \Bigr]\Bigr] = 0,\\
		&\sum_{\sigma \in \fS_{m}}\Bigl[f^{(a)}\bigl(u_{\sigma(1)}\bigr), \Bigl[f^{(a)}\bigl(u_{\sigma(2)}\bigr),\dots \Bigl[f^{(a)}(u_{\sigma(m)}), f^{(b)}(v)\Bigr]\dots \Bigr]\Bigr] = 0,
	\end{aligned}
\end{equation}
where $\CA_{ab}$ is the Cartan matrix of $\fs\fl_{n}$, $m = 1 - \CA_{ab}$, and $a\ne b$.

Unfolding the relations \eqref{Yangian Relations} in modes for $\myY(\fs\fl_{n})$ we get:
\begin{equation}\label{sl(n) Yangian relations in modes}
	\begin{aligned}
		&[e_{n + 1}^{(a)}, e^{(a)}_{k}] - [e^{(a)}_{n}, e^{(a)}_{k + 1}] = \epsilon\{e_{n}^{(a)}, e_{k}^{(a)}\}\,, \\
		&[e_{n + 1}^{(a)}, e_{k}^{(a + 1)}] - [e^{(a)}_{n}, e^{(a + 1)}_{k + 1}] = -\frac{\epsilon}{2}\{e^{(a)}_{n}, e^{(a + 1)}_{k}\}\,, \\
		&[e_{n + 1}^{(a + 1)}, e^{(a)}_{k}] - [e^{(a + 1)}_{n}, e_{k + 1}^{(a)}] = -\frac{\epsilon}{2}\{e^{(a + 1)}_{n}, e^{(a)}_{k}\}\,, \\
		&[\psi_{n + 1}^{(a)}, e_{k}^{(a)}] - [\psi^{(a)}_{n}, e^{(a)}_{k + 1}] = \epsilon\{\psi_{n}^{(a)}, e_{k}^{(a)}\}\,,\\
		&[\psi_{n + 1}^{(a)}, e_{k}^{(a + 1)}] - [\psi^{(a)}_{n}, e_{k + 1}^{(a + 1)}] = -\frac{\epsilon}{2}\{\psi_{n}^{(a)}, e_{k}^{(a + 1)}\}\,,\\
		&[\psi^{(a + 1)}_{n + 1}, e_{k}^{(a)}] - [\psi_{n}^{(a + 1)}, e_{k + 1}^{(a)}] = -\frac{\epsilon}{2}\{\psi_{n}^{(a + 1)}, e_{k}^{(a)}\}\,,\\
		&[f_{n + 1}^{(a)}, f^{(a)}_{k}] - [f^{(a)}_{n}, f^{(a)}_{k + 1}] = -\epsilon\{f_{n}^{(a)}, f_{k}^{(a)}\}\,, \\
		&[f_{n + 1}^{(a)}, f_{k}^{(a + 1)}] - [f^{(a)}_{n}, f^{(a + 1)}_{k + 1}] = \frac{\epsilon}{2}\{f^{(a)}_{n}, f^{(a + 1)}_{k}\}\,, \\
		&[f_{n + 1}^{(a + 1)}, f^{(a)}_{k}] - [f^{(a + 1)}_{n}, f_{k + 1}^{(a)}] = \frac{\epsilon}{2}\{f^{(a + 1)}_{n}, f^{(a)}_{k}\}\,, \\
		&[\psi_{n + 1}^{(a)}, f_{k}^{(a)}] - [\psi^{(a)}_{n}, f^{(a)}_{k + 1}] = -\epsilon\{\psi_{n}^{(a)}, f_{k}^{(a)}\}\,,\\
		&[\psi_{n + 1}^{(a)}, f_{k}^{(a + 1)}] - [\psi^{(a)}_{n}, f_{k + 1}^{(a + 1)}] = \frac{\epsilon}{2}\{\psi_{n}^{(a)}, f_{k}^{(a + 1)}\}\,,\\
		&[\psi^{(a + 1)}_{n + 1}, f_{k}^{(a)}] - [\psi_{n}^{(a + 1)}, f_{k + 1}^{(a)}] = \frac{\epsilon}{2}\{\psi_{n}^{(a + 1)}, f_{k}^{(a)}\}\,,\\
		&[\psi_{n}^{(a)}, \psi_{k}^{(b)}] = 0\,,\\
		&[e_{n}^{(a)}, f_{k}^{(b)}] = \delta_{ab}\psi^{(a)}_{n + k}\,,
	\end{aligned}
\end{equation}
with ``boundary conditions":
\begin{equation}\label{sl(n) boundary relations}
	\begin{aligned}
		&[\psi_{0}^{(a)}, e_{k}^{(a)}] = 2 e_{k}^{(a)}\,,\quad &[\psi_{0}^{(a)}, f_{k}^{(a)}] = -2 f_{k}^{(a)}\,, \\
		&[\psi^{(a + 1)}_{0}, e_{k}^{(a)}] = -e_{k}^{(a)}\,,\quad &[\psi_{0}^{(a + 1)}, f_{k}^{(a)}] = f_{k}^{(a)}\,, \\
		&[\psi_{0}^{(a)}, e^{(a + 1)}_{k}] = -e_{k}^{(a + 1)}\,,\quad &[\psi_{0}^{(a)}, f^{(a + 1)}_{k}] = f_{k}^{(a + 1)}\,.
	\end{aligned}
\end{equation}
The Serre relations \eqref{sl(n) serre relations} in modes take the form introduced in \cite{Drinfeld:1987sy}:
\begin{equation}\label{sl(n) serre relations in modes}
	\begin{aligned}
		&\sum_{\sigma \in \fS_{m}}\Bigl[e^{(a)}_{n_{\sigma(1)}}, \Bigl[e^{(a)}_{n_{\sigma(2)}}, \dots, \Bigl[e^{(a)}_{n_{\sigma(m)}}, e^{(b)}_{m}\Bigl]\dots \Bigr]\Bigr] = 0\,,\\
		&\sum_{\sigma \in \fS_{m}}\Bigl[f^{(a)}_{n_{\sigma(1)}}, \Bigl[f^{(a)}_{n_{\sigma(2)}}, \dots, \Bigl[f^{(a)}_{n_{\sigma(m)}}, f^{(b)}_{m}\Bigl]\dots \Bigr]\Bigr] = 0\,.
	\end{aligned}
\end{equation}
Strictly speaking, we have not included the relations between the generators for which $|a - b| \geqslant 2$ in \eqref{sl(n) Yangian relations in modes}, \eqref{sl(n) boundary relations}. They are trivial, and we can include them by rewriting the relations more compactly using the Cartan matrix:
\begin{equation}\label{sl(n) relations in modes via Cartan matrix}
	\begin{aligned}
		&[e_{n + 1}^{(a)}, e^{(b)}_{k}] - [e^{(a)}_{n}, e^{(b)}_{k + 1}] = \frac{\epsilon}{2}\CA_{ab}\{e_{n}^{(a)}, e_{k}^{(b)}\}\,, \\
		&[\psi_{n + 1}^{(a)}, e_{k}^{(b)}] - [\psi^{(a)}_{n}, e^{(b)}_{k + 1}] = \frac{\epsilon}{2}\CA_{ab}\{\psi_{n}^{(a)}, e_{k}^{(b)}\}\,,\\
		&[f_{n + 1}^{(a)}, f^{(b)}_{k}] - [f^{(a)}_{n}, f^{(b)}_{k + 1}] = -\frac{\epsilon}{2}\CA_{ab}\{f_{n}^{(a)}, f_{k}^{(b)}\}\,, \\
		&[\psi_{n + 1}^{(a)}, f_{k}^{(b)}] - [\psi^{(a)}_{n}, f^{(b)}_{k + 1}] = -\frac{\epsilon}{2}\CA_{ab}\{\psi_{n}^{(a)}, f_{k}^{(b)}\}\,,\\ \hline
		&[\psi_{0}^{(a)}, e_{k}^{(b)}] = \CA_{ab} e_{k}^{(b)}\,,\qquad [\psi_{0}^{(a)}, f_{k}^{(b)}] = -\CA_{ab} f_{k}^{(b)}\,,
	\end{aligned}
\end{equation}
where the form of the remaining relations is not changed.

\subsubsection{The Representations $\Upsilon_{p,\lambda}$}

Having defined the algebraic relations of $\myY(\fs\fl_{n})$, we proceed to its representation. We fix an arbitrary representation $\Upsilon_{p, \lambda}$ \eqref{gen rectangular reps}, where the parameter $p$ is restricted by \eqref{gen rectangular reps after cutoff}, and $\lambda \geqslant 0$.

\subsubsection*{The States}

Now, we discuss the crystal structure of the states of this representation.
The superpotential \eqref{sl(n) superpotential} defines the F-term relations as follows:
\begin{equation}\label{sl(n) F-terms}
	\begin{aligned}
		&\d_{A_{i}}\CW = C_{i}B_{i} - B_{i}C_{i + 1} = 0\,, & i \in \overline{1, n - 2}\,,\\
		&\d_{B_{i}}\CW = A_{i}C_{i} - C_{i + 1}A_{i} = 0\,, & i \in \overline{1, n - 2}\,,\\
		&\d_{C_{i}}\CW = B_{i}A_{i} - A_{i - 1}B_{i - 1} = 0\,, & i \in \overline{2, n - 2}\,,\\
		&\d_{C_{1}}\CW = B_{1}A_{1} = 0\,,
		&\d_{C_{n - 1}}\CW = A_{n - 2}B_{n - 2} = 0\,,\\
		&\d_{S_{p}}\CW = C_{p}^{\lambda}R_{p} = 0\,.
	\end{aligned}
\end{equation}
The fixed points can be defined according to \eqref{sl(n) weights}:
\begin{equation}\label{sl(n) fixed points}
	\begin{aligned}
		&[\Phi_{i}, C_{i}] = \epsilon C_{i}\,,\quad \Phi_{i + 1}A_{i} - A_{i}\Phi_{i} = \Bigl(-\frac{\epsilon}{2} + h\Bigr)A_{i}\,, \quad
		\Phi_{i}B_{i} - B_{i}\Phi_{i + 1} = \Bigl(-\frac{\epsilon}{2} - h\Bigr)B_{i}\,,\\
		&\Phi_{i}R_{i} = 0\,, \qquad S_{p}\Phi_{p} = -\lambda \epsilon S_{p}\,.
	\end{aligned}
\end{equation}
The first two relations in \eqref{sl(n) F-terms} allow us to interchange a field $C_{i}$ and $C_{j}$ for any $i,~j\in Q_{0}$, just as in the previous examples. Therefore, we use the field $C_{p}$. Again, the term $C_{p}^{\lambda}R_{p} = 0$ serves as a cut-off.

The relations involving the fields $A_{i}$ and $B_{i}$ impose more interesting constraints on a crystal. A state grows from the field $R_{p}$. Acting by the fields $A_{i}, B_{j}$ (and $C_{p}$) on the field $R_{p}$ we get some graph whose 2-dimensional slice $(h, R)$ in the space $\Xi$ we depict in \eqref{sl(n) towards crystal}.
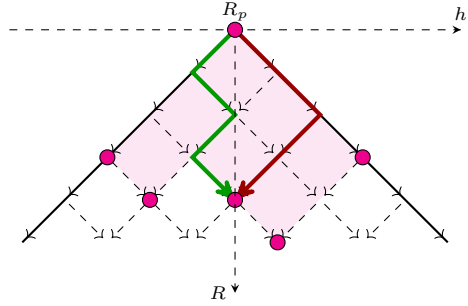
\begin{equation}\label{sl(n) towards crystal}
	\begin{tikzpicture}
		\node[above] at (0, 0) {$\scriptstyle R_{p}$};
		\node[above] at (3, 0) {$\scriptstyle h$};
		\node[left] at (0, -3.5) {$\scriptstyle R$};
		\begin{scope}[shift={(0,0)}, scale=0.8, rotate=-135]
			\draw[white,fill=magenta,opacity=0.1] (0,0) -- (3,0) -- (3,1) -- (2, 1) -- (2,3) -- (0,3) -- cycle;
			\draw[white] (5,0) -- (0,0) -- (0,5) -- cycle;
			\foreach \x in {0,...,4}
			\foreach \y in {0,...,4}
			{
				\pgfmathparse{int(\x+\y-3)}
				\let\r\pgfmathresult
				\ifnum \r > 0
				\breakforeach
				\fi
				\draw[dashed, ->] (\x+0.1,\y)  -- (\x+0.9,\y);
				\draw[dashed, ->] (\x,\y+0.1)  -- (\x,\y+0.9);
			}
			\draw[thick] (0,5) -- (0,0) -- (5,0);
			\draw[ultra thick,->,>=stealth',black!40!green] (0,0) -- (1,0) -- (1,1) -- (2,1) -- (2,2);
			\draw[ultra thick,->,>=stealth',black!40!red] (0,0) -- (0,2) -- (2,2);
			\draw[fill=magenta] (3, 0) circle (0.125) (3, 1) circle (0.125) (2, 3) circle (0.125) (0, 3) circle (0.125) (2, 2) circle (0.125);
		\end{scope}
		\draw[dashed, postaction={decorate},
		decoration={markings, mark= at position 1 with {\arrow{stealth}}}] (-3, 0) -- (3, 0);
		\draw[dashed, postaction={decorate},
		decoration={markings, mark= at position 1 with {\arrow{stealth}}}] (0, 0) -- (0, -3.5);
		\draw[fill=magenta] (0, 0) circle (0.1);
	\end{tikzpicture}
\end{equation}

As we have mentioned, the $F$-terms impose homotopic equivalence on paths of the R-equivariant space $\Xi$. An example of the equivalence is presented in \eqref{sl(n) towards crystal}.

Our aim is to construct the empty box, the set of all the admissible paths on $\Xi$, for the representation $\Upsilon_{p, \lambda}$. One could check that we always can form a rectangular structure \eqref{sl(n) towards crystal 2} that, by construction, is limited by the edges of the quiver \eqref{family of quivers and superpotentials for sl(n)}.
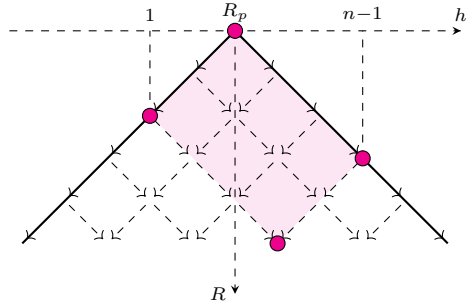
\begin{equation}\label{sl(n) towards crystal 2}
	\begin{tikzpicture}
		\node[above] at (0, 0) {$\scriptstyle R_{p}$};
		\node[above] at (3, 0) {$\scriptstyle h$};
		\node[left] at (0, -3.5) {$\scriptstyle R$};
		\begin{scope}[shift={(0,0)}, scale=0.8, rotate=-135]
			\draw[dashed] (2, 0) -- (1, -1) (0, 3) -- (-1.5, 1.5);
			\draw[white,fill=magenta,opacity=0.1] (0,0) -- (2,0) -- (2,3) -- (0,3) -- cycle;
			\draw[white] (5,0) -- (0,0) -- (0,5) -- cycle;
			\foreach \x in {0,...,4}
			\foreach \y in {0,...,4}
			{
				\pgfmathparse{int(\x+\y-3)}
				\let\r\pgfmathresult
				\ifnum \r > 0
				\breakforeach
				\fi
				\draw[dashed, ->] (\x+0.1,\y)  -- (\x+0.9,\y);
				\draw[dashed, ->] (\x,\y+0.1)  -- (\x,\y+0.9);
			}
			\draw[thick] (0,5) -- (0,0) -- (5,0);
			\draw[fill=magenta] (2, 0) circle (0.125) (2, 3) circle (0.125) (0, 3) circle (0.125);
			\node[above] at (1, -1) {$\scriptstyle 1$};
			\node[above] at (-1.5, 1.5) {$\scriptstyle n - 1$};
		\end{scope}
		\draw[dashed, postaction={decorate},
		decoration={markings, mark= at position 1 with {\arrow{stealth}}}] (-3, 0) -- (3, 0);
		\draw[dashed, postaction={decorate},
		decoration={markings, mark= at position 1 with {\arrow{stealth}}}] (0, 0) -- (0, -3.5);
		\draw[fill=magenta] (0, 0) circle (0.1);
	\end{tikzpicture}
\end{equation}

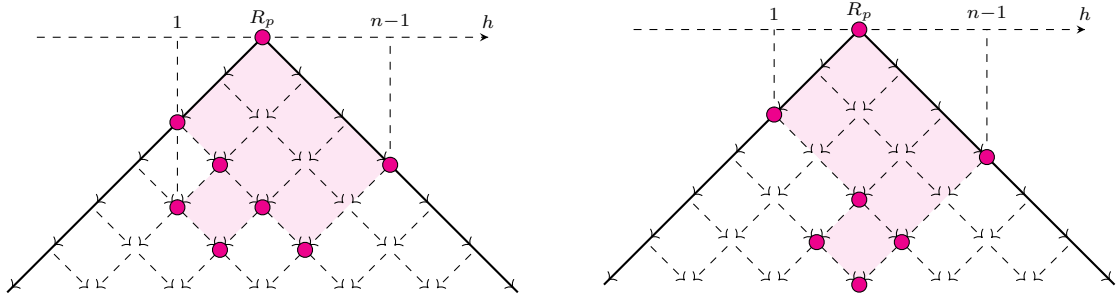
\begin{figure}
	\centering
	\begin{subfigure}{0.4\textwidth}
		\begin{center}
			\begin{tikzpicture}
				\node[above] at (0, 0) {$\scriptstyle R_{p}$};
				\node[above] at (3, 0) {$\scriptstyle h$};
				\begin{scope}[shift={(0,0)}, scale=0.8, rotate=-135]
					\draw[dashed] (3, 1) -- (1, -1) (0, 3) -- (-1.5, 1.5);
					\draw[white,fill=magenta,opacity=0.1] (0,0) -- (2,0) -- (2,1) -- (3,1) -- (3,2) -- (2,2) -- (2,3) -- (0,3) -- cycle;
					\draw[white] (6,0) -- (0,0) -- (0,6) -- cycle;
					\foreach \x in {0,...,5}
					\foreach \y in {0,...,5}
					{
						\pgfmathparse{int(\x+\y-4)}
						\let\r\pgfmathresult
						\ifnum \r > 0
						\breakforeach
						\fi
						\draw[dashed, ->] (\x+0.1,\y)  -- (\x+0.9,\y);
						\draw[dashed, ->] (\x,\y+0.1)  -- (\x,\y+0.9);
					}
					\draw[thick] (0,6) -- (0,0) -- (6,0);
					\draw[fill=magenta] (2, 0) circle (0.125) (2, 3) circle (0.125) (0, 3) circle (0.125) (2, 1) circle (0.125) (3, 1) circle (0.125) (3, 2) circle (0.125) (2, 2) circle (0.125);
					\node[above] at (1, -1) {$\scriptstyle 1$};
					\node[above] at (-1.5, 1.5) {$\scriptstyle n - 1$};
				\end{scope}
				\draw[dashed, postaction={decorate},
				decoration={markings, mark= at position 1 with {\arrow{stealth}}}] (-3, 0) -- (3, 0);
				\draw[fill=magenta] (0, 0) circle (0.1);
			\end{tikzpicture}
		\end{center}
	\end{subfigure}\hfil
	\begin{subfigure}{0.4\textwidth}
		\begin{center}
			\begin{tikzpicture}
				\node[above] at (0, 0) {$\scriptstyle R_{p}$};
				\node[above] at (3, 0) {$\scriptstyle h$};
				\begin{scope}[shift={(0,0)}, scale=0.8, rotate=-135]
					\draw[dashed] (2, 0) -- (1, -1) (0, 3) -- (-1.5, 1.5);
					\draw[white,fill=magenta,opacity=0.1] (0,0) -- (2,0) -- (2,2) -- (3,2) -- (3,3) -- (2,3) -- (0,3) -- cycle;
					\draw[white] (6,0) -- (0,0) -- (0,6) -- cycle;
					\foreach \x in {0,...,5}
					\foreach \y in {0,...,5}
					{
						\pgfmathparse{int(\x+\y-4)}
						\let\r\pgfmathresult
						\ifnum \r > 0
						\breakforeach
						\fi
						\draw[dashed, ->] (\x+0.1,\y)  -- (\x+0.9,\y);
						\draw[dashed, ->] (\x,\y+0.1)  -- (\x,\y+0.9);
					}
					\draw[thick] (0,6) -- (0,0) -- (6,0);
					\draw[fill=magenta] (2, 0) circle (0.125) (2, 3) circle (0.125) (0, 3) circle (0.125) (2, 2) circle (0.125) (3, 2) circle (0.125) (3, 3) circle (0.125);
					\node[above] at (1, -1) {$\scriptstyle 1$};
					\node[above] at (-1.5, 1.5) {$\scriptstyle n - 1$};
				\end{scope}
				\draw[dashed, postaction={decorate},
				decoration={markings, mark= at position 1 with {\arrow{stealth}}}] (-3, 0) -- (3, 0);
				\draw[fill=magenta] (0, 0) circle (0.1);
			\end{tikzpicture}
		\end{center}
	\end{subfigure}
	\caption{Examples of inadmissible empty boxes}\label{sl(n) inadmissible graphs}
\end{figure}

In general, however, we could end up with one of the two cases depicted\footnote{The presented pictures have a very loose scale. We highlight only the edges of the graph. The magenta-colored area represents the shape of the empty box.} in fig. \ref{sl(n) inadmissible graphs} as well.
We claim that none of them is allowed in our case.
\begin{itemize}
	\item In the first case, we have an extra structure near the edge of the quiver. Let us look closer at the graph.
	
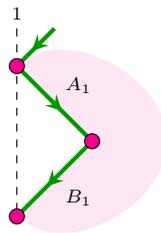
\begin{equation}\label{sl(n) green path on the edge}
		\begin{tikzpicture}
			\draw[dashed] (0, 0.5) -- (0, -2);
			\draw[ultra thick, black!40!green, postaction={decorate},
			decoration={markings, mark= at position 0.6 with {\arrow{stealth}}}] (0.5, 0.5) -- (0, 0);
			\draw[white, fill=magenta,opacity=0.1] (0, 0) to[out=45, in=90] (2, -1) to[out=-90, in=-45] (0, -2) to (1, -1) to (0, 0);
			\draw[ultra thick, black!40!green, postaction={decorate},
			decoration={markings, mark= at position 0.6 with {\arrow{stealth}}}] (0, 0) to node[black,above right] {$\scriptstyle A_{1}$} (1, -1);
			\draw[ultra thick, black!40!green, postaction={decorate},
			decoration={markings, mark=	at position 0.6 with {\arrow{stealth}}}] (1, -1) to node[black, below right] {$\scriptstyle B_{1}$} (0,-2);
			\draw[fill=magenta] (0, 0) circle (0.1) (1, -1) circle (0.1) (0, -2) circle (0.1);
			\node[above] at (0, 0.5) {$\scriptstyle 1$};
		\end{tikzpicture}
	\end{equation}
	We note that the green path in \eqref{sl(n) green path on the edge} is of the form:
	\begin{equation}
		B_{1}A_{1}\cdot\underbrace{\text{combination of the fields }A,\,B}\cdot\, C_{p}^{a}R_{p}\,.
	\end{equation}
	Using the relation:
	\begin{equation}
		B_{1}A_{1} = 0\,,
	\end{equation}
	we verify that this type of growth is indeed prohibited.
	\item In the second case, we have an extra structure ``inside" the quiver. The key fact that we use here is that the path equivalence holds in the whole space $\Xi$, not only the graph we cut off.
	Indeed, let us consider the red path depicted in \eqref{sl(n) red path inside the quiver}.
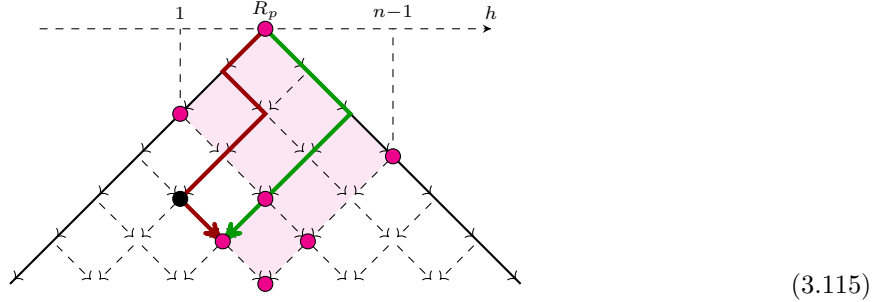
	\begin{equation}\label{sl(n) red path inside the quiver}
		\begin{tikzpicture}
			\node[above] at (0, 0) {$\scriptstyle R_{p}$};
			\node[above] at (3, 0) {$\scriptstyle h$};
			\begin{scope}[shift={(0,0)}, scale=0.8, rotate=-135]
				\draw[dashed] (2, 0) -- (1, -1) (0, 3) -- (-1.5, 1.5);
				\draw[white,fill=magenta,opacity=0.1] (0,0) -- (2,0) -- (2,2) -- (3,2) -- (3,3) -- (2,3) -- (0,3) -- cycle;
				\draw[white] (6,0) -- (0,0) -- (0,6) -- cycle;
				\foreach \x in {0,...,5}
				\foreach \y in {0,...,5}
				{
					\pgfmathparse{int(\x+\y-4)}
					\let\r\pgfmathresult
					\ifnum \r > 0
					\breakforeach
					\fi
					\draw[dashed, ->] (\x+0.1,\y)  -- (\x+0.9,\y);
					\draw[dashed, ->] (\x,\y+0.1)  -- (\x,\y+0.9);
				}
				\draw[thick] (0,6) -- (0,0) -- (6,0);
				\draw[ultra thick,->,>=stealth',black!40!green] (0,0) -- (0,2) -- (3,2);
				\draw[ultra thick,->,>=stealth',black!40!red] (0,0) -- (1,0) -- (1, 1) -- (3,1) -- (3,2);
				\draw[fill=magenta] (2, 0) circle (0.125) (2, 3) circle (0.125) (0, 3) circle (0.125) (2, 2) circle (0.125) (3, 2) circle (0.125) (3, 3) circle (0.125);
				\draw[fill=black] (3, 1) circle (0.125);
				\node[above] at (1, -1) {$\scriptstyle 1$};
				\node[above] at (-1.5, 1.5) {$\scriptstyle n - 1$};
			\end{scope}
			\draw[dashed, postaction={decorate},
			decoration={markings, mark= at position 1 with {\arrow{stealth}}}] (-3, 0) -- (3, 0);
			\draw[fill=magenta] (0, 0) circle (0.1);
		\end{tikzpicture}
	\end{equation}
	The path is equivalent to the green path by the relations $B_{i}A_{i} = A_{i-1}B_{i-1}$; however, it goes through the node that was not included inside the box. Moreover, we can homotopically deform the path to go through the edge nodes and conclude that the vacuum expectation value of this path is zero, as in the previous case.
\end{itemize}
Therefore, we have demonstrated that the $F$-terms limit the growth of a crystal drastically. We solidify our discussion above by providing a concrete example, $\CQ_{5,2,\lambda}$.

In this case, the height on the $\epsilon$-axis is simply restricted by $C_{2}^{\lambda}R_{2} = 0$. Therefore, we are interested in the two-dimensional slice of $\Xi$ that takes the following form:
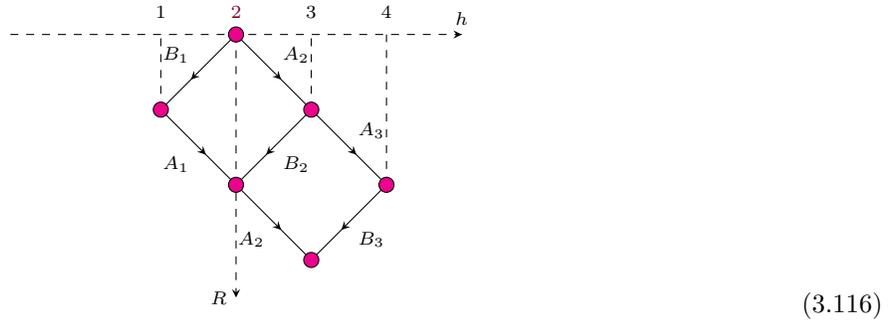
\begin{equation}\label{sl(5) empty box slice}
	\begin{tikzpicture}
		\node[above] at (3, 0) {$\scriptstyle h$};
		\node[left] at (0, -3.5) {$\scriptstyle R$};
		\node[above=0.1] at (-1, 0) {$\scriptstyle 1$};
		\node[above=0.1] at (0, 0) {$\color{burgundy} \scriptstyle 2$};
		\node[above=0.1] at (1, 0) {$\scriptstyle 3$};
		\node[above=0.1] at (2, 0) {$\scriptstyle 4$};
		\draw[postaction={decorate},
		decoration={markings, mark= at position 0.6 with {\arrow{stealth}}}] (0, 0) to node[above left] {$\scriptstyle B_{1}$} (-1, -1);
		\draw[postaction={decorate},
		decoration={markings, mark= at position 0.6 with {\arrow{stealth}}}] (0, 0) to node[above right] {$\scriptstyle A_{2}$} (1, -1);
		\draw[postaction={decorate},
		decoration={markings, mark= at position 0.6 with {\arrow{stealth}}}] (-1, -1) to node[below left] {$\scriptstyle A_{1}$} (0, -2);
		\draw[postaction={decorate},
		decoration={markings, mark= at position 0.6 with {\arrow{stealth}}}] (1, -1) to node[above right] {$\scriptstyle A_{3}$} (2, -2);
		\draw[postaction={decorate},
		decoration={markings, mark= at position 0.6 with {\arrow{stealth}}}] (0, -2) to node[below left] {$\scriptstyle A_{2}$} (1, -3);
		\draw[postaction={decorate},
		decoration={markings, mark= at position 0.6 with {\arrow{stealth}}}] (1, -1) to node[below right] {$\scriptstyle B_{2}$} (0, -2);
		\draw[postaction={decorate},
		decoration={markings, mark= at position 0.6 with {\arrow{stealth}}}] (2, -2) to node[below right] {$\scriptstyle B_{3}$} (1, -3);
		\draw[dashed] (-1, -1) -- (-1, 0) (1, -1) -- (1, 0) (2, -2) -- (2, 0);
		\draw[dashed, postaction={decorate},
		decoration={markings, mark= at position 1 with {\arrow{stealth}}}] (-3, 0) -- (3, 0);
		\draw[dashed, postaction={decorate},
		decoration={markings, mark= at position 1 with {\arrow{stealth}}}] (0, 0) -- (0, -3.5);
		\draw[fill=magenta] (0, 0) circle (0.1) (1, -1) circle (0.1) (-1, -1) circle (0.1) (0, -2) circle (0.1) (1, -3) circle (0.1) (2, -2) circle (0.1);
	\end{tikzpicture}
\end{equation}
The numbers label the nodes of the quiver $\CQ_{5,2,\lambda}$. We focus on a single path depicted in \eqref{sl(5) empty box wrong path}.
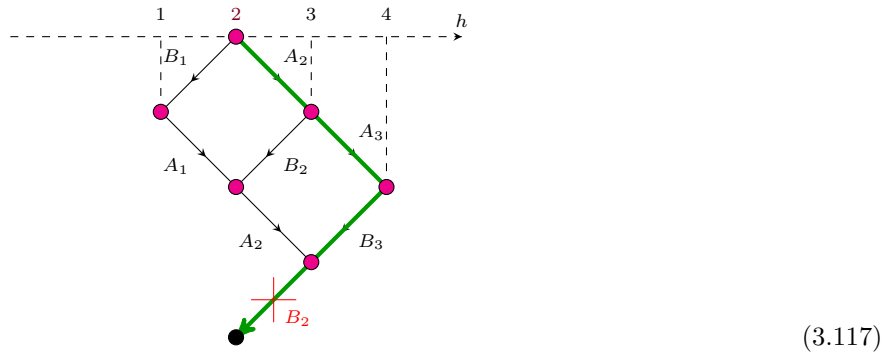
\begin{equation}\label{sl(5) empty box wrong path}
	\begin{tikzpicture}
		\node[above] at (3, 0) {$\scriptstyle h$};
		\node[above=0.1] at (-1, 0) {$\scriptstyle 1$};
		\node[above=0.1] at (0, 0) {$\color{burgundy} \scriptstyle 2$};
		\node[above=0.1] at (1, 0) {$\scriptstyle 3$};
		\node[above=0.1] at (2, 0) {$\scriptstyle 4$};
		\draw[postaction={decorate},
		decoration={markings, mark= at position 0.6 with {\arrow{stealth}}}] (0, 0) to node[above left] {$\scriptstyle B_{1}$} (-1, -1);
		\draw[postaction={decorate},
		decoration={markings, mark= at position 0.6 with {\arrow{stealth}}}] (0, 0) to node[above right] {$\scriptstyle A_{2}$} (1, -1);
		\draw[postaction={decorate},
		decoration={markings, mark= at position 0.6 with {\arrow{stealth}}}] (-1, -1) to node[below left] {$\scriptstyle A_{1}$} (0, -2);
		\draw[postaction={decorate},
		decoration={markings, mark= at position 0.6 with {\arrow{stealth}}}] (1, -1) to node[above right] {$\scriptstyle A_{3}$} (2, -2);
		\draw[postaction={decorate},
		decoration={markings, mark= at position 0.6 with {\arrow{stealth}}}] (0, -2) to node[below left] {$\scriptstyle A_{2}$} (1, -3);
		\draw[postaction={decorate},
		decoration={markings, mark= at position 0.6 with {\arrow{stealth}}}] (1, -1) to node[below right] {$\scriptstyle B_{2}$} (0, -2);
		\draw[postaction={decorate},
		decoration={markings, mark= at position 0.6 with {\arrow{stealth}}}] (2, -2) to node[below right] {$\scriptstyle B_{3}$} (1, -3);
		\draw[ultra thick,->,>=stealth',black!40!green] (0,0) to (2,-2) to (1, -3) to node[red, below right] {$\scriptstyle B_{2}$} (0,-4);
		\draw[rotate=0, red] (0.5-0.3,-3.5) -- (0.5+0.3,-3.5) (0.5,-3.5-0.3) -- (0.5,-3.5+0.3);
		\draw[dashed] (-1, -1) -- (-1, 0) (1, -1) -- (1, 0) (2, -2) -- (2, 0);
		\draw[dashed, postaction={decorate},
		decoration={markings, mark= at position 1 with {\arrow{stealth}}}] (-3, 0) -- (3, 0);
		\draw[fill=magenta] (0, 0) circle (0.1) (1, -1) circle (0.1) (-1, -1) circle (0.1) (0, -2) circle (0.1) (1, -3) circle (0.1) (2, -2) circle (0.1);
		\draw[fill=black] (0, -4) circle (0.1);
	\end{tikzpicture}
\end{equation}
This green path is a path of the form:
\begin{equation}
	B_{2}B_{3}A_{3}A_{2}\,C_{2}^{a}R_{2}\,.
\end{equation}
Now we use the relations \eqref{sl(n) F-terms} and get the result:
\begin{equation}
	B_{2}B_{3}A_{3}A_{2}\,C_{2}^{a}R_{2} = (B_{2}A_{2})(B_{2}A_{2})\,C_{2}^{a}R_{2} = A_{1}(B_{1}A_{1})B_{1}\,C_{2}^{a}R_{2} = 0\,,
\end{equation}
which proves that the path is not admissible.

We can now proceed further and define the empty box for $\Upsilon_{p,\lambda}$. The two-dimensional projection reads as follows:
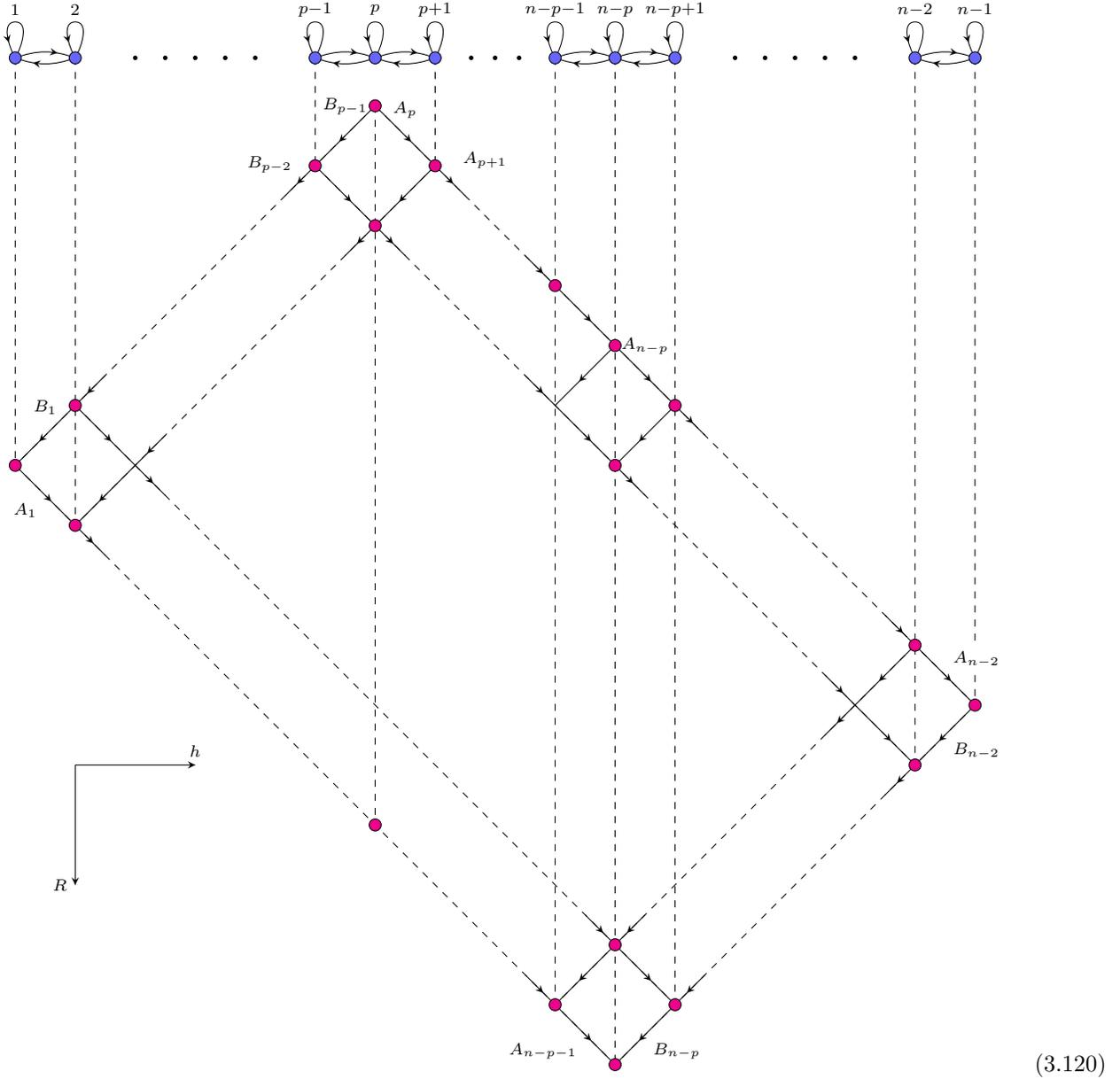
\begin{equation}\label{sl(n) different paths}
	\begin{tikzpicture}[scale=0.9]
		\draw[postaction={decorate},
		decoration={markings, mark= at position 0.6 with {\arrow{stealth}}}] (0, 0) to node[above=0.15] {$\scriptstyle A_{p}$} (1,-1);
		\draw[postaction={decorate},
		decoration={markings, mark= at position 0.6 with {\arrow{stealth}}}] (0, 0) to node[above=0.2] {$\scriptstyle B_{p-1}$} (-1,-1);
		\draw[postaction={decorate},
		decoration={markings, mark= at position 0.6 with {\arrow{stealth}}}] (-1, -1) to node[above left] {$\scriptstyle B_{p-2}$} (-1.5,-1.5);
		\draw[postaction={decorate},
		decoration={markings, mark= at position 0.6 with {\arrow{stealth}}}] (1, -1) to node[above right=0.1] {$\scriptstyle A_{p+1}$} (1.5, -1.5);
		\draw[postaction={decorate},
		decoration={markings, mark= at position 0.6 with {\arrow{stealth}}}] (8.5, -8.5) to (9,-9);
		\draw[postaction={decorate},
		decoration={markings, mark= at position 0.6 with {\arrow{stealth}}}] (9, -9) to node[above right] {$\scriptstyle A_{n-2}$} (10,-10);
		\draw[postaction={decorate},
		decoration={markings, mark= at position 0.6 with {\arrow{stealth}}}] (10, -10) to node[below right] {$\scriptstyle B_{n-2}$} (9,-11);
		\draw[postaction={decorate},
		decoration={markings, mark= at position 0.6 with {\arrow{stealth}}}] (9, -11) to (8.5,-11.5);
		\draw[postaction={decorate},
		decoration={markings, mark= at position 0.6 with {\arrow{stealth}}}] (9, -9) to (8,-10);
		\draw[postaction={decorate},
		decoration={markings, mark= at position 0.6 with {\arrow{stealth}}}] (8, -10) to (9,-11);
		\draw[postaction={decorate},
		decoration={markings, mark= at position 0.6 with {\arrow{stealth}}}] (-1, -1) to (0,-2);
		\draw[postaction={decorate},
		decoration={markings, mark= at position 0.6 with {\arrow{stealth}}}] (1, -1) to (0,-2);
		\draw[postaction={decorate},
		decoration={markings, mark= at position 0.6 with {\arrow{stealth}}}] (0, -2) to (-0.5,-2.5);
		\draw[postaction={decorate},
		decoration={markings, mark= at position 0.6 with {\arrow{stealth}}}] (0, -2) to (0.5,-2.5);
		\draw[postaction={decorate},
		decoration={markings, mark= at position 0.6 with {\arrow{stealth}}}] (3, -3) to (4,-4);
		\draw[postaction={decorate},
		decoration={markings, mark= at position 0.6 with {\arrow{stealth}}}] (4, -4) to node[above=0.2] {$\scriptstyle A_{n-p}$} (5,-5);
		\draw[postaction={decorate},
		decoration={markings, mark= at position 0.6 with {\arrow{stealth}}}] (5, -5) to (5.5,-5.5);
		\draw[postaction={decorate},
		decoration={markings, mark= at position 0.6 with {\arrow{stealth}}}] (2.5, -2.5) to (3,-3);
		\draw[postaction={decorate},
		decoration={markings, mark= at position 0.6 with {\arrow{stealth}}}] (4, -4) to (3,-5);
		\draw[postaction={decorate},
		decoration={markings, mark= at position 0.6 with {\arrow{stealth}}}] (5, -5) to (4,-6);
		\draw[postaction={decorate},
		decoration={markings, mark= at position 0.6 with {\arrow{stealth}}}] (3, -5) to (4,-6);
		\draw[postaction={decorate},
		decoration={markings, mark= at position 0.6 with {\arrow{stealth}}}] (2.5, -4.5) to (3,-5);
		\draw[postaction={decorate},
		decoration={markings, mark= at position 0.6 with {\arrow{stealth}}}] (4, -6) to (4.5,-6.5);
		\draw[postaction={decorate},
		decoration={markings, mark= at position 0.6 with {\arrow{stealth}}}] (7.5, -9.5) to (8,-10);
		\draw[postaction={decorate},
		decoration={markings, mark= at position 0.6 with {\arrow{stealth}}}] (8, -10) to (7.5,-10.5);
		\draw[postaction={decorate}, 
		decoration={markings, mark= at position 0.6 with {\arrow{stealth}}}] (5, -15) to node[below right] {$\scriptstyle B_{n-p}$} (4,-16);
		\draw[postaction={decorate},
		decoration={markings, mark= at position 0.6 with {\arrow{stealth}}}] (4, -14) to (5,-15);
		\draw[postaction={decorate},
		decoration={markings, mark= at position 0.6 with {\arrow{stealth}}}] (4, -14) to (3,-15);
		\draw[postaction={decorate},
		decoration={markings, mark= at position 0.6 with {\arrow{stealth}}}] (3, -15) to node[below left] {$\scriptstyle A_{n-p-1}$} (4,-16);
		\draw[postaction={decorate},
		decoration={markings, mark= at position 0.6 with {\arrow{stealth}}}] (2.5, -14.5) to (3,-15);
		\draw[postaction={decorate},
		decoration={markings, mark= at position 0.6 with {\arrow{stealth}}}] (5.5, -14.5) to (5,-15);
		\draw[postaction={decorate},
		decoration={markings, mark= at position 0.6 with {\arrow{stealth}}}] (4.5, -13.5) to (4,-14);
		\draw[postaction={decorate},
		decoration={markings, mark= at position 0.6 with {\arrow{stealth}}}] (3.5, -13.5) to (4,-14);
		\draw[postaction={decorate},
		decoration={markings, mark= at position 0.6 with {\arrow{stealth}}}] (-5, -5) to node[above=0.2] {$\scriptstyle B_{1}$} (-6,-6);
		\draw[postaction={decorate},
		decoration={markings, mark= at position 0.6 with {\arrow{stealth}}}] (-6, -6) to node[below left] {$\scriptstyle A_{1}$} (-5,-7);
		\draw[postaction={decorate},
		decoration={markings, mark= at position 0.6 with {\arrow{stealth}}}] (-5, -7) to (-4.5,-7.5);
		\draw[postaction={decorate},
		decoration={markings, mark= at position 0.6 with {\arrow{stealth}}}] (-4.5, -4.5) to (-5,-5);
		\draw[postaction={decorate},
		decoration={markings, mark= at position 0.6 with {\arrow{stealth}}}] (-5, -5) to (-4,-6);
		\draw[postaction={decorate},
		decoration={markings, mark= at position 0.6 with {\arrow{stealth}}}] (-4, -6) to (-5,-7);
		\draw[postaction={decorate},
		decoration={markings, mark= at position 0.6 with {\arrow{stealth}}}] (-3.5, -5.5) to (-4,-6);
		\draw[postaction={decorate},
		decoration={markings, mark= at position 0.6 with {\arrow{stealth}}}] (-4, -6) to (-3.5,-6.5);
		\draw[dashed] (-6, 0.8) to (-6, -6);
		\draw[dashed] (-5, 0.8) to (-5, -7);
		\draw[dashed] (10, 0.8) to (10, -9);
		\draw[dashed] (10, -9.5) to (10, -10);
		\draw[dashed] (-1, 0.8) to (-1, -1);
		\draw[dashed] (1, 0.8) to (1, -1);
		\draw[dashed] (3, 0.8) to (3, -15);
		\draw[dashed] (4, 0.8) to (4, -16);
		\draw[dashed] (5, 0.8) to (5, -15);
		\draw[dashed] (9, 0.8) to (9, -11);
		\draw[dashed] (-1, -1) to (9, -11);
		\draw[dashed] (-5, -5) to (5, -15);
		\draw[dashed] (1, -1) to (-5, -7);
		\draw[dashed] (0, 0) to (-6, -6);
		\draw[dashed] (0, 0) to (10, -10);
		\draw[dashed] (-6, -6) to (4, -16);
		\draw[dashed] (10, -10) to (4, -16);
		\draw[dashed] (9, -9) to (3, -15);
		\draw[dashed] (0, 0) to (0, -12);
		\draw[fill=magenta] (0, 0) circle (0.1) (4, -16) circle (0.1) (-6, -6) circle (0.1) (10, -10) circle (0.1) (1, -1) circle (0.1) (-1, -1) circle (0.1) (3, -3) circle (0.1) (4, -4) circle (0.1) (5, -5) circle (0.1) (9, -9) circle (0.1) (9, -11) circle (0.1) (-5, -5) circle (0.1) (-5, -7) circle (0.1) (3, -15) circle (0.1) (5, -15) circle (0.1) (4, -14) circle (0.1) (4, -6) circle (0.1) (0, -2) circle (0.1) (0, -12) circle (0.1);
		\draw[postaction=decorate, decoration={markings, mark= at position 0.8 with {\arrow{stealth}}}] (0,0.8) to[out=60,in=0] (0,1.4) to[out=180,in=120] (0,0.8);
		\draw[postaction=decorate, decoration={markings, mark= at position 0.8 with {\arrow{stealth}}}] (1,0.8) to[out=60,in=0] (1,1.4) to[out=180,in=120] (1,0.8);
		\draw[postaction=decorate, decoration={markings, mark= at position 0.8 with {\arrow{stealth}}}] (-1,0.8) to[out=60,in=0] (-1,1.4) to[out=180,in=120] (-1,0.8);
		\draw[postaction=decorate, decoration={markings, mark= at position 0.8 with {\arrow{stealth}}}] (3,0.8) to[out=60,in=0] (3,1.4) to[out=180,in=120] (3,0.8);
		\draw[postaction=decorate, decoration={markings, mark= at position 0.8 with {\arrow{stealth}}}] (4,0.8) to[out=60,in=0] (4,1.4) to[out=180,in=120] (4,0.8);
		\draw[postaction=decorate, decoration={markings, mark= at position 0.8 with {\arrow{stealth}}}] (5,0.8) to[out=60,in=0] (5,1.4) to[out=180,in=120] (5,0.8);
		\draw[postaction=decorate, decoration={markings, mark= at position 0.8 with {\arrow{stealth}}}] (9,0.8) to[out=60,in=0] (9,1.4) to[out=180,in=120] (9,0.8);
		\draw[postaction=decorate, decoration={markings, mark= at position 0.8 with {\arrow{stealth}}}] (10,0.8) to[out=60,in=0] (10,1.4) to[out=180,in=120] (10,0.8);
		\draw[postaction=decorate, decoration={markings, mark= at position 0.8 with {\arrow{stealth}}}] (-5,0.8) to[out=60,in=0] (-5,1.4) to[out=180,in=120] (-5,0.8);
		\draw[postaction=decorate, decoration={markings, mark= at position 0.8 with {\arrow{stealth}}}] (-6,0.8) to[out=60,in=0] (-6,1.4) to[out=180,in=120] (-6,0.8);
		\draw[postaction=decorate, decoration={markings, mark= at position 0.7 with {\arrow{stealth}}}] (0,0.8) to[out=20,in=160] (1,0.8);
		\draw[postaction=decorate, decoration={markings, mark= at position 0.7 with {\arrow{stealth}}}] (3,0.8) to[out=20,in=160] (4,0.8);
		\draw[postaction=decorate, decoration={markings, mark= at position 0.7 with {\arrow{stealth}}}] (4,0.8) to[out=20,in=160] (5,0.8);
		\draw[postaction=decorate, decoration={markings, mark= at position 0.7 with {\arrow{stealth}}}] (9,0.8) to[out=20,in=160] (10,0.8);
		\draw[postaction=decorate, decoration={markings, mark= at position 0.7 with {\arrow{stealth}}}] (-1,0.8) to[out=20,in=160] (0,0.8);
		\draw[postaction=decorate, decoration={markings, mark= at position 0.7 with {\arrow{stealth}}}] (-6,0.8) to[out=20,in=160] (-5,0.8);
		\draw[postaction=decorate, decoration={markings, mark= at position 0.7 with {\arrow{stealth}}}] (1,0.8) to[out=200,in=340] (0,0.8);
		\draw[postaction=decorate, decoration={markings, mark= at position 0.7 with {\arrow{stealth}}}] (0,0.8) to[out=200,in=340] (-1,0.8);
		\draw[postaction=decorate, decoration={markings, mark= at position 0.7 with {\arrow{stealth}}}] (4,0.8) to[out=200,in=340] (3,0.8);
		\draw[postaction=decorate, decoration={markings, mark= at position 0.7 with {\arrow{stealth}}}] (5,0.8) to[out=200,in=340] (4,0.8);
		\draw[postaction=decorate, decoration={markings, mark= at position 0.7 with {\arrow{stealth}}}] (10,0.8) to[out=200,in=340] (9,0.8);
		\draw[postaction=decorate, decoration={markings, mark= at position 0.7 with {\arrow{stealth}}}] (-5,0.8) to[out=200,in=340] (-6,0.8);
		\draw[fill=\myblue] (0, 0.8) circle (0.1) (-6, 0.8) circle (0.1) (10, 0.8) circle (0.1) (4, 0.8) circle (0.1) (-1, 0.8) circle (0.1) (1, 0.8) circle (0.1) (3, 0.8) circle (0.1) (5, 0.8) circle (0.1) (9, 0.8) circle (0.1) (-5, 0.8) circle (0.1);
		\draw[fill=black] (1.6, 0.8) circle (0.03) (2, 0.8) circle (0.03) (2.4, 0.8) circle (0.03) (6, 0.8) circle (0.03) (6.5, 0.8) circle (0.03) (7, 0.8) circle (0.03) (7.5, 0.8) circle (0.03) (8, 0.8) circle (0.03) (-3, 0.8) circle (0.03) (-2.5, 0.8) circle (0.03) (-2, 0.8) circle (0.03) (-3.5, 0.8) circle (0.03) (-4, 0.8) circle (0.03);
		\node at (0, 1.6) {$\scriptstyle p$};
		\node at (1, 1.6) {$\scriptstyle p+1$};
		\node at (-1, 1.6) {$\scriptstyle p-1$};
		\node at (-6, 1.6) {$\scriptstyle 1$};
		\node at (-5, 1.6) {$\scriptstyle 2$};
		\node at (3, 1.6) {$\scriptstyle n-p-1$};
		\node at (4, 1.6) {$\scriptstyle n-p$};
		\node at (5, 1.6) {$\scriptstyle n-p+1$};
		\node at (10, 1.6) {$\scriptstyle n-1$};
		\node at (9, 1.6) {$\scriptstyle n-2$};
		\begin{scope}[shift={(-5, -11)}]
			\draw[postaction=decorate, decoration={markings, mark= at position 1 with {\arrow{stealth}}}] (0, 0) -- (2, 0);
			\draw[postaction=decorate, decoration={markings, mark= at position 1 with {\arrow{stealth}}}] (0, 0) -- (0, -2);
			\node[above] at (2, 0) {$\scriptstyle h$};
			\node[left] at (0, -2) {$\scriptstyle R$};
		\end{scope}
	\end{tikzpicture}
\end{equation}
As we can see in \eqref{sl(n) different paths}, the crystals that correspond to the quiver $\CQ_{n,p,\lambda}$ contain $p(n - p)$ different types of paths. As in \eqref{sl(4) (0, l, 0) paths on quiver}, we can lift the picture to the whole R-equivariant space, where we can see the fields $C_{i}$. The parameter $\lambda$ becomes the maximal allowed height.

\subsubsection*{Gelfand-Tsetlin Bases}

In this paragraph we lift the correspondence between the crystals and the Gelfand-Tsetlin bases to the representation $\Upsilon_{p,\lambda}$.

We highlight the main ideas of the construction starting with the example $\CQ_{5,2,\lambda}$. We have established earlier that the empty box for the corresponding representation takes the form \eqref{sl(5) paths on quiver example}.
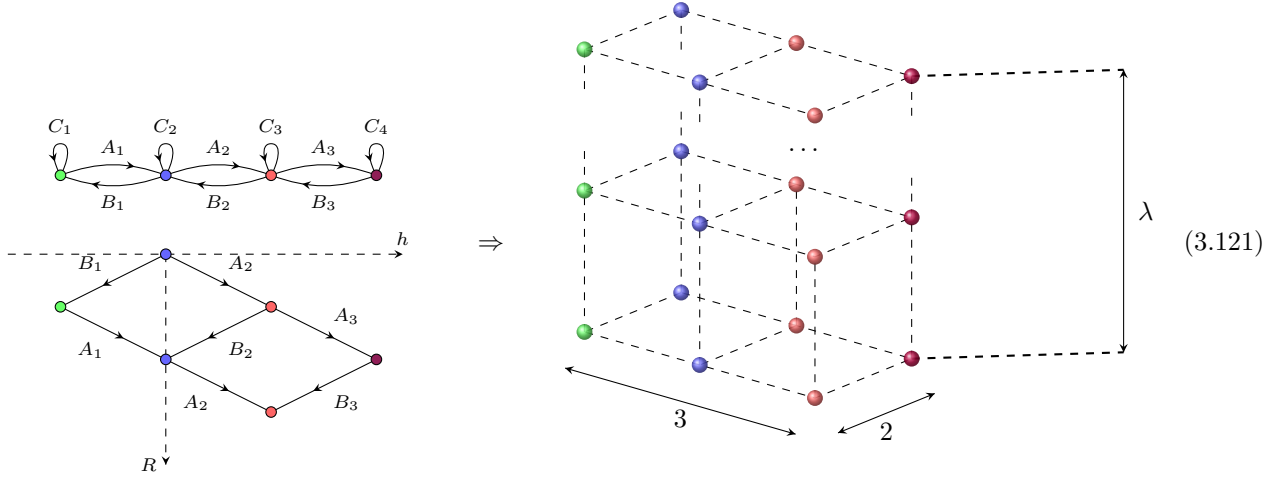
\begin{equation}\label{sl(5) paths on quiver example}
	\begin{tikzpicture}[scale = 0.7, baseline = {(0, -1)}]
		\begin{scope}[shift={(0, -1.5)}]
			\draw[dashed, postaction={decorate},
			decoration={markings, mark= at position 1 with {\arrow{stealth}}}] (-3, 0) -- (4.5, 0);
			\draw[dashed, postaction={decorate},
			decoration={markings, mark= at position 1 with {\arrow{stealth}}}] (0, 0) -- (0, -4);
			\node[above] at (4.5, 0) {$\scriptstyle h$};
			\node[left] at (0, -4) {$\scriptstyle R$};
		\end{scope}
		\begin{scope}[shift = {(0,0)}]
			\draw[postaction=decorate, decoration={markings, mark= at position 0.7 with {\arrow{stealth}}}] (0,0) to[out=20,in=160] node[pos=0.5,above] {$\scriptstyle A_2$} (2,0);
			\draw[postaction=decorate, decoration={markings, mark= at position 0.7 with {\arrow{stealth}}}] (2,0) to[out=200,in=340] node[pos=0.5,below] {$\scriptstyle B_2$} (0,0);
			\begin{scope}[shift={(-2,0)}]
				\draw[postaction=decorate, decoration={markings, mark= at position 0.7 with {\arrow{stealth}}}] (0,0) to[out=20,in=160] node[pos=0.5,above] {$\scriptstyle A_{1}$} (2,0);
				\draw[postaction=decorate, decoration={markings, mark= at position 0.7 with {\arrow{stealth}}}] (2,0) to[out=200,in=340] node[pos=0.5,below] {$\scriptstyle B_{1}$} (0,0);
			\end{scope}
			\begin{scope}[shift={(2,0)}]
				\draw[postaction=decorate, decoration={markings, mark= at position 0.7 with {\arrow{stealth}}}] (0,0) to[out=20,in=160] node[pos=0.5,above] {$\scriptstyle A_{3}$} (2,0);
				\draw[postaction=decorate, decoration={markings, mark= at position 0.7 with {\arrow{stealth}}}] (2,0) to[out=200,in=340] node[pos=0.5,below] {$\scriptstyle B_{3}$} (0,0);
			\end{scope}
			\draw[postaction=decorate, decoration={markings, mark= at position 0.8 with {\arrow{stealth}}}] (0,0) to[out=60,in=0] (0,0.6) to[out=180,in=120] (0,0);
			\node[above] at (0,0.6) {$\scriptstyle C_2$};
			\begin{scope}[shift={(2,0)}]
				\draw[postaction=decorate, decoration={markings, mark= at position 0.8 with {\arrow{stealth}}}] (0,0) to[out=60,in=0] (0,0.6) to[out=180,in=120] (0,0);
				\node[above] at (0,0.6) {$\scriptstyle C_{3}$};
			\end{scope}
			\begin{scope}[shift={(-2,0)}]
				\draw[postaction=decorate, decoration={markings, mark= at position 0.8 with {\arrow{stealth}}}] (0,0) to[out=60,in=0] (0,0.6) to[out=180,in=120] (0,0);
				\node[above] at (0,0.6) {$\scriptstyle C_{1}$};
			\end{scope}
			\begin{scope}[shift={(4,0)}]
				\draw[postaction=decorate, decoration={markings, mark= at position 0.8 with {\arrow{stealth}}}] (0,0) to[out=60,in=0] (0,0.6) to[out=180,in=120] (0,0);
				\node[above] at (0,0.6) {$\scriptstyle C_{4}$};
			\end{scope}
			\draw[fill=purple3] (4,0) circle (0.1);
			\draw[fill=\myred] (2,0) circle (0.1);
			\draw[fill=\mygreen] (-2,0) circle (0.1);
			\draw[fill=\myblue] (0, 0) circle (0.1);
		\end{scope}
		\begin{scope}[shift={(0, -1.5)}]
			\draw[postaction={decorate},
			decoration={markings, mark= at position 0.6 with {\arrow{stealth}}}] (0, 0) to node[above left] {$\scriptstyle B_{1}$} (-2, -1);
			\draw[postaction={decorate},
			decoration={markings, mark= at position 0.6 with {\arrow{stealth}}}] (0, 0) to node[above right] {$\scriptstyle A_{2}$} (2, -1);
			\draw[postaction={decorate},
			decoration={markings, mark= at position 0.6 with {\arrow{stealth}}}] (-2, -1) to node[below left] {$\scriptstyle A_{1}$} (0, -2);
			\draw[postaction={decorate},
			decoration={markings, mark= at position 0.6 with {\arrow{stealth}}}] (2, -1) to node[above right] {$\scriptstyle A_{3}$} (4, -2);
			\draw[postaction={decorate},
			decoration={markings, mark= at position 0.6 with {\arrow{stealth}}}] (0, -2) to node[below left] {$\scriptstyle A_{2}$} (2, -3);
			\draw[postaction={decorate},
			decoration={markings, mark= at position 0.6 with {\arrow{stealth}}}] (2, -1) to node[below right] {$\scriptstyle B_{2}$} (0, -2);
			\draw[postaction={decorate},
			decoration={markings, mark= at position 0.6 with {\arrow{stealth}}}] (4, -2) to node[below right] {$\scriptstyle B_{3}$} (2, -3);
			\draw[fill=\myblue] (0, 0) circle (0.1) (0, -2) circle (0.1);
			\draw[fill=\myred] (2, -1) circle (0.1) (2, -3) circle (0.1);
			\draw[fill=\mygreen] (-2, -1) circle (0.1);
			\draw[fill=purple3] (4, -2) circle (0.1);
		\end{scope}
	\end{tikzpicture}
	\qquad
	\Rightarrow
	\qquad
	\begin{tikzpicture}[tdplot_main_coords, scale=2, baseline={(0, 1)}]
		\coordinate (A) at (0,0,0);
		\coordinate (B) at (1,0,0);
		\coordinate (C) at (1,1,0);
		\coordinate (D) at (0,1,0);
		\coordinate (A') at (0,0,1);
		\coordinate (B') at (1,0,1);
		\coordinate (C') at (1,1,1);
		\coordinate (D') at (0,1,1);
		\coordinate (E) at (0,0,2);
		\coordinate (F) at (1,0,2);
		\coordinate (G) at (1,1,2);
		\coordinate (H) at (0,1,2);
		\coordinate (X) at (0, -1, 0);
		\coordinate (Y) at (1, -1, 0);
		\coordinate (X') at (0, -1, 1);
		\coordinate (Y') at (1, -1, 1);
		\coordinate (L) at (0, -1, 2);
		\coordinate (P) at (1, -1, 2);
		\draw[dashed] (A) -- (B) -- (C) -- (D) -- cycle;
		\draw[dashed] (A') -- (B') -- (C') -- (D') -- cycle;
		\draw[dashed] (E) -- (F) -- (G) -- (H) -- cycle;
		\draw[dashed] (A) -- (A') (B) -- (B') (C) -- (C') (D) -- (D') (X) -- (X') (Y) -- (Y');
		\draw[dashed] (A) -- (X) (B) -- (Y) (X) -- (Y) (A') -- (X') (B') -- (Y') (X') -- (Y') (E) -- (L) (F) -- (P) (L) -- (P);
		\draw[dashed] (B') -- (1,0,1.3);
		\draw[dashed] (D') -- (0,1,1.3);
		\draw[dashed] (F) -- (1,0,1.7);
		\draw[dashed] (H) -- (0,1,1.7);
		\draw[dashed] (Y') -- (1, -1, 1.3);
		\draw[dashed] (P) -- (1, -1, 1.7);
		\draw[dashed] (X') -- (0, -1, 1.3);
		\draw[dashed] (L) -- (0, -1, 1.7);
		\node at (0.5,0.5, 1.5) {$\dots$};
		\draw[postaction={decorate},
		decoration={markings, mark= at position 0 with {\arrowreversed{stealth}}, mark= at position 1 with {\arrow{stealth}}}] (D)++(0,0.2,-0.2) --++ (1,0,0)
		node[midway, below] {$2$};
		\draw[postaction={decorate},
		decoration={markings, mark= at position 0 with {\arrowreversed{stealth}}, mark= at position 1 with {\arrow{stealth}}}] (Y)++(0.2,0,-0.2) --++ (0,2,0)
		node[midway, below] {$3$};
		\draw[thick, dashed] ($(D)+(0,0,0)$) -- ($(D)+(0,0,0)+(-1,1,0)$);
		\draw[thick, dashed] ($(H)+(0,0,0)$) -- ($(H)+(0,0,0)+(-1,1,0)$);
		\draw[postaction={decorate},
		decoration={markings, mark= at position 0 with {\arrowreversed{stealth}}, mark= at position 1 with {\arrow{stealth}}}] ($(D)+(-1,1,0)+(0,0,0)$) -- 
		($(H)+(0,0,0)+(-1,1,0)$)
		node[midway, right=2pt] {$\lambda$};
		\foreach \point in {Y, Y', P}
		{
			\shade[ball color=\mygreen, opacity=0.9] (\point) circle (1.5pt);
		}
		\foreach \point in {X, X', L, B, B', F}
		{
			\shade[ball color=\myblue, opacity=0.9] (\point) circle (1.5pt);
		}
		\foreach \point in {A, A', E, C, C', G}
		{
			\shade[ball color=\myred, opacity=0.9] (\point) circle (1.5pt);
		}
		\foreach \point in {D, D', H}
		{
			\shade[ball color=purple, opacity=0.9] (\point) circle (1.5pt);
		}
	\end{tikzpicture}
\end{equation}
The vector spaces $V_{1}, \dots, V_{4}$ in this case can be represented as:
\begin{equation}\label{sl(5) vector spaces}
	\begin{aligned}
		&V_{1} = \nu_{1, 1} = {\rm Span}\left\{B_{1}C_{2}^{a}\cdot R_{2}\right\}\,,\\
		&V_{2} = \nu_{1, 2} \oplus \nu_{2, 2} = {\rm Span}\left\{C_{2}^{a}\cdot R_{2}\,,~~ A_{1} B_{1} C_{2}^{a}\cdot R_{2} \right\}\,,\\
		&V_{3} = \nu_{2, 3} \oplus \nu_{3, 3} = {\rm Span}\left\{A_{2}C_{2}^{a}\cdot R_{2}\,,~~ B_{3}A_{3}A_{2}C_{2}^{a}\cdot R_{2}\right\}\,,\\
		&V_{4} = \nu_{3, 4} = {\rm Span}\left\{A_{3}A_{2}C_{2}^{a}\cdot R_{2}\right\}\,,
	\end{aligned}
\end{equation}
where the numeration of the subspaces is chosen to match the corresponding elements of the Gelfand-Tsetlin patterns. For example, the second number labels a node:
\begin{equation}
	V_{\color{blue} a} = \bigoplus \nu_{\bullet , \color{blue} a}\,.
\end{equation}
The first number runs through the different types of vectors that live in the corresponding vector space. In order to clarify the range of this parameter, we construct the GT bases next. The result is presented in \eqref{sl(5) gelfand tsetlin and dimensions}.
\begin{equation}\label{sl(5) gelfand tsetlin and dimensions}
	\begin{tikzpicture}[baseline={(0,-0.3)}, every node/.style={inner sep=1pt}]
		\draw[rounded corners=2pt, thick, orange] (-3, 1.3) rectangle (3, 1.7);
		\node[above=2pt] at (0, 1.7) {\footnotesize \textbf{The partition}};
		\GTrow[1.5]{$\lambda$, $\lambda$, $\lambda$, $0$, $0$}
		\GTrow[0.9]{$\lambda$, $\lambda$, \textcolor{purple}{$m_{3, 4}$}, $0$}
		\GTrow[0.3]{$\lambda$, \textcolor{\myred}{$m_{2, 3}$}, \textcolor{\myred}{$m_{3, 3}$}}
		\GTrow[-0.3]{\textcolor{\myblue}{$m_{1, 2}$}, \textcolor{\myblue}{$m_{2, 2}$}}
		\GTrow[-0.9]{\textcolor{\mygreen!50!black}{$m_{1, 1}$}}
		\node[right=10pt] at (3.5, 0.9) {$n_4$};
		\node[right=10pt] at (3.5, 0.3) {$n_3$};
		\node[right=10pt] at (3.5, -0.3) {$n_2$};
		\node[right=10pt] at (3.5, -0.9) {$n_1$};
		
		\draw[postaction={decorate},
		decoration={markings, mark= at position 1 with {\arrow{stealth}}}] (1.9, 0.9) -- (3.2, 0.9);
		\draw[postaction={decorate},
		decoration={markings, mark= at position 1 with {\arrow{stealth}}}] (1.6, 0.3) -- (3.2, 0.3);
		\draw[postaction={decorate},
		decoration={markings, mark= at position 1 with {\arrow{stealth}}}] (1.25, -0.3) -- (3.2, -0.3);
		\draw[postaction={decorate},
		decoration={markings, mark= at position 1 with {\arrow{stealth}}}] (0.7, -0.9) -- (3.2, -0.9);
	\end{tikzpicture}
\end{equation}
At the top we place the partition that corresponds to the Young diagram that labels the representation $\Upsilon_{2, \lambda}$:
\begin{equation}
	\Upsilon_{2, \lambda} \quad \leftrightarrow \quad (0, \lambda, 0, 0) \quad \leftrightarrow \quad [\lambda, \lambda, \lambda, 0, 0]\,.
\end{equation}
The parameters $m_{i, k}$ are the dimensions of the vector spaces \eqref{sl(5) vector spaces}:
\begin{equation}
	m_{i, k} = \dim \nu_{i, k}\,.
\end{equation}
The pattern is organized to satisfy the triangular inequalities:
\begin{equation}\label{sl(5) triangular inequalities}
	\begin{aligned}
		&\lambda \geqslant m_{3, 4} \geqslant 0\,, \quad m_{1, 2} \geqslant m_{1, 1} \geqslant m_{2,2}\,,\\
		&\lambda \geqslant m_{1, 2} \geqslant m_{2, 3}\,, \quad m_{2, 3} \geqslant m_{2, 2} \geqslant m_{3, 3}\,,\\
		&\lambda \geqslant m_{2, 3} \geqslant m_{3, 4}\,, \quad m_{3, 4} \geqslant m_{3, 3} \geqslant 0\,.
	\end{aligned}
\end{equation}
These inequalities define the structure of the GT bases, and from the crystal point of view, they are a direct consequence of the F-relations \eqref{sl(n) F-terms}, which we have seen in the example before in section \ref{subsec: Y(sl(4))}.

Let us prove it in the case $\CQ_{5,2,\lambda}$ as well. In order to get a visualization of a GT base, one could rotate the crystal:
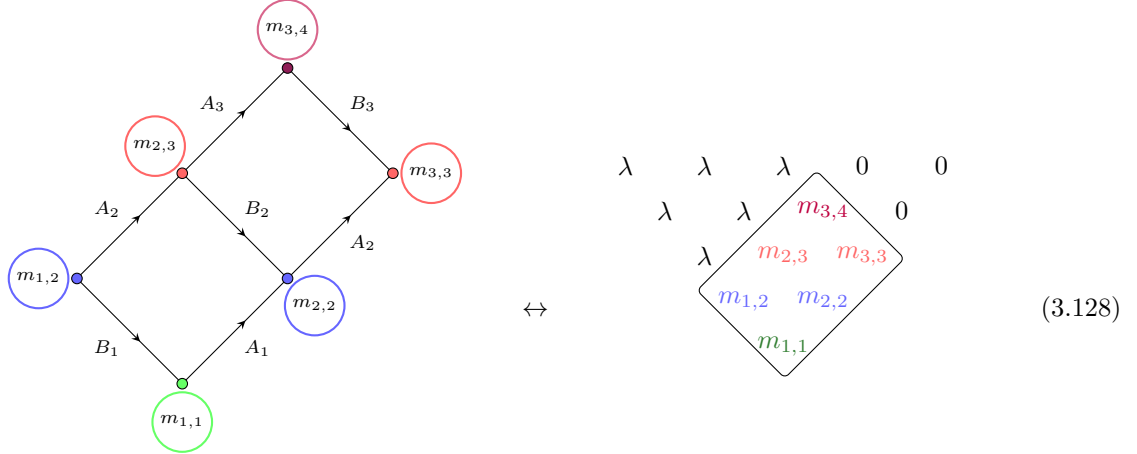
\begin{equation}\label{sl(5) crystal and GT base}
	\begin{tikzpicture}[scale = 0.7, baseline = {(0, -0.5)}]
		\begin{scope}[shift={(0, 0)}, rotate=90]
			\draw[postaction={decorate},
			decoration={markings, mark= at position 0.6 with {\arrow{stealth}}}] (0, 0) to node[below left] {$\scriptstyle B_{1}$} (-2, -2);
			\draw[postaction={decorate},
			decoration={markings, mark= at position 0.6 with {\arrow{stealth}}}] (0, 0) to node[above left] {$\scriptstyle A_{2}$} (2, -2);
			\draw[postaction={decorate},
			decoration={markings, mark= at position 0.6 with {\arrow{stealth}}}] (-2, -2) to node[below right] {$\scriptstyle A_{1}$} (0, -4);
			\draw[postaction={decorate},
			decoration={markings, mark= at position 0.6 with {\arrow{stealth}}}] (2, -2) to node[above left] {$\scriptstyle A_{3}$} (4, -4);
			\draw[postaction={decorate},
			decoration={markings, mark= at position 0.6 with {\arrow{stealth}}}] (0, -4) to node[below right] {$\scriptstyle A_{2}$} (2, -6);
			\draw[postaction={decorate},
			decoration={markings, mark= at position 0.6 with {\arrow{stealth}}}] (2, -2) to node[above right] {$\scriptstyle B_{2}$} (0, -4);
			\draw[postaction={decorate},
			decoration={markings, mark= at position 0.6 with {\arrow{stealth}}}] (4, -4) to node[above right] {$\scriptstyle B_{3}$} (2, -6);
			\draw[fill=\myblue] (0, 0) circle (0.1) (0, -4) circle (0.1);
			\draw[fill=\myred] (2, -2) circle (0.1) (2, -6) circle (0.1);
			\draw[fill=\mygreen] (-2, -2) circle (0.1);
			\draw[fill=purple3] (4, -4) circle (0.1);
			\node[style={circle, draw=\myblue, thick, inner sep=2pt}, left=0.1] at (0, 0) {$\scriptstyle m_{1, 2}$};
			\node[style={circle, draw=\myblue, thick, inner sep=2pt}, below right=0.1] at (0, -4) {$\scriptstyle m_{2, 2}$};
			\node[style={circle, draw=\mygreen, thick, inner sep=2pt}, below=0.1] at (-2, -2) {$\scriptstyle m_{1, 1}$};
			\node[style={circle, draw=\myred, thick, inner sep=2pt}, above left=0.1] at (2, -2) {$\scriptstyle m_{2,3}$};
			\node[style={circle, draw=\myred, thick, inner sep=2pt}, right=0.1] at (2, -6) {$\scriptstyle m_{3, 3}$};
			\node[style={circle, draw=\mypurple, thick, inner sep=2pt}, above=0.1] at (4, -4) {$\scriptstyle m_{3, 4}$};
		\end{scope}
	\end{tikzpicture}
	\qquad
	\leftrightarrow
	\qquad
	\begin{tikzpicture}[scale = 1, baseline = {(0, -0.5)}]
		\begin{scope}[rotate=-45, shift={(-0.7,-0.925)}]
			\draw[rounded corners=2pt, black] (0, 0) -- (1.65, 0) -- (1.65, 2.25) -- (0, 2.25) -- cycle;
		\end{scope}
		\GTrow[1.5]{$\lambda$, $\lambda$, $\lambda$, $0$, $0$}
		\GTrow[0.9]{$\lambda$, $\lambda$, \textcolor{purple}{$m_{3, 4}$}, $0$}
		\GTrow[0.3]{$\lambda$, \textcolor{\myred}{$m_{2, 3}$}, \textcolor{\myred}{$m_{3, 3}$}}
		\GTrow[-0.3]{\textcolor{\myblue}{$m_{1, 2}$}, \textcolor{\myblue}{$m_{2, 2}$}}
		\GTrow[-0.9]{\textcolor{\mygreen!50!black}{$m_{1, 1}$}}
	\end{tikzpicture}
\end{equation}
The picture \eqref{sl(5) crystal and GT base} is in agreement with the decomposition of the vector spaces \eqref{sl(5) vector spaces}. We additionally highlight the dimensions of these vector spaces\footnote{One could depict the dimensions directly as the corresponding heights on the axis $\epsilon$, therefore, had to draw the 3-dimensional space $\Xi$. We believe that our notation is slightly more convenient for our purposes.}.

We claim that there are only two types of inequalities we need to prove, namely $m_{1, 2} \geqslant m_{2, 3}$ and $m_{1, 2}\geqslant m_{1, 1}$. Indeed, every other inequality is of one of two types by construction.

Let us start with the first one,  $m_{1, 2} \geqslant m_{2, 3}$. We assume the contrary, which means that there is a natural number $k$, such that:
\begin{equation}
	C_{2}^{k}\cdot R_{2} = 0\,, \qquad C_{3}^{k}A_{2}\cdot R_{2}\ne 0\,.
\end{equation}
However, applying the F-terms \eqref{sl(n) F-terms}, we get:
\begin{equation}
	0 \ne C_{3}^{k}A_{2}\cdot R_{2} = A_{2}C_{2}^{k}\cdot R_{2} = 0\,,
\end{equation}
which contradicts our initial suggestion. Therefore, $m_{1, 2} \geqslant m_{2, 3}$.

Next, we prove the second inequality, $m_{1, 2}\geqslant m_{1, 1}$. In a similar way, we assume the contrary. Therefore, there is a number $k$ such that:
\begin{equation}
	C_{2}^{k}\cdot R_{2} = 0\,, \qquad C_{1}^{k}B_{1}\cdot R_{2}\ne 0\,.
\end{equation}
Again, we apply the F-terms and end up with the contradiction:
\begin{equation}
	0 \ne C_{1}^{k}B_{1}\cdot R_{2} = B_{1}C_{2}^{k}\cdot R_{2} = 0\,,
\end{equation}
which verifies that $m_{1, 2}\geqslant m_{1, 1}$, and also concludes the proof that the patterns \eqref{sl(5) gelfand tsetlin and dimensions} are indeed the GT bases.

Now, we are ready to present the general construction of the GT bases for $\CQ_{n,p,\lambda}$.
\begin{itemize}
	\item We start with the partition:
	\begin{equation}
		(0, \dots, 0, \underset{p}{\lambda}, 0, \dots, 0) \quad \leftrightarrow \quad [m_{1, n}, \dots, m_{i, n-p}, \dots, m_{n, n}]\,,
	\end{equation}
	where $m_{i,n} = \lambda$ if $i \leqslant n - p$ and $m_{i, n} = 0$ otherwise. The partition is \textit{fixed} by the representation.
	\item We decompose the vector spaces of the quiver nodes into direct sums of linear subspaces:
	\begin{equation}\label{sl(n) vector space decomposition}
		V_{k} = \bigoplus_{i}\nu_{i, k}\,,
	\end{equation}
	and assign to the subspaces their corresponding dimensions:
	\begin{equation}
		m_{i, k} = \dim \nu_{i, k}\,.
	\end{equation}
	\item We form a triangular pattern using these dimensions that reads:
	\begin{equation}\label{sl(n) Gelfand-Tsetlin pattern}
		\mu =\quad
		\begin{tikzpicture}[baseline={(0,-0.1)}]
			\node at (0, 1.2) {$m_{1,n} \quad m_{2,n} \quad \dots \quad m_{n-1,n} \quad m_{n,n}$};
			\node at (0, 0.6) {$m_{1,n-1} \quad \dots \quad m_{n-1,n-1}$};
			\node at (0, 0)     {$\dots \dots \dots \dots \dots$};
			\node at (0, -0.6)  {$m_{1,2} \quad m_{2,2}$};
			\node at (0, -1.2)  {$m_{1,1}$};
		\end{tikzpicture}\,.
	\end{equation}
	We impose the triangular inequalities on the fixed numbers $\lambda,\,0$ in advance. Therefore, the partition cuts off a rectangle of free parameters similar to the case \eqref{sl(5) crystal and GT base}. The shape resembles \eqref{sl(n) different paths}.
	\item The numbers $m_{i, j}$ are dimension parameters. By construction, they satisfy one of two types of inequalities:
	\begin{equation}
		\begin{minipage}{0.4\textwidth}
			\centering
			\begin{tikzpicture}
				\draw[postaction={decorate},
				decoration={markings, mark= at position 0.6 with {\arrow{stealth}}}] (0, 0) to node[below left] {$\scriptstyle B_{i}$} (1,-1);
				\draw[fill=magenta] (0, 0) circle (0.1) (1, -1) circle (0.1);
				\node[above=0.1] at (0, 0) {$\scriptstyle m_{i, j + 1}$};
				\node[below=0.1] at (1, -1) {$\scriptstyle m_{i, j}$};
				\node[draw = black, rounded corners = 0.1cm] at (0.5, -2) {$m_{i, j + 1}\geqslant m_{i, j}$};
			\end{tikzpicture}
		\end{minipage}
		\begin{minipage}{0.4\textwidth}
			\centering
			\begin{tikzpicture}
				\draw[postaction={decorate},
				decoration={markings, mark= at position 0.6 with {\arrow{stealth}}}] (0, 0) to node[below right] {$\scriptstyle A_{i}$} (1,1);
				\draw[fill=magenta] (0, 0) circle (0.1) (1, 1) circle (0.1);
				\node[below=0.1] at (0, 0) {$\scriptstyle m_{i, j}$};
				\node[above=0.1] at (1, 1) {$\scriptstyle m_{i + 1, j + 1}$};
				\node[draw = black, rounded corners = 0.1cm] at (0.5, -1) {$m_{i, j}\geqslant m_{i+1, j+1}$};
			\end{tikzpicture}
		\end{minipage}
	\end{equation}
	These inequalities are, in fact, the triangular inequalities:
	\begin{equation}\label{sl(n) triangual inequalities}
		m_{i, j + 1} \geqslant m_{i, j} \geqslant m_{i + 1, j + 1}\,,
	\end{equation}
	which are a well-known property of the GT bases.
	\item The decompositions \eqref{sl(n) vector space decomposition} can be rewritten in terms of the parameters $m_{i, j}$ as follows:
	\begin{equation}
		\dim V_{k} = n_{k} = \sum_{i = a[k]}^{b[k]}m_{i, k}\,,
	\end{equation}
	where the numbers $a[k],\, b[k]$ can be determined by the structure of \eqref{sl(n) Gelfand-Tsetlin pattern}. In the representation $\Upsilon_{p, \lambda}$ they take the form:
	\begin{equation}
		a[k] = \max\{1, k - p + 1\}\,, \quad b[k] = \min\{n-p, k\}\,.
	\end{equation}
\end{itemize}

\subsubsection*{The Representation}

We proceed now to describe the representation $\Upsilon_{p, \lambda}$ of the algebra $\myY(\fs\fl_{n})$. Again, we adjust the ansatz \eqref{General Ansatz}. The states of the representation are parameterized by the GT bases. As in the example presented in section \ref{subsec: Y(sl(4))}, we denote by $\mu_{k, k + 1}^{i}$ the pattern obtained from $\mu$ by replacing $m_{i, k}$ with $m_{i, k} + 1$. And by $\mu_{k, k - 1}^{i}$ we denote the pattern where $m_{i, k}$ is replaced with $m_{i, k} - 1$.

Implementing the facts above, we end up with the ansatz of the following form:
\begin{equation}\label{sl(n) ansatz}
	\begin{aligned}
		&e^{(k)}_{\Upsilon_{p,\lambda}}(z)|\mu\rangle = \sum_{i = a[k]}^{b[k]}\dfrac{{\bf E}_{\Upsilon_{p, \lambda}}[\mu \to \mu_{k, k + 1}^{i}]}{z - m_{i, k}\epsilon + (i - a[k])\epsilon + |k - p|\frac{\epsilon}{2}}|\mu_{k, k + 1}^{i}\rangle\,,\\
		&f^{(k)}_{\Upsilon_{p, \lambda}}(z)|\mu\rangle = \sum_{i = a[k]}^{b[k]}\dfrac{{\bf F}_{\Upsilon_{p, \lambda}}[\mu \to \mu_{k, k - 1}^{i}]}{z - (m_{i, k} - 1)\epsilon + (i - a[k])\epsilon + |k - p|\frac{\epsilon}{2}}|\mu_{k, k - 1}^{i}\rangle\,,\\
		&\psi^{(k)}_{\Upsilon_{p, \lambda}}(z)|\mu\rangle = \Psi^{(k)}_{\mu, \Upsilon_{p, \lambda}}(z)|\mu\rangle\,,
	\end{aligned}
\end{equation}
where the eigenfunctions are defined by \eqref{general eigenvalues}. The explicit formulas for the $\Upsilon_{p,\lambda}$ are presented below.
The functions on the edges read as follows: 
\begin{equation}\label{sl(n) eigenfunctions on the edge}
	\begin{aligned}
		\Psi^{(1)}_{\mu, \Upsilon_{p, \lambda}}(z) = &-\dfrac{1}{\epsilon}\prod_{l = 1}^{m_{1, 1}}\varphi_{1, 1}\Bigl(z - (l - 1)\epsilon + (i - a[1])\epsilon + |1 - p|\frac{\epsilon}{2}\Bigr)\,\times\\
		&\times\,\prod_{i = a[2]}^{b[2]}\prod_{l = 1}^{m_{i, 2}}\varphi_{1, 2}\Bigl(z - (l - 1)\epsilon + (i - a[2])\epsilon + |2 - p|\frac{\epsilon}{2}\Bigr) =\\
		&= -\frac{1}{\epsilon}\frac{\Bigl(z - m_{1, 2}\epsilon + |1 - p|\frac{\epsilon}{2}\Bigr)\Bigl(z - (m_{2, 2} - 1)\epsilon + |1 - p|\frac{\epsilon}{2}\Bigr)}{\Bigl(z - m_{1, 1}\epsilon + |1 - p|\frac{\epsilon}{2}\Bigr)\Bigl(z - (m_{1, 1} - 1)\epsilon + |1 - p|\frac{\epsilon}{2}\Bigr)}\,,\\
		\Psi_{\mu, \Upsilon_{p,\lambda}}^{(n - 1)} = &-\frac{1}{\epsilon}\prod_{i = a[n - 1]}^{b[n - 1]}\prod_{l = 1}^{m_{i, n - 1}}\varphi_{n - 1, n - 1}\Bigl(z - (l - 1)\epsilon + (i - a[n - 1])\epsilon + |n - 1 - p|\frac{\epsilon}{2}\Bigr)\,\times\\
		&\times\, \prod_{i = a[n - 2]}^{b[n - 2]}\prod_{l = 1}^{m_{i, n - 2}}\varphi_{n - 1, n - 2}\Bigl(z - (l - 1)\epsilon + (i - a[n - 2])\epsilon + |n - 2 - p|\frac{\epsilon}{2}\Bigr) = \\
		&= -\frac{1}{\epsilon}\frac{\Bigl(z - m_{n - p - 1, n - 2}\epsilon + |n - 1 - p|\frac{\epsilon}{2}\Bigr)\Bigl(z - (m_{n - p, n - 2} - 1)\epsilon + |n - 1 - p|\frac{\epsilon}{2}\Bigr)}{\Bigl(z - m_{n - p, n - 1}\epsilon + |n - 1 - p|\frac{\epsilon}{2}\Bigr)\Bigl(z - m_{n - p, n - 1}\epsilon + |n - 1 - p|\frac{\epsilon}{2}\Bigr)}\,,
	\end{aligned}
\end{equation}
where we assumed that $p \ne 1$. If $p = 1$, the eigenfunction $\Psi_{\mu, \Upsilon_{p,\lambda}}^{(1)}(z)$ includes the additional factor:
\begin{equation}\label{sl(n) first eigenfunction when p = 1}
	\begin{aligned}
		&\Psi^{(1)}_{\mu, \Upsilon_{1, \lambda}}(z) \to \frac{z - \lambda \epsilon}{z}\cdot\Psi^{(1)}_{\mu, \Upsilon_{p, \lambda}}(z)\,,\\
		&\Psi^{(1)}_{\mu, \Upsilon_{1, \lambda}}(z) = -\frac{1}{\epsilon}\frac{\bigl(z - \lambda\epsilon\bigr)\bigl(z - (m_{2, 2} - 1)\epsilon\bigr)}{\bigl(z - m_{1, 1}\epsilon\bigr)\bigl(z - (m_{1, 1} - 1)\epsilon\bigr)}\,,
	\end{aligned}
\end{equation}
which in practice substitutes $z$ by $z - \lambda\epsilon$.

The remaining eigenfunctions for $k \ne p\,, n-p$ read as follows:
\begin{equation}\label{sl(n) eigenfunctions}
	\begin{aligned}
		\Psi_{\mu, \Upsilon_{p,\lambda}}^{(k)} = &-\frac{1}{\epsilon}\prod_{i = a[k]}^{b[k]}\prod_{l = 1}^{m_{i, k}}\varphi_{k,k}\Bigl(z - (l - 1)\epsilon + (i - a[k])\epsilon + |k - p|\frac{\epsilon}{2}\Bigr)\,\times\\
		&\times\, \prod_{i = a[k + 1]}^{b[k + 1]}\prod_{l = 1}^{m_{i, k + 1}}\varphi_{k, k + 1}\Bigl(z - (l - 1)\epsilon + (i - a[k + 1])\epsilon + |k + 1 - p|\frac{\epsilon}{2}\Bigr)\,\times\\
		&\prod_{i = a[k - 1]}^{b[k - 1]}\prod_{l = 1}^{m_{i, k - 1}}\varphi_{k, k - 1}\Bigl(z - (l - 1)\epsilon + (i - a[k - 1])\epsilon + |k - 1 - p|\frac{\epsilon}{2}\Bigr) =\\
		&= -\dfrac{1}{\epsilon}\prod_{i = a[k]}^{b[k]}\frac{1}{\Bigl(z - m_{i, k}\epsilon + (i - a[k])\epsilon + |k - p|\frac{\epsilon}{2}\Bigr)\Bigl(z - (m_{i,k} - 1)\epsilon + (i - a[k])\epsilon + |k - p|\frac{\epsilon}{2}\Bigr)}\,\times\\
		&\times\prod_{i = a[k - 1]}^{b[k - 1]}\Bigl(z - m_{i, k - 1}\epsilon + (i - a[k - 1])\epsilon + |k - p|\frac{\epsilon}{2}\Bigr)\prod_{i = a[k + 1]}^{b[k + 1]}\Bigl(z - (m_{i, k + 1} - 1)\epsilon + (i - a[k + 1])\epsilon + |k - p|\frac{\epsilon}{2}\Bigr)\,.
	\end{aligned}
\end{equation}

Finally, we present the explicit form of eigenfunctions for $k = p$ and $k = n - p$:
\begin{equation}\label{sl(n) eigenfunctions k = p}
	\begin{aligned}
		\Psi_{\mu, \Upsilon_{p,\lambda}}^{(p)} = &-\frac{1}{\epsilon}(z - \lambda\epsilon)\prod_{i = 1}^{b[p]}\frac{1}{\Bigl(z - m_{i, p}\epsilon + (i - 1)\epsilon\Bigr)\Bigl(z - (m_{i,p} - 1)\epsilon + (i - 1)\epsilon\Bigr)}\,\times\\
		&\times\prod_{i = a[p - 1]}^{b[p - 1]}\Bigl(z - m_{i, p - 1}\epsilon + (i - a[p - 1])\epsilon\Bigr)\prod_{i = a[p + 1]}^{b[p + 1]}\Bigl(z - (m_{i, p + 1} - 1)\epsilon + (i - a[p + 1])\epsilon\Bigr)\,,\\
		\Psi_{\mu, \Upsilon_{p,\lambda}}^{(n - p)} = &-\frac{1}{\epsilon}\Bigl(z - n\frac{\epsilon}{2}\Bigr)\prod_{i = 1}^{b[p]}\frac{1}{\Bigl(z - m_{i, p}\epsilon + (i - 1)\epsilon\Bigr)\Bigl(z - (m_{i,p} - 1)\epsilon + (i - 1)\epsilon\Bigr)}\,\times\\
		&\times\prod_{i = a[p - 1]}^{b[p - 1]}\Bigl(z - m_{i, p - 1}\epsilon + (i - a[p - 1])\epsilon\Bigr)\prod_{i = a[p + 1]}^{b[p + 1]}\Bigl(z - (m_{i, p + 1} - 1)\epsilon + (i - a[p + 1])\epsilon\Bigr)\,.
	\end{aligned}
\end{equation}
If we have the case when $p = n - p$, then:
\begin{equation}
	\begin{aligned}
		\Psi_{\mu, \Upsilon_{p,\lambda}}^{(p)} = &-\frac{1}{\epsilon}\bigl(z - \lambda\epsilon\bigr)\Bigl(z - n\frac{\epsilon}{2}\Bigr)\prod_{i = 1}^{b[p]}\frac{1}{\Bigl(z - m_{i, p}\epsilon + (i - 1)\epsilon\Bigr)\Bigl(z - (m_{i,p} - 1)\epsilon + (i - 1)\epsilon\Bigr)}\,\times\\
		&\times\prod_{i = a[p - 1]}^{b[p - 1]}\Bigl(z - m_{i, p - 1}\epsilon + (i - a[p - 1])\epsilon\Bigr)\prod_{i = a[p + 1]}^{b[p + 1]}\Bigl(z - (m_{i, p + 1} - 1)\epsilon + (i - a[p + 1])\epsilon\Bigr)\,.
	\end{aligned}
\end{equation}

\subsubsection*{Amplitudes}

Our next step is to apply equivariant techniques and get the matrix coefficients of the raising and lowering operators in the ansatz \eqref{sl(n) ansatz}. In order to evaluate the corresponding Euler classes, we need to construct the vacuum expectation values of the fields in a fixed point $\mu$.

We start with the fields $C_{k}$. These matrices act from $V_{k}$ to $V_{k}$. As we have mentioned, the vector space $V_{k}$ can be decomposed as follows:
\begin{equation}\label{sl(n) V_k decomposition}
	V_{k} = \bigoplus_{i = a[k]}^{b[k]}\nu_{i, k}\,.
\end{equation} 
This allows us to determine the general structure of $C_{k}$ as block matrices. The corresponding vacuum expectation values read:
\begin{equation}
	\begin{aligned}
		&\bar{C}_{1} = C(m_{1, 1}) = C(n_{1})\,, \\
		&\bar{C}_{k} = \left(\begin{array}{c|c|c|c|c}
			\mathbf{C(m_{a[k], k})} & O(m_{a[k], k}, m_{a[k] + 1, k}) & \cdots & O(m_{a[k], k}, m_{b[k] - 1, k}) & O(m_{a[k], k}, m_{b[k], k})\\ \hline
			O(m_{a[k] + 1, k}, m_{a[k], k}) & \mathbf{C(m_{a[k] + 1, k})} & \cdots & O(m_{a[k] + 1, k}, m_{b[k] - 1, k}) & O(m_{a[k] + 1, k}, m_{b[k], k}) \\ \hline
			\vdots & \vdots & \ddots & \vdots & \vdots  \\ \hline
			O(m_{b[k] - 1, k}, m_{a[k], k}) & O(m_{b[k]- 1, k}, m_{a[k] + 1, k}) & \cdots & \mathbf{C(m_{b[k] - 1, k})} & O(m_{b[k] - 1, k}, m_{b[k], k})\\ \hline
			O(m_{b[k], k}, m_{a[k], k}) & O(m_{b[k], k}, m_{a[k] + 1, k}) & \cdots & O(m_{b[k] - 1, k}, m_{b[k], k}) & \mathbf{C(m_{b[k], k})}
		\end{array}\right)\,.
	\end{aligned}
\end{equation}

In order to construct the matrices $\bar{A}_{k}$ and $\bar{B}_{k}$, we need to split the region $k \in [1, n-1]$ into three pieces. They correspond to the three different regions in \eqref{sl(n) different paths}.

\begin{enumerate}
	\item {\bf The First Case} $k \leqslant p - 1$:
	This area corresponds to the left triangle in \eqref{sl(n) different paths}. The ranges read:
	\begin{equation}
		\begin{aligned}
			&a[k] = 1\,, \quad b[k] = k\,,\\
			&a[k + 1] = 1 \quad b[k + 1] = k + 1\,.
		\end{aligned}
	\end{equation}
	In this case the fields acquire the following vacuum expectation values:
	\begin{equation}
		\begin{aligned}
			&\bar{A}_{k} = \left(\begin{array}{c|c|c|c|c}
				O(m_{1, k + 1}, m_{1, k}) & O(m_{1, k+1}, m_{2, k}) & \cdots & O(m_{1, k + 1}, m_{k - 1, k}) & O(m_{1, k + 1}, m_{k, k}) \\ \hline
				\mathbf{I(m_{2, k + 1}, m_{1, k})} & O(m_{2, k + 1}, m_{2, k}) & \cdots & O(m_{2, k + 1}, m_{k - 1, k}) & O(m_{2, k + 1}, m_{k, k}) \\ \hline
				\vdots & \ddots & \ddots & \vdots & \vdots \\ \hline
				O(m_{k -1, k + 1}, m_{1, k}) & O(m_{k - 1, k + 1}, m_{2, k}) & \cdots & O(m_{k - 1, k + 1}, m_{k - 1, k}) & O(m_{k - 1, k + 1}, m_{k, k}) \\ \hline
				O(m_{k, k + 1}, m_{1, k}) & O(m_{k, k + 1}, m_{2, k}) & \cdots & \mathbf{I(m_{k, k + 1}, m_{k - 1, k})} & O(m_{k, k + 1}, m_{k, k})\\ \hline
				O(m_{k + 1, k + 1}, m_{1, k}) & O(m_{k + 1, k + 1}, m_{2, k}) & \cdots & O(m_{k + 1, k + 1}, m_{k-1, k}) & \mathbf{I(m_{k + 1, k + 1}, m_{k, k})}
			\end{array}\right)\,,\\ \\
			&\bar{B}_{k} = \left(\begin{array}{c|c|c|c|c}
				\mathbf{I(m_{1, k}, m_{1, k + 1})} & O(m_{1, k}, m_{2, k + 1}) & \cdots & O(m_{1, k}, m_{k, k + 1}) & O(m_{1, k}, m_{k + 1, k + 1})\\ \hline
				O(m_{2, k}, m_{1, k + 1}) & \mathbf{I(m_{2, k}, m_{2, k + 1})} & \cdots & O(m_{2, k}, m_{k, k + 1}) & O(m_{2, k}, m_{k + 1, k + 1})\\ \hline
				\vdots & \ddots & \ddots & \vdots & \vdots \\ \hline
				O(m_{k, k}, m_{1, k + 1}) & O(m_{k, k}, m_{2, k + 1}) & \cdots & \mathbf{I(m_{k, k}, m_{k, k + 1})} & O(m_{k, k}, m_{k + 1, k + 1})
			\end{array}\right)\,.
		\end{aligned}
	\end{equation}
	\item {\bf The Second Case} $p \leqslant k \leqslant n - p - 1$:
	This area corresponds to the parallelogram in the center. The ranges read:
	\begin{equation}
		\begin{aligned}
			&a[k] = k - p + 1\,, \quad b[k] = k\,,\\
			&a[k + 1] = k - p + 2\,, \quad b[k + 1] = k + 1\,.
		\end{aligned}
	\end{equation}
	We can notice that the number of the vertical blocks is equal to the number of the horizontal ones:
	\begin{equation}
		b[k] - a[k] + 1 = b[k + 1] - a[k + 1] + 1 = p
	\end{equation}
	It means that the blocks of the matrices $\bar{A}_{k}$ and $\bar{B}_{k}$ form a ``square" matrix themselves. Therefore, the vacuum expectation values of the fields read as follows:
	\begin{equation}
		\resizebox{\textwidth}{!}{$
		\begin{aligned}
			&\bar{A}_{k} = \left(\begin{array}{c|c|c|c}
				\mathbf{I(m_{a[k + 1], k + 1}, m_{a[k], k})} & O(m_{a[k + 1], k + 1}, m_{a[k] + 1, k}) & \cdots & O(m_{a[k + 1], k + 1}, m_{k, k})\\ \hline
				O(m_{a[k + 1] + 1, k + 1}, m_{a[k], k}) & \mathbf{I(m_{a[k + 1] + 1, k + 1}, m_{a[k] + 1, k})} & \cdots & O(m_{a[k + 1] + 1, k + 1}, m_{k, k}) \\ \hline
				\vdots & \vdots & \ddots & \vdots \\ \hline
				O(m_{k + 1, k + 1}, m_{a[k], k}) & O(m_{k + 1, k + 1}, m_{a[k] + 1, k}) & \cdots & \mathbf{I(m_{k + 1, k + 1}, m_{k, k})}
			\end{array}\right)\,,\\ \\
			&\bar{B}_{k} = \left(\begin{array}{c|c|c|c|c}
				O(m_{a[k], k}, m_{a[k + 1], k + 1}) & O(m_{a[k], k}, m_{a[k + 1] + 1, k + 1}) & \cdots & O(m_{a[k], k}, m_{k, k + 1}) & O(m_{a[k], k}, m_{k + 1, k + 1}) \\ \hline
				\mathbf{I(m_{a[k] + 1, k}, m_{a[k + 1], k + 1})} & O(m_{a[k] + 1, k}, m_{a[k + 1] +1, k + 1}) & \cdots & O(m_{a[k] + 1, k}, m_{k, k + 1}) & O(m_{a[k] + 1, k}, m_{k + 1, k + 1}) \\ \hline
				\vdots & \ddots & \ddots & \vdots & \vdots \\ \hline
				O(m_{k - 1, k}, m_{a[k + 1], k + 1}) & O(m_{k - 1, k}, m_{a[k + 1] + 1, k + 1}) & \cdots & O(m_{k - 1, k}, m_{k, k + 1}) & O(m_{k - 1, k}, m_{k + 1, k + 1}) \\ \hline
				O(m_{k, k}, m_{a[k + 1], k + 1}) & O(m_{k, k}, m_{a[k + 1] + 1, k + 1}) & \cdots & \mathbf{I(m_{k, k}, m_{k - 1, k})} & O(m_{k, k}, m_{k + 1, k + 1})\\
			\end{array}\right)\,.
		\end{aligned}
		$}
	\end{equation}
	\item {\bf The Third Case} $k \geqslant n - p$:
	This is the final area in \eqref{sl(n) different paths}. The ranges read:
	\begin{equation}
		\begin{aligned}
			&a[k] = k - p + 1, \quad b[k] = n - p = \tilde{p},\\
			&a[k + 1] = k - p + 2 \quad b[k + 1] = n - p = \tilde{p}.
		\end{aligned}
	\end{equation}
	This sector is connected to the first sector via quiver symmetry \eqref{sl(n) quiver symmetry Z_2}. This means that we can get the vacuum expectation values by interchanging the fields $A_{k}$ and $B_{k}$ from the first area and applying the corresponding change of the dimension parameters $m_{i, k}$.
	
	Finally, we end up with the result:
	\begin{equation}
		\resizebox{\textwidth}{!}{$
		\begin{aligned}
			&\bar{A}_{k} = \left(\begin{array}{c|c|c|c|c}
				\mathbf{I(m_{a[k + 1], k + 1}, m_{a[k], k})} & O(m_{a[k + 1], k + 1}, m_{a[k] + 1, k}) & \cdots & O(m_{a[k + 1], k + 1}, m_{\tilde{p} - 1, k}) & O(m_{a[k + 1], k + 1}, m_{\tilde{p}, k})\\ \hline
				O(m_{a[k + 1] + 1, k + 1}, m_{a[k], k}) & \mathbf{I(m_{a[k + 1] + 1, k + 1}, m_{a[k] + 1, k})} & \cdots & O(m_{a[k + 1] + 1, k + 1}, m_{\tilde{p} - 1, k}) & O(m_{a[k + 1] + 1, k + 1}, m_{\tilde{p}, k})\\ \hline
				\vdots & \ddots & \ddots & \vdots & \vdots \\ \hline
				O(m_{\tilde{p}, k + 1}, m_{a[k], k}) & O(m_{\tilde{p}, k + 1}, m_{a[k] + 1, k}) & \cdots & \mathbf{I(m_{\tilde{p}, k + 1}, m_{\tilde{p} - 1, k})} & O(m_{\tilde{p}, k + 1}, m_{\tilde{p}, k})
			\end{array}\right)\,,\\ \\
			&\bar{B}_{k} = \left(\begin{array}{c|c|c|c|c}
				O(m_{a[k], k}, m_{a[k + 1], k + 1}) & O(m_{a[k], k}, m_{a[k + 1] + 1, k + 1}) & \cdots & O(m_{a[k], k}, m_{\tilde{p} - 1, k + 1}) & O(m_{a[k], k}, m_{\tilde{p}, k + 1}) \\ \hline
				\mathbf{I(m_{a[k] + 1, k}, m_{a[k + 1], k + 1})} & O(m_{a[k] + 1, k}, m_{a[k + 1] + 1, k + 1}) & \cdots & O(m_{a[k] + 1, k}, m_{\tilde{p} - 1, k + 1}) & O(m_{a[k] + 1, k}, m_{\tilde{p}, k + 1}) \\ \hline
				\vdots & \ddots & \ddots & \vdots & \vdots \\ \hline
				O(m_{\tilde{p} - 2, k}, m_{a[k + 1], k + 1}) & O(m_{\tilde{p} - 2, k}, m_{a[k + 1] + 1, k + 1}) & \cdots & O(m_{\tilde{p} - 2, k}, m_{\tilde{p} - 1, k + 1}) & O(m_{\tilde{p} - 2, k}, m_{\tilde{p}, k + 1}) \\ \hline
				O(m_{\tilde{p} - 1, k}, m_{a[k + 1], k + 1}) & O(m_{\tilde{p} - 1, k}, m_{a[k + 1] + 1, k + 1}) & \cdots & \mathbf{I(m_{\tilde{p} - 1, k}, m_{\tilde{p} - 1, k + 1})} & O(m_{\tilde{p} - 1, k}, m_{\tilde{p}, k + 1})\\ \hline
				O(m_{\tilde{p}, k}, m_{a[k + 1], k + 1}) & O(m_{\tilde{p}, k}, m_{a[k + 1] + 1, k + 1}) & \cdots & O(m_{\tilde{p}, k}, m_{\tilde{p} - 1, k + 1}) & \mathbf{I(m_{\tilde{p}, k}, m_{\tilde{p}, k + 1})}
			\end{array}\right)\,.
		\end{aligned}
		$}
	\end{equation}
\end{enumerate}

Next, we apply the algorithm described in section \ref{subsec: Equivariant Matrix Coefficients} using the vacuum expectation values of the fields at the fixed points defined above.

First, we point out that:
\begin{equation}
	{\bf E}_{\Upsilon_{p,\lambda}}[\mu \to \mu_{p, p + 1}^{1}] = -\frac{1}{\epsilon}\,,
\end{equation}
which is a result of the choice of the normalization that is provided by equivariant integration. For later convenience we introduce shifted dimensions following \cite{Gelfand:1950ihs}:
\begin{equation}
	l_{i, k} = m_{i, k} - i\,.
\end{equation}
Then, after calculating the Euler classes, we end up with the following matrix coefficients for $k \ne p$:
\begin{equation}\label{sl(n) equiv coefficients k != p}
	\begin{aligned}
		&{\bf E}_{\Upsilon_{p,\lambda}}[\mu \to \mu_{k, k + 1}^{j}] = \frac{\prod\lm_{i = 1}^{j}(l_{i, k + 1} - l_{j, k})\prod\lm_{i = 1}^{j - 1}(l_{i, k - 1} - l_{j, k} - 1)}{\prod\lm_{i = 1}^{j - 1}(l_{i, k} - l_{j, k})(l_{i,k} - l_{j, k} - 1)}\frac{1}{l_{j, k}\epsilon + a[k]\epsilon - |k - p|\frac{\epsilon}{2}}\,,\\
		&{\bf F}_{\Upsilon_{p,\lambda}}[\mu \to \mu_{k, k - 1}^{j}] = \frac{(-1)\prod\lm_{i = j + 1}^{k + 1}(l_{i, k + 1} - l_{j, k} + 1)\prod\lm_{i = j}^{k - 1}(l_{i, k - 1} - l_{j, k})}{\prod\lm_{i = j + 1}^{k}(l_{i, k} - l_{j, k} + 1)(l_{i, k} - l_{j, k})}\Bigl((l_{j, k} - 1)\epsilon + a[k]\epsilon - |k - p|\frac{\epsilon}{2}\Bigr)\,.
	\end{aligned}
\end{equation}
For the case $k = p$ we have:
\begin{equation}\label{sl(n) equiv coefficients k = p}
	\begin{aligned}
		&{\bf E}_{\Upsilon_{p, \lambda}}[\mu \to \mu_{p, p + 1}^{j}] = -\frac{\prod\lm_{i = 2}^{j}(l_{i, p + 1} - l_{j, p})\prod\lm_{i = 1}^{j - 1}(l_{i, p - 1} - l_{j, p} - 1)}{\prod\lm_{i = 1}^{j - 1}(l_{i, p} - l_{j, p})(l_{i, p} - l_{j, p} - 1)}\frac{1}{\epsilon}\,,\\
		&{\bf F}_{\Upsilon_{p, \lambda}}[\mu \to \mu_{p, p - 1}^{j}] = \frac{(l_{1, p + 1} - l_{j, p})\prod\lm_{i = j + 1}^{p + 1}(l_{i, p + 1} - l_{j, p} + 1)\prod\lm_{i = j}^{p - 1}(l_{i, p - 1} - l_{j, p})}{\prod\lm_{i = j + 1}^{p}(l_{i, p} - l_{j, p} + 1)(l_{i, p} - l_{j, p})}\epsilon\,.
	\end{aligned}
\end{equation}

Again, we can normalize the coefficients using the formula \eqref{Rescaling Equiv-to-Root}. At the zero level of the algebra we get:
\begin{equation}
	\begin{aligned}
		&e^{(k)}_{0}|\mu\rangle = \sum_{j = a[k]}^{b[k]}{\bf E}^{(root)}_{\Upsilon_{p,\lambda}}[\mu \to \mu_{k, k + 1}^{j}]\,|\mu_{k, k + 1}^{j}\rangle =  \sum_{j = a[k]}^{b[k]}a_{k, k + 1}^{j}\,|\mu_{k, k + 1}^{j}\rangle\,,\\
		&f^{(k)}_{0}|\mu\rangle = \sum_{j = a[k]}^{b[k]}{\bf F}^{(root)}_{\Upsilon_{p, \lambda}}[\mu \to \mu_{k, k - 1}^{j}]\,|\mu_{k, k - 1}^{j}\rangle = \sum_{j = a[k]}^{b[k]}b_{k, k - 1}^{j}\,|\mu_{k, k - 1}^{j}\rangle\,,
	\end{aligned}
\end{equation}
where:
\begin{equation}\label{gelfand-coefficients}
	\begin{aligned}
		&a_{k, k + 1}^{j} = \left[(-1)\frac{\prod\lm_{i = 1}^{k+1}(l_{i, k+1} - l_{j, k})\prod\lm_{i = 1}^{k-1}(l_{i, k-1} - l_{j, k} - 1)}{\prod\lm_{i\ne j}^{k}(l_{i, k} - l_{j, k})(l_{i, k} - l_{j, k} - 1)}\right]^{\frac{1}{2}}\,, \\
		&b_{k, k - 1}^{j} = \left[(-1)\frac{\prod\lm_{i = 1}^{k + 1}(l_{i, k + 1} - l_{j, k} + 1)\prod\lm_{i = 1}^{k-1}(l_{i, k-1} - l_{j, k})}{\prod\lm_{i\ne j}^{k}(l_{i, k} - l_{j, k} + 1)(l_{i, k} - l_{j, k})}\right]^{\frac{1}{2}}\,.
	\end{aligned}
\end{equation}
These formulas are precisely the formulas introduced by Gelfand in \cite{Gelfand:1950ihs} for the representations of $\fs\fl_{n}$ algebras.

We still need to check the hysteresis relations \eqref{general hysteresis} to prove that the coefficients \eqref{sl(n) equiv coefficients k != p} and \eqref{sl(n) equiv coefficients k = p} satisfy the Yangian algebraic relations \eqref{sl(n) Yangian relations in modes}. Empirically, these relations are satisfied in concrete examples for low $n$ and take the form:
\begin{equation}\label{sl(n) hysteresis}
	\begin{aligned}
		{\bf E}_{\Upsilon_{p,\lambda}}[\mu \to \mu_{k, k + 1}^{j}]{\bf F}_{\Upsilon_{p,\lambda}}[\mu_{k, k + 1}^{j} \to (\mu_{k, k + 1}^{j})^{j}_{k, k - 1}] &= (-1)\frac{\prod\lm_{i = 1}^{k+1}(l_{i, k+1} - l_{j, k})\prod\lm_{i = 1}^{k-1}(l_{i, k-1} - l_{j, k} - 1)}{\prod\lm_{i\ne j}^{k}(l_{i, k} - l_{j, k})(l_{i, k} - l_{j, k} - 1)} = \\
		&= \mathop{\rm res}_{z = l_{j, k}\epsilon + a[k]\epsilon - |k - p|\frac{\epsilon}{2}}\Psi_{\mu, \Upsilon_{p,\lambda}}^{(k)}(z)\,,\\
		\frac{{\bf E}_{\Upsilon_{p, \lambda}}\bigl[\mu \to \mu_{k, k + 1}^{i}\bigr]{\bf E}_{\Upsilon_{p,\lambda}}\bigl[\mu_{k, k + 1}^{i} \to (\mu_{k, k + 1}^{i})_{k, k + 1}^{j}\bigr]}{{\bf E}_{\Upsilon_{p,\lambda}}\bigl[\mu \to \mu_{k, k + 1}^{j}\bigr]{\bf E}_{\Upsilon_{p, \lambda}}\bigl[\mu_{k, k + 1}^{j} \to (\mu_{k,k + 1}^{j})_{k, k + 1}^{i}\bigr]} &= \frac{l_{i, k} - l_{j, k} + 1}{l_{i, k} - l_{j, k} - 1} = \\
		&= \varphi_{k, k}\Bigl(l_{j, k}\epsilon + a[k]\epsilon - |k - p|\frac{\epsilon}{2} - l_{i, k}\epsilon - a[k]\epsilon + |k - p|\frac{\epsilon}{2}\Bigr)\,,\\
		\frac{{\bf E}_{\Upsilon_{p, \lambda}}\bigl[\mu \to \mu_{k, k + 1}^{i}\bigr]{\bf E}_{\Upsilon_{p,\lambda}}\bigl[\mu_{k, k + 1}^{i} \to (\mu_{k, k + 1}^{i})_{t, t + 1}^{j}\bigr]}{{\bf E}_{\Upsilon_{p,\lambda}}\bigl[\mu \to \mu_{t, t + 1}^{j}\bigr]{\bf E}_{\Upsilon_{p, \lambda}}\bigl[\mu_{t, t + 1}^{j} \to (\mu_{t, t + 1}^{j})_{k, k + 1}^{i}\bigr]} &= \frac{l_{j, t}  + a[t] - |t - p|\frac{1}{2} - l_{i, k} - a[k] + |k - p|\frac{1}{2} + \frac{1}{2}}{l_{j, t} + a[t] - |t - p|\frac{1}{2} - l_{i, k} - a[k] + |k - p|\frac{1}{2} - \frac{1}{2}} = \\
		&= \varphi_{t, k}\Bigl(l_{j, t}\epsilon + a[t]\epsilon - |t - p|\frac{\epsilon}{2} - l_{i, k}\epsilon - a[k]\epsilon + |k - p|\frac{\epsilon}{2}\Bigr)\,,\\
	\end{aligned}
\end{equation}
where $k,\,t\in\, Q_{0}$, $k \ne t$. The relations involving ${\bf F}_{\Upsilon_{p,\lambda}}[\mu \to \mu_{k, k - 1}^{j}]$ take a similar form; therefore, we do not present them here.

\subsection{General Comments}\label{subsec: Comments}

In this section, we give a few closing comments regarding the algebras $\myY(\fs\fl_{n})$ and their representations.

\subsubsection{Connection to Drinfeld Yangians}

Having constructed the algebras $\myY(\fs\fl_{n})$ combining the quiver approach and equivariant integration, we are now interested in how these Yangians are connected to the ones introduced by Drinfeld in \cite{Drinfeld:1985rx, Drinfeld:1987sy}.

One might notice that our algebras are one-parametric: $\myY(\fs\fl_{n}) = \myY_{\epsilon}(\fs\fl_{n})$. We claim, however, that one could scale away this parameter completely. Indeed, let us consider the algebra isomorphism \eqref{Z symmetry} and choose the parameter of the symmetry $\sigma$ to be:
\begin{equation}
	\begin{aligned}
		&\sigma = \frac{\epsilon}{\epsilon'}\,,\\
		&e^{(a)}_{n} \to \left(\dfrac{\epsilon}{\epsilon'}\right)^{n}e^{(a)}_{n}\,,\\
		&f^{(a)}_{n} \to \left(\dfrac{\epsilon}{\epsilon'}\right)^{n}f^{(a)}_{n}\,,\\
		&\psi^{(a)}_{n} \to \left(\dfrac{\epsilon}{\epsilon'}\right)^{n}\psi^{(a)}_{n}\,.
	\end{aligned}
\end{equation}
This transformation is an isomorphism of the algebras $\myY_{\epsilon}(\fs\fl_{n})$ and $\myY_{\epsilon'}(\fs\fl_{n})$. Therefore, one can always fix the parameter to be $\epsilon' = 1$, provided that $\epsilon \ne 0$. This also proves that the Yangians we have constructed are isomorphic to Drinfeld Yangians for the algebras $\fs\fl_{n}$.

In fact, one could include the deformation parameter into the Drinfeld definitions directly. In the same way, the resulting algebra is essentially independent of this parameter \cite{chari1994guide}. In the physical literature, however, the equivariant parameters are considered as complex mass (flavor fugacity) parameters, which means $\epsilon$ is a dimensional parameter; see, for example, 
\cite{Galakhov:2021xum, Chen:2025xoe}. Therefore, one should in general keep $\epsilon \ne 1$.

Now, let us return to the parameters of the algebra $\myY(\fs\fl_{n})$. To preserve the crystal structure of the states, we have included the second parameter $h$. We, however, treated it as an effective additional parameter and always set it to be $0$ in the calculations. One could ask if this parameter adds a new structure to the algebra.

Indeed, the algebra relations \eqref{sl(n) Yangian relations in modes} modify. We present a few examples of these relations:
\begin{equation}
	\begin{aligned}
		&[e_{n + 1}^{(a)}, e^{(a)}_{m}] - [e^{(a)}_{n}, e^{(a)}_{m + 1}] = \epsilon\{e_{n}^{(a)}, e_{m}^{(a)}\}\,,\\
		&[e_{n + 1}^{(a)}, e_{m}^{(a + 1)}] - [e^{(a)}_{n}, e^{(a + 1)}_{m + 1}] = -\frac{\epsilon}{2}\{e^{(a)}_{n}, e^{(a + 1)}_{m}\} - h[e_{n}^{(a)}, e^{(a + 1)}_{m}]\,,\\
		&[e_{n + 1}^{(a + 1)}, e^{(a)}_{m}] - [e^{(a + 1)}_{n}, e_{m + 1}^{(a)}] = -\frac{\epsilon}{2}\{e^{(a + 1)}_{n}, e^{(a)}_{m}\} + h[e^{(a + 1)}_{n}, e_{m}^{(a)}]\,,
	\end{aligned}
\end{equation}
where the similar behavior manifests in other relations as well. One could notice that these relations are not exactly of Drinfeld type and involve the commutator term. The commutator terms arise from the bond functions:
\begin{equation}
	\varphi_{a, a + 1}(z - w) = \varphi_{a, a - 1}(z - w) = \frac{z - w - \frac{\epsilon}{2} + h}{z - w + \frac{\epsilon}{2} + h}\,.
\end{equation}
One, however, could apply the spectral shift \eqref{spectral shift} in the position $z \to z - h$ at a vertex $a$ and completely eliminate the parameter $h$ from these factors:
\begin{equation}
	\varphi_{a, a + 1}(z - w) = \varphi_{a, a - 1}(z - w) = \frac{z - w - \frac{\epsilon}{2}}{z - w + \frac{\epsilon}{2}}\,.
\end{equation}
In fact, these functions are exactly the ones we have used in our calculations \eqref{sl(n) bonding factors} to construct the algebras $\myY_{\epsilon}(\fs\fl_{n})$. This proves that the algebras $\myY_{\epsilon, h}(\fs\fl_{n})$ and $\myY_{\epsilon}(\fs\fl_{n})$ are isomorphic. This essentially means that the parameter $h$ is indeed effective and does not contribute to the algebraic calculations, as was suggested in section \ref{subsec: Quiver Dynkin Diagrams}.

\subsubsection{Embedding structure}

Mathematically, the triangular form of GT bases represents the hidden structure of an algebra. Namely, one should consider a given Lie algebra $\fA_{n}$ not as a single object but as a part of a chain of subalgebras:
\begin{equation}
	\fA_{1} \subset \fA_{2} \subset \dots \subset \fA_{n}\,.
\end{equation}
This idea can also be applied to Yangian algebras \cite{molev2002gelfandtsetlinbasesclassicallie, nazarov2000representationsyangiansgelfandzetlinbases, Cherednik1987ANI}. Since we have introduced the Gelfand-Tsetlin bases in our construction as well, we expect that they capture the following embedding structure:
\begin{equation}\label{sl(n) embeddings}
	\myY(\fs\fl_{2}) \subset \myY(\fs\fl_{3}) \subset \myY(\fs\fl_{4}) \subset \dots \subset \myY(\fs\fl_{n}) \subset \dots\,.
\end{equation}
The explicit construction of the embeddings in the general case is quite bulky. Therefore, we limit ourselves in this text to the simplest example:
\begin{equation}
	\myY(\fs\fl_{2}) \subset \myY(\fs\fl_{3})\,,
\end{equation}
where we demonstrate how the embedding structure manifests in terms of the quiver diagrams and the corresponding GT bases.

\begin{figure}[h]
	\begin{center}
		\begin{tikzpicture}
			\draw[postaction=decorate, decoration={markings, mark= at position 0.7 with {\arrow{stealth}}}] (0,0) to[out=20,in=160] node[pos=0.5,above] {$\scriptstyle A$} (1.5,0);
			\draw[postaction=decorate, decoration={markings, mark= at position 0.7 with {\arrow{stealth}}}] (1.5,0) to[out=200,in=340] node[pos=0.5,below] {$\scriptstyle B$} (0,0);
			\draw[postaction=decorate, decoration={markings, mark= at position 0.8 with {\arrow{stealth}}}] (0,0) to[out=60,in=0] (0,0.6) to[out=180,in=120] (0,0);
			\node[above] at (0,0.6) {$\scriptstyle C_1$};
			\begin{scope}[shift={(1.5,0)}]
				\draw[postaction=decorate, decoration={markings, mark= at position 0.8 with {\arrow{stealth}}}] (0,0) to[out=60,in=0] (0,0.6) to[out=180,in=120] (0,0);
				\node[above] at (0,0.6) {$\scriptstyle C_2$};
			\end{scope}
			\begin{scope}[shift = {(0,0)}]
				\draw[postaction=decorate, decoration={markings, mark= at position 0.7 with {\arrow{stealth}}}] (0,-1.2) to[out=100,in=260] node[pos=0.3, left] {$\scriptstyle R_{1}$} (0,0);
				\draw[postaction=decorate, decoration={markings, mark= at position 0.7 with {\arrow{stealth}}}] (0,0) to[out=280,in=80] node[pos=0.7, right] {$\scriptstyle S_{1}$} (0,-1.2);
				\begin{scope}[shift={(0,-1.2)}]
					\draw[fill=burgundy] (-0.08,-0.08) -- (-0.08,0.08) -- (0.08,0.08) -- (0.08,-0.08) -- cycle;
				\end{scope}
			\end{scope}
			\begin{scope}[shift={(1.5, 0)}]
				\draw[postaction=decorate, decoration={markings, mark= at position 0.7 with {\arrow{stealth}}}, dashed] (0,-1.2) to[out=100,in=260] (0,0);
				\draw[postaction=decorate, decoration={markings, mark= at position 0.7 with {\arrow{stealth}}}, dashed] (0,0) to[out=280,in=80] (0,-1.2);
				\begin{scope}[shift={(0,-1.2)}, dashed]
					\draw[fill=gray!30!white] (-0.08,-0.08) -- (-0.08,0.08) -- (0.08,0.08) -- (0.08,-0.08) -- cycle;
				\end{scope}
			\end{scope}
			\draw[fill=\myblue] (0,0) circle (0.08);
			\draw[fill=\mygray] (1.5,0) circle (0.08);
			\node[left] at (0, 0) {$\scriptstyle n_{1}$};
			\node[right] at (1.5, 0) {$\scriptstyle n_{2}$};
			\begin{scope}[shift={(1.5,0)}]
				\draw[rotate=45, red, thick] (0.0-0.3,0.0) -- (0.0+0.3,0.0) (0.0,0.0-0.3) -- (0.0,0.0+0.3);
			\end{scope}
			\begin{scope}[shift={(6, 0)}, rotate=0]
				\draw[postaction=decorate, decoration={markings, mark= at position 0.8 with {\arrow{stealth}}}] (0,0) to[out=60,in=0] (0,0.6) to[out=180,in=120] (0,0);
				\begin{scope}[shift = {(0,0)}]
					\draw[postaction=decorate, decoration={markings, mark= at position 0.7 with {\arrow{stealth}}}] (0,-1.2) to[out=100,in=260] node[pos=0.3, left] {$\scriptstyle R$} (0,0);
					\draw[postaction=decorate, decoration={markings, mark= at position 0.7 with {\arrow{stealth}}}] (0,0) to[out=280,in=80] node[pos=0.7, right] {$\scriptstyle S$} (0,-1.2);
					\begin{scope}[shift={(0,-1.2)}]
						\draw[fill=burgundy] (-0.08,-0.08) -- (-0.08,0.08) -- (0.08,0.08) -- (0.08,-0.08) -- cycle;
					\end{scope}
				\end{scope}
				\draw[fill=\myblue] (0,0) circle (0.08);
				\node[above] at (0,0.6) {$\scriptstyle C$};
				\node[left] at (0, 0) {$\scriptstyle n$};
			\end{scope}
			\draw[postaction=decorate, decoration={markings, mark= at position 1 with {\arrow{stealth}}}] (3, 0) to (4.5, 0);
		\end{tikzpicture}
	\end{center}
	\caption{The reduction of the quiver $\CQ_{3, 1, \lambda}$ to the Jacobian quiver}\label{sl(3) to sl(2) quiver reduction}
\end{figure}
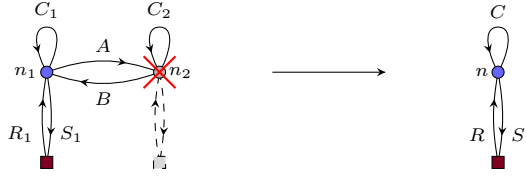

Let us consider the quiver in fig. \ref{sl(3) quiver} that we used to study the representations $\Upsilon_{1, \lambda}$ of the algebra $\myY(\fs\fl_{3})$. Next, we remove the second vertex from this quiver, or more formally, delete the vertex and all the arrows incident to it. We depict the procedure in fig. \ref{sl(3) to sl(2) quiver reduction}, where we have erased the indices for the resulting quiver.
In fact, we end up with the Jacobian quiver \cite{Ginzburg_lect}, which indeed is used to describe the representations of the weight $\lambda$ of the algebra $\myY(\fs\fl_{2})$; see, for example, \cite{Galakhov_2024, Yang:2024ubh, Bykov:2019cst}.

Moreover, this representation can be derived from the representation $\Upsilon_{1,\lambda}$ of $\myY(\fs\fl_{3})$. The crystals \eqref{sl(3) random fixed point} transform into:
\begin{equation}
	\begin{array}{c}
		\begin{tikzpicture}[rotate = 0, scale=0.8]
			\draw[postaction={decorate},
			decoration={markings, mark= at position 0.6 with {\arrow{stealth}}}] (0, 0) -- (2,0);
			\draw[postaction={decorate},
			decoration={markings, mark= at position 0.6 with {\arrow{stealth}}}] (2, 0) -- (4,0);
			\draw[postaction={decorate},
			decoration={markings, mark= at position 0.6 with {\arrow{stealth}}}] (4, 0) -- (6,0);
			\draw[postaction={decorate},
			decoration={markings, mark= at position 0.6 with {\arrow{stealth}}}] (6, 0) -- (8,0);
			\draw[postaction={decorate},
			decoration={markings, mark= at position 0.6 with {\arrow{stealth}}}] (8, 0) -- (9,0);
			\draw[fill=white, draw=burgundy] circle (0.4);
			\draw[fill=\myblue] (0, 0) circle (0.2);
			\draw[fill=\myblue] (2, 0) circle (0.2);
			\draw[fill=\myblue] (4, 0) circle (0.2);
			\draw[fill = \myblue] (6, 0) circle (0.2);
			\draw[fill = \myblue] (8, 0) circle (0.2);
			\draw (0, 0.6) node {$\scriptstyle R$};
			\draw (2, 0.6) node {$\scriptstyle CR$};
			\draw (4, 0.6) node {$\scriptstyle C^{2}R$};
			\draw (6, 0.6) node {$\scriptstyle C^{3}R$};
			\draw (8, 0.6) node {$\scriptstyle \dots$};
		\end{tikzpicture}
	\end{array}
\end{equation}
We can acquire the corresponding GT bases by removing the last row in $\myY(\fs\fl_{3})$ GT bases \eqref{sl(3) states picture} and setting $n_{2} = 0$:
\begin{equation}\label{sl(3) to sl(2) GT reduction}
	|n_{1}, n_{2}\rangle = \quad
	\begin{tikzpicture}[baseline={(0,-0.0)}, every node/.style={inner sep=1pt}]
		\begin{scope}[shift={(0,0)}]
			\GTrow[0.6]{$\lambda$, $\lambda$, $0$}
			\GTrow[0]{$\lambda$, $n_{2}$}
			\GTrow[-0.6]{$n_{1}$}
		\end{scope}
	\end{tikzpicture}
	\qquad \to \qquad
	|n\rangle = \quad 
	\begin{tikzpicture}[baseline={(0,-0.0)}, every node/.style={inner sep=1pt}]
		\begin{scope}[shift={(0,0)}]
			\GTrow[0.3]{$\lambda$, $0$}
			\GTrow[-0.3]{$n$}
		\end{scope}
	\end{tikzpicture}\,.
\end{equation}
The ansatz \eqref{Y(sl(3)) anzatz} reduces to:
\begin{equation}\label{Y(sl(2)) Ansatz}
	\begin{aligned}
		&e^{[\lambda]}(z)|n\rangle = \dfrac{{\bf E}_{n, n+1}^{[\lambda]}}{z - n\epsilon}|n+1\rangle\,,\\
		&f^{[\lambda]}(z)|n\rangle = \dfrac{{\bf F}_{n,n-1}^{[\lambda]}}{z - (n-1)\epsilon}|n-1\rangle\,,\\
		&\psi^{[\lambda]}(z)|n\rangle = \bPsi_{n}^{[\lambda]}(z)|n\rangle\,,
	\end{aligned}
\end{equation}
where $\bPsi_{n}^{[\lambda]}(z)$ takes the form:
\begin{equation}\label{sl(2) eigenvalues}
	\bPsi_{n}^{[\lambda]}(z) = \frac{(z-\lambda\epsilon)(z+\epsilon)}{(z-n\epsilon)(z-(n-1)\epsilon)}\,.
\end{equation}
Notably, the Euler classes \eqref{sl(3) Euler classes} also allow the reduction:
\begin{equation}
	\begin{aligned}
		&\Eul_{(n_{1}, n_{2})}^{[\lambda, 0]} \quad \to \quad \Eul_n^{[\lambda]} = \prod_{k = 1}^{n}(\lambda + 1 - k)\epsilon (n + 1 - k)\epsilon = \dfrac{n!\lambda!}{(\lambda - n)!}\epsilon^{2n}\,,\\
		&\Eul_{(n_{1}, n_{2}) \to (n_{1} + 1, n_{2})}^{[\lambda, 0]} \quad \to \quad \Eul_{n,n+1}^{[\lambda]} =-\epsilon \prod_{k = 1}^{n}(k - \lambda - 1)\epsilon (k - n - 1)\epsilon = \dfrac{n!\lambda!}{(\lambda - n)!}(-\epsilon)^{2n + 1}\,,\\
		&\Eul_{(n_{1}, n_{2}) \to (n_{1}, n_{2} + 1)} \quad \to \quad 0\,,
	\end{aligned}
\end{equation}
which also affects the matrix coefficients \eqref{sl(3) yangian coefficients}:
\begin{equation}\label{sl(2) coefficients}
	{\bf E}_{n, n+1}^{[\lambda]} = -\dfrac{1}{\epsilon}\,, \quad {\bf F}_{n, n - 1}^{[\lambda]} = -n(\lambda - n +1)\epsilon\,.
\end{equation}

In general, a similar procedure can be implemented to acquire the following reduction:
\begin{equation}
	\begin{aligned}
		\CQ_{n,p,\lambda} \quad &\to \quad \CQ_{n-1,p,\lambda}\,,\\
		\Upsilon_{p,\lambda}|_{\fs\fl_{n}} \quad &\to \quad \Upsilon_{p,\lambda}|_{\fs\fl_{n - 1}}\,.
	\end{aligned}
\end{equation}


\section{Conclusion}

Simple Lie algebras remain one of the central objects in discussions of symmetries from various physical systems to mathematical objects. One of the natural generalizations of these algebras is Yangian algebras.

In this paper, we implement the quiver approach to investigate Yangians associated to Dynkin diagrams of A-type and their representations. While the approach reconstructs the known algebras \cite{Drinfeld:1987sy, chari1994guide}, it has its own specifics. For instance, we consider only the rectangular representations. This is tied to the fact that in order to construct more general representations of $\myY(\fs\fl_{n})$, it seems that one should work with more complicated framings of the quivers. This naturally brings tensor product structure to the appearing representations \cite{Galakhov:2022uyu, Galakhov:2021xum}. It results in more complicated representation theory, where one could get reducible representations instead of irreducible ones. Therefore, one needs to split them into a sum of the irreducible representations, and only after that could extract a required representation, which drastically complicates the analysis. The existence of the choice of framing suitable to describe an arbitrary irreducible representation remains unknown.

The quiver varieties for $\CQ_{n,p,\lambda}$ are still toric, which allowed us to use powerful equivariant localization. While the construction is quite similar to the cases $\myY(\widehat{\fg\fl}_{n})$ related to toric ${\bf CY}_{3}$, we end up with only one equivariant parameter $\epsilon$. However, to preserve the crystal structure of the states, we had to violate the gauge-fixing condition, namely, vertex constraints, prioritizing the no-overlap condition instead \cite{Bao_2025}. The generalization of the construction for other types of Dynkin diagrams is not straightforward. The quiver varieties, even for the graphs of D and E types, are no longer toric \cite{Bao:2023ece, Bao_2025, Li:2023zub}. Therefore, the applicability of the approach and its generalization to the non-simply laced cases requires more careful consideration.

Notably, we showed that the quiver approach seems to inherit important information about the corresponding algebras $\fs\fl_{n}$. For example, the crystal states resemble the famous Gelfand-Tsetlin bases. Namely, the GT patterns match the structure of different elements of Jacobian algebras of the quivers $\CQ_{n,p,\lambda}$, which we have shown on various examples. Moreover, these bases highlight the embedding structure of the algebras and its natural correspondence to the quiver diagrams. As expected, the zero level of our Yangians restores the representations of the corresponding Lie algebras $\fs\fl_{n}$. We also constructed the representations of the Yangian algebras themselves, providing the calculations with explicit examples.

This work is another step towards understanding the connection between quiver diagrams and Yangian algebras. While the question about what types of algebras can be constructed from quiver diagrams remains open, there are quite a few facts that fuel the investigation of this area:
\begin{itemize}
	\item The quivers seem to be closely related to various Dynkin diagrams and encode the information about the Yangian algebras associated to the corresponding Lie algebras, as we have seen in this text;
	\item Drinfeld's second realization of the Yangian algebras offers the description that resembles the Chevalley basis in simple Lie algebras;
	\item The quivers can be constructed for toric diagrams of ${\bf CY}_{3}$ \cite{Li:2020rij}, which are already richer than Dynkin classification;
	\item There are generalizations of the quiver diagrams for ${\bf CY}_{4}$ \cite{Bao_2024, Galakhov_2024csd, Bao_2025}. Moreover, \cite{Bao_2025} develops a generalization of Yangian algebras to double Yangians for these cases; see also the review \cite{Bao:2025dqs}.
\end{itemize}
That all gives us hope that the classification similar to Cartan exists for Yangian algebras as well. The roles of the effective gauge theories and their BPS states, Calabi-Yau manifolds, and the structure of simple Lie algebras, remain unclear in the general picture. However, it was proposed in \cite{Li:2023zub} that the quiver Yangians for any quivers should still give rise to
the BPS algebras.

There are still quite a few more open questions considering the algebras $\myY(\fs\fl_{n})$ aside from the ones mentioned above. For example, the lift to quantum toroidal and elliptic analogues of Yangians. Although the procedure is seemingly transparent \cite{Galakhov:2021vbo, Noshita:2021ldl}, the behavior of crystal representations with different values of the parameters requires an additional investigation.

Finally, we have not raised the questions about various stability chambers of the moduli space of quiver representations \cite{Bryan_2010,  Galakhov:2024foa, Nakajimarnq195, Yamazaki:2010fz} and limited ourselves only to the discussion of a particular chamber where we get crystal states. The questions about the structure of the states in various chambers, as well as the wall-crossing between different chambers, remain relevant.

\section*{Acknowledgments}\label{sec: Acknowledgments}

We would like to express our gratitude to D. Galakhov, A. Morozov and N. Tselousov for enlightening discussions and suggestions regarding our work. We are indebted to D. Galakhov for providing valuable comments on the draft.

Our work is funded within the state assignment of NRC Kurchatov institute.


\bibliographystyle{utphys}
\bibliography{biblio}

\providecommand{\href}[2]{#2}\begingroup\raggedright\begin{thebibliography}{10}

\bibitem{Drinfeld:1985rx}
V.~G. Drinfeld, ``{Hopf algebras and the quantum Yang-Baxter equation},'' {\em
  Sov. Math. Dokl.} {\bfseries 32} (1985) 254--258.

\bibitem{Drinfeld:1987sy}
V.~G. Drinfeld, ``{A New realization of Yangians and quantized affine
  algebras},'' {\em Sov. Math. Dokl.} {\bfseries 36} (1988) 212--216.

\bibitem{chari1994guide}
V.~Chari and A.~Pressley, {\em {A Guide to Quantum Groups}}.
\newblock Cambridge University Press, 1994.
\newblock \url{https://books.google.ru/books?id=_fRhQgAACAAJ}.

\bibitem{Harvey:1995fq}
J.~A. Harvey and G.~W. Moore, ``{Algebras, BPS states, and strings},''
  \href{http://dx.doi.org/10.1016/0550-3213(95)00605-2}{{\em Nucl. Phys. B}
  {\bfseries 463} (1996) 315--368},
  \href{http://arxiv.org/abs/hep-th/9510182}{{\ttfamily arXiv:hep-th/9510182}}.

\bibitem{Harvey:1996gc}
J.~A. Harvey and G.~W. Moore, ``{On the algebras of BPS states},''
  \href{http://dx.doi.org/10.1007/s002200050461}{{\em Commun. Math. Phys.}
  {\bfseries 197} (1998) 489--519},
  \href{http://arxiv.org/abs/hep-th/9609017}{{\ttfamily arXiv:hep-th/9609017}}.

\bibitem{Kontsevich:2008fj}
M.~Kontsevich and Y.~Soibelman, ``{Stability structures, motivic
  Donaldson-Thomas invariants and cluster transformations},''
  \href{http://arxiv.org/abs/0811.2435}{{\ttfamily arXiv:0811.2435 [math.AG]}}.

\bibitem{Kontsevich:2010px}
M.~Kontsevich and Y.~Soibelman, ``{Cohomological Hall algebra, exponential
  Hodge structures and motivic Donaldson-Thomas invariants},''
  \href{http://dx.doi.org/10.4310/CNTP.2011.v5.n2.a1}{{\em Commun. Num. Theor.
  Phys.} {\bfseries 5} (2011) 231--352},
  \href{http://arxiv.org/abs/1006.2706}{{\ttfamily arXiv:1006.2706 [math.AG]}}.

\bibitem{Galakhov:2018lta}
D.~Galakhov, ``{BPS Hall Algebra of Scattering Hall States},''
  \href{http://dx.doi.org/10.1016/j.nuclphysb.2019.114693}{{\em Nucl. Phys. B}
  {\bfseries 946} (2019) 114693},
  \href{http://arxiv.org/abs/1812.05801}{{\ttfamily arXiv:1812.05801
  [hep-th]}}.

\bibitem{rapcak2021branesquiversbpsalgebras}
M.~Rapcak, ``{Branes, Quivers and BPS Algebras},''
  \href{http://arxiv.org/abs/2112.13878}{{\ttfamily arXiv:2112.13878
  [hep-th]}}.

\bibitem{Yamazaki:2008bt}
M.~Yamazaki, ``{Brane Tilings and Their Applications},''
  \href{http://dx.doi.org/10.1002/prop.200810536}{{\em Fortsch. Phys.}
  {\bfseries 56} (2008) 555--686},
  \href{http://arxiv.org/abs/0803.4474}{{\ttfamily arXiv:0803.4474 [hep-th]}}.

\bibitem{Li:2020rij}
W.~Li and M.~Yamazaki, ``{Quiver Yangian from Crystal Melting},''
  \href{http://dx.doi.org/10.1007/JHEP11(2020)035}{{\em JHEP} {\bfseries 11}
  (2020) 035}, \href{http://arxiv.org/abs/2003.08909}{{\ttfamily
  arXiv:2003.08909 [hep-th]}}.

\bibitem{Ooguri:2009ijd}
H.~Ooguri and M.~Yamazaki, ``{Crystal Melting and Toric Calabi-Yau
  Manifolds},'' \href{http://dx.doi.org/10.1007/s00220-009-0836-y}{{\em Commun.
  Math. Phys.} {\bfseries 292} (2009) 179--199},
  \href{http://arxiv.org/abs/0811.2801}{{\ttfamily arXiv:0811.2801 [hep-th]}}.

\bibitem{Aganagic:2010qr}
M.~Aganagic and K.~Schaeffer, ``{Wall Crossing, Quivers and Crystals},''
  \href{http://dx.doi.org/10.1007/JHEP10(2012)153}{{\em JHEP} {\bfseries 10}
  (2012) 153}, \href{http://arxiv.org/abs/1006.2113}{{\ttfamily arXiv:1006.2113
  [hep-th]}}.

\bibitem{Yamazaki:2010fz}
M.~Yamazaki, ``{Crystal Melting and Wall Crossing Phenomena},''
  \href{http://dx.doi.org/10.1142/S0217751X11051482}{{\em Int. J. Mod. Phys. A}
  {\bfseries 26} (2011) 1097--1228},
  \href{http://arxiv.org/abs/1002.1709}{{\ttfamily arXiv:1002.1709 [hep-th]}}.

\bibitem{Rapcak:2018nsl}
M.~Rapcak, Y.~Soibelman, Y.~Yang, and G.~Zhao, ``{Cohomological Hall algebras,
  vertex algebras and instantons},''
  \href{http://dx.doi.org/10.1007/s00220-019-03575-5}{{\em Commun. Math. Phys.}
  {\bfseries 376} no.~3, (2019) 1803--1873},
  \href{http://arxiv.org/abs/1810.10402}{{\ttfamily arXiv:1810.10402
  [math.QA]}}.

\bibitem{Galakhov:2020vyb}
D.~Galakhov and M.~Yamazaki, ``{Quiver Yangian and Supersymmetric Quantum
  Mechanics},'' \href{http://dx.doi.org/10.1007/s00220-022-04490-y}{{\em
  Commun. Math. Phys.} {\bfseries 396} no.~2, (2022) 713--785},
  \href{http://arxiv.org/abs/2008.07006}{{\ttfamily arXiv:2008.07006
  [hep-th]}}.

\bibitem{NakajimaALE}
H.~Nakajima, ``{Instantons on ALE spaces, quiver varieties, and Kac-Moody
  algebras},'' {\em Duke Mathematical Journal} {\bfseries 76} no.~2, (1994) 365
  -- 416.

\bibitem{Pestun:2016qko}
V.~Pestun, ``{Review of localization in geometry},''
  \href{http://dx.doi.org/10.1088/1751-8121/aa6161}{{\em J. Phys. A} {\bfseries
  50} no.~44, (2017) 443002}, \href{http://arxiv.org/abs/1608.02954}{{\ttfamily
  arXiv:1608.02954 [hep-th]}}.

\bibitem{guillemin2013supersymmetry}
{Guillemin, Victor W and Sternberg, Shlomo}, {\em {Supersymmetry and
  equivariant de Rham theory}}.
\newblock Springer Science \& Business Media, 2013.

\bibitem{Prochazka:2015deb}
T.~Proch\'azka, ``{$\mathcal{W}$ -symmetry, topological vertex and affine
  Yangian},'' \href{http://dx.doi.org/10.1007/JHEP10(2016)077}{{\em JHEP}
  {\bfseries 10} (2016) 077}, \href{http://arxiv.org/abs/1512.07178}{{\ttfamily
  arXiv:1512.07178 [hep-th]}}.

\bibitem{Tsymbaliuk_2017}
A.~Tsymbaliuk, ``{The affine Yangian of $\mathfrak{gl}_1$ revisited},''
  \href{http://dx.doi.org/10.1016/j.aim.2016.08.041}{{\em Adv. Math.}
  {\bfseries 304} (2017) 583--645},
  \href{http://arxiv.org/abs/1404.5240}{{\ttfamily arXiv:1404.5240 [math.RT]}}.

\bibitem{maulik2018quantumgroupsquantumcohomology}
{Davesh Maulik and Andrei Okounkov}, ``{Quantum Groups and Quantum
  Cohomology},'' \href{http://arxiv.org/abs/1211.1287}{{\ttfamily
  arXiv:1211.1287 [math.AG]}}.

\bibitem{calogero1971solution}
F.~Calogero, ``{Solution of the one-dimensional N-body problems with quadratic
  and/or inversely quadratic pair potentials},'' {\em Journal of Mathematical
  Physics} {\bfseries 12} no.~3, (1971) 419--436.

\bibitem{MOSER1975197}
J.~Moser, ``{Three integrable Hamiltonian systems connected with isospectral
  deformations},'' {\em Advances in Mathematics} {\bfseries 16} no.~2, (1975)
  197--220.

\bibitem{Wang:2022fxr}
R.~Wang, F.~Liu, C.-H. Zhang, and W.-Z. Zhao, ``{Superintegrability for ($\beta
  $-deformed) partition function hierarchies with W-representations},''
  \href{http://dx.doi.org/10.1140/epjc/s10052-022-10875-z}{{\em Eur. Phys. J.
  C} {\bfseries 82} no.~10, (2022) 902},
  \href{http://arxiv.org/abs/2206.13038}{{\ttfamily arXiv:2206.13038
  [hep-th]}}.

\bibitem{Mironov:2023pnd}
A.~Mironov, V.~Mishnyakov, A.~Morozov, A.~Popolitov, R.~Wang, and W.-Z. Zhao,
  ``{Interpolating matrix models for WLZZ series},''
  \href{http://dx.doi.org/10.1140/epjc/s10052-023-11549-0}{{\em Eur. Phys. J.
  C} {\bfseries 83} no.~5, (2023) 377},
  \href{http://arxiv.org/abs/2301.04107}{{\ttfamily arXiv:2301.04107
  [hep-th]}}.

\bibitem{Mironov:2023zwi}
A.~Mironov and A.~Morozov, ``{Many-body integrable systems implied by WLZZ
  models},'' \href{http://dx.doi.org/10.1016/j.physletb.2023.137964}{{\em Phys.
  Lett. B} {\bfseries 842} (2023) 137964},
  \href{http://arxiv.org/abs/2303.05273}{{\ttfamily arXiv:2303.05273
  [hep-th]}}.

\bibitem{Galakhov_2024}
D.~Galakhov, A.~Gavshin, A.~Morozov, and N.~Tselousov, ``{Algorithms for
  representations of quiver Yangian algebras},''
  \href{http://dx.doi.org/10.1007/JHEP08(2024)209}{{\em JHEP} {\bfseries 08}
  (2024) 209}, \href{http://arxiv.org/abs/2406.20074}{{\ttfamily
  arXiv:2406.20074 [hep-th]}}.

\bibitem{nakajima1996jackpolynomialshilbertschemes}
H.~Nakajima, ``{Jack polynomials and Hilbert schemes of points on surfaces},''
  \href{http://arxiv.org/abs/alg-geom/9610021}{{\ttfamily
  arXiv:alg-geom/9610021 [alg-geom]}}.

\bibitem{van2012calogero}
J.~van Diejen and L.~Vinet, {\em {Calogero—Moser— Sutherland Models}}.
\newblock CRM Series in Mathematical Physics. Springer New York, 2012.
\newblock \url{https://books.google.com.vn/books?id=sBXSBwAAQBAJ}.

\bibitem{Mishnyakov:2024cgl}
V.~Mishnyakov and I.~Myakutin, ``{Superintegrability of the monomial Uglov
  matrix model},'' \href{http://arxiv.org/abs/2403.19538}{{\ttfamily
  arXiv:2403.19538 [hep-th]}}.

\bibitem{Mironov_2022}
A.~Mironov and A.~Morozov, ``{Superintegrability summary},''
  \href{http://dx.doi.org/10.1016/j.physletb.2022.137573}{{\em Phys. Lett. B}
  {\bfseries 835} (2022) 137573},
  \href{http://arxiv.org/abs/2201.12917}{{\ttfamily arXiv:2201.12917
  [hep-th]}}.

\bibitem{Azheev:2025wti}
B.~Azheev and N.~Tselousov, ``{Towards construction of superintegrable basis in
  matrix models},''
  \href{http://dx.doi.org/10.1016/j.nuclphysb.2025.116975}{{\em Nucl. Phys. B}
  {\bfseries 1018} (2025) 116975},
  \href{http://arxiv.org/abs/2503.07583}{{\ttfamily arXiv:2503.07583
  [hep-th]}}.

\bibitem{Morozov_2023}
A.~Morozov and N.~Tselousov, ``{3-Schurs from explicit representation of
  Yangian $ \textrm{Y}\left({\hat{\mathfrak{gl}}}_1\right) $. Levels
  1{\textendash}5},'' \href{http://dx.doi.org/10.1007/JHEP11(2023)165}{{\em
  JHEP} {\bfseries 11} (2023) 165},
  \href{http://arxiv.org/abs/2305.12282}{{\ttfamily arXiv:2305.12282
  [hep-th]}}.

\bibitem{Galakhov:2024mbz}
D.~Galakhov, A.~Morozov, and N.~Tselousov, ``{Simple Representations of BPS
  Algebras: the case of $Y(\widehat{\mathfrak{gl}}_2)$},''
  \href{http://arxiv.org/abs/2402.05920}{{\ttfamily arXiv:2402.05920
  [hep-th]}}.

\bibitem{Galakhov:2023mak}
D.~Galakhov, A.~Morozov, and N.~Tselousov, ``{Super-Schur polynomials for
  Affine Super Yangian Y($ \hat{\mathfrak{gl}} _{1|1}$)},''
  \href{http://dx.doi.org/10.1007/JHEP08(2023)049}{{\em JHEP} {\bfseries 08}
  (2023) 049}, \href{http://arxiv.org/abs/2307.03150}{{\ttfamily
  arXiv:2307.03150 [hep-th]}}.

\bibitem{Galakhov:2024foa}
D.~Galakhov, A.~Morozov, and N.~Tselousov, ``{Wall-crossing effects on quiver
  BPS algebras},'' \href{http://dx.doi.org/10.1007/JHEP05(2024)118}{{\em JHEP}
  {\bfseries 05} (2024) 118}, \href{http://arxiv.org/abs/2403.14600}{{\ttfamily
  arXiv:2403.14600 [hep-th]}}.

\bibitem{Galakhov_2024mpsp}
D.~Galakhov, A.~Morozov, and N.~Tselousov, ``{Macdonald polynomials for
  super-partitions},''
  \href{http://dx.doi.org/10.1016/j.physletb.2024.138911}{{\em Phys. Lett. B}
  {\bfseries 856} (2024) 138911},
  \href{http://arxiv.org/abs/2407.03301}{{\ttfamily arXiv:2407.03301
  [hep-th]}}.

\bibitem{Bao:2023ece}
J.~Bao, ``{More on affine Dynkin quiver Yangians},''
  \href{http://dx.doi.org/10.1007/JHEP07(2023)153}{{\em JHEP} {\bfseries 07}
  (2023) 153}, \href{http://arxiv.org/abs/2304.00767}{{\ttfamily
  arXiv:2304.00767 [hep-th]}}.

\bibitem{Li:2023zub}
W.~Li, ``{Quiver algebras and their representations for arbitrary quivers},''
  \href{http://dx.doi.org/10.1007/JHEP12(2024)089}{{\em JHEP} {\bfseries 12}
  (2024) 089}, \href{http://arxiv.org/abs/2303.05521}{{\ttfamily
  arXiv:2303.05521 [hep-th]}}.

\bibitem{Bykov:2019cst}
D.~Bykov and P.~Zinn-Justin, ``{Higher spin ${{\mathfrak {s}}}{{\mathfrak
  {l}}}_2$R-matrix from equivariant (co)homology},'' {\em Lett. Math. Phys.}
  {\bfseries 110} no.~9, (2020) 2435--2470,
  \href{http://arxiv.org/abs/1904.11107}{{\ttfamily arXiv:1904.11107
  [math-ph]}}.

\bibitem{Yang:2024ubh}
Y.~Yang and P.~Zinn-Justin, ``{Higher spin representations of the Yangian of
  $\mathfrak{sl}_2$ and R-matrices},''
  \href{http://arxiv.org/abs/2403.17433}{{\ttfamily arXiv:2403.17433
  [math.RT]}}.

\bibitem{molev2002gelfandtsetlinbasesclassicallie}
A.~I. Molev, ``{Gelfand-Tsetlin bases for classical Lie algebras},''
  \href{http://arxiv.org/abs/math/0211289}{{\ttfamily arXiv:math/0211289
  [math.RT]}}.

\bibitem{Gelfand:1950ihs}
I.~M. Gelfand and M.~L. Tsetlin, ``{Finite-dimensional representations of the
  group of unimodular matrices},'' {\em Dokl. Akad. Nauk SSSR} {\bfseries 71}
  no.~5, (1950) 825--828.

\bibitem{Galakhov:2021xum}
D.~Galakhov, W.~Li, and M.~Yamazaki, ``{Shifted quiver Yangians and
  representations from BPS crystals},''
  \href{http://dx.doi.org/10.1007/JHEP08(2021)146}{{\em JHEP} {\bfseries 08}
  (2021) 146}, \href{http://arxiv.org/abs/2106.01230}{{\ttfamily
  arXiv:2106.01230 [hep-th]}}.

\bibitem{Galakhov:2021vbo}
D.~Galakhov, W.~Li, and M.~Yamazaki, ``{Toroidal and elliptic quiver BPS
  algebras and beyond},'' \href{http://dx.doi.org/10.1007/JHEP02(2022)024}{{\em
  JHEP} {\bfseries 02} (2022) 024},
  \href{http://arxiv.org/abs/2108.10286}{{\ttfamily arXiv:2108.10286
  [hep-th]}}.

\bibitem{Noshita:2021dgj}
G.~Noshita and A.~Watanabe, ``{Shifted quiver quantum toroidal algebra and
  subcrystal representations},''
  \href{http://dx.doi.org/10.1007/JHEP05(2022)122}{{\em JHEP} {\bfseries 05}
  (2022) 122}, \href{http://arxiv.org/abs/2109.02045}{{\ttfamily
  arXiv:2109.02045 [hep-th]}}.

\bibitem{Kodera:2016faj}
R.~Kodera and H.~Nakajima, ``{Quantized Coulomb branches of Jordan quiver gauge
  theories and cyclotomic rational Cherednik algebras},'' {\em Proc. Symp. Pure
  Math.} {\bfseries 98} (2018) 49--78,
  \href{http://arxiv.org/abs/1608.00875}{{\ttfamily arXiv:1608.00875
  [math.RT]}}.

\bibitem{Negut:2022pka}
A.~Negu\c{t}, ``{Quantum loop groups for arbitrary quivers},''
  \href{http://arxiv.org/abs/2209.09089}{{\ttfamily arXiv:2209.09089
  [math.RT]}}.

\bibitem{Bezerra:2019kfl}
L.~Bezerra and E.~Mukhin, ``{Braid actions on quantum toroidal
  superalgebras},''
  \href{http://dx.doi.org/10.1016/j.jalgebra.2021.06.012}{{\em J. Algebra}
  {\bfseries 585} (2021) 338--369},
  \href{http://arxiv.org/abs/1912.08729}{{\ttfamily arXiv:1912.08729
  [math.QA]}}.

\bibitem{Bezerra_2020}
L.~Bezerra and E.~Mukhin, ``{Quantum Toroidal Algebra Associated with
  $\mathfrak {gl}_{m|n}$},''
  \href{http://dx.doi.org/10.1007/s10468-020-09959-9}{{\em Algebr. Represent.
  Theory} {\bfseries 24} no.~2, (2021) 541--564},
  \href{http://arxiv.org/abs/1904.07297}{{\ttfamily arXiv:1904.07297
  [math.QA]}}.

\bibitem{Bao_2025}
J.~Bao and M.~Yamazaki, ``{Crystals and double quiver algebras from
  Jeffrey-Kirwan residues},''
  \href{http://dx.doi.org/10.21468/SciPostPhys.18.4.143}{{\em SciPost Phys.}
  {\bfseries 18} no.~4, (2025) 143},
  \href{http://arxiv.org/abs/2501.03365}{{\ttfamily arXiv:2501.03365
  [hep-th]}}.

\bibitem{chuang2022hilbertschemesnonreduceddivisors}
W.-y. Chuang, T.~Creutzig, D.~E. Diaconescu, and Y.~Soibelman, ``{Hilbert
  schemes of nonreduced divisors in Calabi-Yau threefolds and W-algebras},''
  \href{http://arxiv.org/abs/1907.13005}{{\ttfamily arXiv:1907.13005
  [math.AG]}}.

\bibitem{Szendroi:2007nu}
B.~Szendroi, ``{Non-commutative Donaldson\textendash{}Thomas invariants and the
  conifold},'' \href{http://dx.doi.org/10.2140/gt.2008.12.1171}{{\em Geom.
  Topol.} {\bfseries 12} no.~2, (2008) 1171--1202},
  \href{http://arxiv.org/abs/0705.3419}{{\ttfamily arXiv:0705.3419 [math.AG]}}.

\bibitem{kirillov2016quiver}
A.~Kirillov and A.~Kirillov, {\em {Quiver Representations and Quiver
  Varieties}}.
\newblock Graduate studies in mathematics. American Mathematical Society, 2016.
\newblock \url{https://books.google.ru/books?id=p8zeswEACAAJ}.

\bibitem{Galakhov_2024csd}
D.~Galakhov and W.~Li, ``{Charging solid partitions},''
  \href{http://dx.doi.org/10.1007/JHEP01(2024)043}{{\em JHEP} {\bfseries 01}
  (2024) 043}, \href{http://arxiv.org/abs/2311.02751}{{\ttfamily
  arXiv:2311.02751 [hep-th]}}.

\bibitem{franco20234dcrystalmeltingtoric}
S.~Franco, ``{4d crystal melting, toric Calabi-Yau 4-folds and brane brick
  models},'' \href{http://dx.doi.org/10.1007/JHEP03(2024)091}{{\em JHEP}
  {\bfseries 03} (2024) 091}, \href{http://arxiv.org/abs/2311.04404}{{\ttfamily
  arXiv:2311.04404 [hep-th]}}.

\bibitem{Bao_2024}
J.~Bao, R.-K. Seong, and M.~Yamazaki, ``{The origin of Calabi-Yau crystals in
  BPS states counting},'' \href{http://dx.doi.org/10.1007/JHEP03(2024)140}{{\em
  JHEP} {\bfseries 03} (2024) 140},
  \href{http://arxiv.org/abs/2401.02792}{{\ttfamily arXiv:2401.02792
  [hep-th]}}.

\bibitem{jeffrey1994localizationnonabeliangroupactions}
L.~C. Jeffrey and F.~C. Kirwan, ``{Localization for nonabelian group
  actions},'' \href{http://arxiv.org/abs/alg-geom/9307001}{{\ttfamily
  arXiv:alg-geom/9307001}}.

\bibitem{Galakhov:2022uyu}
D.~Galakhov, W.~Li, and M.~Yamazaki, ``{Gauge/Bethe correspondence from quiver
  BPS algebras},'' \href{http://dx.doi.org/10.1007/JHEP11(2022)119}{{\em JHEP}
  {\bfseries 11} (2022) 119}, \href{http://arxiv.org/abs/2206.13340}{{\ttfamily
  arXiv:2206.13340 [hep-th]}}.

\bibitem{Chen:2025xoe}
T.~Chen and W.~Li, ``{Quiver Yangians as Coulomb branch algebras},''
  \href{http://arxiv.org/abs/2502.01323}{{\ttfamily arXiv:2502.01323
  [hep-th]}}.

\bibitem{Rapcak:2020ueh}
M.~Rapcak, Y.~Soibelman, Y.~Yang, and G.~Zhao, ``{Cohomological Hall algebras
  and perverse coherent sheaves on toric Calabi\textendash{}Yau $3$-folds},''
  \href{http://dx.doi.org/10.4310/CNTP.2023.v17.n4.a2}{{\em Commun. Num. Theor.
  Phys.} {\bfseries 17} no.~4, (2023) 847--939},
  \href{http://arxiv.org/abs/2007.13365}{{\ttfamily arXiv:2007.13365
  [math.QA]}}.

\bibitem{king1994moduli}
A.~D. King, ``{Moduli of representations of finite dimensional algebras},''
  {\em The Quarterly Journal of Mathematics} {\bfseries 45} no.~4, (1994)
  515--530.

\bibitem{Nakajima98}
H.~Nakajima, ``{Quiver varieties and Kac-Moody algebras},'' {\em Duke
  Mathematical Journal} {\bfseries 91} no.~3, (1998) 515.

\bibitem{Nakajima_lect}
H.~{Nakajima}, ``{More lectures on Hilbert schemes of points on surfaces},''
  \href{http://arxiv.org/abs/1401.6782}{{\ttfamily arXiv:1401.6782 [math.RT]}}.

\bibitem{Witten:1982im}
E.~Witten, ``{Supersymmetry and Morse theory},'' {\em J. Diff. Geom.}
  {\bfseries 17} no.~4, (1982) 661--692.

\bibitem{Cordes:1994fc}
S.~Cordes, G.~W. Moore, and S.~Ramgoolam, ``{Lectures on 2-d Yang-Mills theory,
  equivariant cohomology and topological field theories},''
  \href{http://dx.doi.org/10.1016/0920-5632(95)00434-B}{{\em Nucl. Phys. B
  Proc. Suppl.} {\bfseries 41} (1995) 184--244},
  \href{http://arxiv.org/abs/hep-th/9411210}{{\ttfamily arXiv:hep-th/9411210}}.

\bibitem{huybrechts2006fourier}
D.~Huybrechts, {\em Fourier-Mukai transforms in algebraic geometry}.
\newblock Clarendon Press, 2006.

\bibitem{zbMATH03814600}
N.~Berline and M.~Vergne, ``{Zeros d'un champ de vecteurs et classes
  characteristiques {\'e}quivariantes},''
  \href{http://dx.doi.org/10.1215/S0012-7094-83-05024-X}{{\em Duke Math. J.}
  {\bfseries 50} (1983) 539--549}.

\bibitem{Atiyah:1984px}
M.~F. Atiyah and R.~Bott, ``{The moment map and equivariant cohomology},''
  \href{http://dx.doi.org/10.1016/0040-9383(84)90021-1}{{\em Topology}
  {\bfseries 23} (1984) 1--28}.

\bibitem{Alekseev:2000fe}
A.~Alekseev, ``{Notes on equivariant localization},''
  \href{http://dx.doi.org/10.1007/3-540-46552-9_1}{{\em Lect. Notes Phys.}
  {\bfseries 543} (2000) 1--24}.

\bibitem{McKay:1981aap}
J.~McKay, {\em {Graphs, singularities, and finite groups}}.
\newblock 1981.
\newblock \url{https://books.google.com.vn/books?id=_9sDCAAAQBAJ}.

\bibitem{Bridgeland:2001xf}
T.~Bridgeland, A.~King, and M.~Reid, ``{The McKay correspondence as an
  equivalence of derived categories},''
  \href{http://dx.doi.org/10.1090/S0894-0347-01-00368-X}{{\em J. Am. Math.
  Soc.} {\bfseries 14} (2001) 535--554}.

\bibitem{mclean2023mckay}
M.~McLean and A.~F. Ritter, ``{The McKay correspondence for isolated
  singularities via Floer theory},'' {\em Journal of Differential Geometry}
  {\bfseries 124} no.~1, (2023) 113--168.

\bibitem{cirafici2013curvecountinginstantonsmckay}
M.~Cirafici and R.~J. Szabo, ``{Curve counting, instantons and McKay
  correspondences},''
  \href{http://dx.doi.org/10.1016/j.geomphys.2013.03.020}{{\em J. Geom. Phys.}
  {\bfseries 72} (2013) 54--109},
  \href{http://arxiv.org/abs/1209.1486}{{\ttfamily arXiv:1209.1486 [hep-th]}}.

\bibitem{Slodowy_1980}
P.~Slodowy, \href{http://dx.doi.org/10.1007/bfb0090294}{{\em {Simple
  Singularities and Simple Algebraic Groups}}}.
\newblock Springer Berlin Heidelberg, 1980.

\bibitem{stekolshchik2005notescoxetertransformationsmckay}
R.~Stekolshchik, ``{Notes on Coxeter Transformations and the McKay
  correspondence},'' \href{http://arxiv.org/abs/math/0510216}{{\ttfamily
  arXiv:math/0510216 [math.RT]}}.

\bibitem{Bershadsky_1996}
M.~Bershadsky, K.~A. Intriligator, S.~Kachru, D.~R. Morrison, V.~Sadov, and
  C.~Vafa, ``{Geometric singularities and enhanced gauge symmetries},''
  \href{http://dx.doi.org/10.1016/S0550-3213(96)90131-5}{{\em Nucl. Phys. B}
  {\bfseries 481} (1996) 215--252},
  \href{http://arxiv.org/abs/hep-th/9605200}{{\ttfamily arXiv:hep-th/9605200}}.

\bibitem{Cecotti:2012gh}
S.~Cecotti and M.~Del~Zotto, ``{4d N=2 Gauge Theories and Quivers: the
  Non-Simply Laced Case},''
  \href{http://dx.doi.org/10.1007/JHEP10(2012)190}{{\em JHEP} {\bfseries 10}
  (2012) 190}, \href{http://arxiv.org/abs/1207.7205}{{\ttfamily arXiv:1207.7205
  [hep-th]}}.

\bibitem{Ginzburg:2006fu}
V.~Ginzburg, ``{Calabi-Yau algebras},''
  \href{http://arxiv.org/abs/math/0612139}{{\ttfamily arXiv:math/0612139}}.

\bibitem{Guay:2018wel}
N.~Guay, V.~Regelskis, and C.~Wendlandt, ``{Vertex Representations for Yangians
  of Kac-Moody algebras},'' \href{http://arxiv.org/abs/1804.04081}{{\ttfamily
  arXiv:1804.04081 [math.RT]}}.

\bibitem{nazarov2000representationsyangiansgelfandzetlinbases}
M.~Nazarov and V.~Tarasov, ``{Representations of Yangians with Gelfand-Zetlin
  Bases},'' \href{http://arxiv.org/abs/q-alg/9502008}{{\ttfamily
  arXiv:q-alg/9502008 [q-alg]}}.

\bibitem{Cherednik1987ANI}
I.~Cherednik, ``{A new interpretation of Gelfand-Tzetlin bases},'' {\em Duke
  Mathematical Journal} {\bfseries 54} (1987) 563--577.
  \url{https://api.semanticscholar.org/CorpusID:121056379}.

\bibitem{Ginzburg_lect}
{{Ginzburg}, Victor}, ``{Lectures on Nakajima's Quiver Varieties},''
  \href{http://arxiv.org/abs/0905.0686}{{\ttfamily arXiv:0905.0686 [math.RT]}}.

\bibitem{Bao:2025dqs}
J.~Bao, ``{An Overview of Crystals and Double Quiver Yangians},''
  \href{http://arxiv.org/abs/2509.16918}{{\ttfamily arXiv:2509.16918
  [hep-th]}}.

\bibitem{Noshita:2021ldl}
G.~Noshita and A.~Watanabe, ``{A note on quiver quantum toroidal algebra},''
  \href{http://dx.doi.org/10.1007/JHEP05(2022)011}{{\em JHEP} {\bfseries 05}
  (2022) 011}, \href{http://arxiv.org/abs/2108.07104}{{\ttfamily
  arXiv:2108.07104 [hep-th]}}.

\bibitem{Bryan_2010}
B.~Young and J.~Bryan, ``{Generating functions for colored 3D Young diagrams
  and the Donaldson-Thomas invariants of orbifolds},''
  \href{http://dx.doi.org/10.1215/00127094-2010-009}{{\em Duke Math. J.}
  {\bfseries 152} (2010) 115--153},
  \href{http://arxiv.org/abs/0802.3948}{{\ttfamily arXiv:0802.3948 [math.CO]}}.

\bibitem{Nakajimarnq195}
K.~Nagao and H.~Nakajima, ``{Counting invariant of perverse coherent sheaves
  and its wall-crossing},'' \href{http://arxiv.org/abs/0809.2992}{{\ttfamily
  arXiv:0809.2992 [math.AG]}}.

\end{thebibliography}\endgroup

\end{document}